\input amstex
\documentstyle{amsppt}
\magnification=\magstep 1
\document
\vsize6.7in


\chardef\oldatsign=\catcode`\@
\catcode`\@=11
\newif\ifdraftmode			
\global\draftmodefalse



%
\font@\twelverm=cmr12 
\font@\twelvei=cmmi12 \skewchar\twelvei='177 
\font@\twelvesy=cmsy10 scaled\magstep1 \skewchar\twelvesy='060 
\font@\twelveex=cmex10 scaled\magstep1 
\font@\twelvemsa=msam10 scaled\magstep1 
\font@\twelvemsb=msbm10 scaled\magstep1 
\font@\twelvebf=cmbx12 
\font@\twelvett=cmtt12 
\font@\twelvesl=cmsl12 
\font@\twelveit=cmti12 
\font@\twelvesmc=cmcsc10 scaled\magstep1 
%
%
\font@\ninerm=cmr9 
\font@\ninei=cmmi9 \skewchar\ninei='177 
\font@\ninesy=cmsy9 \skewchar\ninesy='60 
\font@\ninemsa=msam9
\font@\ninemsb=msbm9
\font@\ninebf=cmbx9
%
%
%
\font@\ttlrm=cmbx12 scaled \magstep2 
\font@\ttlsy=cmsy10 scaled \magstep3 
\font@\tensmc=cmcsc10 
%
%
\def\normaltype{
	\def\pointsize@{12}%
	\abovedisplayskip18\p@ plus5\p@ minus9\p@
	\belowdisplayskip18\p@ plus5\p@ minus9\p@
	\abovedisplayshortskip1\p@ plus3\p@
	\belowdisplayshortskip9\p@ plus3\p@ minus4\p@
	\textonlyfont@\rm\twelverm
	\textonlyfont@\it\twelveit
	\textonlyfont@\sl\twelvesl
	\textonlyfont@\bf\twelvebf
	\textonlyfont@\smc\twelvesmc
	\ifsyntax@
		\def\big##1{{\hbox{$\left##1\right.$}}}%
	\else
		\let\big\twelvebig@
 \textfont0=\twelverm \scriptfont0=\ninerm \scriptscriptfont0=\sevenrm
 \textfont1=\twelvei  \scriptfont1=\ninei  \scriptscriptfont1=\seveni
 \textfont2=\twelvesy \scriptfont2=\ninesy \scriptscriptfont2=\sevensy
 \textfont3=\twelveex \scriptfont3=\twelveex  \scriptscriptfont3=\twelveex
 \textfont\itfam=\twelveit \def\it{\fam\itfam\twelveit}%
 \textfont\slfam=\twelvesl \def\sl{\fam\slfam\twelvesl}%
 \textfont\bffam=\twelvebf \def\bf{\fam\bffam\twelvebf}%
 \scriptfont\bffam=\ninebf \scriptscriptfont\bffam=\sevenbf
 \textfont\ttfam=\twelvett \def\tt{\fam\ttfam\twelvett}%
 \textfont\msafam=\twelvemsa \scriptfont\msafam=\ninemsa
 \scriptscriptfont\msafam=\sevenmsa
 \textfont\msbfam=\twelvemsb \scriptfont\msbfam=\ninemsb
 \scriptscriptfont\msbfam=\sevenmsb
	\fi
 \normalbaselineskip=\twelvebaselineskip
 \setbox\strutbox=\hbox{\vrule height12\p@ depth6\p@
      width0\p@}%
 \normalbaselines\rm \ex@=.2326ex%
}
%
%
%
\def\smalltype{
	\def\pointsize@{10}%
	\abovedisplayskip12\p@ plus3\p@ minus9\p@
	\belowdisplayskip12\p@ plus3\p@ minus9\p@
	\abovedisplayshortskip\z@ plus3\p@
	\belowdisplayshortskip7\p@ plus3\p@ minus4\p@
	\textonlyfont@\rm\tenrm
	\textonlyfont@\it\tenit
	\textonlyfont@\sl\tensl
	\textonlyfont@\bf\tenbf
	\textonlyfont@\smc\tensmc
	\ifsyntax@
		\def\big##1{{\hbox{$\left##1\right.$}}}%
	\else
		\let\big\tenbig@
	\textfont0=\tenrm \scriptfont0=\sevenrm \scriptscriptfont0=\fiverm 
	\textfont1=\teni  \scriptfont1=\seveni  \scriptscriptfont1=\fivei
	\textfont2=\tensy \scriptfont2=\sevensy \scriptscriptfont2=\fivesy 
	\textfont3=\tenex \scriptfont3=\tenex \scriptscriptfont3=\tenex
	\textfont\itfam=\tenit \def\it{\fam\itfam\tenit}%
	\textfont\slfam=\tensl \def\sl{\fam\slfam\tensl}%
	\textfont\bffam=\tenbf \def\bf{\fam\bffam\tenbf}%
	\scriptfont\bffam=\sevenbf \scriptscriptfont\bffam=\fivebf
	\textfont\msafam=\tenmsa
	\scriptfont\msafam=\sevenmsa
	\scriptscriptfont\msafam=\fivemsa
	\textfont\msbfam=\tenmsb
	\scriptfont\msbfam=\sevenmsb
	\scriptscriptfont\msbfam=\fivemsb
		\textfont\ttfam=\tentt \def\tt{\fam\ttfam\tentt}%
	\fi
 \normalbaselineskip 14\p@
 \setbox\strutbox=\hbox{\vrule height10\p@ depth4\p@ width0\p@}%
 \normalbaselines\rm \ex@=.2326ex%
}

\def\titletype{
	\def\pointsize@{17}%
	\textonlyfont@\rm\ttlrm
	\ifsyntax@
		\def\big##1{{\hbox{$\left##1\right.$}}}%
	\else
		\let\big\twelvebig@
		\textfont0=\ttlrm \scriptfont0=\twelverm
		\scriptscriptfont0=\tenrm
		\textfont2=\ttlsy \scriptfont2=\twelvesy
		\scriptscriptfont2=\tensy
	\fi
	\normalbaselineskip 25\p@
	\setbox\strutbox=\hbox{\vrule height17\p@ depth8\p@ width0\p@}%
	\normalbaselines
	\rm
	\ex@=.2326ex%
}

\def\tenbig@#1{
	{%
		\hbox{%
			$%
			\left
			#1%
			\vbox to8.5\p@{}%
			\right.%
			\n@space
			$%
		}%
	}%
}

\def\twelvebig@#1{%
	{%
		\hbox{%
			$%
			\left
			#1%
			\vbox to10.2\p@{}
			\right.%
			\n@space
			$%
		}%
	}%
}

%
%
%
%
%
\newif\ifl@beloutopen
\newwrite\l@belout
\newread\l@belin

\global\let\currentfile=\jobname

\def\getfile#1{%
	\immediate\closeout\l@belout
	\global\l@beloutopenfalse
	\gdef\currentfile{#1}%
	\input #1%
	\par
	\newpage
}

\def\getxrefs#1{%
	\bgroup
		\def\gobble##1{}
		\edef\list@{#1,}%
		\def\gr@boff##1,##2\end{
			\openin\l@belin=##1.xref
			\ifeof\l@belin
			\else
				\closein\l@belin
				\input ##1.xref
			\fi
			\def\list@{##2}%
			\ifx\list@\empty
				\let\next=\gobble
			\else
				\let\next=\gr@boff
			\fi
			\expandafter\next\list@\end
		}%
		\expandafter\gr@boff\list@\end
	\egroup
}

\def\testdefined#1#2#3{%
	\expandafter\ifx
	\csname #1\endcsname
	\relax
	#3%
	\else #2\fi
}

\def\document{%
	\minaw@11.11128\ex@ 
	\def\alloclist@{\empty}%
	\def\fontlist@{\empty}%
	\openin\l@belin=\jobname.xref	
	\ifeof\l@belin\else
		\closein\l@belin
		\input \jobname.xref
	\fi
}

\def\getst@te#1#2{%
	\edef\st@te{\csname #1s!#2\endcsname}%
	\expandafter\ifx\st@te\relax
		\def\st@te{0}%
	\fi
}

\def\setst@te#1#2#3{%
	\expandafter
	\gdef\csname #1s!#2\endcsname{#3}%
}

\outer\def\setupautolabel#1#2{%
	\def\newcount@{\global\alloc@0\count\countdef\insc@unt}	
	\def\newtoks@{\global\alloc@5\toks\toksdef\@cclvi}
	\expandafter\newcount@\csname #1Number\endcsname
	\expandafter\global\csname #1Number\endcsname=1%
	\expandafter\newtoks@\csname #1l@bel\endcsname
	\expandafter\global\csname #1l@bel\endcsname={#2}%
}

\def\reflabel#1#2{%
	\testdefined{#1l@bel}
	{
		\getst@te{#1}{#2}%
		\ifcase\st@te
			???
			\message{Unresolved forward reference to
				label #2. Use another pass.}%
		\or	
			\setst@te{#1}{#2}2
			\csname #1l!#2\endcsname 
		\or	
			\csname #1l!#2\endcsname 
		\or	
			\csname #1l!#2\endcsname 
		\fi
	}{
		{\escapechar=-1 
		\errmessage{You haven't done a
			\string\\setupautolabel\space for type #1!}%
		}%
	}%
}

{\catcode`\{=12 \catcode`\}=12
	\catcode`\[=1 \catcode`\]=2
	\xdef\Lbrace[{]
	\xdef\Rbrace[}]%
]%

\def\setlabel#1#2{%
	\testdefined{#1l@bel}
	{
		\edef\templ@bel@{\expandafter\the
			\csname #1l@bel\endcsname}%
		\def\@rgtwo{#2}%
		\ifx\@rgtwo\empty
		\else
			\ifl@beloutopen\else
				\immediate\openout\l@belout=\currentfile.xref
				\global\l@beloutopentrue
			\fi
			\getst@te{#1}{#2}%
			\ifcase\st@te
			\or	
			\or	
				\edef\oldnumber@{\csname #1l!#2\endcsname}%
				\edef\newnumber@{\templ@bel@}%
				\ifx\newnumber@\oldnumber@
				\else
					\message{A forward reference to label 
						#2 has been resolved
						incorrectly.  Use another
						pass.}%
				\fi
			\or	
				\errmessage{Same label #2 used in two
					\string\setlabel s!}%
			\fi
			\expandafter\xdef\csname #1l!#2\endcsname
				{\templ@bel@}
			\setst@te{#1}{#2}3%
			\immediate\write\l@belout 
				{\string\expandafter\string\gdef
				\string\csname\space #1l!#2%
				\string\endcsname
				\Lbrace\templ@bel@\Rbrace
				}%
			\immediate\write\l@belout 
				{\string\expandafter\string\gdef
				\string\csname\space #1s!#2%
				\string\endcsname
				\Lbrace 1\Rbrace
				}%
		\fi
		\templ@bel@	
		\expandafter\ifx\envir@end\endref 
			\gdef\marginalhook@{\marginal{#2}}%
		\else
			\marginal{#2}
		\fi
		\expandafter\global\expandafter\advance	
			\csname #1Number\endcsname
			by 1 %
	}{
		{\escapechar=-1
		\errmessage{You haven't done a \string\\setupautolabel\space
			for type #1!}%
		}%
	}%
}


\newcount\SectionNumber
\setupautolabel{t}{\number\SectionNumber.\number\tNumber}
\setupautolabel{r}{\number\rNumber}
\setupautolabel{T}{\number\TNumber}

\define\rref{\reflabel{r}}
\define\tref{\reflabel{t}}

\define\tnum{\setlabel{t}}
\define\rnum{\setlabel{r}}

%
\def\strutdepth{\dp\strutbox}%
\def\strutheight{\ht\strutbox}%

\newif\iftagmode
\tagmodefalse

\let\old@tagform@=\tagform@
\def\tagform@{\tagmodetrue\old@tagform@}

\def\marginal#1{%
	\ifvmode
	\else
		\strut
	\fi
	\ifdraftmode
		\ifmmode
			\ifinner
				\let\Vorvadjust=\Vadjust
			\else
				\let\Vorvadjust=\vadjust
			\fi
		\else
			\let\Vorvadjust=\Vadjust
		\fi
		\iftagmode	
			\llap{%
				\smalltype
				\vtop to 0pt{%
					\pretolerance=2000
					\tolerance=5000
					\raggedright
					\hsize=.72in
					\parindent=0pt
					\strut
					#1%
					\vss
				}%
				\kern.08in
				\iftagsleft@
				\else
					\kern\hsize
				\fi
			}%
		\else
			\Vorvadjust{%
				\kern-\strutdepth 
				{%
					\smalltype
					\kern-\strutheight 
					\llap{%
						\vtop to 0pt{%
							\kern0pt
							\pretolerance=2000
							\tolerance=5000
							\raggedright
							\hsize=.5in
							\parindent=0pt
							\strut
							#1%
							\vss
						}%
						\kern.08in
					}%
					\kern\strutheight
				}%
				\kern\strutdepth
			}
		\fi
	\fi
}


\newbox\Vadjustbox

\def\Vadjust#1{
	\global\setbox\Vadjustbox=\vbox{#1}%
	\ifmmode
		\ifinner
			\innerVadjust
		\fi		
	\else
		\innerVadjust
	\fi
}

\def\innerVadjust{%
	\def\nexti{\aftergroup\innerVadjust}%
	\def\nextii{%
		\ifvmode
			\hrule height 0pt 
			\box\Vadjustbox
		\else
			\vadjust{\box\Vadjustbox}%
		\fi
	}%
	\ifinner
		\let\next=\nexti
	\else
		\let\next=\nextii
	\fi
	\next
}%

\global\let\marginalhook@\empty

\def\endref{%
\setbox\tw@\box\thr@@
\makerefbox?\thr@@{\endgraf\egroup}%
  \endref@
  \endgraf
  \endgroup
  \keyhook@
  \marginalhook@
  \global\let\keyhook@\empty 
  \global\let\marginalhook@\empty 
}

\catcode`\@=\oldatsign

\nologo
\input xy
\xyoption{all}
\vsize6.7in
\topmatter

\title A study of singularities on rational curves via syzygies\endtitle
  \leftheadtext{Cox, Kustin, Polini, and Ulrich}
\rightheadtext{Singularities via syzygies}
\author David Cox,  
Andrew R. Kustin\footnote{Supported in part by the National Security Agency.\hphantom{the National Science Foundation and XXX}}, Claudia Polini\footnote{Supported in part by the National Science Foundation and the National Security Agency.\hphantom{XXX}}, and Bernd  Ulrich\footnote{Supported in part by the National Science Foundation.\hphantom{and the National Security Agency XXX}}\endauthor
\address
Department of Mathematics, 
Amherst College, 
Amherst,  MA 01002-5000
\endaddress
\email dac\@math.amherst.edu \endemail
 \address
Mathematics Department,
University of South Carolina,
Columbia, SC 29208\endaddress
\email kustin\@math.sc.edu \endemail
\address Mathematics Department,
University of Notre Dame,
Notre Dame, IN 46556\endaddress 
\email     cpolini\@nd.edu \endemail
\address 
Department of Mathematics, 
Purdue University,
West Lafayette, IN 47907
\endaddress
\email          ulrich\@math.purdue.edu \endemail


\endtopmatter

\document

{\eightpoint \flushpar{\bf Abstract.}
Consider a rational projective curve $\Cal C$ of degree $d$ over an algebraically closed field $\pmb k$.  There are $n$ homogeneous forms $g_1,\dots,g_n$ of degree $d$ in $B=\pmb k[x,y]$ which parameterize  $\Cal C$ in a birational, base point free, manner. We study the singularities of $\Cal C$ by studying a  Hilbert-Burch matrix $\varphi$ for the row vector $[g_1,\dots,g_n]$. In the ``General Lemma'' we use   the generalized row ideals of $\varphi$ to identify the singular points on $\Cal C$, their multiplicities, the number of branches at each singular point, and   the multiplicity of each branch.

Let $p$ be a singular point on the parameterized planar curve  $\Cal C$ which corresponds to a
generalized zero of $\varphi$.
In the ``Triple Lemma''  we give a matrix $\varphi'$ whose
 maximal minors parameterize the closure, in $\Bbb P^2$, of the
blow-up at $p$ of $\Cal C$ in a neighborhood of $p$. We  apply the General Lemma  to $\varphi'$ in order to learn about the singularities of $\Cal C$ in the first neighborhood of $p$. If $\Cal C$ has even degree $d=2c$ and the multiplicity of $\Cal C$ at $p$ is equal to $c$, then we apply the Triple Lemma   again to learn about the singularities of $\Cal  C$ in the second neighborhood   of $p$.

Consider  rational plane curves $\Cal C$ of even degree $d=2c$. We classify curves according to the configuration of multiplicity $c$ singularities on or infinitely near $\Cal C$. There are $7$ possible configurations of such singularities. We classify the Hilbert-Burch matrix which corresponds to each configuration. 
The study of multiplicity $c$ singularities on, or infinitely near, a fixed rational plane curve $\Cal C$ of degree $2c$ is equivalent to the study of the scheme of generalized zeros of the fixed balanced Hilbert-Burch matrix $\varphi$ for a parameterization of $\Cal C$. Let   $$\operatorname{BalH}_d=\left\{ \varphi\left\vert \matrix\format\l\\\text{$\varphi$ is a $3\times 2$ matrix; each entry in $\varphi$ is a homogeneous} \\\text{form of degree $c$ from $B$; and  $\operatorname{ht} I_2(\varphi)=2$} \endmatrix\right.\right\}.$$
The group $G=\operatorname{GL}_3(\pmb k)\times \operatorname{GL}_2(\pmb k)$ acts on $\operatorname{BalH}_d$ by way of $(\chi,\xi)\cdot \varphi=\chi\varphi\xi^{-1}$. We decompose $\operatorname{BalH}_d$ into a disjoint union of 11 orbits. Each orbit has the form $G\cdot M$, where $M$ is the closed irreducible subspace of affine space defined by the maximal order minors of a generic matrix.

We introduce the parameter space $\Bbb A_d$. Each element $\pmb g$ of $\Bbb A_d$ is an ordered triple of $d$-forms from $B$, and each $\pmb g\in \Bbb A_d$ induces a rational map $\Psi_{\pmb g}\: \xymatrix{\Bbb P^1\ar@{-->}[r]&\Bbb P^{2}}$. We define 
$$\split \Bbb T_d&{}=\left \{\pmb g\in \Bbb A_d\left\vert \text{$\Psi_{\pmb g}$ is  birational onto its image  without base points}\right.\right\}\ \text{and}\\
\Bbb B_d&{}=\left\{\pmb g\in \Bbb T_d\left\vert \matrix\format\l\\\text{every entry in the corresponding  homogeneous  Hilbert-Burch}\\\text{matrix   has degree $d/2$}\endmatrix\right.\right\}.\endsplit $$ 
In practice we are only interested in the subset $\Bbb T_d$ of $\Bbb A_d$. Every element of $\Bbb A_d$ which is not in $\Bbb T_d$ corresponds to an unsuitable parameterization of a curve. We prove that if  there is a multiplicity $c$ singularity on or infinitely near a point $p$ on a curve $\Cal C$ of degree $d=2c$, then the parameterization of $\Cal C$ is an element of $\Bbb B_d$. We prove that $\Bbb B_d\subseteq \Bbb T_d$ are open subsets of $\Bbb A_d$. We   identify an open cover $\cup \Bbb B_d^{(i)}$ of $\Bbb B_d$ and for each $\Bbb B_d^{(i)}$ in this open cover, we identify a generic Hilbert-Burch matrix which specializes to give a Hilbert-Burch matrix for $\pmb g$ for each $\pmb g\in \Bbb B_d^{(i)}$. As an application of this result, we  identify a universal projective resolution  for the graded Betti numbers
$$0\to B(-3c)^2 \to B(-2c)^3\to B. $$ 

We decompose the space $\Bbb B_d$ of balanced  triples into strata. Each stratum consists of those triples $\pmb g$ in $\Bbb B_d$ for which the corresponding curve $\Cal C_{\pmb g}$ exhibits one particular configuration of multiplicity $c$ singularities. 

We   use the Jacobian matrix associated to the parameterization to identify the non-smooth branches of the curve as well as the multiplicity of each branch. The starting point for this line of reasoning  is the result  that if $D$ is an algebra which is essentially of finite type over the ring $C$, then the ramification locus of $D$ over $C$ is equal to the support of the module of K\"ahler differentials $\Omega_{D/C}$. The General Lemma  is a  local result. Once one knows the singularities $\{p_i\}$ on a parameterized curve $\Cal C$, then the General Lemma     shows how to read the multiplicity and number of branches at each point $p_i$   from the Hilbert-Burch matrix of the parameterization. The result about the Jacobian matrix is a global result. It describes, in terms of the parameterization,  all of the points $p$ on $\Cal C$ and all of the branches of $\Cal C$ at $p$ for which the multiplicity  is at least two.  In contrast to General Lemma  one may apply the Jacobian matrix technique  before one knows the singularities on $\Cal C$. 

  We use conductor techniques to study the singularity degree  $\delta$. 
Let $\pmb g=(g_1,g_2,g_3)$ be an element of $\Bbb T_d$ and $\Cal C_{\pmb g}$ be the corresponding parameterized plane curve.
We produce a polynomial $c_{\pmb g}$ whose factorization into linear factors gives the value of the invariant $\delta$ at each singular point of   $\Cal C_{\pmb g}$. The polynomial $c_{\pmb g}$ is obtained  in a polynomial manner from the coefficients of the entries of a Hilbert-Burch matrix for $\pmb g$.
We use these ideas to produce closed sets in $\Bbb B_d$ which separate various configurations of singularities. To create $c_{\pmb g}$, we start with the coordinate ring $A_{\pmb g}=\pmb k[g_1,g_2,g_3]\subseteq B$. If $\pmb V$ is the $d^{\text{th}}$ Veronese subring of $B$ and $\frak c_{\pmb g}$ is the conductor $A_{\pmb g}\:\! \pmb V$, then $c_{\pmb g}$ generates the saturation of the extension of $\frak c_{\pmb g}$ to $B$.

In the final section of the paper, we apply our results to rational plane quartics.  We exhibit a stratification of $\Bbb B_4$ in which every curve associated to a given stratum has the same configuration of singularities and we compute the dimension of each stratum.}

\heading  Table of Contents \endheading


\halign{
#\hfil&\quad#\hfil&\quad#\hfil&\quad#\hfil&\quad#\hfil\cr
0.&Introduction, terminology, and preliminary results.\cr
1.&The General Lemma.\cr
2.&The Triple Lemma.\cr
3.&The BiProj Lemma.\cr 
4.&Singularities of multiplicity equal to degree divided by two.\cr
5.&The space of true triples of forms of degree $d$: the base point free locus, the\cr
 &birational locus, and the generic   Hilbert-Burch matrix.\cr
6.& Decomposition of the space of true triples.\cr
7.&The Jacobian matrix and the ramification locus.\cr
8.&The conductor and the branches of a rational plane curve.\cr
9.&Rational plane quartics: a stratification   and the  correspondence between\cr&the Hilbert-Burch matrices and the configuration of singularities.\cr
\cr
}

\SectionNumber=0\tNumber=1
\heading Section \number\SectionNumber. \quad Introduction, terminology, and preliminary results.
\endheading 

\bigskip

\heading Subsection \number\SectionNumber.A \quad Introduction.
\endheading
\bigskip
Consider a rational projective curve $\Cal C$ of degree $d$ over an algebraically closed field $\pmb k$. This curve, by definition,  is the closure of the image of a rational map 
$\Psi\: \xymatrix{\Bbb P^1\ar@{-->}[r]&\Bbb P^{n-1}}$. One can arrange the parameterization so that $\Psi$ has no base points and is birational onto its image. Thus, there are $n$ homogeneous forms $g_1,\dots,g_n$ of degree $d$ in $B=\pmb k[x,y]$ so that $\Psi(q)=[g_1(q):\cdots:g_n(q)]$ for all $q\in \Bbb P^1$ and $$\Cal C=\{\Psi(q)\mid q\in \Bbb P^1\}.$$ One traditional approach to the study of $\Cal C$ involves finding the polynomial relations $F\in\pmb k[T_1,\dots,T_n]$ on $g_1,\dots,g_n$, with $F(g_1,\dots,g_n)=0$. If one follows this line of investigation, then the defining equations of $\Cal C$ form the ideal
$$I(\Cal C)=\{F\in\pmb k[T_1,\dots,T_n]\mid F(g_1,\dots,g_n)=0\}.$$ Our approach instead is to look at the linear relations on $g_1,\dots,g_n$ with coefficients in $B$. We study the syzygy module
$$\left\{\left.\bmatrix c_1\\\vdots\\c_n\endbmatrix \in B^n\right \vert \sum_{i=1}^c g_ic_i=0\right\}.$$ In particular, we focus on a Hilbert-Burch matrix $\varphi$  so that 
$$0\to \bigoplus_{i=1}^{n-1} B(-d-d_i)@> \varphi >> B(-d)^n @> [g_1,\dots,g_n] >> B\tag\tnum{*8}$$is an exact sequence. It is plausible to learn about $\Cal C$ from $\varphi$ because, for example, when $\Cal C$ is a plane curve (that is, $n=3$), then the ideal $I(\Cal C)$ is generated by the resultant of the two polynomials $\sum_{i=1}^3 T_i\varphi_{i,j}$, where $\varphi=(\varphi_{i,j})$ and the two listed polynomials ($1\le j\le 2$) are viewed as homogeneous polynomials of degree $d_1$ and $d_2$ in the variables $x$ and $y$ with coefficients being linear forms in $k[T_1,T_2,T_3]$. 

Section 1 is concerned with the General Lemma. We use  the generalized row ideals of $\varphi$ to identify the singular points on $\Cal C$, their multiplicities, and the number of branches at each singular point. Indeed, the parameterization $\Psi$ automatically parameterizes the branches of $\Cal C$; see Observation \tref{.O34.2}. The General Lemma also uses the generalized row ideals of $\varphi$ to identify   the multiplicity of each branch of $\Cal C$. 

Let $\Cal C$ be a parameterized planar curve and let $d_1\le d_2$, as shown in (\tref{*8}), with $n=3$, be the shifts for a homogeneous Hilbert-Burch matrix for a parameterization of $\Cal C$.
One consequence of the General Lemma  is the observation that if $p$ is a   singularity on $\Cal C$, then  the multiplicity of $p$ is either equal to $d_2$ or less than or equal to $d_1$; furthermore, there is a one-to-one correspondence between the generalized zeros of $\varphi$ and the singularities on $\Cal C$ of multiplicity $d_1$ and $d_2$. (This observation was already known by the Geometric Modeling Community \cite{\rref{SCG}, \rref{CWL}}, where these singularities are called axial singularities. Our contribution in this context is a complete description of all singularities on $\Cal C$ (not just axial singularities) in terms of information which may be read from $\varphi$. We study space curves as well as planar curves and we calculate the information about  the branches and the multiplicity.)

Let $p$ be a singular point on the parameterized planar curve  $\Cal C$ which corresponds to a
generalized zero of $\varphi$.
The Triple Lemma in Section 2 concerns the blow-up $\Cal C'$ of $\Cal C$ centered at   $p$.   Theorem \tref{TPL1} gives a matrix $\varphi'$ whose
 maximal minors parameterize the closure, in $\Bbb P^2$, of the
blow-up at $p$ of $\Cal C$ in a neighborhood of $p$. We are able to apply the General Lemma  to $\varphi'$ in order to learn about the singularities of $\Cal C$ in the first neighborhood of $p$. If $\Cal C$ has even degree $d=2c$ and the multiplicity of $\Cal C$ at $p$ is equal to $c$, then we are able to apply Theorem \tref{TPL1} again and learn about the singularities of $\Cal  C$ in the second neighborhood (that is, after two blow-ups) of $p$. 
In section 2 we also prove that if $p$ is a point on $\Cal C$ and $q$ is
a singular point of multiplicity $c$ infinitely near $p$, then the
multiplicity of $p$ is also equal to $c$. It follows that if there is a
multiplicity $c$ singularity on or infinitely near $\Cal C$, then every
entry in a homogeneous Hilbert-Burch matrix for a parameterization of
$\Cal C$ is a homogeneous form of degree $c$. 

The rest of the paper, with the exception of Section 7, is  about rational plane curves of even degree $d=2c$. We classify curves according to the configuration of multiplicity $c$ singularities on or infinitely near $\Cal C$. There are $7$ possible configurations of such singularities. In Section 4, we classify the Hilbert-Burch matrix which corresponds to each configuration. For example, one assertion of Theorem \tref{d=2c} is as follows. If $\Cal C$ is a rational plane curve of  degree $d=2c$ and there is a singularity  on $\Cal C$ of multiplicity $c$ such that after one blow-up the singularity still has multiplicity $c$ and after a second blow-up the singularity still has multiplicity $c$, then there exists a linear automorphism $\Lambda$ of $\Bbb P^2$ and linearly independent forms $Q_1,Q_2,Q_3$ in $B_c$ such that $\Lambda \Cal C$ is parameterized by the maximal order minors of $$\bmatrix Q_1&Q_2\\Q_3&Q_1\\0&Q_3\endbmatrix.$$ The multiplicity $c$ singularity on $\Lambda \Cal C$ occurs at the point $[0:0:1]$. The repeated entry $Q_3$ ensures that there is a multiplicity $c$ singularity in the first neighborhood of $p$. The repeated entry $Q_1$ ensures that there is a multiplicity $c$ entry in the second neighborhood of $c$. 

Section three addresses the following question. 
If one has an arbitrary Hilbert-Burch matrix (or, equivalently, an arbitrary parameterization of a curve), how can one determine what the Hilbert-Burch matrix will be once enough row and column operations have been applied to transform it into the pretty form which is promised by Theorem \tref{d=2c}? 
 The other motivation for Section 3 is our desire   to separate the set of parameterizations into strata where each 
parameterization in a given stratum gives rise to a curve with a predetermined   configuration of singularities, and the closure of any given stratum is the union of  all strata less than or equal to the given stratum. This stratification is carried out in Section 6. The results of Section 3 are used to determine the appropriate stratum for each parameterization.  

At any rate, the study of multiplicity $c$ singularities on, or infinitely near, a fixed rational plane curve $\Cal C$ of degree $2c$ is equivalent to the study of the scheme of generalized zeros of the fixed balanced Hilbert-Burch matrix $\varphi$ for a parameterization of $\Cal C$. We prove that the projection maps
$$\eightpoint  \xymatrix{ &\{(p,q)\in \Bbb P^2\times \Bbb P^1\mid p\varphi q^{\text{\rm T}}=0\}\ar[ld]_{\pi_1}^{\cong}\ar[rd]_{\cong}^{\pi_2}&\\
\left\{p\in \Bbb P^2\left\vert {\matrix\format\l\\\text{$p$ is a}\\\text{multiplicity $c$}\\\text{singularity}\\\text{of $\Cal C$}\endmatrix}\right.\right \}&&\left\{q\in \Bbb P^1\left\vert{\matrix\format\l\\ \exists p\in \Bbb P^2\\\text{with}\\p\varphi q^{\text{\rm T}}=0\endmatrix}\right.\right\}}$$ are isomorphisms and we identify the defining equations for each of the three schemes. The scheme on the lower left is the scheme of interest. The scheme on the lower right does not appear to have any intrinsic significance; however, it is a subscheme of $\Bbb P^1$, so it is easy to make computations concerning this scheme. In particular, our ultimate criteria in Theorem \tref{TBA} for determining the configuration of multiplicity $c$ singularities on or infinitely near $\Cal C$ amount to looking at the factorization of the greatest common divisor of the $3\times 3$  minors of a matrix $A$ whose entries are linear forms in $\pmb k[u_1,u_2]$. The matrix $A$ is created quickly from $\varphi$ by way of extracting the variables $x$ and $y$ from the critical equation $[T_1,T_2,T_3] \varphi\left [\smallmatrix u_1\\u_2\endsmallmatrix \right ]=0$.   We write $[T_1,T_2,T_3] \varphi=[y^c,\dots, x^c]C$, where $C$ is a matrix of linear forms in $\pmb k[T_1,T_2,T_3]$ and we write $C\left [\smallmatrix u_1\\u_2\endsmallmatrix\right ]=A\left [\smallmatrix T_1\\T_2\\T_3\endsmallmatrix\right ]$. The ideal generated by the entries of the product matrix on either side of the most recent equation defines the subscheme of $\Bbb P^2\times \Bbb P^1$ on the top of the above picture.  Of course, Theorem \tref{TBA} is phrased in terms of multiplicity and therefore we are able to apply the Generic Freeness Lemma, by way of Theorem \tref{L37.7}, to families of parameterizations. 

 As previously mentioned, Section 4 contains a classification of Hilbert-Burch matrices and how these matrices correspond to a particular configuration of multiplicity $c$ singularities. This classification is the starting point for our study of families of curves. We use it to prove that our strata are irreducible and to calculate their dimensions. This classification is also an appropriate place to begin a study of the other singularities on, or infinitely near, the curves $\Cal C$ those with multiplicity not equal to $c$. We carry out this study in the classical situation of quartics in Section 9. We have begun a study of sextics beginning from exactly this point. 

It is possible that Theorem \tref{orbitz} could be of independent interest. There is no explicit mention of curves in the statement of this result. Let $d=2c$ be an even integer and $$\operatorname{BalH}_d=\left\{ \varphi\left\vert \matrix\format\l\\\text{$\varphi$ is a $3\times 2$ matrix with entries in $B_c$ such that $\operatorname{ht} I_2(\varphi)=2$} \endmatrix\right.\right\},$$where, as always, $B=\pmb k[x,y]$. The group $G=\operatorname{GL}_3(\pmb k)\times \operatorname{GL}_2(\pmb k)$ acts on $\operatorname{BalH}_d$ by way of $(\chi,\xi)\cdot \varphi=\chi\varphi\xi^{-1}$. In Theorem \tref{orbitz} we decompose $\operatorname{BalH}_d$ into a disjoint union of 11 orbits. Each orbit has the form $G\cdot M$, where $M$ is the closed irreducible subspace of affine space defined by the maximal order minors of a generic matrix.

In Section 5 we introduce the parameter space $\Bbb A_d=B_d\times B_d\times B_d$. Each element $\pmb g$ of $\Bbb A_d$ is an ordered triple of $d$-forms from $B=\pmb k[x,y]$ (so, $\Bbb A_d$ is an affine space of dimension $3d+3$), and each $\pmb g\in \Bbb A_d$ induces a rational map $\Psi_{\pmb g}\: \xymatrix{\Bbb P^1\ar@{-->}[r]&\Bbb P^{2}}$. We define 
$$\split \Bbb T_d&{}=\left \{\pmb g\in \Bbb A_d\left\vert \text{$\Psi_{\pmb g}$ is  birational onto its image  without base points}\right.\right\}\ \text{and}\\
\Bbb B_d&{}=\left\{\pmb g\in \Bbb T_d\left\vert \matrix\format\l\\\text{every entry in a homogeneous  Hilbert-Burch}\\\text{matrix for $\pmb d_1(\pmb g)$ has degree $d/2$}\endmatrix\right.\right\}.\endsplit $$(If $\pmb g$ is the ordered triple $(g_1,g_2,g_3)$ of $\Bbb A^d$, then $\pmb d_1(\pmb g)$ is the row vector $[g_1,g_2,g_3]$.)
We call $\Bbb T_d$ the space of {\bf true} triples of forms of degree $d$ and $\Bbb B_d$ the space of  {\bf balanced} true triples of  forms of degree $d$.

In practice we are only interested in the subset $\Bbb T_d$ of $\Bbb A_d$. Every element of $\Bbb A_d$ which is not in $\Bbb T_d$ corresponds to an unsuitable parameterization of a curve. As we have already seen, the results of Section 2 show that if there is a multiplicity $c$ singularity on or infinitely near a point $p$ on a curve $\Cal C$ of degree $d=2c$, then the parameterization of $\Cal C$ is an element of $\Bbb B_d$. In Section 5 we prove that $\Bbb B_d\subseteq \Bbb T_d$ are open subsets of $\Bbb A_d$. We also identify an open cover $\cup \Bbb B_d^{(i)}$ of $\Bbb B_d$ and for each $\Bbb B_d^{(i)}$ in this open cover, we identify a generic Hilbert-Burch matrix which specializes to give a Hilbert-Burch matrix for $\pmb g$ for each $\pmb g\in \Bbb B_d^{(i)}$. 

Corollary \tref{UPR} is    another result which   could be of independent interest.  We  identify a universal projective resolution $\Bbb U\Bbb P\Bbb R_{\Bbb Z}$ for the graded Betti numbers
$$0\to B(-3c)^2 \to B(-2c)^3\to B.\tag\tnum{0.*}$$The resolution $\Bbb U\Bbb P\Bbb R_{\Bbb Z}$ is built over the ring $(\pmb R_{\Bbb Z})_{w_{\Bbb Z}}[x,y]$, where $\pmb R_{\Bbb Z}$ is a polynomial ring over $\Bbb Z$,   $w_{\Bbb Z}$ is a non-zero homogeneous element in $\pmb R_{\Bbb Z}$, and $(\pmb R_{\Bbb Z})_{w_{\Bbb Z}}$ is the localization of $\pmb R_{\Bbb Z}$ at the multiplicatively closed set $\{1, w_{\Bbb Z},w_{\Bbb Z}^2,\dots\}$. If  (\tref{0.*}) is a minimal homogeneous resolution of $B/I$ over $B=\pmb k[x,y]$, then there exists a homomorphism $(\pmb R_{\Bbb Z})_{w_{\Bbb Z}}\to \pmb k$ so that $\Bbb U\Bbb P\Bbb R_{\Bbb Z}\otimes_{(\pmb R_{\Bbb Z})_{w_{\Bbb Z}}}\pmb k$ is a minimal homogeneous resolution of $B/I$ over $B$. 

In Section 6 we decompose the space $\Bbb B_d$ of balanced  triples into strata. Each stratum consists of those triples $\pmb g$ in $\Bbb B_d$ for which the corresponding curve $\Cal C_{\pmb g}$ exhibits one particular configuration of multiplicity $c$ singularities. If $\Cal C$ is a parameterized plane curve of degree $d=2c$, then there are $7$ possible configurations of multiplicity $c$ singularities on or infinitely near $\Cal C$:
$$\operatorname{CP}=\{\emptyset,\ \{c\},\ \{c,c\},\ \{c,c,c\},\ \{c:c\},\ \{c:c,c\},\   \{c:c:c\}\},$$ where a colon indicates an infinitely near singularity, a comma indicates a different singularity on the curve, and $\emptyset$ indicates that there are no singularities of multiplicity $c$ on, or infinitely near, $\Cal C$. For example, if $\pmb g\in S_{c:c:c}$, then there is exactly one multiplicity $c$ singularity $p$ on the corresponding curve $\Cal C_{\pmb g}$ and there are two multiplicity $c$ singularities infinitely near to $p$: one in the first neighborhood of $p$ and one in the second neighborhood of $p$; and if $\pmb g\in S_{c:c,c}$, then there are exactly two multiplicity $c$ singularities  on   $\Cal C_{\pmb g}$ and exactly one of these has a multiplicity $c$ singularity infinitely near it. We decompose $\Bbb B_d$ as  $$\xymatrix{&&S_{c:c}\ar[dr]\\
S_{c:c:c} \ar[r]  &S_{c:c,c} \ar[ru] \ar[rd]& &S_{c,c} \ar[r] &S_{c}\ar[r] &S_{\emptyset}\cap \Bbb B_d,\\&&S_{c,c,c}\ar[ur]}$$where $S_{\#'}\to S_{\#}$ means that $S_{\#'}$ is contained in the closure of $S_{\#}$. Think of $\operatorname{CP}$ as a poset with $\#'\le \#$ whenever we have drawn $S_{\#'}\to S_{\#}$. We prove that each set $\Bbb T_{\#}=\cup_{\#'\le \#}S_{\#'}$ is a closed irreducible subset of $\Bbb B_d$ and we compute the dimension of each stratum $S_{\#}$. The proof that $T_{\#}$ is closed uses Corollary \tref{CorND} to associate numerical data to each $S_{\#}$, Theorem \tref{L37.7} to convert this numerical data into closed conditions on families of curves as calculated by way of the Hilbert-Burch matrix, and Theorem \tref{T1} to create a generic Hilbert-Burch matrix from the coefficients of the parameterization. The calculation of the dimension of $T_{\#}$ is based on the decomposition of balanced Hilbert-Burch matrices which takes place in Section 4.

In Section 7 we return to the hypotheses of Section 1; that is, $\Cal C$ is a parameterized space curve of arbitrary degree and we learn about infinitely near singularities of arbitrary multiplicity. 
Remark \tref{brn} provides a method of parameterizing the branches of a parameterized curve. Theorem \tref{,T4.1} shows that the Jacobian matrix associated to the parameterization identifies the non-smooth branches of the curve as well as the multiplicity of each branch. The starting point for this line of reasoning  is the result  that if $D$ is an algebra which is essentially of finite type over the ring $C$, then the ramification locus of $D$ over $C$ is equal to the support of the module of K\"ahler differentials $\Omega_{D/C}$. It is important to compare Theorems  \tref{Cor1} and  \tref{Cor2'}. 
Theorem \tref{Cor1} is a  local result. Once one knows the singularities $\{p_i\}$ on a parameterized curve $\Cal C$, then this   result shows how to read $m_{\Cal C,p_i}$ and $s_{\Cal C,p_i}$, for each  $p_i$,  from the Hilbert-Burch matrix of the parameterization. Theorem \tref{Cor2'} is a global result. It describes, in terms of the parameterization,  all of the points $p$ on $\Cal C$ and all of the branches of $\Cal C$ at $p$ for which the multiplicity  is at least two.  In contrast to Theorem \tref{Cor1}, one may apply Theorem \tref{Cor2'} before one knows the singularities on $\Cal C$. 

  Sections 1, 2,  3, and 7 are largely concerned with the multiplicity and the number of branches at each singularity on a curve.    Section 8 deals with the singularity degree $\delta$. 
Let $\pmb g=(g_1,g_2,g_3)$ be an element of $\Bbb T_d$ and $\Cal C_{\pmb g}$ be the corresponding parameterized plane curve.
We produce a polynomial $c_{\pmb g}$ whose factorization into linear factors gives the value of the invariant $\delta$ at each singular point of   $\Cal C_{\pmb g}$. The polynomial $c_{\pmb g}$ is obtained  in a polynomial manner from the coefficients of the entries of a Hilbert-Burch matrix for $\pmb g$.
We use these ideas to produce closed sets in $\Bbb B_d$ which separate various configurations of singularities. To create $c_{\pmb g}$ we start with the coordinate ring $A_{\pmb g}=\pmb k[g_1,g_2,g_3]\subseteq B=\pmb k[x,y]$. If $\pmb V$ is the $d^{\text{th}}$ Veronese subring of $B$ and $\frak c_{\pmb g}$ is the conductor $A_{\pmb g}\:\! \pmb V$, then $c_{\pmb g}$ generates the saturation of the extension of $\frak c_{\pmb g}$ to $B$.

In Section 9 we apply our results to rational plane quartics. It has been known for well over one hundred years (see, for example, Basset \cite{\rref{Bas}}, Hilton \cite{\rref{Hil}}, Namba \cite{\rref{Nam}}, or Wall \cite{\rref{Wa04}}) that there are 13 possible configurations of singularities on a rational plane quartic:

$$\matrix \format \l&\quad\ \l&\quad\ \l&\quad\ \l\\
(3:1,1,1)&2:2:2:1,1&(2:2:1,1),(2:1,1)&(2:1,1)^3\\
(3:1,1)&2:2:2:1&(2:2:1,1),(2:1)&(2:1,1)^2,(2:1)\\
3:1&&(2:2:1),(2:1,1)&(2:1,1),(2:1)^2\\
&&(2:1)^3&(2:2:1),(2:1).\endmatrix\tag\tnum{intro}$$
(The notation is fully defined around (\tref{osc}) and the full set of names is recorded in Section 9.) These 13 singularity configurations have been named and stratified many times already; see, for example,  \cite{\rref{BG}} and \cite{\rref{Wa}}. In Section 9, we show what our approach produces in this context. The techniques from Section 6 separate 
$$3:*\qquad 2:2:2:*\qquad (2:2:*),(2:*),\qquad\text{and}\qquad(2:*)^3.$$The finer techniques of Section 7 separate the parameterizations according to the value of the invariant
$$\sum\limits_{p\in \operatorname{Sing} \Cal C}m_p-s_p,$$ where $p$ roams over all of the singular points on $\Cal C$, $m_p$ is the multiplicity of $\Cal C$ at $p$ and $s_p$ is the number of branches of $\Cal C$ at $p$. The conductor technique of Section 8 distinguishes the configuration of singularities 
  $(2:2:1), (2:1,1)$ from $(2:2:1,1), (2:1)$. In Corollary \tref{MC9}  we exhibit a stratification of $\Bbb B_4$ in which every curve associated to a given stratum has the same configuration of singularities.   We compute the dimension of each stratum. 

\bigskip We obtain a rational plane curve $d$ from each true triple of $d$-forms $\pmb g$ in $\Bbb T_d$. The space $\Bbb T_d$ is large and provides us with plenty of room for maneuvering. Nonetheless, it is important that we exhibit the precise relationship between the space $\Bbb T_d$ and the space $\Bbb R\Bbb P\Bbb C^r_d$ of {\bf r}ight equivalence classes of {\bf R}ational {\bf P}lane {\bf C}urves of degree $d$. Recall that the plane curves $\Cal C$ and $\Cal C'$ are {\it right equivalent} if there exists a linear automorphism $\Lambda:\Bbb P^2\to \Bbb P^2$ with $\Lambda \Cal C=\Cal C'$. We have surjective morphisms 
$$\Bbb T_d@> \alpha_1 >> \Bbb P\Bbb T_d@> \alpha_2>> \Bbb R\Bbb P\Bbb C_d@> \alpha_3 >> \Bbb R\Bbb P\Bbb C^r_d, \tag\tnum{sm}$$
 where $\Bbb P\Bbb T_d$ is the space of degree $d$ base point free morphisms from $\Bbb P^1$ to $\Bbb P^2$ which are birational onto their image and $\Bbb R\Bbb P\Bbb C_d$ is the space of Rational Plane Curves of degree $d$. If $\Psi\:\Bbb P^1\to \Bbb P^2$ is an element of $\Bbb P\Bbb T_d$, then the fiber of $\alpha_1$ over $\Psi$ is isomorphic to $\operatorname{GL}_1(\pmb k)$. (Indeed, if $\pmb g=(g_1,g_2,g_3)$ is in $\Bbb T_d$ and $u$ is a non-zero element of $\pmb k$, then $\pmb g$ and $u\pmb g=(ug_1,ug_2,ug_3)$ give rise to the same morphism $\Bbb P^1\to \Bbb P^2$.) If $\Cal C$ is an element of $\Bbb R\Bbb P\Bbb C_d$, then the fiber of $\alpha_2$ over $\Cal C$ is isomorphic to $\operatorname{SL}_2(\pmb k)$. Indeed, if $\Psi\:\Bbb P^1\to \Bbb P^2$ is a morphism and $\xi\:\Bbb P^1\to \Bbb P^1$ is a linear automorphism, then $\Psi$ and $\Psi\circ \xi$ give rise to same curve. The definition of right equivalence yields that the fiber of $\alpha_3$ over each equivalence class of  $\Bbb R\Bbb P\Bbb C^r_d$ is isomorphic to $\operatorname{SL}_3(\pmb k)$. All together the fiber of $\alpha_3\circ \alpha_2\circ \alpha_1$ over any equivalence class of  $\Bbb R\Bbb P\Bbb C^r_d$ is isomorphic to $\operatorname{GL}_1(\pmb k)\times \operatorname{SL}_2(\pmb k)\times \operatorname{SL}_3(\pmb k)$ and has dimension $12$.

\bigskip

\heading Subsection \number\SectionNumber.B \quad Terminology.
\endheading
\bigskip
All of our work takes place over a field $\pmb k$. We write $\pmb k^*$ for the multiplicative group $\pmb k\setminus\{0\}$. Often (but not always) the field is algebraically closed. Sometimes the field $\pmb k$ has characteristic zero. Projective space $\Bbb P^n$ means projective space $\Bbb P^n_{\pmb k}$ over $\pmb k$. Let $\Cal O_{\Cal C,p}$ represent the local ring of the curve $\Cal C$ at the rational point $p$.  The {\it multiplicity} of $\Cal C$ at $p$, denoted $m_p$ (or $m_{\Cal C,p}$ if there is any ambiguity about what $\Cal C$ is), is the multiplicity of the local ring $\Cal O_{\Cal C,p}$. (Recall that the multiplicity of the $d$-dimensional Notherian local ring $(A,\frak m)$ is $e(A)=\lim\limits_{n\to \infty}\frac {d!\lambda_A(A/\frak m^n)}{n^d}$, where $\lambda_A(M)$ is the length of the $A$-module $M$.) Let $\widehat{\Cal O_{\Cal C,p}}$ be 
the completion  of $\Cal O_{\Cal C,p}$ with respect to its maximal ideal $\frak m_{\Cal C,p}$.
 Each minimal prime ideal of $\widehat{\Cal O_{\Cal C,p}}$ is called a {\it branch} of $\Cal C$ at $p$. If $\Cal J$ is a minimal prime of $\widehat {\Cal O_{\Cal C,p}}$, then the multiplicity of the local ring $\widehat {\Cal O_{\Cal C,p}}/\Cal J$  is called the {\it multiplicity of  the branch $\Cal J$.} The {\it number of branches} of $\Cal C$ at $p$ is denoted $s_p$ (or $s_{\Cal C,p}$).
The {\it singularity degree} of $\Cal C$ at the point $p$, denoted $\delta_p$ (or $\delta_{\Cal C,p}$), is  
$\lambda\left (\overline{\Cal O_{\Cal C,p}}/\Cal O_{\Cal C,p}\right )$, where $\overline{\Cal O_{\Cal C,p}}$ is the integral closure of 
$\Cal O_{\Cal C,p}$ and $\lambda$ means $\lambda_{\Cal O_{\Cal C,p}}$. The invariant $\delta_p$ may also be realized as
$\delta_p=\sum_q\binom{m_q}{2}$, where $q$ varies over $p$ together with ``all singularities infinitely near to $p$''. This notion is described below.

Let $p$ be a point on the curve $\Cal C$ and let $\Cal C'@>\sigma >> \Cal C$ be the blow-up of $\Cal C$ at $p$. The
points in   the {\it first neighborhood} of $p$ are the
points on $\Cal C'$ which map to $p$. There are finitely many such points. The local rings of the points in
the first neighborhood of $p$ are exactly the local rings of $\operatorname{Proj}(\Cal
O[\frak m t])$, where $(\Cal O,\frak m)=(\Cal O_{\Cal C,p},\frak m_p)$. As
$\dim \Cal O=1$, and $\pmb k$ is infinite, there exists $f\in \frak m$ so
that $\frak m^{r+1}=f\frak m^r$ for some $r\ge 0$; hence,
$\operatorname{Proj}(\Cal O[\frak m t])=\operatorname{Spec}(\Cal O[\frac{\frak m}f])$ and therefore the
local rings of the points in the first neighborhood of $p$ are exactly
the localizations of $\Cal O[\frac{\frak m}f]$ at one of its finitely many
maximal ideals. Apply this process to each point in the first
neighborhood of $p$ to obtain the points in the {\it second neighborhood} of $p$. Repeat the process until all of the points in the $r^{\text{th}}$ neighborhood of $p$ are smooth points, for some $r$. The union of the singular points in the $i^{\text{th}}$ neighborhood of $p$, for $1\le i\le r$, is called the set of {\it singularities of $\Cal C$ infinitely near to $p$}. All of the local rings for the singularities of $\Cal C$ infinitely near to $p$ occur between $\Cal O$ and its integral closure $\overline{\Cal O}$, which is finitely generated as a $\Cal O$-module.

Notice that we reserve the phrase ``singularity of $\Cal C$'' to refer to points on $\Cal C$ which are singular points and we reserve the phrase ``infinitely near singularity of $\Cal C$'' to refer to singular points which are not on $\Cal C$, but  which  appear after a  finite sequence of blow-ups has occurred.

When we write that the {\it multiplicity sequence} for the oscnode $q_0$ on the curve $\Cal C_0$ is $2:2:2:1,1$, we mean that there is a sequence of blow-ups:
$$ \Cal C_0@<\sigma_1 <<\Cal C_1@<\sigma_2 <<\Cal C_{2}@<\sigma_{3} <<\Cal C_{3}\tag\tnum{osc}$$ and a sequence of points 
$q_0$ on $ \Cal C_0$, $q_1$ on $ \Cal C_1$, $q_2$ on $ \Cal C_2$, and $q_3\neq q_3'$ on $ \Cal C_3$, such that $\sigma_i$ is the blowup of $\Cal C_{i-1}$ centered at $q_{i-1}$ for $1\le i\le 3$, 
$$\sigma_i^{-1}(q_{i-1})=\cases \phantom{,q_3}q_i&\text{if $1\le i\le 2$}\\ \{q_3,q_3'\}&\text{if $i=3$},\endcases$$
$$m_{\Cal C_0,q_0}= m_{\Cal C_1,q_1}=m_{\Cal C_2,q_2}=2\quad\text{and}\quad m_{\Cal C_3,q_3}= m_{\Cal C_3,q_3'}=1.$$So, in particular, 
$$\delta_{q_0}=\sum\limits_{i=0}^2\binom{m_{q_i}}2=3\quad\text{and}\quad s_{q_0}=2$$ because there are two smooth points on $\Cal C_3$ which lie over $q_0$. When we write that the configuration of singularities of the curve $\Cal C$ is described by $(2:2:1,1),(2:1)$, we mean that there are exactly 2 singularities  on $\Cal C$: one is a tacnode and the other is an ordinary cusp. There is exactly one singular point in the first neighborhood of the tacnode. This point has multiplicity 2. There are exactly two points in the second neighborhood of the tacnode; both of these points are smooth points. There is exactly one point in the first neighborhood of the cusp and this point is a smooth point. 

All rings in this paper are commutative and Noetherian. The {\it grade} of a proper ideal $I$ in a ring $R$ is the length of the
longest regular sequence on $R$ in $I$. The expression ``let $(A,\frak m,\pmb k)$ be a local ring'' means that $A$ is a local ring with unique maximal ideal $\frak m$ and $\pmb k$ is the residue class field $A/\frak m$. 
 If $(A,\frak m)$ is a local ring and $z$ is an indeterminate over $A$, then $$\text{$A(z)$ means the ring $A[z]_{\frak mA[z]}$}.\tag\tnum{Aofz}$$ 
In other words, $A(z)$ is the localization of the polynomial ring $A[z]$ at the prime ideal $\frak mA[z]$.

If $M$ is a matrix, then $I_r(M)$ is the ideal generated by the $r\times r$ minors of $M$ and $M^{\text{\rm T}}$ is the transpose of $M$. 
If $L$ is a graded module over a graded ring $B=B_0\oplus B_1\oplus \dots$, then the $B_0$-module $L_d$ is the homogeneous component of $L$ of degree $d$. 
We write $\mu(L)$ to denote the minimal number of generators of the module $L$. This concept makes sense whenever $L$ is a finitely generated module over a local ring or $L$ is a finitely generated graded module over a graded ring with a unique maximal homogeneous ideal. 
If $\pmb T=(T_{i,j})$ is a   matrix of indeterminates and $R$ is a ring,  then $R[\pmb T]$ is the polynomial ring $R[\{T_{i,j}\}]$.  
If $I$ is an ideal of a ring $R$, then we write $\dim I$ to represent the Krull dimension of the quotient ring $R/I$. In particular, if $I$ is a homogeneous ideal in a graded ring $R$ with unique maximal homogeneous ideal $\frak m$, then $I$ is a zero-dimensional ideal if and only if $I$ is an $\frak m$-primary ideal.

If $v_1,\dots,v_r$ are elements in a vector space $V$, then we write ${<}v_1,\dots,v_r{>}$ to mean the subspace of $V$ spanned by $v_1,\dots,v_r$.

If $R$ is a ring, then we write $\operatorname{Quot}(R)$ for the {\it total quotient ring} of $R$; that is, $\operatorname{Quot}(R)=U^{-1}R$, where $U$ is the set of non zerodivisors on $R$. 
If $R$ is a domain, then the total quotient ring of $R$ is usually called the {\it quotient field} of $R$. The {\it integral closure} $\overline{R}$  of a ring $R$ is    the integral closure of 
$R$ in the $\operatorname{Quot} (R)$.

If $g_1,\dots,g_n$ are elements of $B=\pmb k[x,y]$ which generate an ideal of height two, then the Hilbert-Burch Theorem asserts that the relations on the row vector $[g_1,\dots,g_n]$ fit into an exact sequence
$$0\to B^{n-1}@> \varphi>> B^n@> [g_1,\dots,g_n]>> B,$$ and that there is a unit $u$ in $B$ so that $g_i$ is equal to $(-1)^{i+1}u$ times the determinant of $\varphi$ with row $i$ removed. We call $\varphi$ a {\it Hilbert-Burch matrix}   for $[g_1,\dots,g_n]$.
If $I$ and $J$ are ideals of a ring $R$, then the {\it saturation of $I$ with respect to $J$} is $I\:\! J^{\infty}=\bigcup\limits_{i=1}^{\infty}(I\:\! J^i)$. We write $\gcd$ to mean {\it greatest common divisor}. If $I$ is a homogeneous ideal in $\pmb k[x,y]$, then 
we denote the $\gcd$ of a generating set of $I$ by $\gcd I$; notice that this polynomial generates the saturation   $I\:\! (x,y)^{\infty}$.

\remark{Remark \tnum{gri}} The concept of  ``generalized row ideals'' (as well as other related concepts) appears widely in the literature; see, for example, \cite{\rref{G},\rref{EHU},\rref{EU}}.
  Let $M$ be a matrix with entries in a $\pmb k$-algebra, where $\pmb k$ is a field. A {\it generalized row} of $M$ is the product $pM$, a  {\it generalized column}  of $M$ is the product $Mq^{\text{\rm T}}$, and a {\it generalized entry} of $M$ is the product $pMq^{\text{\rm T}}$, where $p$ and $q$ are non-zero row vectors with entries from $\pmb k$. A {\it generalized zero} of $M$ is a generalized entry of $M$ which is zero. 
If $pM$ is a generalized row of $M$, then the ideal $I_1(pM)$ is called a {\it generalized row ideal} of $M$.
If $Mq^{\text{\rm T}}$ is a generalized column of $M$, then the ideal $I_1(Mq^{\text{\rm T}})$ is called a {\it generalized column ideal} of $M$. \endremark

\remark{Remark \tnum{Unam}}Often,  we  think of a point $p=[a_1\:\cdots \:a_n]\in \Bbb P^{n-1}$ as a row vector $[a_1,\dots,a_n]$, where $[a_1,\dots,a_n]$ is any one of the row vectors that represent the point $p$ in $\Bbb P^{n-1}$. For example, if $M$ is a matrix with $n$ rows, then the point $p$ in $\Bbb P^{n-1}$ gives rise to a unique generalized row ideal $I_1([a_1,\dots,a_n]M)$ of $M$ since the ideals 
$I_1([a_1,\dots,a_n]M)$ and $I_1([\lambda a_1,\dots,\lambda a_n]M)$ are equal for all non-zero constants $\lambda$ of $\pmb k$. It is convenient, and unambiguous, to write $I_1(pM)$, with $p\in \Bbb P^{n-1}$, instead of $I_1([a_1,\dots,a_n]M)$. \endremark

\remark{Remark \tnum{sigh}} If  $g_1,\dots,g_n$ are linearly independent homogeneous polynomials of the same  degree  in $\pmb k[x,y]$, then we consider  the rational map
$\Psi\: \xymatrix{\Bbb P^1\ar@{-->}[r]&\Bbb P^{n-1}}$ which sends  the point $q$ of $\Bbb P^1\setminus  V(g_1,\dots,g_n)$ to $[g_1(q):\dots:g_n(q)]\in \Bbb P^{n-1}$. The image of $\Psi$ is a curve $\Cal C$ in  $\Bbb P^{n-1}$. 
Let $\varphi$ be a matrix   of relations on the row vector $g=[g_1,\dots,g_n]$.
Row operations on $\varphi$
correspond to replacing $\varphi$ by $\varphi'=\chi \varphi$, where
  $\chi$ is an invertible $n\times n$ matrix with entries in $\pmb k$. Now $\varphi'$ is a matrix of relations on the row vector $g'=g\chi^{-1}$.  Let $\Lambda\: \Bbb P^{n-1}\to \Bbb P^{n-1}$ be the linear automorphism sending $p$ to $p\chi^{-1}$. When we replace $g$ by $g'$, the rational map $\Psi$ and the curve $\Cal C$ become $\Lambda\circ \Psi$ and  
$\Lambda(\Cal C)$, respectively. Column operations on $\varphi$ have no effect on  the row vector $g$, the rational map $\Psi$, or the curve $\Cal C$. 
\endremark

\bigskip

\heading Subsection \number\SectionNumber.C \quad Preliminary results.
\endheading
\bigskip
The following result is established in \cite{\rref{KPU-B}}; it interprets the birationality of a parameterization of a curve in terms of the generalized Hilbert-Burch matrix which corresponds to this parameterization. Recall that $R$ is called a {\it standard graded algebra} over the ring $R_0$ if $R$ is a graded ring $\bigoplus_{0\le i}R_i$, with  $R_1$ finitely generated as a module over $R_0$, and $R$   generated as an algebra over $R_0$ by $R_1$. If $R$ is a standard graded algebra over a field, then the multiplicity of $R$ is given by ${e(R)=\lim\limits_{n\to \infty}\frac {d!\lambda_R(R/\frak m^n)}{n^d}}$, where $d$ is the Krull dimension of $R$ and $\frak m$ is the maximal homogeneous ideal $\bigoplus_{1\le i}R_i$.
\proclaim{Theorem \tnum{2}}  Let  $B$ be the standard graded polynomial ring $\pmb k[x,y]$, with $\pmb k$ an infinite field,  
 $I$ be a height two  ideal of  $B$ generated by forms $g_1,\dots,g_n$
     of degree $d$, and $\varphi$   be a homogeneous Hilbert-Burch matrix for the row vector  $[g_1,\dots,g_n]$.
       If $C$ and $D$  are  the standard graded $\pmb k$-algebras $C=\pmb k[I_d]$ and  
$D=\pmb k[B_d]$, $r$ is the degree of the field extension $[\operatorname{Quot}(D):\operatorname{Quot}(C)]$, and $e$ is the multiplicity of  $C$,  then the following statements hold.
\roster
\item $re=d$.
\item The morphism $\Bbb P^1\to \Bbb P^{n-1}$, which is given by $q\mapsto  [g_1(q):\dots:g_n(q)]$, is birational onto its image if and only if $r=1$.
\item There exist 
forms $f_1$ and $f_2$  of degree $r$ in $B$ such that $C\subseteq \pmb k[f_1,f_2]$. 
 In particular, $I$ is extended from an ideal in
$\pmb k[f_1,f_2]$ in the sense that ${I=(I\cap \pmb k[f_1,f_2])B}$.    Furthermore,   $f_i=\gcd(I_1(p_i\varphi))$ 
  for general points ${p_1}$ and  ${p_2}$  in
$\Bbb P^n$.
\endroster\endproclaim
In the context of Theorem \tref{2} it is clear that some homogeneous Hilbert-Burch matrix for $[g_1,\dots,g_n]$ has entries in $\pmb k[f_1,f_2]$. We identify a sufficient condition which guarantees that every Hilbert-Burch matrix for $[g_1,\dots,g_n]$ has entries in $\pmb k[f_1,f_2]$.
\proclaim{Observation \tnum{2.1}} Retain the notation and hypotheses of Theorem {\rm\tref{2}}.  Assume that every entry of $\varphi$ has the same degree $c$. The following statements hold.
\roster \item Every entry of $\varphi$ is in the ring $\pmb k[f_1,f_2]$.\item If $c$ is prime, then the morphism of Theorem {\rm\tref{2}} is birational if and only if $\mu(I_1(\varphi))\ge 3$.\endroster
\endproclaim
\demo{Proof} (1) Let $\pmb k[z_1,z_2]$ be a new polynomial ring and $\theta\:\pmb k[z_1,z_2]\to \pmb k[f_1,f_2]$ be the $\pmb k$-algebra homomorphism which is defined by $\theta(z_i)=f_i$. Identify $g_i'$ in $\pmb k[z_1,z_2]$ with $\theta(g_i')=g_i$ for all $i$. It is clear that the ideal $(g_1',\dots,g_n')$ of $\pmb k[z_1,z_2]$ has height two. Let $\varphi'$, with entries in $\pmb k[z_1,z_2]$, be any homogeneous Hilbert-Burch matrix for $[g_1',\dots,g_n']$. It follows that $\theta(\varphi')$ is a Hilbert-Burch matrix for $[g_1,\dots,g_n]=\theta([g_1',\dots,g_n'])$. Thus, there exist homogeneous invertible matrices $\chi$ and $\xi$ with entries in $B$ so that $\chi\theta(\varphi')\xi=\varphi$. Degree considerations show that the entries of $\chi$ and $\xi$ are actually in $\pmb k$; and therefore, $\varphi=\chi\theta(\varphi')\xi=\theta(\chi\varphi'\xi)$. 

(2) The direction $(\Rightarrow)$ is obvious. We prove $(\Leftarrow)$. Assume that $\mu(I_1(\varphi))\ge 3$. Every entry of $\varphi$ is a homogeneous form of degree, say  $s$, in the ring $\pmb k[f_1,f_2]$ and the homogeneous forms $f_i$ have degree $r$. So, the prime $c$ is equal to $rs$. The hypothesis $\mu(I_1(\varphi))\ge 3$ shows that   $r\neq c$. Thus,    $r$ must equal $1$. \qed \enddemo

\bigpagebreak
\SectionNumber=1\tNumber=1
\heading Section \number\SectionNumber. \quad The General Lemma.
\endheading 

In this section we prove Lemma \tref{.gl}, which we call the General Lemma. Let  $\Cal C$ be   the curve  which is parameterized by  the homogeneous polynomials $g_1,\dots,g_n$  of the same degree and let  $p$  be the point $[a_1:\dots:a_n]$  in $\Bbb P^{n-1}$. 
If $\varphi$ is a Hilbert-Burch matrix for $[g_1,\dots,g_n]$, then 
the generalized row ideal  $I_1([a_1,\dots,a_n]\varphi)$  captures a significant amount of geometric information about the behavior of $\Cal C$ at $p$. Indeed, from this ideal, one may determine whether or not $p$ is on $\Cal C$ , the multiplicity $m_{\Cal C,p}$ of $\Cal C$ at $p$, the number of branches $s_{\Cal C,p}$ of $\Cal C$ at $p$, and the multiplicity of  each  branch.

The basic set-up is given in Data \tref{34.1} and explained in Remark \tref{34.2}. 
 Observation \tref{O-1-18} is a small calculation, given originally in \cite{\rref{EU}}, which shows how the generic row ideal $I_1(p\varphi)$ can be used to determine whether or not $p$ is on $\Cal C$.   After stating Lemma \tref{.gl}, we immediately describe the first round of consequences; further consequences are found throughout the paper.
Theorem \tref{Cor1} covers  essentially the same ground as Lemma \tref{.gl}; but it might be easier to apply because the statement of Theorem \tref{Cor1} is less ambitious.  Corollary \tref{ctgl}, which is a result due to Song, Chen, and Goldman \cite{\rref{SCG}}, may be deduced from Theorem \tref{Cor1}.

After giving the proof of Corollary \tref{ctgl}, we  lay the foundations for the proof of Lemma \tref{.gl}.   
Lemma \tref{.O33.10'} concerns the relationship between a one-dimensional local domain and its integral closure; in particular, this result is about the structure of the extension of the maximal ideal of the smaller ring.  Proposition \tref{!32.1} is the statement of the well-known   correspondence between the minimal prime ideals of the completion $\widehat{\Cal O_{\Cal C,p}}$ of the local ring $\Cal O_{\Cal C,p}$ and the maximal ideals of the integral closure $\overline{\Cal O_{\Cal C,p}}$ of $\Cal O_{\Cal C,p}$.  
Various comments about the statement of Lemma \tref{.gl} are collected in Remarks \tref{many}. The proof of Lemma \tref{.gl} is given in (\tref{p.gl}).

We close this section with a second application of the General Lemma: every parameterization of a curve leads to a parameterization of the branches of the curve. In particular, in Observation \tref{.O34.2} and Remark \tref{brn}, we show that if Data \tref{34.1} is in effect, with $\pmb k$ an algebraically closed field, then there is a one-to-one correspondence between the height one linear ideals  of $B$ and the branches of $\Cal C$. Indeed, it makes sense to talk about the branch ``$\Cal C(\ell)$'' of $\Cal C$ that corresponds to the ideal $(\ell)$ of $B$, where $\ell$ is a non-zero homogeneous linear form of $B$; furthermore, it makes sense to talk about the multiplicity of the branch $\Cal C(\ell)$. We return to this idea in Section 7, where we calculate the multiplicity of the branch $\Cal C(\ell)$ in terms of data found in the Jacobian ideal of the parameterization of $\Cal C$.

\definition{Data \tnum{34.1}} Let $\pmb k$ be a  field, $g_1,\dots,g_n$ be homogeneous forms of degree $d$ in  the polynomial ring $B=\pmb k[x,y]$, 
$\Psi\:\Bbb P^1\to \Bbb P^{n-1}$ be the morphism    which is given by $[g_1:\dots:g_n]$, $\Cal C$ be the image of $\Psi$,
and  $I$ be the ideal $(g_1,\dots,g_n)B$ of $B$. 
Assume that 
\roster
\item $I$ is minimally generated by $g_1,\dots,g_n$,
\item $I$ has height two, and
\item $\Psi\:\Bbb P^1\to \Cal C$ is a birational morphism.
\endroster
Let $\varphi$ be a homogeneous Hilbert-Burch matrix for $[g_1,\dots,g_n]$ and define the integers $d_1\le \dots \le d_{n-1}$ by requiring that 
$$0\to B(-d_1-d)\oplus \dots\oplus B(-d_{n-1}-d)@> \varphi>> B(-d)^n@> [g_1,\dots,g_n]>> I\to 0$$
is a homogeneous resolution of $I$. 
\enddefinition

\remark{Remark \tnum{34.2}}The hypotheses imposed on the parameterization $\Psi$ in Data \tref{34.1} are fairly mild. 
Furthermore, if a given parameterization of a rational curve $\Cal C$ fails to satisfy these hypotheses, one can reparameterize and obtain a parameterization of $\Cal C$ which does satisfy the hypotheses.  Hypothesis (1) is equivalent to the statement that $\Cal C\subseteq \Bbb P^{n-1}$ is ``non-degenerate''; that is, $\Cal C$ is not contained in any hyperplane in $\Bbb P^{n-1}$. Hypothesis (2) is equivalent to the statement ``the rational map $\Psi$ has no base points''. 
The homogeneous coordinate ring of $\Cal C$ is $\pmb k[g_1,\dots,g_n]=\pmb k[I_d]$ and the homomorphism $\pmb k[T_1,\dots, T_n]\twoheadrightarrow \pmb k[I_d]$, which sends $T_i$ to  $g_i$ for each $i$, induces the isomorphism $$\frac{\pmb k[T_1,\dots, T_n]}{I(\Cal C)}\cong \pmb k[I_d],\tag\tnum{Ts}$$ where $I(\Cal C)$ is the ideal generated by the homogeneous  polynomials which vanish on $\Cal C$. The homogeneous coordinate ring for the Veronese curve of degree $d$ is $\pmb k[x^d,x^{d-1}y,\dots,y^d]=\pmb k[B_d]$. Hypothesis  (3) is equivalent to the statement ``the domains $\pmb k[I_d]\subseteq \pmb k[B_d]$ have the same quotient field''; see Theorem \tref{2}. Theorem \tref{2} also provides a procedure for producing a birational reparameterization of $\Cal C$; an alternative procedure may be found in Section 6.1 of \cite{\rref{SWP}}. The integers $\{d_i\}$ could also be defined by insisting that each element of column $i$ of $\varphi$ be a homogeneous form of degree $d_i$. The hypotheses already in place guarantee that $1\le d_1$ and $\sum_{i=1}^{n-1}d_i=d$.\endremark

\bigskip
Observation \tref{O-1-18} is a small calculation, given originally in \cite{\rref{EU}}, which shows how the generalized row ideal $I_1(p\varphi)$ can be used to determine whether or not $p$ is on $\Cal C$. It serves as a good introduction to the Eisenbud-Ulrich interpretation of the fibers of the morphism $\Psi$ in terms of generalized rows of the syzygy matrix $\varphi$. This calculation is expanded and fine-tuned in the proof of the General Lemma.

\proclaim{Observation \tnum{O-1-18}} Adopt Data \tref{34.1} with $\pmb k$ algebraically closed. If $p$ is a point of $\Bbb P^{n-1}$, then   $p$ is on $\Cal C$ if and only if  $\operatorname{ht} (I_1(p\varphi))=1$.
\endproclaim
\remark{Remark} Notice that we use the symbol ``$p$'' to represent   both a point  in $\Bbb P^{n-1}$ and a row vector with $n$-entries from $\pmb k$; see Remark \tref{Unam} for more explanation.\endremark

\demo{Proof}     Assume first that $p$ is on $\Cal C$. In this case, $p=[g_1(q):\dots:g_n(q)]$ for some  point $q\in \Bbb P^1$. The equation $[g_1,\dots,g_n]\varphi=0$ holds by hypothesis. Therefore, the ideal 
$$I_1(p\varphi)=I_1([g_1(q),\cdots,g_n(q)]\varphi)$$ vanishes at $q$, and $I_1(p\varphi)$ is contained in the height one ideal $I(q)$ in $B$. This completes the proof of the direction $(\Rightarrow)$.

   Write $p$ as $[a_1:\dots:a_n]$ in $\Bbb P^{n-1}$.
Consider the ideal $$\frak I=I_2\left (\bmatrix a_1&\dots&a_n\\g_1&\dots&g_n\endbmatrix \right )$$of $B$. Observe that 
the zero locus of $\frak I$ is
$$V(\frak I)=\{q\in \Bbb P^1\mid \Psi(q)=p\},$$ which is the fiber $\Psi^{-1}(p)$ of the morphism $\Psi$ over $p$.

Let $\chi$ be an invertible matrix with entries in $\pmb k$ and $(0,\dots,0,1)\chi=(a_1,\dots,a_n)$. 
Define $[g_1',\dots,g_n']=g'$ by
$$g'=[g_1,\dots,g_n]\chi^{-1}.\tag\tnum{.gi0}$$The entries of $g'$ generate $I$ and the Hilbert-Burch matrix for $g'$ is $\varphi'=\chi \varphi$. One consequence of this last statement is the fact that the bottom row of $\varphi'$ generates the ideal 
$(g_1',\dots,g_{n-1}')\:\! g_n'$.
  Observe  that
$$\bmatrix a_1&\dots&a_n\\g_1&\dots&g_n\endbmatrix \chi^{-1}= \bmatrix 0&\dots&0&1\\g_1'&\dots&g_{n-1}'&g_n'\endbmatrix;$$ and therefore, $\frak I$ is generated by $(g_1',\dots,g_{n-1}')$. We now see that
$$I_1(p\varphi)=I_1([0,\dots,1]\varphi')=(g_1',\dots,g_{n-1}')\:\! g_n'= \frak I\:\! I.\tag\tnum{BBB}$$

Assume $p\notin \Cal C$. In this case, the zero locus $V(\frak I)$  of $\frak I$ is empty. The field $\pmb k$ is algebraically closed; so Hilbert's Nullstellensatz ensures that $\frak I$ is an $(x,y)$-primary ideal of $B$. We saw in (\tref{BBB}) that $\frak I\subseteq I_1(p\varphi)$; so, $I_1(p\varphi)$ is also an $(x,y)$-primary ideal. We conclude that $\operatorname{ht} (I_1(p\varphi))=2$. \qed \enddemo

\proclaim{Lemma \tnum{.gl}. (The General Lemma)} Adopt Data \tref{34.1}.
 Fix the point $p$ on $\Cal C$.   Assume that the fiber $\Psi^{-1}(p)$ is equal to the fiber $\Psi_{\bar {\pmb k}}^{-1}(p)$, where $\Psi_{\bar {\pmb k}}\:\Bbb P_{\bar {\pmb k}}^1 \to \Bbb P_{\bar {\pmb k}}^{n-1}$ is the extension of $\Psi$ to a morphism over the algebraic closure $\bar {\pmb k}$ of $\pmb k$. Let $q_1,\dots,q_s$ be the  $s$ distinct points  in $\Bbb P^1$
which comprise 
the fiber $\Psi^{-1}(p)$; $\frak p$ be the prime ideal in 
$\pmb k[I_d]$ which corresponds to the point $p$ on $\Cal C$;  $\frak q_1,\dots,\frak q_s$ be the prime ideals in $B$ which correspond to the points $q_1,\dots,q_s$ in $\Bbb P^1$,   and  
 $\Delta$ be the greatest common divisor of the entries of the row vector $p\varphi$. Write $\Delta=\ell_1^{c_1}\cdots \ell_t^{c_t}$, where $\ell_1,\dots,\ell_t$ are pairwise non-associate irreducible homogeneous forms in $B$. Let $R\subseteq S$ be the rings $\pmb k[I_d]_{\frak p}\subseteq B_{\frak p}$, $\hat R\subseteq\hat S$ be the completions of $R\subseteq S$ in the $\frak m_R$-adic topology, $J_1,\dots,J_u$ be the minimal prime ideals of $\hat R$, and $\frak M_1,\dots\frak M_v$ be the maximal ideals of $\hat S$. Then the parameters $s$, $t$, $u$, $v$, are all equal  to $s_{\Cal C,p}$ {\rm(}which is the number of branches of $\Cal C$ at $p${\rm)},  and, after re-numbering, $$q_i=V(\ell_i),\quad \frak q_i=\ell_i B,\quad \frak M_i=\frak q_i \hat S,\quad J_i=\ker (\hat R\to \hat S_{\frak M_i}),\quad\text{and}\quad c_i=e(\hat R/J_i),$$for  $1\le i\le s$. In particular, $\deg \Delta=e(R)=e(\Cal O_{\Cal C,p})=m_{\Cal C,p}$, and $c_1,\dots,c_s$ are the multiplicities of the branches of $\Cal C$ at $p$. 
\endproclaim
Some remarks pertaining to the statement of Lemma \tref{.gl} may be found in Remarks \tref{many}. The proof of Lemma \tref{.gl} is given in (\tref{p.gl}). We proceed immediately to the first round of applications. 

\proclaim{Theorem \tnum{Cor1}} Adopt Data \tref{34.1} with $\pmb k$ an algebraically closed field. 
Let  $p_1,\dots,p_z$ be the   singularities of $\Cal C$.
For each singular point $p_j$, let  $m_j$ be the multiplicity of $\Cal C$ at $p_j$  and  
$s_j$ be the number of branches of $\Cal C$ at $p_j$.
Then the following statements hold\,{\rm :}
\roster
\item The polynomial $\gcd I_1( p_j\varphi)$ in $\pmb k[x,y]$ has  degree equal to $m_j$ and has $s_j$ pairwise non-associated linear factors.
\item The polynomials $\gcd I_1( p_i\varphi)$ and $\gcd I_1( p_j\varphi)$ are   relatively prime for $i\neq j$.
\endroster\endproclaim

\demo{Proof}Assertion (1)  is explicitly stated as part of Lemma \tref{.gl}. We prove (2). Lemma \tref{.gl} shows that $\gcd I_1( p_j\varphi)=\prod_{u=1}^{s_j}\ell_{u,j}^{e_{u,j}}$, where the linear factors $\ell_{u,j}$ correspond to the points in the fiber $\Psi^{-1}(p_j)$. If $i\neq j$, then the fibers $\Psi^{-1}(p_i)$ and $\Psi^{-1}(p_j)$ are disjoint, so; the polynomials $\gcd I_1( p_i\varphi)$ and $\gcd I_1( p_j\varphi)$ are relatively prime. \qed\enddemo

As our first application of Theorem \tref{Cor1} we recover results of Song, Chen, and Goldman. The language and techniques of \cite{\rref{SCG}} are much different from ours; nonetheless, \cite{\rref{SCG}, Theorem 3} is essentially the same as the next result. Indeed, the General Lemma is our first attempt to interpret the results of \cite{\rref{SCG}} in terms of generalized row ideals. 
\proclaim {Corollary \tnum{ctgl}}Adopt Data \tref{34.1} with $\pmb k$ an algebraically closed field and $n=3$, and 
let  $p$ be a point on $\Cal C$. The following statements hold. 
\roster
\item The multiplicity of $\Cal C$ at $p$ satisfies $m_p\le d_2$.
\item If $m_p<d_2$, then $m_p\le d_1$.
\item The curve $\Cal C$ has a  singularity of multiplicity $d-1$ if and only if $d_1=1$.

\item The multiplicity of $\Cal C$ at $p$ is equal to $d_j$, for $j$ equal to $1$ or $2$, if and only if there exist homogeneous invertible matrices $\chi$ and $\xi$ {\rm(}the entries of $\chi$ are in $\pmb k$, the entries of $\xi$ are in $B${\rm)}  such that $p=[0:0:1]\chi$ and $\chi\varphi \xi$ is equal to 
$$\varphi'=\bmatrix P_1&Q_1\\P_2&Q_2\\0&Q_3\endbmatrix, $$ where the $P_i$ and the $Q_i$ are homogeneous forms from $B$ with $\deg P_i=d-d_j$, and $\deg Q_i=d_j$. Furthermore, in this situation,   the matrices $\chi$ and $\xi$   may be chosen  so that $\gcd (P_1,Q_1)=1$ and  $\gcd (P_2,Q_2)=1$. 
\endroster\endproclaim

\remark{Remark} The  invertible matrix $\chi$ from (4) gives rise to a linear automorphism of $\Bbb P^2$ which sends $p$ to $[0:0:1]$ as described in Remark \tref{sigh}.\endremark 
\demo{Proof} Let $p\varphi$ be the row vector $[a_1,a_2]$. The hypotheses tell us that
$$\deg a_1=d_1\le d_2=\deg a_2.$$ Theorem \tref{Cor1} now gives
 $$m_p=\deg \gcd (a_1,a_2)\le d_2$$ and  (1) is established. If $m_p<d_2$, then $a_1\neq 0$,   
and  $$m_p=\deg \gcd (a_1,a_2)\le \deg a_1= d_1,$$  and (2) is established.
We now prove (3). If $p'$ is a point on $\Cal C$ with $m_{p'}=d-1$, then (1) shows that
$$d-1=m_{p'}\le d_2=d-d_1\le d-1.$$  Thus, $m_{p'}=d-1\implies d_1=1$.
Conversely, if $d_1=1$, then the entries of column $1$ of $\varphi$ span the two-dimensional vector space $B_1$; hence, there exists a point $p'$ in $\Bbb P^2$ so that the left entry of $p'\varphi$ is zero. The ideal $(g_1,g_2,g_3)$ is $3$-generated; so the ideal $I_1(p'\varphi)$ is non-zero and is generated by a homogeneous form of degree $d-1$. Lemma \tref{.gl} (together with Observation \tref{O-1-18}) show that $p'$ is a singularity on $\Cal C$ of multiplicity $d-1$.

 Finally, we prove (4). Assume that $m_p=d_j$. It follows that 
$$\deg a_j=d_j=m_p=\deg \gcd(a_1,a_2)$$ and the ideal $(a_1,a_2)$ of $B$ must be principal. 
Thus, there are homogeneous invertible matrices $\chi$ and $\xi$  so that $(0,0,1)\chi=p$ and $\chi\varphi \xi$ has the form 
$$\varphi'=\bmatrix P_1&Q_1\\P_2&Q_2\\0&Q_3\endbmatrix.   $$The maximal minors of $\varphi'$ generate an ideal of height two; hence, $P_1,Q_1,Q_3$ have no factor in common. One may add a scalar multiple of row three to row one in order to ensure that the entries of row one are relatively prime. The analogous argument works for row two.  The later rounds of row operations do not change the bottom row of $\chi$. \qed\enddemo 
 
\remark{Remark} The method of \cite{\rref{SCG}} (see also \cite{\rref{CWL}, Theorem 1}) shows that  a rational plane curve   of degree $d$
can have at most one singularity of multiplicity greater than $d/2$. This is clear by B\'ezout's Theorem and is also implied by our method since one can not create two zeros in a column of a Hilbert-Burch matrix without violating the height condition. \endremark 

We next collect the preliminary results (Lemmas \tref{.O33.10'} and \tref{hat}, and Proposition \tref{!32.1})   which are used in the proof of Lemma \tref{.gl}. 
\proclaim{Lemma \tnum{.O33.10'}}Let    $(A,\frak m_A,\pmb k_A)\subseteq (B,\frak m_B,\pmb k_B)$ be one-dimensional local domains. Assume that $\pmb k_A$ is infinite. If $B$ is the integral closure of $A$ and $B$ is  finitely generated as an $A$-module,
then $\frak m_AB=\frak m_B^{e/r}$, where $e=e(A)$ is the multiplicity of the local ring $A$ and $r=[\pmb k_B \:\! \pmb k_A]$ is the degree of the field extension $\pmb k_A\subseteq \pmb k_B$.\endproclaim

\demo{Proof} The hypothesis about $\pmb k_A$ ensures that there exists a minimal reduction 
  of $\frak m_A$ generated by one element $z$ and that $e=\lambda_{A}(A/zA)$; see, for example, \cite{\rref{SH}, Prop.~11.2.2}. The domain $B$ is  normal, local, and one dimensional; so, $B$ is a Principal Ideal Domain and the equation $z\frak m_A^uB=\frak m_A^{u+1}B$, for some $u$, tells us that $zB=\frak m_AB$.
We compute $$e=\lambda_{A}(A/zA)=\lambda_{A}(B/zB)=\lambda_{B}(B/zB)\,r\,.\tag\tnum{str}$$
The middle equality is obtained from the picture $$\xymatrix{ &B\ar@{-}[d]\\ A  \ar@{-}[ur]\ar@{-}[d]&
 zB\\
zA\ar@{-}[ur] }$$
All lengths are finite; in particular, $\lambda_{A}(B/A)$ is finite because $A$ is a one-dimensional domain and $B$ is a module-finite extension of $A$ with $B\subseteq \operatorname{Quot}(A)$.  Multiplication by $z$ gives an isomorphism of $A$-modules $B/A\cong zB/zA$; therefore the $A$-modules $A/zA$ and $B/zB$ have the same length.  The equality on the right is due to the fact that every factor in a composition series for the $B$-module $B/zB$ is $\pmb k_B$.  The formulas of (\tref{str}) have been established and $\lambda_{B}(B/zB)=e/r$. 
The only quotient of $B$ with length $e/r$, as a $B$-module, is $B/\frak m_B^{e/r}$. 
We conclude that $\frak m_AB=zB=\frak m_B^{e/r}$. \qed \enddemo

 Proposition \tref{!32.1} is the well-known   correspondence between the minimal prime ideals of the completion $\widehat{\Cal O_{\Cal C,p}}$ of the local ring $\Cal O_{\Cal C,p}$ and the maximal ideals of the integral closure $\overline{\Cal O_{\Cal C,p}}$ of $\Cal O_{\Cal C,p}$. This result is stated as Exercise 1 on page 122 in Nagata \cite{\rref{Na}};
see \cite{\rref{Ka}, Corollary 5} and  \cite{\rref{Kunz}, Thm. 16.14} for a proof.

 \proclaim{Proposition \tnum{!32.1}} Let $T$ be the integral closure of the  one-dimensional local domain $(R,\frak m)$. Suppose $T$ is finitely generated as a module over $R$.  Let $\hat R$ and $\hat T$ represent the completions of $R$ and $T$ in the $\frak m$-adic topology. 
Then  there is a   one-to-one correspondence between the minimal prime ideals of $\hat R$ and  the maximal ideals of $T$. If $\Cal M$ is a maximal ideal of $T$, then the corresponding 
minimal prime ideal of $\hat R$ is
$\ker (\hat R\to \hat T_{\Cal M\hat T})$.
\endproclaim

\remark{\bf Remark \tnum{back}} When the hypotheses of Proposition \tref{!32.1} are in effect then $\hat T$ is the integral closure of $\hat R$. Indeed, $\hat T$ is a normal ring which is finitely generated as a module over $\hat R$. Lemma \tref{!.p171} shows that $\hat T\subseteq \operatorname{Quot} \hat R$. \endremark 

\proclaim{Lemma \tnum{!.p171}} Let $R\subseteq T\subseteq \operatorname{Quot}(R)$ be rings with $(R, \frak m)$ a local domain and $T$ a finitely generated $R$-module. Let $\hat R$ and $\hat T$ be the completions of $R$ and $T$ in the $\frak m$-adic topology. Then $\hat R$ and $\hat T$ have the same total quotient ring.\endproclaim 

\demo{Proof} Let $U=R\setminus \{0\}$. Notice first that $\operatorname{Quot}(R)=U^{-1}R=U^{-1}T$. 
Every element of $U$ is regular on both rings $\hat R\subseteq\hat T$; so, $U^{-1}(\hat R)\subseteq U^{-1}(\hat T)\subseteq \operatorname{Quot}(\hat T)$. The  rings $U^{-1}(\hat R)\subseteq U^{-1}(\hat T)$ are equal because
$$\split U^{-1}(\hat T)&{}=U^{-1}(T\otimes_R\hat R)=U^{-1}(T)\otimes_{U^{-1}(R)} U^{-1}(\hat R)=U^{-1}(R)\otimes_{U^{-1}(R)} U^{-1}(\hat R)\\&{}= U^{-1}(\hat R).\endsplit $$
We have $U^{-1}\hat T=U^{-1}\hat R\subseteq \operatorname{Quot}\hat T$. A typical element of $\operatorname{Quot}\hat T$ is $z/w$, where $z$ and $w$ are in $\hat T$ with $w$ regular on $\hat T$. There exists $u\in U$ with $uz,uw$ in $\hat R$. Of course $uw$ is regular on $\hat R\subseteq \hat T$. So, $z/w=uz/uw\in \operatorname{Quot}\hat R$.  
\qed \enddemo

In Lemma \tref{hat} we show that if $A$ is a complete local domain and $z$ is an indeterminate over $A$, then the local ring $A(z)$, see (\tref{Aofz}), is analytically irreducible. Our proof uses the fact that the completion of an excellent normal ring is a normal ring. We will apply Lemma \tref{hat} in the proof of Lemma \tref{.gl} to calculate the multiplicity of a branch of a curve.

\proclaim{Observation \tnum{O-apr-2}} If $A\subseteq C$ is a module finite extension of local rings and $z$ is an indeterminate over $C$, then $A(z)\subseteq C(z)$ is module finite extension of local rings. \endproclaim

\demo{Proof} Let $\frak m_A$ and  $\frak m_C$   be the maximal ideals  of $A$ and $C$, respectively. The Cohen-Seidenberg Theorems ensure that $\frak m_C$ is the only prime ideal of $C$ which contracts to $\frak m_A$; so, in particular, the radical of the ideal $\frak m_A C$ of $C$  is $\operatorname{rad} \frak m_A C=\frak m_C$, and the ideal $\frak m_CC[z]$ contracts to $\frak m_AA[z]$. 

Consider the multiplicatively closed subsets
$$S=A[z]\setminus \frak m_AA[z]\quad \subseteq \quad T= C[z]\setminus \frak m_CC[z]$$ of $C[z]$. We know that
$S^{-1}(A[z])\subseteq S^{-1}(C[z])$ is a module finite extension, with $A(z)=S^{-1}(A[z])$. Moreover,  $C(z)=T^{-1}(C[z])$ is a further localization of $S^{-1}(C[z])$ at a maximal ideal of $S^{-1}(C[z])$. To prove the result, it suffices to show that $C(z)$ is already equal to $S^{-1}(C[z])$; and therefore, it suffices to show that $S^{-1}(C[z])$  is already a local ring. 

The inclusion $A(z)\subseteq S^{-1}(C[z])$ is a module finite extension and $A(z)$ is a local ring with maximal ideal $\frak m_AA(z)$. According to the Cohen-Seidenberg Theorems, the maximal spectrum of $S^{-1}(C[z])$ is contained in
$$ \text{$\{PS^{-1}(C[z])\mid   P$ is a prime ideal of $C[z]$, minimal over $\frak m_AC[z]\}$}.$$The above set consists of exactly one ideal because $\operatorname{rad}(\frak m_AC)C[z]$ is a prime ideal of $C[z]$ since 
$$\tsize \frac{C[z]}{\operatorname{rad}(\frak m_AC[z])}=\frac{C }{\operatorname{rad}(\frak m_AC )}[z]=\frac{C }{ \frak m_C }[z],$$ which is a domain. Thus, 
$S^{-1}(C[z])$ is a local ring and the proof is complete. \qed \enddemo

\smallskip
\proclaim{Lemma \tnum{hat}} If $A$ is a complete local domain and $z$ is an indeterminate over $A$, then the completion  ${A(z)}\!\widehat{\phantom{x}}$ of the local ring $A(z)$ is also a domain.\endproclaim

\demo{Proof} Let $V$ be the integral closure of $A$. Our proof amounts to showing that each vertical map in the commutative diagram of rings:
$$ \xymatrix{ V&\subseteq &V(z)&\subseteq &V(z)\!\widehat{\phantom{x}}\\
A\ar[u]&\subseteq &A(z)\ar[u]&\subseteq &A(z)\!\widehat{\phantom{x}}\ar[u]}
$$
is a module-finite inclusion.
The ring $(A,\frak m_A)$ is  local, Noetherian, and complete; hence excellent; see, for example, \cite{\rref{M80}, pg.~260}. It follows that the integral closure $V$ of $A$ is finitely generated as an $A$-module. The ring $V$ is necessarily local. Indeed, $V$ is semi-local and complete in the $\frak m_A$-adic topology; hence, $V$ is a direct product of local rings. On the other hand, $V$ is a domain since it is contained in the quotient field of $A$ and so only one local ring appears in the previously mentioned direct product decomposition of $V$. Apply Observation \tref{O-apr-2} to see that $V(z)$ is finitely generated as an $A(z)$-module. Let $V(z)\!\widehat{\phantom{x}}$ be the completion of the local ring $V(z)$. We now see that the natural map ${A(z)}\!\widehat{\phantom{x}}\to V(z)\!\widehat{\phantom{x}}$ is an injection; since 
$V(z)\!\widehat{\phantom{x}}=V(z)\otimes_{A(z)}{A(z)}\!\widehat{\phantom{x}}$ and the completion map $A(z)\to {A(z)}\!\widehat{\phantom{x}}$ is flat. The rings $V$ and $V(z)$ are excellent normal local rings and therefore the completion $V(z)\!\widehat{\phantom{x}}$ is also an excellent normal local ring; see, for example, \cite{\rref{M80}, Thm 79}. The local normal ring $V(z)\!\widehat{\phantom{x}}$ is necessarily a normal domain and the subring ${A(z)}\!\widehat{\phantom{x}}$ of $V(z)\!\widehat{\phantom{x}}$ is also a domain. \qed \enddemo

The following remarks pertain to the statement of Lemma \tref{.gl}. 
\remark{\bf Remarks \tnum{many}} {(a)} We use the symbol $p$ to represent a point $[a_1:\dots:a_n]$ in $\Bbb P^{n-1}$ as well as a row vector $[a_1,\dots,a_n]$. The meaning will be clear from context; see Remark \tref{Unam}.

\smallskip\flushpar 
{(b)} The important conclusions are $s=t$, $\frak q_i=\ell_i B$, and $c_1,\dots,c_s$ are the multiplicities of the branches of $\Cal C$ at $p$.  The other statements are well-known, or follow easily from these, as we explain below.

\smallskip\flushpar{(c)}  The prime ideal in $\pmb k[T_1,\dots, T_n]$ which corresponds to the point $p$ in $\Bbb P^{n-1}$ is $I_2(M)$ for 
$$M=\pmatrix a_1&\dots&a_n\\
T_1&\dots&T_{n}\endpmatrix.$$ The prime ideal $\frak p$ in $\pmb k[I_d]$ which corresponds to $p$ on $\Cal C$ is $I_2(M)\pmb k[I_d]$,  see (\tref{Ts}); thus, $$\frak p=I_2\pmatrix a_1&\dots&a_n\\
g_1&\dots&g_{n}\endpmatrix \pmb k[I_d].$$

\smallskip\flushpar{(d)} Recall that $\frak p$ is a prime ideal of the ring $\pmb k[I_d]$, $R=\pmb k[I_d]_{\frak p}$, and $S=B_{\frak p}$. Our notation  means that $S$ is the localization $U^{-1}B$ of $B$ at the multiplicatively closed set $U=\pmb k[I_d]\setminus \frak p$. Let $T=\pmb k[B_d]_{\frak p}$. (In other words, $T$ is equal to $U^{-1}(\pmb k[B_d])$.) The ring inclusions $\pmb k[I_d]\subseteq \pmb k[B_d]\subseteq B$ are integral extensions; so the ring inclusions $R\subseteq T\subseteq S$ are integral extensions. The Veronese ring $\pmb k[B_d]$ is a normal domain and the domains $\pmb k[I_d]\subseteq \pmb k[B_d]$  have the same quotient field by hypothesis (3) of Data \tref{34.1}; hence, $\pmb k[B_d]$ is the integral closure of $\pmb k[I_d]$ and $T$ is the integral closure of $R$. 

\smallskip\flushpar{(e)} The ring inclusions $R\subseteq T\subseteq S$ are module finite extensions and $R$ is a local ring with maximal ideal $\frak m_R=\frak p\pmb k[I_d]_{\frak p}$. It follows from the Cohen-Seidenberg Theorems that $T$ and $S$ are semi-local rings. Moreover, $\operatorname{Proj} \pmb k[B_d]=\operatorname{Proj} B$; so the function $M\mapsto M\cap T$ gives a one-to-one correspondence between the maximal ideals of $S$ and the maximal ideals of $T$.

\smallskip\flushpar{(f)} The $\frak m_R$-adic topology on $S$ is equivalent to the $J$-adic topology on $S$, where $J$ is the Jacobson radical of $S$. It is well known (see, for example, \cite{\rref{M89}, Thm. 8.15} or \cite{\rref{Kunz}, Thm. K.11}) that the natural map $$\hat S\to \hat S_{\frak M_1}\times \cdots \times \hat S_{\frak M_v}$$ is an isomorphism; furthermore, the local ring $\hat S_{\frak M_i}$ is complete for each $i$ and $\hat S_{\frak M_i}$ is equal to the completion of the local ring $S_{\frak M_i\cap S}$ in the $(\frak M_i\cap S)$-adic topology.
Each ring $S_{\frak M_i\cap S}$ is a one-dimensional regular ring; hence, each $\hat S_{\frak M_i}$ is a complete DVR.
Furthermore, the maximal ideals of $S$ are $\{\frak M_i\cap S\mid 1\le i\le v\}$. 

\smallskip\flushpar{(g)} The statements of Remark (f) also apply to $T$. So there is a commutative diagram
$$ \xymatrix{ \hat S&=&\hat S_{\frak M_1}&\times&\dots&\times&\hat S_{\frak M_v}\\
\hat T\ar@{^{(}->}[u]&=&\hat T_{\frak M_1\cap \hat T}\ar@{^{(}->}[u]&\times&\dots&\times&\hat T_{\frak M_v\cap \hat T},\ar@{^{(}->}[u]}$$
where $\hat T_{\frak M_i\cap \hat T}$ is the localization of the ring $ \hat T$ at the maximal ideal $\frak M_i\cap \hat T$ and also is the completion of the local ring $T_{\frak M_i\cap T}$ at the maximal ideal $\frak M_i\cap T$. Each map $\hat T_{\frak M_i\cap \hat T}\hookrightarrow \hat S_{\frak M_i}$ is an integral extension.

\smallskip\flushpar{(h)} Proposition \tref{!32.1} may be applied to the rings $R\subseteq T$ in order to see that
   $u=v$, and, after renumbering, $J_i=\ker (\hat R\to \hat T_{\frak M_i\cap \hat T})$. 
Also, $R$ is a reduced local ring which is the localization of a finitely generated algebra over a field; so $\hat R$ is reduced; see, \cite{\rref{ZS}, Vol. 2, Chapt\. 8, Sect\. 13} or \cite{\rref{SH}, page 177, item (1)}.
We enlarge the commutative diagram of (g) to obtain the  commutative diagram:
$$ \xymatrix{ \hat S&=&\hat S_{\frak M_1}&\times&\dots&\times&\hat S_{\frak M_v}\\
\hat T\ar@{^{(}->}[u]&=&\hat T_{\frak M_1\cap \hat T}\ar@{^{(}->}[u]&\times&\dots&\times&\hat T_{\frak M_v\cap \hat T}\ar@{^{(}->}[u]\\
\hat R\ar@{^{(}->}[u]&\hookrightarrow&\frac{\hat R}{J_1}\ar@{^{(}->}[u]&\times&\dots&\times&\frac{\hat R}{J_v}.\ar@{^{(}->}[u]}$$
The ring $\hat T$ is the integral closure of $\hat R$ by Remark \tref{back}; and therefore, for each $i$, the ring $\hat T_{\frak M_i\cap \hat T}$ is the integral closure of $\frac{\hat R}{J_i}$.
\endremark

\demo{\bf (\tnum{p.gl}) Proof of Lemma \tref{.gl}}The ring extension $\pmb k[I_d]\subseteq B$ is integral; so every maximal ideal of $S$ has the form $\frak q S$, where $\frak q$ is a height one homogeneous prime ideal of $B$ which is minimal over $\frak p B$ and which satisfies $\frak q \cap \pmb k[I_d]=\frak p$. Let $\frak q$ be a homogeneous prime ideal in $B$ for which $\frak qS$ is a maximal ideal of $S$. The ideal $\frak q$ is principal and is generated by some homogeneous form $f\in B$. Let $q\in \Bbb P^1_{\bar {\pmb k}}$ be a root of $f$. The generators of $\frak p$, which may be found in  Remark 
\tref{many}~(c), are in the ideal $\frak q=(f)$; and therefore, the generators of $\frak p$ all vanish at $q$. It follows that $\Psi_{\bar {\pmb k}}(q)=p$. The hypothesis $\Psi_{\bar {\pmb k}}^{-1}(p)=\Psi^{-1}(p)$ ensures that $q$ is already in $\Bbb P^1$; and therefore, $q\in \{q_1,\dots,q_s\}$, $f$ is a linear polynomial, and $\frak q\in \{\frak q_1,\dots,\frak q_s\}$.  We conclude that $$\text{the maximal ideals of $S$ are $\{\frak q_1S,\dots,\frak q_sS\}$.}\tag\tnum{maxs}$$ It follows from Remark \tref{many}~(f) that $v=s$ and the maximal ideals of $\hat S$ are $\frak M_i=\frak q_i\hat S$, for $1\le i\le s$. 

We saw in (\tref{BBB}) that  $$\frak pB\:\! I=I_1(p\varphi).\tag\tnum{.rowop}$$ 
 Fix $g'=g\chi$, as described in (\tref{.gi0}), with $\frak p=(g_1',\dots,g_{n-1}')$. 
We compute the saturation $\frak pB\:\! (x,y)^{\infty}$ two different ways. On the one hand, $\frak pB\:\! (x,y)^{\infty}$ is equal to the intersection of the $\frak q$-primary components of $\frak pB$ as $\frak q$ roams over all of the height one prime ideals of $B$ in $\operatorname{Ass} B/\frak pB$. For each such $\frak q$, the $\frak q$-primary component of $\frak pB$ is $\frak pB_{\frak q}\cap B$ and the ring $B_{\frak q}$ is a DVR; so, $\frak pB_{\frak q}=\frak q^wB_{\frak q}$ for some exponent $w$. Thus,
$$\frak pB\:\! (x,y)^{\infty}=\frak q_1^{(w_1)}\cap \dots \cap\frak q_s^{(w_s)}=\frak q_1^{ w_1}\cap \dots \cap\frak q_s^{w_s}=\frak q_1^{ w_1}\cdots \frak q_s^{w_s}.$$ We have taken advantage of the fact that each $\frak q_i$ is principal in the  Unique Factorization Domain $B$. On the other hand, the ideal $I$ is $(x,y)$-primary and ``$\:\!$'' is associative, so we use (\tref{.rowop}) to see that
$$\split \frak pB\:\! (x,y)^{\infty}={}&\frak pB\:\! I^{\infty}=(\frak pB\:\! I)\:\! I^{\infty}=I_1(p\varphi)\:\! I^{\infty}
=I_1(p\varphi)\:\! (x,y)^{\infty}=\Delta B\\{}={}&\ell_1^{c_1}\cdots \ell_t^{c_t}B.\endsplit$$We have two factorizations of $\Delta$ into non-associate irreducible factors. We conclude that $s=t$, and after renumbering, $\ell_iB=\frak q_i$ and $w_i=c_i$. It follows that each $\ell_i$ is a linear form and  $q_i=V(\ell_i)$. 

We next calculate the multiplicity $e(\hat R/J_i)$ for each $i$. Recall that the maximal ideals of $T$ are
$$(\frak M_i\cap \hat T)\cap T=(\frak q_i\cap \pmb k[B_d])T\quad\text{for $1\le i\le s$}$$ and that the completion $(T_{(\frak q_i\cap k[B_d])T})\!\widehat{\phantom{x}}$ of the local ring $T_{(\frak q_i\cap \pmb k[B_d])T}$ at the maximal ideal ${(\frak q_i\cap \pmb k[B_d])}T_{(\frak q_i\cap \pmb k[B_d])T}$ is equal to the localization $\hat T_{\frak M_i\cap \hat T}$ of the complete ring $\hat T$ at the maximal ideal $\frak M_i\cap \hat T$. We simplify the notation by letting $$R_i=\hat R/J_i,\  T_i=(T_{(\frak q_i\cap \pmb k[B_d])T})\!\widehat{\phantom{x}}=\hat T_{\frak M_i\cap \hat T},\ \text{and}\ \frak m_{T_i}=(\frak q_i\cap \pmb k[B_d])T_i=(\frak M_i\cap \hat T)T_i.$$
Recall from Remark \tref{many}~(h) that $T_i$ is the integral closure $\bar R_i$ of $R_i$. The maximal ideal of $R_i$ is 
$$\frak m_{R_i}=\frak m_{\hat{R}/J_i}=\frak m_{\hat{R}}/J_i=(\frak m_R\hat R)/J_i = ((\frak pR)\hat R)/J_i= (\frak p\hat R)/J_i.$$We know (see Remark \tref{many}~(h)) that $J_i=\ker(\hat R\to T_i)$; so, $$\frak m_{R_i}T_i=\frak pT_i.\tag\tnum{.mrt}$$ The ideal $
\frak pT_i$ is completely determined by the $(\frak q_i\cap \pmb k[B_d])$-primary component of $\frak p\pmb k[B_d]$. We have already computed the primary components of $\frak pB$ corresponding to the prime ideals minimal in the support of $B/\frak pB$:
$$\frak pB\:\!(x,y)^{\infty}=\frak q_1^{c_1}\cap\dots\cap\frak q_s^{c_s}.$$ 
The rings  $\pmb k[B_d]_{\frak q_i\cap \pmb k[B_d]}\subseteq B_{\frak q_i}$ are DVRs. One can choose the same generator for the maximal ideal of these two rings. There is no difficulty in seeing that if $f$ is an element of $\pmb k[B_d]$ and $r$ is arbitrary, then 
$$f\in (\frak q_i\cap \pmb k[B_d])^r\pmb k[B_d]_{\frak q_i\cap \pmb k[B_d]}\iff f\in \frak q_i^r B_{\frak q_i};$$ and therefore, there is no difficulty in seeing that the 
$(\frak q_i\cap \pmb k[B_d])$-primary component of $\frak p\pmb k[B_d]$ is $(\frak q_i\cap \pmb k[B_d])^{c_i}$. It follows from  (\tref{.mrt}) that
$$\frak m_{R_i}T_i=\frak pT_i=\frak m_{T_i}^{c_i}T_i.$$
On the other hand, Lemma \tref{.O33.10'} shows that $\frak m_{R_i}T_i=\frak m_{T_i}^{e(R_i)}T_i$ because $R_i\subseteq T_i$ are local one-dimensional domains with common residue class field $\pmb k(g_n')\subseteq R_i$, for $g_n'$ as found in (\tref{.gi0}); furthermore, $T_i$ is the integral closure of $R_i$ and $T_i$ is finitely generated as an $R_i$-module. Thus,
$$\frak m_{T_i}^{c_i}T_i=\frak m_{R_i}T_i=\frak m_{T_i}^{e(R_i)}T_i,$$ and $c_i=e(R_i)$.

Next, we relate the degree of the polynomial $\Delta$ to the multiplicity of the local ring $R$. The modules $R/\frak m_R^r$ and $\hat R/\frak m_R^r\hat R$ are equal for all $r$; so $e(R)=e(\hat R)$. The associativity formula for multiplicities yields 
$$e(\hat R)=\sum_{i=1}^se(\hat R/J_i).$$
Thus,
$$e(R)=e(\hat R)= \sum_{i=1}^se(\hat R/J_i)=\sum_{i=1}^sc_i=\deg (\ell_1^{c_1}\cdots \ell_s^{c_s})=\deg \Delta.$$

We translate the information we have collected about the rings $R\subseteq T$ and $\hat R \subseteq \hat T$ to information about the rings
$\Cal O_{\Cal C,p}\subseteq \overline{\Cal O_{\Cal C,p}}$ and $\widehat{\Cal O_{\Cal C,p}}\subseteq \left (\overline{\Cal O_{\Cal C,p}}\right )\!\widehat{\phantom{X}}$. Recall first that all four rings are  subrings  of $\operatorname{Quot}(\pmb k[I_d])=\operatorname{Quot}(\pmb k[B_d])$:   
$$\allowdisplaybreaks \alignat 3 \Cal O_{\Cal C,p}={}&\{\tsize\frac fg\mid \text{ $f\in \pmb k[I_d]$ and $g\in \pmb k[I_d]\setminus \frak p$ are homogeneous of the same degree}\},\\
\overline{\Cal O_{\Cal C,p}}={}&\{\tsize\frac fg\mid \text{ $f\in \pmb k[B_d]$ and $g\in \pmb k[I_d]\setminus \frak p$ are homogeneous of the same degree}\},\\ R={}&\pmb k[I_d]_{\frak p}=\{\tsize\frac fg\mid \text{ $f\in \pmb k[I_d]$ and $g\in \pmb k[I_d]\setminus \frak p$} \},\text{ and} \\ T={}&\pmb k[B_d]_{\frak p}=\{\tsize\frac fg\mid \text{ $f\in \pmb k[B_d]$ and $g\in \pmb k[I_d]\setminus \frak p$}\}.\endalignat $$ 
Observe that $g_n'$ is a unit of $R$ which is transcendental over $\Cal O_{\Cal C,p}$.
It is not difficult to check that the two subrings $\Cal O_{\Cal C,p}(g_n')$, see (\tref{Aofz}), and $R$ of $\operatorname{Quot} (\pmb k[I_d])$ are equal and therefore, $R$ and $\Cal O_{\Cal C,p}$ have the same multiplicity.  The ring $\overline{\Cal O_{\Cal C,p}}$ is a subring 
of  $T$. Dehomogenization and homogenization
$$\matrix \frak a&\to& \frak a T,\\\frak A\cap \overline{\Cal O_{\Cal C,p}}&\leftarrow &\frak A,\endmatrix\tag\tnum{hd}$$ 
 provide   a one-to-one correspondence between the maximal ideals $\frak a$ of $\overline{\Cal O_{\Cal C,p}}$ and the maximal ideals of $\frak A\in\{\frak M_i\cap T\mid 1\le i\le v\}$ of $T$. Recall that the invariant $s_{\Cal C,p}$ is equal to the number of branches of $\Cal C$ at $p$. Thus,
$$\split s_{\Cal C,p}&{}=\text{ the number of minimal prime ideals of $\widehat{\Cal O_{\Cal C,p}}$}\\
&{}=\text{ the number of maximal prime ideals of $\overline{\Cal O_{\Cal C,p}}$}\\
&{}=\text{ the number of maximal prime ideals of $T$}=v
.\endsplit$$

Finally, we compare the multiplicities $e(\hat R/J_i)$ and $e(\widehat{\Cal O_{\Cal C,p}}/\Cal J_i)$, where $\Cal J_i$ is the minimal prime of $\widehat{\Cal O_{\Cal C,p}}$ which corresponds to the minimal prime ideal $J_i$ of $\hat R$. We have already shown that $e(\hat R/J_i)=c_i$. We complete the proof of the Lemma by showing that 
$$e(\widehat{\Cal O_{\Cal C,p}}/\Cal J_i)=e(\hat R/J_i).\tag\tnum{shw}$$
  As before, we have $R=\Cal O_{\Cal C,p}(g_n')$. It is well known, and easy to show, that if $z$ is an indeterminate over a local ring $A$, then
$$\widehat{A(z)}=\widehat{\hat A(z)},$$where the meaning of $A(z)$ is given in (\tref{Aofz}) and each completion is the completion of a local ring in its maximal ideal adic topology.  It follows that $$\hat R=\left (\widehat{\Cal O_{\Cal C,p}}(g_n')\right )\!\widehat{\phantom{X}}.$$The ideal $\Cal J_i$ of $\widehat{\Cal O_{\Cal C,p}}$  is $J_i\cap \widehat{\Cal O_{\Cal C,p}}$. We first prove 
$$ \Cal J_i\hat R=J_i.\tag\tnum{ce}$$ 
The ring ${\widehat{\Cal O_{\Cal C,p}}}/{\Cal J_i}$ is a complete local domain and $g_n'$ is an indeterminate over this ring; so, Lemma \tref{hat} yields that $\left (\left (\frac{\widehat{\Cal O_{\Cal C,p}}}{\Cal J_i}\right )(g_n')\right )\!\widehat{\phantom{X}}$
  is a domain. On the other hand,
$$\tsize \left (\left (\frac{\widehat{\Cal O_{\Cal C,p}}}{\Cal J_i}\right )(g_n')\right )\!\widehat{\phantom{X}}= 
\left (\frac{\widehat{\Cal O_{\Cal C,p}}(g_n')}{\Cal J_i\widehat{\Cal O_{\Cal C,p}}(g_n')}\right )\!\widehat{\phantom{X}}
= \frac{\left (\widehat{\Cal O_{\Cal C,p}}(g_n')\right )\!\widehat{\phantom{X}}} {\Cal J_i\left (\widehat{\Cal O_{\Cal C,p}}(g_n')\right )\!\widehat{\phantom{X}}}
=\frac{\hat R}{\Cal J_i\hat R}.
$$ Thus, ${\Cal J_i}\hat R$ is a prime ideal of $\hat R$ with ${\Cal J_i}\hat R\subseteq J_i$. The ideal $J_i$ is a minimal prime ideal of $\hat R$; so  (\tref{ce}) is established and we have
$$\tsize e\left (\frac{\hat R}{J_i}\right )=e\left (\frac{\hat R}{\Cal J_i\hat R}\right )=e\left (\left (\left (\frac{\widehat{\Cal O_{\Cal C,p}}}{\Cal J_i}\right )(g_n')\right )\!\widehat{\phantom{X}}\right )=e \left (\left (\frac{\widehat{\Cal O_{\Cal C,p}}}{\Cal J_i}\right )(g_n')\right ) 
 =e \left (\frac{\widehat{\Cal O_{\Cal C,p}}}{\Cal J_i}\right ),   $$which is (\tref{shw}).
\qed \enddemo

We close this section with the observation that every parameterization of a curve leads to a parameterization of the branches of the curve.

\proclaim{Observation \tnum{.O34.2}} Adopt the Data of \tref{34.1}, with $\pmb k$ algebraically closed. Then there is a one-to-one correspondence between the points of $\Bbb P^1$ and the branches of $\Cal C$. \endproclaim

\demo{Proof}
Fix  a point  $q$  in $\Bbb P^1$. Let $p$ be the point $\Psi(q)$ on $\Cal C$. Form the ideal $\frak p$ of $\pmb k[I_d]$ as described in Remark \tref{many}~(c) and form the rings $$R=\pmb k[I_d]_{\frak p}\subseteq T=\pmb k[B_d]_{\frak p}\subseteq S=B_{\frak p}$$as described in Remark \tref{many}~(d). There are explicit, well-defined, one-to-one correspondences between each of the following sets:\phantom{\tnum{c1},\tnum{c2},\tnum{c3},\tnum{c4},\tnum{c5}}
$$\xymatrix{
\Psi^{-1}(p)\ar@{<->}[r]_{(\tref{c1})}
&\operatorname{MaxSpec}(S)\ar@{<->}[r]_{(\tref{c2})}
&\operatorname{MaxSpec}(T)\ar@{<->}[d]_{(\tref{c3})}\\
\text{The branches of $\Cal C$ at $p$}\ar@{<->}[r]_{(\tref{c5})}& \text{The Min Primes of $\widehat{\Cal O_{\Cal C,p}}$}\ar@{<->}[r]_{\ \ \ \ \ (\tref{c4})} &\operatorname{MaxSpec}(\overline{\Cal O_{\Cal C,p}})
}$$
If $q_1$ is a point in $\Psi^{-1}(p)$ and $\frak q_1$ is the homogeneous prime ideal of $B$ which corresponds to $q_1$, then the correspondence (\tref{c1}) sends $q_1$ to $\frak q_1S$ as shown in (\tref{maxs}). The correspondence (\tref{c2}) is described in Remark \tref{many}~(e). The correspondence (\tref{c3}) is explained in (\tref{hd}). Proposition \tref{!32.1} accounts for   (\tref{c4}), and (\tref{c5}) is the definition of branch. \qed
\enddemo
\remark{\bf Remark \tnum{brn}} We say that an ideal of $B=\pmb k[x,y]$ is a {\it height one linear ideal} if it is generated by one non-zero linear form. There is a one-to-one correspondence between the height one linear ideals of $B$ and the points of $\Bbb P^1$. Thus, Observation \tref{.O34.2} gives a one-to-one correspondence between the height one linear ideals of $B$ and the branches of $\Cal C$. If $(\ell)$ is a height one linear ideal of $B$, then let $\Cal C(\ell)$ be the corresponding branch of $\Cal C$. It makes sense to speak about the {\it multiplicity of the branch $\Cal C(\ell)$} because $\ell$ corresponds to a point $q$ in $\Bbb P^1$, $\Psi(q)=p$ is a point on $\Cal C$, $\Cal C(\ell)$ is a minimal prime ideal of $\widehat{\Cal O_{\Cal C,p}}$, and the multiplicity of the local ring $\widehat{\Cal O_{\Cal C,p}}/\Cal C(\ell)$ is called the multiplicity of  the branch $\Cal C(\ell)$.\endremark  

\bigskip
We return to the ideas of Remark \tref{brn} in Section 7, where we calculate the multiplicity of the branch $\Cal C(\ell)$ in terms of data found in the Jacobian ideal of the parameterization of $\Cal C$.

\bigpagebreak
\SectionNumber=2\tNumber=1
\heading Section \number\SectionNumber. \quad The Triple Lemma.
\endheading 

The General Lemma in Section 1 shows how the generalized row ideals of the Hilbert-Burch matrix of a parameterization of a rational    curve $\Cal C$ encode information about the singularities on $\Cal C$. In the present section we describe how one can read information about the infinitely near singularities of $\Cal C$ from the Hilbert-Burch matrix. Theorem \tref{TPL1} deals with infinitely near singularities in the first neighborhood of $\Cal C$. Corollary \tref{CTLL} is mainly concerned with infinitely near singularities in the second neighborhood of $\Cal C$. We used our first version of Corollary \tref{CTLL} to classify singularities of multiplicity three (triple points) on rational plane curves of degree six and for that reason we call this result the Triple Lemma. We now are able to classify singularities of degree $d/2$ on rational plane curves of even degree $d$; see Theorem \tref{d=2c}.

Adopt Data \tref{34.1} with $n=3$, and $\pmb k$ an algebraically
closed field. Assume that $p$ is a singular point on the curve $\Cal C$ of multiplicity $d_j$, for $j$ equal to $1$ or $2$. In Theorem \tref{TPL1} we describe how to read the infinitely near singularities in the first neighborhood of $p$ from the Hilbert-Burch matrix $\varphi$ for the parameterization of $\Cal C$. Corollary \tref{ctgl} shows how to re-arrange the data so that $p$ becomes the point $[0:0:1]$ and $\varphi$ has the form described in Theorem \tref{TPL1}. 

\proclaim{Theorem \tnum{TPL1}}
Adopt Data \tref{34.1} with $n=3$, and $\pmb k$ an algebraically
closed field.  Assume that $p=[0:0:1]$ is a  singular point on the curve $\Cal C$ of multiplicity $d_j$, for $j$ equal to $1$ or $2$. Assume further that 
$$\varphi=\bmatrix P_1&Q_1\\P_2&Q_2\\0&Q_3\endbmatrix, $$ where the $P_i$ and the $Q_i$ are homogeneous forms from $B$ with $\deg P_i=d-d_j$ and  $\deg Q_i=d_j$.
Then,
a point $P= (p, [a:b]) \in {\Bbb P}^2 \times {\Bbb P}^1$ lies on the
blowup of $\Cal C$ at $p$ if and only if $\gcd(Q_3, aP_1+bP_2)$ is
not a constant. In this case, the multiplicity of $P$ on the blowup
is the degree of $\gcd(Q_3, aP_1+bP_2)$.
\endproclaim

\demo{Proof} Write $\Delta=P_1Q_2-P_2Q_1$. Notice that $\gcd(P_1, P_2)=1$
and $\gcd(Q_3, \Delta)=1$  since $I_2(\varphi)$ has height $2$.
Corollary \tref{ctgl} shows how to modify $Q_1$ and $Q_2$ in order to have
$\gcd (P_1,Q_1)=1$ and  $\gcd (P_2,Q_2)=1$. 
Passing to an affine chart we may assume that $p$ is the origin on
the affine curve parametrized by  $( \frac{g_1}{g_3},
\frac{g_2}{g_3} )$. The blowup of this curve at the origin has two
charts, parametrized by $(\frac{g_1}{g_3}, \frac{g_2}{g_1})$ and
$(\frac{g_2}{g_3}, \frac{g_1}{g_2})$, respectively. Homogenizing we
obtain two curves $\Cal C'$ and $\Cal C''$ in $\Bbb P^2$  parametrized by
$[g_1^2: g_2g_3: g_1g_3]$ and $[g_2^2: g_1g_3: g_2g_3]$,
respectively. After dividing by a common factor these
parameterizations become $$[P_2^2Q_3:-P_1\Delta:P_2\Delta]\quad\text{and}\quad
[P_1^2Q_3:P_2\Delta:-P_1\Delta],$$ respectively. One easily sees that the
entries of either vector generate an ideal of height $2$ and the
corresponding Hilbert Burch matrices are
$$
\bmatrix
             0 & -\Delta\\
             P_2 & 0 \\
             P_1 & P_2Q_3 \\
           \endbmatrix \quad\text{and}\quad
\bmatrix
             0 & \Delta \\
             P_1 & 0 \\
             P_2 & P_1Q_3 \\
           \endbmatrix.
\tag\tnum{Dsp}$$
Notice that the exceptional fiber of the blowup consists of  the
points on $\Cal C'$ or $\Cal C''$ that are of the form $[0:c:1]$. Observation \tref{O-1-18} shows that a point $[0:c:1]$ lies on $\Cal C'$ if and only
if $\gcd(Q_3, P_1+cP_2)$ is not a constant, and that it is on
$\Cal C''$ if and only if $\gcd(Q_3, cP_1+P_2)$ is not a constant.
Furthermore, Theorem \tref{Cor1} shows that the degree of this gcd gives the multiplicity of the
point on the blowup. \qed\enddemo

\proclaim{Corollary \tnum{CTL}}
Adopt Data \tref{34.1} with $n=3$ and $\pmb k$ an algebraically
closed field.  Then the  infinitely near singularities of $\Cal C$ have
multiplicity at most $d_1$.
\endproclaim
\demo{Proof} It suffices to see this for the multiplicity of the singular
points in the first neighborhood of a point $p$ on $\Cal C$.  If $m_p
\leq d_1$ then the assertion is clear.
Otherwise, Corollary \tref{ctgl} shows that $m_p =d_2$. Notice that in the setting of Theorem \tref{TPL1} the
form $aP_1+bP_2$ cannot be zero and hence $\deg(\gcd(Q_3,
aP_1+P_2))\leq d_1$. Now the present assertion follows from the
theorem. \qed \enddemo

\remark{Remark \tnum{R1.3}} Let  $q$ be a root of $Q_3$. Notice that the
point $$P=(p, [a:b])=(\Psi(q), [-P_2(q): P_1(q)])$$ is a point in the
blowup of $\Cal C$ at $p$ and conversely  each point in the blowup of
$\Cal C$ at $p$ is obtain in this manner.
\endremark

Corollary \tref{CTLL} is the ``Triple Lemma''; it concerns singularities of multiplicity $c=d/2$, on, or infinitely near, $\Cal C$, where $d$ is the degree of $\Cal C$. We first show that if $q$ is a singularity of multiplicity $c$ infinitely near to $\Cal C$, then $q$ is infinitely near to a singularity $p$ on $\Cal C$ of multiplicity $c$.   In particular, it is not possible for $c=m_{q}<m_p$. The most significant part of the result is assertion (4), where we study singular points of multiplicity $c$ which are in the second neighborhood of $\Cal C$. A complete description of all multiplicity $c$ singularities which appear on or infinitely near a rational plane curve $\Cal C$ of degree $d=2c$ appears in Theorem \tref{d=2c}.

\remark{Remark} Assertion (1) in the next result is a purely geometric statement. We wonder if there is a geometric argument for it.\endremark

\proclaim{Corollary \tnum{CTLL}}
Adopt Data \tref{34.1} with $n=3$, $\pmb k$ an algebraically
closed field, and $d$ equal to the even number $2c$. 

\roster
\item"{(1)}" If $p$ is a point on $\Cal C$ and $q$ is a singularity of multiplicity $c$ infinitely near to $p$,
then the multiplicity of $p$ is also equal to  $c$.   
\item"{(2)}" If $p$ is a singularity on $\Cal C$ of multiplicity $c$, then the parameters $d_1$ and $d_2$ are both equal to $c$.
\item"{(3)}" Suppose that $p$ is a singularity on $\Cal C$ of multiplicity $c$. Then there is a singularity of $\Cal C$, infinitely near to $p$,  of multiplicity $c$ if and only if there exist invertible matrices $\chi$ and $\xi$, with entries in $\pmb k$, such that  $(0,0,1)=p\chi$ and 
$$
\chi \varphi \xi=\bmatrix P_1&Q_1\\Q_3&Q_2\\0&Q_3\endbmatrix,$$where the homogeneous forms $P_1$ and  $Q_i$ all have degree $c$.
\item"{(4)}" Suppose that $p$ is a singularity on $\Cal C$ of multiplicity $c$. Then there  are two singularities of $\Cal C$, infinitely near to $p$, of multiplicity $c$  if and only if there exist invertible matrices $\chi$ and $\xi$, with entries in $\pmb k$, such that  $(0,0,1)=p\chi$ and 
$$
\chi \varphi \xi=\bmatrix Q_2&Q_1\\Q_3&Q_2\\0&Q_3\endbmatrix,$$where the homogeneous forms    $Q_i$ all have degree $c$.
\endroster
\endproclaim

\demo{Proof}First apply Corollary \tref{CTL}. If some infinitely near singularity of $\Cal C$ has multiplicity $c$ or higher, then $c\le d_1\le (1/2)d=c$. Assertions (1) and (2) have been established. Henceforth, we assume that $p$ is a singularity of multiplicity $m_p$ on $\Cal C$. Assertion (3) is an immediate consequence of Theorem \tref{TPL1}. We prove (4). In light of (3), we start with the singularity $p=[0:0:1]$ on the curve $\Cal C$ whose Hilbert-Burch matrix is equal to
$$\varphi'=\bmatrix P_1&Q_1\\Q_3&Q_2\\0&Q_3\endbmatrix,\tag\tnum{phi'}$$where the homogeneous forms $P_1$ and  $Q_i$ all have degree $c$ and
$$\gcd(Q_2,Q_3)=\gcd(P_1,Q_3)=\gcd(P_1,Q_1)=1.$$ Let $p'$   be the singular point of multiplicity $c$ in the first  neighborhood of $p$. Remark \tref{R1.3} shows that $p'$ is the point $(p,[0:1])$ on the blowup of $\Cal C$. Thus, in the language of the proof of Theorem \tref{TPL1}, $p'$ is on $\Cal C''$. Keep in mind that the polynomial $P_2$ of the original Hilbert-Burch matrix $\varphi$ has been replaced by $Q_3$ in the Hilbert-Burch matrix $\varphi'$ of (\tref{phi'}) and apply  column operations to the matrix on the right side of (\tref{Dsp}),  to see that one Hilbert-Burch matrix for the parameterization of $\Cal C''$ is
$$\bmatrix \Delta&0\\-P_1^2&P_1\\0&Q_3\endbmatrix.$$ 
Let $p''=(p',[a:b])\in \Bbb P^2\times \Bbb P^1$ be a point on the blowup of $\Cal C''$ at $p'$. 
Theorem \tref{TPL1} now shows that $p''$ is a singular point of multiplicity $c$ in the first neighborhood of $p'$ if and only if 
$$\deg\gcd(Q_3,a\Delta-bP_1^2)=c.$$
Recall that $\Delta=P_1Q_2-Q_1Q_3$, $\gcd(P_1,Q_3)=1$, and $Q_2$, $Q_3$, and $P_1$ all are homogeneous forms of   degree $c$. It follows that $p''$ is a singular point of multiplicity $c$ in the first neighborhood of $p'$ if and only if $$\gamma Q_3=aQ_2-bP_1\tag\tnum{bnot0}$$ for some constant $\gamma\in \pmb k$. The equation (\tref{bnot0}) cannot hold if $b=0$, since $[a:b]$ is a point in $\Bbb P^1$ and $Q_2$ and $Q_3$ are relatively prime. Therefore, 
$$\hskip-3pt\matrix\format\l\\ 
\phantom{{}\iff {}}\text{there exists a singular point of multiplicity $c$ in the first neighborhood of $p'$}\\ 
{}\iff \text {$P_1$ is in the vector space generated by $Q_2$ and $Q_3$}\\
{}\iff \text {there exist invertible matrices $\chi'$ and $\xi'$ so that $\chi \varphi'\xi$ has the form of the}\\\phantom{{}\iff{}} \text {matrix of (4).} \qed\endmatrix$$
   \enddemo

\bigpagebreak
\SectionNumber=3\tNumber=1
\heading Section \number\SectionNumber. \quad The BiProj Lemma. 
\endheading

Let $\varphi$ be a Hilbert-Burch matrix which corresponds to a parameterized plane curve $\Cal C$ of  degree $d$.
 Points of $\Bbb P^2$ give rise to generalized row ideals of the matrix $\varphi$. 
Thus, features of the generalized row ideals reflect properties of the corresponding points.
For example, a generalized row ideal of $\varphi$ encodes information that can be used to determine if the corresponding point $p$ is on $\Cal C$  and, if so,  what type of singularity occurs at $p$.

In this section, we focus on the situation where the degree of $\Cal C$ is even (so $d=2c$) and we describe the singular points on, or infinitely near,  $\Cal C$ of multiplicity $c$.
We know from Corollary \tref{CTLL} that such a point exists if and only if all entries of $\varphi$ have the same degree $c$ and  the corresponding  generalized row has a generalized zero.  This leads us to consider column operations on $\varphi$, which we     identify   with points in $\Bbb P^1$. Thus, inside $\Bbb P^2\times \Bbb P^1$ we consider the closed subset consisting of pairs $$\text{(row operation, column operation)}$$  that lead to a generalized zero of $\varphi$. Projection onto the first factor gives the singular points of multiplicity $c$. On the other hand, projection onto the second factor yields  a finite set of points in $\Bbb P^1$ that is easier to study yet reflects properties of the set of singular points of multiplicity $c$ on the plane curve $\Cal C$.

Fix a Hilbert-Burch matrix $\varphi$ in which every entry is a homogeneous form of degree $c$.
To find a singularity of multiplicity c on $\Cal C$ we need to describe a generalized zero of $\varphi$. In other words, we look for $(p,q)$ in $\Bbb P^2\times \Bbb P^1$ such that $p\varphi q^{\text{\rm T}}=0$. Consider the polynomial 
 $\pmb T\varphi \pmb u^{\text{\rm T}}\in \pmb k[\pmb T,\pmb u,x,y]$, where $\pmb T=[T_1,T_2,T_3]$ and $\pmb u=[u_1,u_2]$ are matrices of indeterminates. We extract the variables $x$ and $y$ from the critical polynomial  $\pmb T\varphi \pmb u^{\text{\rm T}}$. For each positive integer $i$, let $\rho^{(i)}$ be the following $1\times (i+1)$ row vector of all monomials of degree $i$ in $\pmb k[x,y]$: $$\rho^{(i)}=[y^i,xy^{i-1},x^2y^{i-2},\dots,x^{i-1}y,x^i].\tag\tnum{rho}$$ Define $C$ and $A$ to be the matrices with
$$\pmb T \varphi=\rho^{(c)}C\quad\text{and}\quad C\pmb u^{\text{\rm T}}=A\pmb T^{\text{\rm T}},$$
so that the entries of $C$ are linear forms in $\pmb k[\pmb T]$ and the entries   of $A$ are linear forms in $\pmb k[\pmb u]$. One now has 
$$\pmb T\varphi \pmb u^{\text{\rm T}}= \rho^{(c)} C\pmb u^{\text{\rm T}}= \rho^{(c)} A\pmb T^{\text{\rm T}}.$$
Thus, 
$$\matrix \text {the set of $(p,q)$ in $\Bbb P^2\times \Bbb P^1$ such that $p\varphi q^{\text{\rm T}}=0$ is 
the zero set, in $\Bbb P^2 \times \Bbb P^1$,}\hfill\\\text{of the bihomogeneous ideal $I_1(C\pmb u^{\text{\rm T}})=I_1(A\pmb T^{\text{\rm T}})$.}\hfill\endmatrix \tag\tnum{mot}$$
In Theorem \tref{XXX}, we obtain two isomorphisms of schemes
$$ \xymatrix{ &\operatorname{BiProj}\left (\pmb k[\pmb T,\pmb u]/I_1(C\pmb u^{\text{\rm T}})\right )\ar[ld]_{\pi_1}^{\cong}\ar[rd]_{\cong}^{\pi_2}&\\
\operatorname{Proj}\left (\pmb k[\pmb T]/I_2(C)\right ) &&\operatorname{Proj}\left (\pmb k[\pmb u]/I_3(A)\right )} \tag\tnum{pict}$$
and we exploit these isomorphisms in Theorem \tref{TBA} to describe the singularities on $\Cal C$ of multiplicity $c$. The equation $$C\pmb u^{\text{\rm T}}=A\pmb T^{\text{\rm T}}\tag\tnum{CE}$$ provides symmetry. Theorem \tref{bipro} is used twice to produce the isomorphisms of (\tref{pict})

\proclaim{Theorem \tnum{bipro}} Let $S=\pmb k[x_1,\dots,x_m, y_1,\dots,y_n]$ be a bi-graded polynomial ring with $\deg x_i=(1,0)$ and $\deg y_i=(0,1)$, and $R$ be the sub-algebra $\pmb k[x_1,\dots,x_m]$ of $S$. Let $J$ be an  $S$-ideal   generated by bi-homogeneous forms which are linear in the $y$'s.  Write $J=I_1(\phi \pmb y)$ where $\pmb y= [y_1,\dots,y_n]^{\text{\rm T}} $ and $\phi$ is a matrix with entries in $R$. The entries in each row of $\phi$ are homogeneous of the same degree. Consider the natural projection map $\pi\:\operatorname{BiProj}(S/J)\to \operatorname{Proj} (R)$. If the ideal $I_{n-1}(\phi)$ is zero-dimensional in $R$, then $\pi$ is an isomorphism onto its image and this image is defined scheme-theoretically by the $R$-ideal $I_n(\phi)$. \endproclaim

Notice  that $\operatorname{im} \pi=\operatorname {Proj} \left (R/\left (\vphantom{E^{E^{E}}_{E_{E}}} J\:\! (I_1(\pmb x)I_1(\pmb y))^{\infty}\right )\cap R\right )\subseteq \operatorname{Proj} (R)$. The theorem means that $\pi$ gives a bijection 
$\operatorname{BiProj}(S/J)\to \operatorname{Proj} (R/I_n(\phi))$ which induces  isomorphisms   at the level of local rings.

\demo{Proof} 
We prove the result after localizing at any fixed one-dimensional homogeneous prime ideal $\frak m$ of $R$. The ideal $I_{n-1}(\phi)R_\frak m$ is the unit ideal by hypothesis; therefore,  there are invertible matrices $U$ and $V$ with $$U\phi V=\bmatrix I_{n-1}&0\\0&\phi'\endbmatrix$$ for some column vector  $\phi'$ with entries in $R_\frak m$. Write $S_\frak m$ to mean $(R\setminus \frak m)^{-1}S$. We have $S_\frak m=R_\frak m[y_1',\dots,y_n']$, where $\pmb y'=V^{-1}\pmb y$. The ideal $JS_\frak m$ is equal to $I_1(\phi\pmb y)S_{\frak m}=(y_1',\dots,y_{n-1}',I_1(\phi')y_n')$.

The ring map $S/J\to S_\frak m/JS_\frak m$ induces the inclusion  $$\operatorname{Proj} (S_\frak m/JS_\frak m)\subseteq \operatorname{BiProj}(S/J),$$ and the   natural projection map $\pi\:\operatorname{BiProj}(S/J)\to \operatorname{Proj} (R)$ restricts to become 
$$\operatorname{Proj}(S_\frak m/JS_\frak m)@> \pi|>> \operatorname{Spec} (R_\frak m).$$  The image of $\pi|$ is the subscheme $\operatorname{Spec} (\frac{R_\frak m}{(JS_\frak m\:I_1(\pmb y)^{\infty})\cap R_\frak m})$ of $\operatorname{Proj} (R)$. We must show that

\flushpar{\bf (1)} $\pi|$ is an isomorphism onto its image, and 

\flushpar{\bf (2)} the image of $\pi|$ is defined scheme-theoretically by the ideal $I_n(\phi)R_\frak m$. 

Notice  that 
$$\split &JS_\frak m\:I_1(\pmb y)^{\infty}= \left (y_1',\dots,y_{n-1}',I_1(\phi')y_n'\right )S_\frak m\:\left (y_1',\dots,y_{n-1}',y_n'\right )^{\infty}\\&{}= \left (y_1',\dots,y_{n-1}',I_1(\phi')y_n'\right )S_\frak m\:(y_n')^{\infty}= \left (y_1',\dots,y_{n-1}',I_1(\phi')\right )S_\frak m.\endsplit $$
Therefore, the image of $\pi|$ is defined scheme-theoretically by
$$(JS_\frak m\:I_1(\pmb y)^{\infty})\cap R_\frak m= I_1(\phi')R_\frak m=I_n(\phi)R_\frak m,$$ which proves the second assertion. To show the first claim, notice that the natural map 
$$R_\frak m\to \frac{R_\frak m[y_n']}{I_1(\phi')y_n'R_\frak m[y_n']}=\frac{S_\frak m}{JS_\frak m}$$
induces $\pi|\: \operatorname{Proj}(S_\frak m/JS_\frak m)\to \operatorname{Spec} (R_\frak m)$; furthermore,
$$\split &{}\operatorname{Proj}\left (\frac{S_\frak m}{JS_\frak m}\right )=\operatorname{Proj} \left (\frac{R_\frak m[y_n']}{I_1(\phi')y_n'R_\frak m[y_n']}\right )=\operatorname{Proj} \left (\frac{R_\frak m[y_n']}{I_1(\phi')y_n'R_\frak m[y_n']\:(y_n')^{\infty}}\right )\\&{}=\operatorname{Proj} \left (\frac{R_\frak m[y_n']}{I_1(\phi')R_\frak m[y_n']}\right )
=\operatorname{Proj} \left (\frac{R_\frak m}{I_1(\phi')R_\frak m}[y_n']\right )=\operatorname{Spec} \left (\frac{R_\frak m}{I_1(\phi')R_\frak m}\right )\\&{}=\operatorname{Spec} \left (\frac{R_\frak m}{I_n(\phi)R_\frak m}\right ). \qed \endsplit$$\enddemo

The next Lemma analyzes a general version of (\tref{CE}). This Lemma is a step toward verifying that the hypothesis of Theorem \tref{bipro} are satisfied twice in (\tref{pict}). It is fortuitous that the information about the minors of $A$ translates into information about the columns of $C$ and vice versa. 
 In this result $C$ is a matrix with $s$ columns and $A$ is a matrix with entries from $\pmb k[u_1,\dots,u_s]$. 

\remark{Remark \tnum{301}} If $\lambda\in \Bbb A^s$, then we abbreviate the generalized column $C\lambda^{\text{\rm T}}$ of $C$  as $C_{\lambda}$; also, $A(\lambda)$ is the matrix of constants which is obtained by evaluating each entry of $A$ at the point $\lambda$.\endremark

\proclaim{Lemma \tnum{1-8}} Let $\pmb u=[u_1,\dots,u_s]$ and  $\pmb T=[T_1,\dots,T_t]$ be matrices of indeterminates  over an algebraically closed field $\pmb k$, 
$C$ be an ${n\times s}$ matrix of linear forms from $\pmb k[\pmb T]$, and  $A$ be an $n\times t$ matrix of linear forms from $\pmb k[\pmb u]$. If 
$C\pmb u^{\text{\rm T}} = A \pmb T^{\text{\rm T}}$, then the following statements hold:\phantom{\tnum{topf},\tnum{botf}}
$$\alignat2
I_i(A)\neq 0&{}\iff \operatorname{ht} I_1(C_{\lambda})\ge i \text{ for     general }\lambda\in \Bbb A^{s}, \text{ and}\tag{\tref{topf}}\\
\operatorname{ht} I_i(A)\ge s &{}\iff \operatorname{ht} I_1(C_{\lambda})\ge i\text{  for all non-zero }\lambda \in \Bbb A^{s}.\tag{\tref{botf}}\endalignat $$
 \endproclaim
\demo{Proof} We have $C_\lambda=C\lambda^{\text{\rm T}} =A(\lambda)\pmb T^{\text{\rm T}}$.
The right hand side of (\tref{topf}) holds if and only if
$$\split &\text{$\operatorname{rank} A(\lambda)\ge i$\ for  general $\lambda \in \Bbb A^{s}$}\\
\iff {}& \text{$I_i(A(\lambda))\neq 0$\ for  general $\lambda \in \Bbb A^{s}\iff I_i(A)\neq 0$}.\endsplit$$
The right side of (\tref{botf}) holds if and only if 
$$\split &\text{$\operatorname{rank} A(\lambda)\ge i$\ for all non-zero  $\lambda \in \Bbb A^{s}$}\\
\iff {}& \text{$I_i(A(\lambda))\neq 0$\ for all non-zero $\lambda \in \Bbb A^{s}\iff \operatorname{ht} I_i(A)\ge s$.}\endsplit$$
The last equivalence  is due to Hilbert's Nullstellensatz   which applies since $\pmb k$ is algebraically closed. \qed
\enddemo   

\definition{Data \tnum{biData}} Adopt Data \tref{34.1} with $n=3$, $d=2c$, $d_1=d_2=c$, and $\pmb k$ an algebraically closed field. Let $\pmb T$ and $\pmb u$ be the row vectors $\pmb T=[T_1,T_2,T_3]$ and $\pmb u=[u_1,u_2]$ of indeterminates. Recall the row vector $\rho^{(c)}$ from (\tref{rho}). Define the matrices
$C$ and $A$ by\phantom{\tnum{no1}\tnum{no2}}
$$\alignat2 \pmb T \varphi&{}={}\rho^{(c)}C\ \text{and} \tag{\tref{no1}}\\C\pmb u^{\text{\rm T}}&{}={}A\pmb T^{\text{\rm T}},\tag{\tref{no2}}\endalignat$$
so that the entries of $C$ are linear forms in $\pmb k[\pmb T]$ and the entries entries of $A$ are linear forms in $\pmb k[\pmb u]$. \enddefinition 

\proclaim{Lemma \tnum{1-8.1}} Adopt Data \tref{biData}.  The following statements hold{\rm:}
\roster
\item  $\operatorname{ht} I_1(C)=3$,
\item $\operatorname{ht} I_1(C_\lambda)\ge 2$ for all non-zero $\lambda\in\Bbb A^2$, and 
\item $\operatorname{ht} I_2(A)=2$.\endroster\endproclaim

\demo{Proof} (1) The height of $I_1(C)$  is three otherwise, since the entries of $C$ are linear, $I_1(C)$ is contained in an ideal generated by two linear  forms. Equation (\tref{no1}) shows that, after row operations, $\varphi$ has a row of zeros yielding $\operatorname{ht}(I)=1$. Now we prove (2). Equation (\tref{no1}) gives $\rho^{(c)}C_{\lambda}=\pmb T\varphi_{\lambda}$. If $\operatorname{ht} (I_1(C_{\lambda}))\le 1$, then $\mu(I_1(C_{\lambda}))\le 1$ because the entries of $C$ are linear. Thus, $\mu(\pmb T\varphi_{\lambda}))\le 1$, which shows that after row operations $\varphi_{\lambda}$ has at most one non-zero entry. This would imply that  $\operatorname{ht}(I)\le 1$. Assertion (3) follows from (2) and Lemma \tref{1-8}. 
\qed
\enddemo

\proclaim{Theorem \tnum{XXX}} Adopt Data \tref{biData}. The following statements hold.
\roster
\item"{(1)}" The projections
$$ \xymatrix{ &\operatorname{BiProj}\left (\pmb k[\pmb T,\pmb u]\right )\ar[ld]\ar[rd]&\\
\operatorname{Proj}\left (\pmb k[\pmb T]\right ) &&\operatorname{Proj}\left (\pmb k[\pmb u]\right )} $$
induce isomorphisms
$$ \xymatrix{ &\operatorname{BiProj}\left (\pmb k[\pmb T,\pmb u]/I_1(C\pmb u^{\text{\rm T}})\right )\ar[ld]_{\pi_1}^{\cong}\ar[rd]_{\cong}^{\pi_2}&\\
\operatorname{Proj}\left (\pmb k[\pmb T]/I_2(C)\right ) &&\operatorname{Proj}\left (\pmb k[\pmb u]/I_3(A)\right );} $$
 in particular, 
the schemes $\operatorname{Proj}(\frac {\pmb k[\pmb T]}{I_2(C)})$ and $\operatorname{Proj}(\frac {\pmb k[\pmb u]}{I_3(A)})$ are isomorphic.

\item"{(2)}" As a subset of $\Bbb P^2$,  $\operatorname{Proj}(\frac {\pmb k[\pmb T]}{I_2(C)})$ is equal to $\{p\in \Cal C\mid m_p=c\}$.

\endroster
\endproclaim

\demo{Proof}  We apply Theorem \tref{bipro} twice. Each time $S=\pmb k[\pmb T,\pmb u]$ and $J=I_1(C\pmb u^{\text{\rm T}})=I_1(A\pmb T^{\text{\rm T}})$. In the first application $R=\pmb k[\pmb T]$.  Lemma \tref{1-8.1} ensures that $I_1(C)$ is zero-dimensional in $\pmb k[\pmb T]$.
In the second application $R=\pmb k[\pmb u]$. Again Lemma \tref{1-8.1} ensures that $I_2(A)$ is zero-dimensional in $\pmb k[\pmb u]$. Assertion (1) is established.

Return to the first setting. Theorem \tref{bipro} also yields that the image of the map \phantom{\tnum{t1}\tnum{t2}}
$$\Bbb P^2\times \Bbb P^1 \supseteq \operatorname{BiProj}(S/J)@> \pi >> \operatorname{Proj}(\pmb k[\pmb T])$$
is $\operatorname{Proj} (\pmb k[\pmb T]/I_2(C))\subseteq \operatorname{Proj} (\pmb k[\pmb T])$. On the other hand, as a set $$\alignat 2\operatorname{im} \pi&{}=\{p\in \Bbb P^2\mid \exists q\in \Bbb P^1 \text{ with } (p,q)\in V(J)\}\\
&{}=\{p\in \Bbb P^2\mid \exists q\in \Bbb P^1 \text{ with } p\varphi q^{\text{\rm T}} =0\}\tag{\tref{t1}}\\
&{}=\{p\in \Cal C\mid m_p=c\}.\tag{\tref{t2}}\endalignat
$$ 
The equality (\tref{t1}) is explained   in (\tref{mot}) and the equality (\tref{t2}) is established in Corollary \tref{ctgl}. \qed
\enddemo
Theorem \tref{XXX} produces a morphism (even an isomorphism) 
$$\pi_2\circ \pi_1^{-1}\:\operatorname{Proj} (\pmb k[\pmb T]/I_2(C)) \to \operatorname{Proj}\left (\pmb k[\pmb u]/I_3(A)\right );$$the graph of this morphism is 
$\operatorname{BiProj} (\pmb k[\pmb u,\pmb T]/I_1(C\pmb u^{\text{\rm T}}))$, as can be seen from diagram (\tref{pict}). 

\remark{Remark \tnum{CCC}}The theorem also    says, in particular, that if  $p$ and $q$ are non-zero row vectors in $\pmb k^3$ and $\pmb k^2$, respectively, with $p\varphi q^{\text{\rm T}}=0$, then, up to multiplication by non-zero scalars, $q$ is determined by $p$ and vice versa. \endremark

We next isolate one quick result about how the Data of \tref{biData} can be used to determine the configuration of multiplicity $c$ singularities on a curve $\Cal C$. The complete story is contained in Theorem \tref{TBA}.

\proclaim{Corollary \tnum{mt}} Adopt Data \tref{biData}. The curve $\Cal C$ has only singularities of multiplicity at most $c-1$ if and only if every generalized column ideal of $C$ has height three.\endproclaim
\demo{Proof} All singularities of $\Cal C$ have  multiplicity at most $c-1$ if and only if $\Cal C$ has no singularity of multiplicity $c$ by Corollary \tref{ctgl}. This in turn means $\operatorname{ht} I_3(A)=2$ according to Theorem \tref{XXX}. Finally, by Lemma \tref{1-8}, $\operatorname{ht} I_3(A)=2$ if and only if every generalized column ideal of $C$ has height three. \qed \enddemo

\proclaim{Corollary \tnum{blah}}Adopt Data \tref{biData} with $c\ge 2$.  
The following statements hold:
\roster
\item  $\operatorname{ht} I_2(C)\ge 2$,
\item $\operatorname{ht} I_3(A)\ge 1$, and
\item $\operatorname{ht} I_1(C_\lambda)=3$ for general $\lambda\in\Bbb A^2$ .\endroster\endproclaim

\demo{Proof}Theorem \tref{XXX} shows that $\operatorname{Proj}(\pmb k[\pmb T]/I_2(C))$ is either empty or is a finite set; therefore, (1) and (2) hold. Assertion (3) follows from (2) and Lemma \tref{1-8}. \qed\enddemo

\proclaim{Lemma \tnum{e+mu}} Let $C$ be  matrix of linear forms from a polynomial ring $R$ in three variables over a field $\pmb k$.   If $C$ has $2$ columns,  some generalized column ideal of $C$ is a zero-dimensional ideal of $R$, and   $\operatorname{ht}(I_2(C))=2$, then 
$$e(R/I_2(C))+\mu(I_2(C))=6.$$\endproclaim

\demo{Proof}  Let $R=\pmb k[T_1,T_2,T_3]$. 
Since $I_2(C)$ is an ideal generated by quadrics in a polynomial in three variables, 
 $\mu(I_2(C))\le 5$; for otherwise, $I_2(C)=(T_1,T_2,T_3)^2$, which would contradict the hypothesis that $\operatorname{ht} I_2(C)=2$.

We are allowed to apply 
 row and column operations, to perform linear changes of variables, and to suppress zero rows. Using these operations, we transform   $C$ into one of the following four forms:
$$M_1= \bmatrix T_1&*\\T_2&*\\T_3&*\endbmatrix,\quad M_2= \bmatrix T_1&f_1\\T_2&f_2\\T_3&f_3\\0&T_1\endbmatrix,\quad 
M_3= \bmatrix T_1&g_1\\T_2&g_2\\T_3&g_3\\0&T_1\\0&T_2\endbmatrix,\quad M_4= \bmatrix T_1&0\\T_2&0\\T_3&0\\0&T_1\\0&T_2\\0&T_3\endbmatrix,
$$
where the $f_i$ are linear forms in $\pmb k[T_2,T_3]$ and the $g_i$ are linear forms in $R=\pmb k[T_3]$.

 In the case $C=M_1$, the ideal $I_2(C)$ is a perfect ideal of height $2$ with a $2$-linear resolution; hence $e(R/I_2(C))=3=\mu(I_2(C))$. 

In the case $C=M_2$, we have $I_2(C)=T_1(T_1,T_2,T_3)+JR$, where $J$ is a non-zero ideal of $\pmb k[T_2,T_3]$ generated by quadrics. Write $\Delta=\gcd(J)\in \pmb k[T_2,T_3]$. Notice that $T_1$ is in $I_2(C)^{\text{unm}}$, the unmixed part of $I_2(C)$. Therefore, $I_2(C)^{\text{unm}}=(T_1,J)^{\text{unm}}=(T_1,\Delta)$, and $e(R/I_2(C))=e(R/I_2(C)^{\text{unm}})=\deg \Delta\le 2$.  Observe that $\deg \Delta=2$ if and only if $\mu(J)=1$;
thus,  $e(R/I_2(C))=2$ if and only if $\mu(I_2(C))=4$. It follows that $e(R/I_2(C))=1$ if and only if $\mu(I_2(C))\ge 5$, which means that $\mu(I_2(C))= 5$.

In the case $C=M_3$, we have $(T_1,T_2)(T_1,T_2,T_3)\subseteq I_2(C)$ and these two ideals are equal because 
$I_2(C)$ is generated by at most $5$ quadrics. Therefore, $\mu(I_2(C))= 5$ and $e(R/I_2(C))=e(R/(T_1,T_2))=1$.

Finally, if $C=M_4$, then $(T_1,T_2,T_3)^2= I_2(C)$, which is impossible since  $I_2(C)$ has height two.  
\qed
\enddemo

Before proceeding, 
we describe the effect on Data \tref{biData} of applying row and column operations to $\varphi$. 
\remark{Remark \tnum{305}}If we replace $\varphi$ by $\varphi'=\chi \varphi$ for some invertible matrix $\chi$ with entries in $\pmb k$, then the matrices $C$ and $A$ of Data \tref{biData} become
$C'=C(\pmb T\chi)$ and $A'=A\chi^{\text{\rm T}}$. Recall that changing $\varphi$ by $\varphi'$ amounts to  applying to the curve $\Cal C$ the linear automorphism of $\Bbb P^2$ defined by the matrix $\chi^{-1}$; see Remark \tref{sigh}.  If we replace $\varphi$ by $\varphi'=\varphi\xi$ for some invertible matrix $\xi$ with entries in $\pmb k$, then the matrices $C$ and $A$  become
$C'= C\xi$ and $A'=A(\pmb u\xi^{\text{\rm T}})$, for $A(\lambda)$ as defined in Remark \tref{301}. Notice that column operations on $\varphi$ have no effect on the $g$'s or on the curve $\Cal C$.\endremark

Ultimately, in Corollary \tref{CorND}, we use the Jacobian ideal of $\pmb k[\pmb u]/I_3(A)$ over $\pmb k$ to count the number of infinitely near singularities of multiplicity $c$. Section 7 contains much information about this object. In the mean time, we recall that if $D$ is an algebra essentially of finite type over a Noetherian ring $L$, then  
the module of  K\"ahler differentials of $D$ over $L$ is  denoted by $\Omega_{D/L}$. The zeroth Fitting ideal of the {$D$-module} $\Omega_{D/L}$ is called the Jacobian ideal of $D$ over $L$ and is denoted by $\operatorname{Jac}(D/L)$. The Jacobian ideal $\operatorname{Jac}(D/L)$ is an ideal of $D$.  These objects may be computed using any presentation of the $L$-algebra $D$ and their formation is compatible with localization. Notice that the subsets $\operatorname{Supp}(\Omega_{D/L})$ and $V(\operatorname{Jac}(D/L))$ of $\operatorname{Spec} D$ are equal. In Section  7 we will use the fact that this common subset of $\operatorname{Spec}(D)$ is the ramification locus of $D$ over $L$.

In \cite{\rref{KPU-08}} we studied  Rees algebras of ideals generated by $3$ forms of degree $6$ in $\pmb k[x,y]$. We found that there are exactly seven families of Rees algebras for such ideals. Each member of a given family has the same graded Betti numbers. Four of these family arise when  every entry in the syzygy matrix has degree 3. The distinguishing invariant turned out to be $\mu (I_2(C))$. 
Ad hoc methods  led to this classification. We are very pleased  to now see  in Theorem \tref{TBA} that this invariant arises naturally and reflects the geometry of the curve.

\proclaim{Theorem \tnum{TBA}} Adopt Data \tref{biData}.

\roster 

\item The ideal $I_2(C)$ is zero-dimensional  if and only if   $I_3(A)$ is zero-dimensional if and only if $ \gcd I_3 (A)$ is a unit; otherwise
$$e(\pmb k[\pmb T]/I_2(C))=e(\pmb k[\pmb u]/I_3(A))=\deg\gcd I_3(A)=6-\mu(I_2(C)).$$

\item The non-associate linear factors of $\gcd I_3(A)$ correspond to the distinct singular points on $\Cal C$ of multiplicity c.

\item Write $\gcd(I_3(A))=\prod \ell_i^{e_i}$, where the $\ell_i$ are non-associate linear factors and $e_i\ge 1$. 
Then $e_i-1$ is the number of singular points of multiplicity $c$ infinitely near to the point on $\Cal C$ corresponding to $\ell_i$. 

\item
The degree of  $\gcd I_3(A)$ is equal to the number of distinct singular points of multiplicity  $c$ that are either on $\Cal C$ or infinitely near to $\Cal C$.

\endroster
\endproclaim 

Before proving Theorem \tref{TBA}, we reformulate part of it in the exact form we will use in Section 6 where we 
decompose the space of parameterizations  of rational plane curves of even degree $d=2c$  into strata which reflect   the configuration of multiplicity $c$ singularities on or infinitely near the corresponding curve.

 \proclaim{Corollary \tnum{CorND}} Retain the notation of Theorem \tref{TBA} and let $D=\pmb k[\pmb u]/I_3(A)$.
For assertions {\rm(3b)} and {\rm(3c)} assume that $\operatorname{char} \pmb k$ is not equal to $2$ or $3$.  The following assertions hold.
\roster
\item The dimension of $D$ is either zero or one.
\item There are no singularities of multiplicity $c$ on or infinitely near $\Cal C$ if and only if $\dim D=0$.
\item If $D$ has   dimension one, then 
\itemitem{\rm(a)} $e(D)$ is the number of distinct singular points of multiplicity  $c$ that are either on $\Cal C$ or infinitely near to $\Cal C$,  
\itemitem {\rm(b)} if $D/\operatorname{Jac} (D/\pmb k)$ is a ring of dimension zero then all of the multiplicity $c$ singularities on or infinitely near $\Cal C$ are actually on $\Cal C$, and 
\itemitem{\rm(c)} if $D/\operatorname{Jac} (D/\pmb k)$ is a ring of dimension one,  $e(D)- e(D/\operatorname{Jac} (D/\pmb k))$ is equal to the number of distinct singular points of multiplicity  $c$ that are   on $\Cal C$.
 \endroster
\endproclaim

\demo{Proof of Corollary \tref{CorND}} 
Assertion (2) follows from items (2) of Theorem \tref{TBA} and (1) of Corollary \tref{CTLL}. Now we prove (3b) and (3c). Observe that $\deg \gcd I_3(A)\le 3$; hence each $e_i\le 3$ and $e_i\neq \operatorname{char} \pmb k$. If $\dim D=1$, then the minimal primes of $D$ are $\ell_iD$. Notice that $D_{\ell_iD}\cong \pmb k[\pmb u]_{\ell_i\pmb k[\pmb u]}/(\ell_i^{e_i})$. Hence,
$$\operatorname{Jac}(D/\pmb k)_{\ell_iD}=\operatorname{Jac}(D_{\ell_i}/\pmb k)=\ell_i^{e_i-1}D_{\ell_i D}.$$
It follows that

$$\split \dim D/\operatorname{Jac} (D/\pmb k)=0&{}\iff \operatorname{Jac} (D/\pmb k)=D\\&{}\iff \text{the roots of $\gcd I_3(A)$ are distinct}\endsplit$$ and
$$e(D/\operatorname{Jac} (D/\pmb k))=\sum (e_i-1)$$ and this is equal to 
the degree of $\deg\gcd I_3(A)$ minus the number of distinct linear factors of $\gcd I_3(A)$. \qed \enddemo

\demo{Proof of Theorem \tref{TBA}} The schemes  $\operatorname{Proj}(  {\pmb k[\pmb T]}/{I_2(C)})$ and $\operatorname{Proj}(  {\pmb k[\pmb u]}/{I_3(A)})$ are isomorphic  and are either empty or zero-dimensional by Theorem \tref{XXX} and Corollary  \tref{blah}. It follows that the rings ${\pmb k[\pmb T]}/{I_2(C)}$ and ${\pmb k[\pmb u]}/{I_3(A)}$ have the same dimension and this common dimension is either $0$ or $1$; furthermore, these two rings have the same multiplicity. The height of  $I_3(A)$ is $2$ if and only if the  $\gcd$ of $I_3(A)$ is a unit. Otherwise,  $I_3(A)$ has height $1$ and  $(I_3(A))^{\text{unm}}= (\gcd I_3(A))$. In this case,    $$e( {\pmb k[\pmb u]}/{I_3(A)})=e({\pmb k[\pmb u]}/{(\gcd I_3(A))})= \deg \gcd I_3(A).$$ The equality $e(\pmb k[\pmb T]/I_2(C))=6-\mu(I_2(C))$, when $\operatorname{ht} I_2(C)=2$, is proven in Lemma \tref{e+mu}, because the hypothesis that some generalized column ideal of $C$ is zero-dimensional is established in part (3)  of Corollary \tref{blah}. This completes the proof of (1).

Assertion (2) follows from Theorem \tref{XXX}: if  $\ell$ is a linear factor of $\gcd I_3(A)$, then $V(\ell)$  is a point in $\operatorname{Proj}(\pmb k[\pmb u]/I_3(A))$ and the corresponding singular point on $\Cal C$ of multiplicity $c$ is $\pi_1\circ \pi_2^{-1}(V(\ell))$ as described in Remark \tref{CCC}.

Assertion (4) is an immediate consequence of (2) and (3). We now prove (3).  Consider the factor $\ell_1$ of $\gcd I_3(A)$. Let $p$ be the point on $\Cal C$ which corresponds to $\ell_1$ as described in  the proof of (2). Now apply a linear automorphism of $\Bbb P^2$  to move the point $p$ to $(0,0,1)$. According to the General Lemma, the entries of the bottom row of the new $\varphi$ are linearly dependent. 
Now apply a column operation to make the bottom row become $[0,*]$. Under the correspondence of Remark \tref{CCC} the point $(0,0,1)$ is paired with $(1,0)$ and henceforth    $\ell_1$ has moved to  $u_2$. Furthermore,   $\varphi$ has become
$$\varphi =\bmatrix Q_1&Q_3\\Q_2&Q_4\\0&Q_5\endbmatrix,$$ where the $Q_i$ are forms of degree $c$. The fact that $\operatorname{ht}(I)=2$ forces $Q_1$ and $Q_2$ to be linearly independent over $\pmb k$. Write $Q_i=\rho^{(c)}\gamma_i$, where $\gamma_i$ is a column vector with entries in $\pmb k$. A straightforward calculation yields $$A=\bmatrix u_1\gamma_1+u_2\gamma_3& u_1\gamma_2+u_2\gamma_4&u_2\gamma_5\endbmatrix.$$It follows that $$I_3(A)=u_2I_3(A'),$$for 
$$A'=\bmatrix u_1\gamma_1+u_2\gamma_3& u_1\gamma_2+u_2\gamma_4&\gamma_5\endbmatrix.$$
Notice that $$I_3(A')\cong u_1^2 I_3(\bmatrix \gamma_1&\gamma_2&\gamma_5\endbmatrix)\mod u_2;$$hence,   $u_2$ divides every 
element of $I_3(A')$ if and only if $\gamma_1,\gamma_2,\gamma_5$ are linearly dependent. 

The fact that $\gamma_1$ and $\gamma_2$ are linearly independent guarantees that
$$\split &\text{$\gamma_1,\gamma_2,\gamma_5$ are linearly dependent} \\{}\iff& \text{$\gamma_5$ is in the vector space spanned by $\gamma_1$ and $\gamma_2$}\\{}\iff& \text{$Q_5$ is in the vector space spanned by $Q_1$ and $Q_2$.}\endsplit$$
 By the Triple Lemma \tref{TPL1} this is equivalent to the  existence of at least one singular point  of multiplicity $c$   infinitely near to  $(0,0,1)$. Thus such a point exists if and only if $e_1\ge 2$. 

Assume now that there does exist a singular point  of multiplicity $c$   infinitely near to  $(0,0,1)$; in particular $Q_5$ is in the span of $Q_1$ and $Q_2$ and therefore applying row operations on $\varphi$ only involving the first two rows we may assume that $\gamma_5=\gamma_2$. Notice that $I_3(A)$ is unchanged. We have
$I_3(A)=u_2^2I_3(A'')$ for $A''= \bmatrix u_1\gamma_1+u_2\gamma_3& \gamma_4&\gamma_2\endbmatrix$. Now one sees as before that 
$u_2$ divides every element of $I_3(A'')$ if and only if $\gamma_1,\gamma_4,\gamma_2$ are linearly dependent. Therefore,
$u_2^3$ divides every element of $I_3(A)$ if and  only if $\gamma_1,\gamma_4,\gamma_2$ are linearly dependent. The fact that $\gamma_1$ and $\gamma_2$ are linearly independent guarantees that
$$\split &\text{$\gamma_1,\gamma_2,\gamma_4$ are linearly dependent} \\{}\iff& \text{$\gamma_4$ is in the vector space spanned by $\gamma_1$ and $\gamma_2$}\\{}\iff& \text{$Q_4$ is in the vector space spanned by $Q_1$ and $Q_2$.}\endsplit$$
By the Triple Lemma \tref{TPL1} this is equivalent to the  existence of exactly two singular points  of multiplicity $c$   infinitely near to  $(0,0,1)$. Thus such points exist if and only if $e_1=3$. \qed
\enddemo

\bigpagebreak

\SectionNumber=4\tNumber=1
\heading Section \number\SectionNumber. \quad Singularities of multiplicity equal to degree divided by two.
\endheading 

Throughout this section $B$ is the polynomial ring $\pmb k[x,y]$, where $\pmb k$ is a field, and $d=2c$ is an even integer. 

When $\pmb k$ is algebraically closed, Theorem \tref{d=2c} completely classifies the parameterizations of rational plane curves $\Cal C$ of degree $d$ as a function of the configuration of  multiplicity $c$ singularities which appear on, or infinitely near,   $\Cal C$. Theorem \tref{d=2c} is typical of classification theorems in general in the sense that a classification theorem is always the culmination of one project and is often the starting point of new projects. Indeed, Theorem \tref{d=2c} is the starting point of the decomposition of $\Bbb T_d$ into strata which is carried out in Section 6. Also, we anticipate that Theorem \tref{d=2c} will eventually lead to a better understanding of the other singularities on and infinitely near the curve $\Cal C$. An analysis of this sort for quartics is carried out in Section 8. A future paper will contain our analysis of all singularities on sextics.

Let  $\Cal C$ be a rational plane curve of degree $d=2c$.  Recall, from Corollary \tref{CTLL}, that if 
there is a multiplicity $c$ singularity on, or infinitely near, 
$\Cal C$, then every entry in a homogeneous Hilbert-Burch matrix for a parameterization of $\Cal C$ is a form of degree $c$. In Theorem \tref{orbitz} we decompose the space of all such Hilbert-Burch matrices (called $\operatorname{BalH}_d$, see Definition \tref{H}) into $11$ disjoint orbits under the action of the group $G=\operatorname{GL}_3(\pmb k)\times \operatorname{GL}_2(\pmb k)$. The orbits are parameterized by the poset   $\operatorname{ECP}$ of Definition \tref{ECP}. If $\natural$ is in $\operatorname{ECP}$, then the corresponding orbit is called $\operatorname{DO}^{\text{\rm Bal}}_{\natural}$. Each orbit $\operatorname{DO}^{\text{\rm Bal}}_{\natural}$ has the form $G\cdot M^{\text{\rm Bal}}_{\natural}$, where $M^{\text{\rm Bal}}_{\natural}$ is 
the intersection of $\operatorname{BalH}_d$ and 
an open subset of some affine space.  Let $\varphi$ be a Hilbert-Burch matrix in $M^{\text{\rm Bal}}_{\natural}$, for some $\natural$, and let $\Cal C$ be the corresponding curve. In Lemma \tref{d=2c'} we apply Theorem \tref{TBA} and Corollary \tref{CTLL} to explicitly identify the multiplicity $c$ singularities  on $\Cal C$, together with all infinitely near multiplicity $c$ singularities.

We notice that Theorem \tref{TBA} guarantees that the number of   distinct singular points of multiplicity  $c$ that are either on $\Cal C$ or infinitely near to $\Cal C$ is at most $3$ since this number is $\deg\gcd I_3(A)$ and $I_3(A)$ is generated by cubic forms. Of course, this bound is also implied by Max Noether's formula
$$g=\binom {d-1}2-\sum_{q} \binom{m_q}2\tag\tnum{MNF}$$ which   gives the genus $g$ of the irreducible plane curve  $\Cal C$ of degree $d$, where $q$ varies over all singularities  and infinitely near  singularities of $\Cal C$, and  $m_q$ is the multiplicity at $q$. (See, for example, chapter 4 exercise 1.8 and chapter 5 examples 3.9.2 and 3.9.3 in \cite{\rref{Ha}}.) 
 At any rate, there are seven possible configurations of multiplicity $c$ singularities. The curve itself might have $0$, $1$, $2$, or $3$ singularities of multiplicity $c$ and any one of these singularities   has $0$,  $1$, or $2$ infinitely near singularities of multiplicity $c$, provided that the total number   is at most $3$. The seven possibilities are: $\emptyset$, $\{c\}$, $\{c,c\}$, $\{c,c,c\}$, $\{c:c\}$, $\{c:c,c\}$, and $\{c:c:c\}$, where a colon indicates an infinitely near singularity, a comma indicates a different singularity on the curve, and $\emptyset$ indicates that there are no singularities of multiplicity $c$ on, or infinitely near, $\Cal C$. Henceforth, we refer to the set of seven possible configurations of multiplicity $c$ singularities as $\operatorname{CP}$. The order that we impose on this set is dictated by the decomposition of $\Bbb T_d$ into strata which takes place in Section 6. 

\definition{Definition \tnum{po}} Let $(\operatorname{CP},\le)$ be the {\it Configuration Poset}.  The elements of $\operatorname{CP}$ are the seven possible configurations for multiplicity $c$ singularities on or infinitely near a rational plane curve of degree $d=2c$. We read ``$\#$'' as the sharp symbol  and we write $\#'\to \#$, for $\#'$ and $\#$ in $\operatorname{CP}$,  to mean $\#'\le \#$. The poset $\operatorname{CP}$ is:
$$
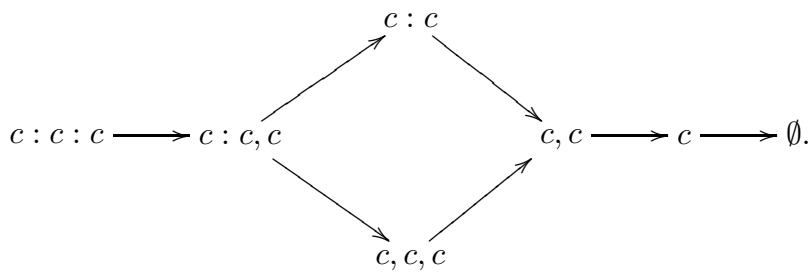
\xymatrix{&&c:c\ar[dr]\\
c:c:c \ar[r]  &c:c,c \ar[ru] \ar[rd]& &c,c \ar[r] &c\ar[r] &\emptyset.\\&&c,c,c\ar[ur]}$$ 
\enddefinition

\definition{Definition \tnum{H}} Let $B$ be the polynomial ring $\pmb k[x,y]$, over the field  $\pmb k$, and   $d$ be the even integer $2c$.  
\smallskip\flushpar{\bf(1)} Define $\Bbb H_d$ to be the space of $3\times 2$ matrices with entries from $B_c$.
\smallskip\flushpar{\bf(2)} Define $\operatorname{BalH}_d=\{\varphi\in \Bbb H_d\mid \operatorname{ht} I_2(\varphi)=2\}$.
\smallskip\flushpar{\bf(3)} Define
$$\Bbb B\Bbb H_d=\left\{ \varphi\in \operatorname{BalH}_d\left\vert \matrix\format\l\\\text{the morphism $\Psi\:\Bbb P^1\to \Cal C$ determined by the signed}\\\text{maximal order minors of $\varphi$ is birational}\endmatrix\right.\right\}.$$
\smallskip\flushpar{\bf(4)}  The group $G=\operatorname{GL}_3(\pmb k)\times \operatorname{GL}_2(\pmb k)$ acts on $\Bbb H_d$. If $g=(\chi,\xi)\in G$ and $\varphi\in \Bbb H_d$, then ${g\varphi=\chi\varphi\xi^{-1}}$. 
\enddefinition

\remark{\bf Remarks \tnum{401}} {\bf (1)} The variety $\Bbb H_d$ is isomorphic to affine space $\Bbb A^{6c+6}$. The subsets $\Bbb B\Bbb H_d$ and $\operatorname{BalH}_d$ of $\Bbb H_d$ are both open, see Observation \tref{o5}. 

\smallskip\flushpar{\bf (2)} If 
$$\varphi=\bmatrix Q_{1,1}&Q_{1,2}\\Q_{2,1}&Q_{2,2}\\Q_{3,1}&Q_{3,2}\endbmatrix$$
is a $3\times 2$ matrix $\varphi$
with entries in $B$, then 
the ordered triple of ``signed maximal order minors'' of $\varphi$  is the following   ordered triple $\Phi(\varphi)$ of polynomials from $B$:
$$ \Phi(\varphi)=
\left (\left\vert \matrix  Q_{2,1}&Q_{2,2}\\Q_{3,1}&Q_{3,2}\endmatrix\right\vert, - \left\vert \matrix Q_{1,1}&Q_{1,2}\\ Q_{3,1}&Q_{3,2}\endmatrix\right\vert, \left\vert \matrix Q_{1,1}&Q_{1,2}\\Q_{2,1}&Q_{2,2}\endmatrix\right\vert \right ).\tag\tnum{PHI}$$

\smallskip\flushpar{\bf (3)} If $\varphi$ is in $\operatorname{BalH}_d$, then $\varphi$ is a Hilbert-Burch matrix for the 
row vector determined by the signed maximal order minors of $\varphi$; furthermore, $\varphi$ is {\bf balanced} in the sense that every element of $\varphi$ is a homogeneous element of the same degree. We refer to $\operatorname{BalH}_d$ as the space of ``Balanced Hilbert-Burch matrices for triples of homogeneous $d$-forms''. 

\smallskip\flushpar{\bf (4)} If $\varphi$ is in $\Bbb H_d$, then 
$$\varphi\in \Bbb B\Bbb H_d \iff  \matrix\format\l\\\text{$\varphi$ is a {\bf B}alanced {\bf H}ilbert-Burch matrix and the morphism}\\\text{determined by the signed maximal order minors of $\varphi$ is}\\   \text{{\bf B}irational and  {\bf B}ase point free}.\endmatrix$$
We refer to $\Bbb B\Bbb H_d$ as the space of ``Balanced Hilbert-Burch matrices for true triples of homogeneous $d$-forms''.
The notion of ``true triples'' is introduced in Remark \tref{UB}. In practice we are interested in the geometry which corresponds to $\Bbb B\Bbb H_d$. On the other hand, Theorem \tref{orbitz} and two thirds of Lemma \tref{d=2c'} do not require the birationality hypothesis. These results make sense in $\operatorname{BalH}_d$. 

\smallskip\flushpar {\bf (5)} The action of $G$ on $\Bbb H_d$ restricts to given  actions of $G$ on $\operatorname{BalH}_d$ and also on $\Bbb B\Bbb H_d$. Indeed,    $I_2(\varphi)$ and $I_2(g\varphi)$ are equal for all $g\in G$ and $\varphi\in \Bbb H_d$ and the curve parameterized by the signed maximal order minors of $g\varphi$ is the image, under a linear automorphism of $\Bbb P^2$ of the curve parameterized by the signed maximal order minors of $\varphi$; see Remark \tref{sigh}.

\smallskip\flushpar{\bf (6)}
The well known formula $\xi^{-1}=(\det \xi)^{-1}\operatorname{Adj} \xi$, where $\operatorname{Adj} \xi$ is the classical adjoint of $\xi$,  expresses the inverse of the matrix  as a rational function in the entries of $\xi$. Thus,
the function $\Upsilon\:G\times \Bbb H_d \to \Bbb H_d$, which is defined  by $\Upsilon(g,\varphi)=g\varphi$, is a morphism of varieties. 
\endremark 

There are $7$ configurations of multiplicity $c$ singularities in $\operatorname{CP}$, but $11$ disjoint orbits in our decomposition of $\operatorname{BalH}_d$. We form the Total Configuration Poset ($\operatorname{TCP}$) by adjoining $6$ new elements to $\operatorname{CP}$ and the Extended Configuration Poset ($\operatorname{ECP}$) by removing $2$ elements of $\operatorname{CP}$ from $\operatorname{TCP}$. The order in $\operatorname{TCP}$ is used when we combine some of the disjoint orbits $\{\operatorname{DO}_{\natural}|\natural\in \operatorname{ECP}\}$ to form $\{\operatorname{CO}_{\#}|\#\in \operatorname{CP}\}$ in Definition \tref{CO}.

\definition{Definition \tnum{ECP}}As a set, $\operatorname{TCP}$ consists of $\operatorname{CP}$ together with $6$ new elements:
$$\mu_2,\quad (c,\mu_4),\quad (c,\mu_5),\quad (\emptyset,\mu_4),\quad (\emptyset, \mu_5),\quad\text{and}\quad (\emptyset,\mu_6).$$
The order in $\operatorname{CP}$ is extended to give the order in $\operatorname{TCP}$:
$$\hskip-2.4pt\eightpoint
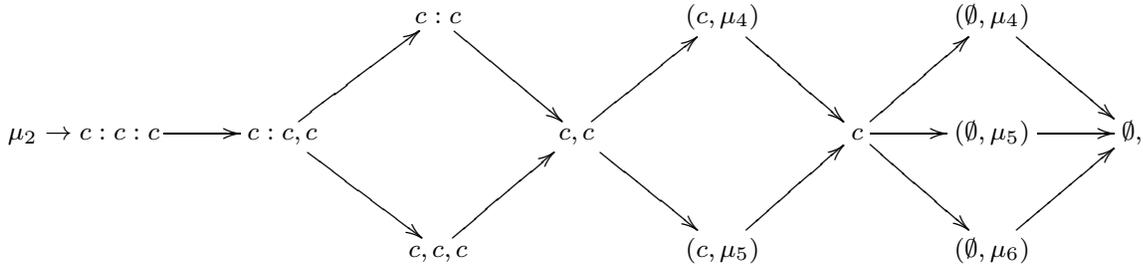
\xymatrix{&&c:c\ar[dr]&&(c,\mu_4)\ar[rd]&&(\emptyset,\mu_4)\ar[rd]\\
\mu_2\to c:c:c \hskip-1.02pt\ar[r]  &c:c,c \ar[ru] \ar[rd]& &c,c \ar[ru]\ar[rd] & &c\ar[ru]\ar[r]\ar[rd]&(\emptyset,\mu_5)\ar[r]&\emptyset,\\&&c,c,c\ar[ur]&&(c,\mu_5)\ar[ru]&&(\emptyset,\mu_6)\ar[ru]}$$ where we write $\natural'\to\natural$, for $\natural'$ and $\natural$ in $\operatorname{TCP}$ to mean $\natural'\le \natural$. The symbol ``$\natural$'' is read as ``natural''. The poset $\operatorname{ECP}$ is $\operatorname{TCP}$ with $c$ and $\emptyset$ removed. So $\operatorname{ECP}$ is 
$$
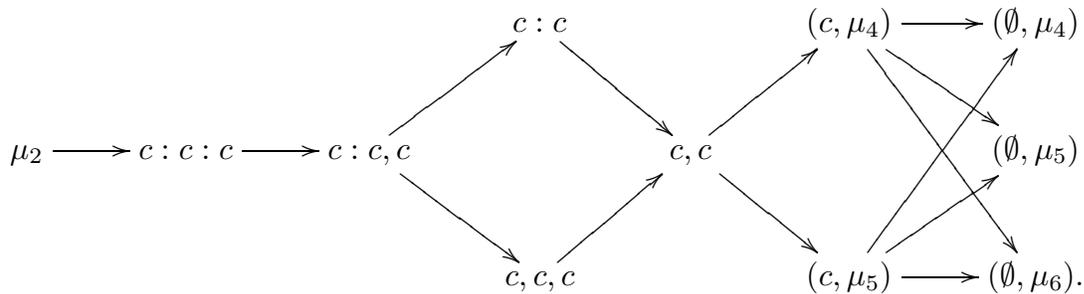
\xymatrix{&&&c:c\ar[dr]&&(c,\mu_4)\ar[r]\ar[rd]\ar[rdd]&(\emptyset,\mu_4)\\
\mu_2\ar[r]&c:c:c \ar[r]  &c:c,c \ar[ru] \ar[rd]& &c,c \ar[ru]\ar[rd] & &(\emptyset,\mu_5)\\&&&c,c,c\ar[ur]&&(c,\mu_5)\ar[r]\ar[ru]\ar[ruu]&(\emptyset,\mu_6).}$$\enddefinition

\definition{Definition \tnum{Mnat}} {\bf (1)} For each $\natural\in \operatorname{ECP}$, let   $M^{\text{\rm Bal}}_{\natural}$ be the following subset of  $\operatorname{BalH}_d$: 
$$\allowdisplaybreaks\alignat 1 M^{\text{\rm Bal}}_{\mu_2}&{}=\left\{\left.\bmatrix Q_1&0\\Q_2&Q_1\\0&Q_2\endbmatrix\in \operatorname{BalH}_d\right\vert \dim_{\pmb k}{<}Q_1,Q_2{>}=2\right\},\\
M^{\text{\rm Bal}}_{c:c:c}&{}=\left\{\left.\bmatrix Q_1&Q_2\\Q_3&Q_1\\0&Q_3\endbmatrix\in \operatorname{BalH}_d\right\vert \dim_{\pmb k}{<}Q_1,Q_2,Q_3{>}=3\right\},\\
M^{\text{\rm Bal}}_{c:c,c}&{}=\left\{\left.\bmatrix Q_1&0\\Q_2&Q_3\\0&Q_2\endbmatrix\in \operatorname{BalH}_d\right\vert \dim_{\pmb k}{<}Q_1,Q_2,Q_3{>}=3\right\},\\
M^{\text{\rm Bal}}_{c,c,c}&{}=\left\{\left.\bmatrix Q_1&Q_1\\Q_2&0\\0&Q_3\endbmatrix\in \operatorname{BalH}_d\right\vert \dim_{\pmb k}{<}Q_1,Q_2,Q_3{>}=3\right\},\\
M^{\text{\rm Bal}}_{c:c}&{}=\left\{\left.\bmatrix Q_1&Q_2\\Q_3&Q_4\\0&Q_3\endbmatrix\in \operatorname{BalH}_d\right\vert \dim_{\pmb k}{<}Q_1,Q_2,Q_3,Q_4{>}=4\right\},\\
M^{\text{\rm Bal}}_{c,c}&{}=\left\{\left.\bmatrix Q_1&Q_2\\Q_3&Q_3\\0&Q_4\endbmatrix\in \operatorname{BalH}_d\right\vert \dim_{\pmb k}{<}Q_1,Q_2,Q_3,Q_4{>}=4\right\},\\
M^{\text{\rm Bal}}_{(c,\mu_4)}&{}=\left\{\left.\bmatrix Q_1&Q_2\\Q_3&Q_1\\0&Q_4\endbmatrix\in \operatorname{BalH}_d\right\vert \dim_{\pmb k}{<}Q_1,Q_2,Q_3,Q_4{>}=4\right\},\\
M^{\text{\rm Bal}}_{(c,\mu_5)}&{}=\left\{\left.\bmatrix Q_1&Q_2\\Q_3&Q_4\\0&Q_5\endbmatrix\in \operatorname{BalH}_d\right\vert \dim_{\pmb k}{<}Q_1,Q_2,Q_3,Q_4,Q_5{>}=5\right\},\\
M^{\text{\rm Bal}}_{(\emptyset,\mu_4)}&{}=\left\{\left.\bmatrix Q_1&Q_2\\Q_2&Q_3\\Q_3&Q_4\endbmatrix\in \operatorname{BalH}_d\right\vert \dim_{\pmb k}{<}Q_1,Q_2,Q_3,Q_4{>}=4\right\},\\
M^{\text{\rm Bal}}_{(\emptyset,\mu_5)}&{}=\left\{\left.\bmatrix Q_1&Q_2\\Q_3&Q_4\\Q_5&Q_1\endbmatrix\in \operatorname{BalH}_d\right\vert \dim_{\pmb k}{<}Q_1,Q_2,Q_3,Q_4,Q_5{>}=5\right\},\ \text{and}\\M^{\text{\rm Bal}}_{(\emptyset,\mu_6)}&{}=\left\{\left.\bmatrix Q_1&Q_2\\Q_3&Q_4\\Q_5&Q_6\endbmatrix\in \operatorname{BalH}_d\right\vert \dim_{\pmb k}{<}Q_1,Q_2,Q_3,Q_4,Q_5,Q_6{>}=6\right\}.\endalignat$$

\medskip\flushpar{\bf (2)} For each $\natural\in \operatorname{ECP}$, let $\operatorname{DO}^{\text{\rm Bal}}_{\natural}$  be the subset $G\cdot M^{\text{\rm Bal}}_{\natural}=\{g\varphi\mid g\in G, \varphi\in M^{\text{\rm Bal}}_{\natural}\}$ of $\operatorname{BalH}_d$.

\medskip\flushpar{\bf (3)} For each $\natural\in \operatorname{ECP}$, let $M_{\natural}=M^{\text{\rm Bal}}_{\natural}\cap \Bbb B\Bbb H_d$ and $\operatorname{DO}_{\natural}=\operatorname{DO}^{\text{\rm Bal}}_{\natural}{}\cap \Bbb B\Bbb H_d$.  

\medskip\flushpar{\bf (4)} Define $M_c=M_{(c,\mu_4)}\cup M_{(c,\mu_5)}$.

\enddefinition 

\remark{Remarks} {\bf(1)} Notice that each set $M^{\text{\rm Bal}}_{\natural}$ and $M_{\natural}$ has the form $\operatorname{BalH}_d{} \cap \, U$ or $\Bbb B\Bbb H_d{}\cap{} U$, where $U$ is an open subset of some affine space. Indeed, for example the condition ``$\dim_{\pmb k}{<}Q_1,Q_2,Q_3{>}=3$'' is equivalent to the statement that the coefficients of $Q_1,Q_2,Q_3$ are not solutions of the  of the maximal order minors of a generic $3\times (c+1)$ matrix.

\flushpar{\bf (2)} We think of (4) from the above definition as an abbreviation. We don't use this abbreviation until (\tref{abr}) where its use does simplify the exposition.\endremark

\proclaim{Theorem \tnum{d=2c}} Let $\Cal C$ be a rational  plane curve of even degree $d=2c$ over an algebraically closed field $\pmb k$.  Assume that  there exists at least one singularity of multiplicity $c$ on or infinitely near $\Cal C$. Then there exists a linear automorphism 
$\Lambda$ of $\Bbb P^2$ and a matrix $\varphi$ in $M_{\natural}$, for some $\natural\in \operatorname{ECP}\setminus\{\mu_2,(\emptyset,\mu_4),(\emptyset,\mu_5),(\emptyset,\mu_6)\}$,
such that $\Lambda \Cal C$ is parameterized by the signed maximal order minors $\varphi$.
Furthermore, the following statements hold.

\smallskip\flushpar{\bf(1)} If the configuration of multiplicity $c$ singularities on or infinitely near $\Cal C$ is described by $\{c\}$, then $\Lambda \Cal C$ is parameterized by the signed maximal order minors of $\varphi$ in $M_{(c,\mu_4)}$ or $M_{(c,\mu_5)}$. In this case, $[0:0:1]$ is the singularity of multiplicity $c$. 
 
\smallskip\flushpar{\bf(2)} If the configuration of multiplicity $c$ singularities on or infinitely near $\Cal C$ is described by $\{c,c\}$, then $\Lambda \Cal C$ is parameterized by the signed maximal order minors of $\varphi\in M_{c,c}$. In this case, $p_1=[0:0:1]$ and $p_2=[0:1:0]$ are the singularities of multiplicity $c$. 
 
\smallskip\flushpar{\bf(3)}  If the configuration of multiplicity $c$ singularities on or infinitely near $\Cal C$ is described by $\{c,c,c\}$, then $\Lambda \Cal C$ is parameterized by the signed maximal order minors of $\varphi\in M_{c,c,c}$.  In this case, $p_1=[0:0:1]$, $p_2=[0:1:0]$, and  $p_3=[1:0:0]$,  are the singularities on $\Cal C$ of multiplicity $c$. 
 
\smallskip\flushpar{\bf(4)}  If the configuration of multiplicity $c$ singularities on or infinitely near $\Cal C$ is described by $\{c:c\}$, then $\Lambda \Cal C$ is parameterized by the signed maximal order minors of $\varphi\in M_{c:c}$.   In this case, $p=[0:0:1]$ is the singularity on $\Cal C$ of  multiplicity $c$ and there is one singularity of multiplicity $c$  infinitely near to $p$. 
 
\smallskip\flushpar{\bf(5)}  If the configuration of multiplicity $c$ singularities on or infinitely near $\Cal C$ is described by $\{c:c,c\}$, then $\Lambda \Cal C$ is parameterized by the signed maximal order minors of $\varphi\in M_{c:c,c}$.  In this case, $p_1=[0:0:1]$ and $p_2=[1:0:0]$  are the singularities on $\Cal C$ of multiplicity $c$ and there is one  singularity of multiplicity $c$ infinitely near to $p_1$. 
 
\smallskip\flushpar{\bf(6)}  If the configuration of multiplicity $c$ singularities on or infinitely near $\Cal C$ is described by $\{c:c:c\}$, then $\Lambda \Cal C$ is parameterized by the signed maximal order minors of $\varphi\in M_{c:c:c}$.  In this case, $p_1=[0:0:1]$ is the singularity  on $\Cal C$ of multiplicity $c$ and there are two  singularities of  multiplicity $c$ infinitely near to $p_1$. 
  
\endproclaim

\demo{Proof} Recall from Corollary \tref{CTLL}, that if there is a singularity of multiplicity $c$ on or infinitely near $\Cal C$, then there is a singularity of multiplicity $c$ on $\Cal C$ and  every entry in a homogeneous Hilbert-Burch matrix for $\Cal C$ is a form of degree $c$. Thus, $\Cal C$ is birationally parameterized by the signed maximal order minors of some matrix $\varphi$ in $\Bbb B\Bbb H_d$.   Theorem \tref{orbitz} 
shows that there exists $g\in G$ and $\varphi_{\natural}\in M_{\natural}$ for some $\natural\in\operatorname{ECP}$ with $\varphi=g\varphi_{\natural}$. Thus, there is a linear automorphism of $\Bbb P^2$ which carries $\Cal C$ to the curve $\Cal C_{\natural}$ which is parameterized by the signed maximal order minors of $\varphi_{\natural}$, see Remark \tref{sigh}. Part (3) of Lemma \tref{d=2c'} shows that $\natural\notin\{\mu_2,(\emptyset,\mu_6),(\emptyset,\mu_5),(\emptyset,\mu_4)\}$.
Part (3) of Lemma \tref{d=2c'} also records the multiplicity $c$ singularities on or infinitely near $\Cal C_{\natural}$. \qed
\enddemo

\proclaim{Theorem \tnum{orbitz}} If 
 every polynomial in $\pmb k[x]$ of degree $2$ or $3$   has a root in the field $\pmb k$, then
the space of balanced Hilbert-Burch matrices $\operatorname{BalH}_d$ is the disjoint union of the orbits $\operatorname{DO}^{\text{\rm Bal}}_{\natural}$ as $\natural$ varies over $\operatorname{ECP}$. \endproclaim

\demo{Proof}  Fix $\varphi\in \operatorname{BalH}_d$. We first prove that there exists $g\in G$ with $g\varphi\in M^{\text{\rm Bal}}_{\natural}$ for some $\natural\in \operatorname{ECP}$. Consider the parameter $\mu=\mu(I_1(\varphi))$. The matrix $\varphi$ has six homogeneous entries of degree $c$, so $\mu \le 6$. On the other hand, the hypothesis that $\operatorname{ht} I_2(\varphi)=2$ guarantees that $2\le \mu$. Thus, $2\le \mu \le 6$.   We treat each possible value for $\mu $ separately. If $\mu=6$, then the entries of $\varphi$ are linearly independent and $\varphi\in M^{\text{\rm Bal}}_{(\emptyset,\mu_6)}$. Suppose now that $\mu=5$. If $\varphi$ has a generalized zero (see Remark \tref{gri}), then, after row and column operations, $\varphi$ is transformed into $g\varphi\in M^{\text{\rm Bal}}_{(c,\mu_5)}$. If $\varphi$ does not have a generalized zero, one may apply row and column operations to put 
$\varphi$ in the form $$\bmatrix Q_1&Q_4\\Q_2&Q_5\\Q_3&\sum\limits_{i=1}^5\alpha_iQ_i\endbmatrix,$$ where the $\alpha_i\in \pmb k$ are constants. Further row and column operations (and re-naming the entries of $\varphi$) put $\varphi$ in the form 
$$\bmatrix Q_1&Q_4\\Q_2&Q_5\\Q_3&\sum\limits_{i=1}^2\alpha_iQ_i\endbmatrix,$$
and ultimately one finds $g\in G$ with $g\varphi\in M^{\text{\rm Bal}}_{(\emptyset,\mu_5)}$.  Lemmas \tref{rest'} and \tref{L4q} show that if $\mu$ is equal to $4$ or $3$,   then there exists $g\in G$ with $g\varphi\in M^{\text{\rm Bal}}_{\natural}$ for some $\natural\in \operatorname{ECP}$.   Finally, if $\mu=2$, then one quickly puts $\varphi$ in the form
$$\bmatrix *&*\\*&*\\0&Q_2\endbmatrix.$$ Further elementary row and column operations transform $\varphi$ into
$$\bmatrix Q_1&\alpha Q_1\\Q_2&\beta Q_1\\0&Q_2\endbmatrix,\text{ then }\bmatrix Q_1& 0\\Q_2&\beta Q_1\\0&Q_2\endbmatrix,$$for some constants $\alpha$ and $\beta$. The constant $\beta$ is non-zero since $\operatorname{ht} I_2(\varphi)=2$ and there exists $g\in G$ with $g\varphi\in M^{\text{\rm Bal}}_{\mu_2}$.  We have shown that the $\operatorname{BalH}_d=\cup_{\natural\in\operatorname{ECP}}\operatorname{DO}^{\text{\rm Bal}}_{\natural}$. 
The chart of invariants in (2) of Lemma \tref{d=2c'} shows that the orbits $\operatorname{DO}^{\text{\rm Bal}}_{\natural}$, with $\natural\in \operatorname{ECP}$, are disjoint.
\qed 
\enddemo

\proclaim {Lemma \tnum{d=2c'}} Let $\varphi_{\natural}$, $C_{\natural}$, and $A_{\natural}$ be matrices which satisfy {\rm(\tref{no1})} and {\rm(\tref{no2})} with $\varphi_{\natural}\in \operatorname{DO}^{\text{\rm Bal}}_{\natural}$.  

\smallskip\flushpar{\bf (1)} One may 
transform the matrices $(C_{\natural},A_{\natural})$, using elementary operations and the suppression of zero rows, into the matrices $(C_{\natural}',A_{\natural}')$, which are given by{\rm:}
$$\eightpoint  C_{(\emptyset,\mu_6)}'=\left [\smallmatrix T_1&0\\T_2&0\\T_3&0\\0&T_1\\0&T_2\\0&T_3\endsmallmatrix\right ],\  
A_{(\emptyset,\mu_6)}'=\left [\smallmatrix u_1&0&0\\0&u_1&0\\0&0&u_1\\u_2&0&0\\0&u_2&0\\0&0&u_2\endsmallmatrix\right ],\ 
C_{(\emptyset,\mu_5)}'=\left [\smallmatrix T_1&T_3\\T_2&0\\T_3&0\\0&T_1\\0&T_2\endsmallmatrix\right ],\ 
A_{(\emptyset,\mu_5)}'=\left [\smallmatrix u_1&0&u_2\\0&u_1&0\\0&0&u_1\\u_2&0&0\\0&u_2&0\endsmallmatrix\right ],$$$$
\eightpoint  
C_{(c,\mu_5)}'=\left [\smallmatrix T_1&0\\T_2&0\\0&T_1\\0&T_2\\0&T_3\endsmallmatrix\right ],\  
A_{(c,\mu_5)}'=\left [\smallmatrix u_1&0&0\\0&u_1&0\\u_2&0&0\\0&u_2&0\\0&0&u_2\endsmallmatrix\right ], 
C_{(\emptyset,\mu_4)}'=\left [\smallmatrix T_1&0\\T_2&T_1\\T_3&T_2\\0&T_3\endsmallmatrix\right ],\ 
A_{(\emptyset,\mu_4)}'=\left [\smallmatrix u_1&0&0\\u_2&u_1&0\\0&u_2&u_1\\0&0&u_2\endsmallmatrix\right ],$$$$\eightpoint 
C_{(c,\mu_4)}'=\left [\smallmatrix T_1&T_2\\T_2&0\\0&T_1\\0&T_3\endsmallmatrix\right ],\ 
A_{(c,\mu_4)}'=\left [\smallmatrix u_1&u_2&0\\0&u_1&0\\u_2&0&0\\0&0&u_2\endsmallmatrix\right ],\ 
C_{c,c}'=\left [\smallmatrix T_1&0\\0&T_1\\T_2&T_2\\0&T_3\endsmallmatrix\right ],\ 
A_{c,c}'=\left [\smallmatrix u_1&0&0\\u_2&0&0\\0&(u_1+u_2)&0\\0&0&u_2\endsmallmatrix\right ],$$$$
\eightpoint C_{c:c}'=\left [\smallmatrix T_1&0\\T_2&T_3\\0&T_1\\0&T_2\endsmallmatrix\right ],\ 
A_{c:c}'=\left [\smallmatrix u_1&0&0\\0&u_1&u_2\\u_2&0&0\\0&u_2&0\endsmallmatrix\right ],\ 
C_{c:c:c}'=\left [\smallmatrix T_1&T_2\\T_2&T_3\\0&T_1\endsmallmatrix\right ],\ 
A_{c:c:c}'=\left [\smallmatrix u_1&u_2&0\\0&u_1&u_2\\u_2&0&0\endsmallmatrix\right ],
$$$$\eightpoint 
C_{c:c,c}'=\left [\smallmatrix T_1&0\\T_2&T_3\\0&T_2\endsmallmatrix\right ],\ 
A_{c:c,c}'=\left [\smallmatrix u_1&0&0\\0&u_1&u_2\\0&u_2&0\endsmallmatrix\right ],\ 
C_{c,c,c}'=\left [\smallmatrix T_1&T_1\\T_2&0\\0&T_3\endsmallmatrix\right ],\  
A_{c,c,c}'=\left [\smallmatrix u_1+u_2&0&0\\0&u_1&0\\0&0&u_2\endsmallmatrix\right ],$$$$ \eightpoint 
C_{\mu_2}'=\left [\smallmatrix T_1&T_2\\T_2&T_3\endsmallmatrix\right ],
\  A_{\mu_2}'=\left [\smallmatrix u_1&u_2&0\\0&u_1&u_2\endsmallmatrix\right ].\hphantom{\ 
C_{c,c,c}'=\left [\smallmatrix T_1&T_1\\T_2&0\\0&T_3\endsmallmatrix\right ],\  
A_{c,c,c}'=\left [\smallmatrix u_1+u_2&0&0\\0&u_1&0\\0&0&u_2\endsmallmatrix\right ],}$$

\smallskip\flushpar{\bf (2)} The matrices $(\varphi_{\natural}, C_{\natural},A_{\natural})$ satisfy
$$\matrix \format \c&\qquad \c&\qquad \c&\qquad \c\\ 
\natural&\mu(I_1(\varphi_{\natural}))&\mu(I_2(C_{\natural}))&\gcd(I_3(A_{\natural}))\\
(\emptyset,\mu_6)&6&6&1\\
(\emptyset,\mu_5)&5&6&1\\
(c,\mu_5)&5&5&\ell_1\\
(\emptyset,\mu_4)&4&6&1\\
(c,\mu_4)&4&5&\ell_1\\
c,c&4&4&\ell_1\ell_2\\
c:c&4&4&\ell_1^2\\
c,c,c&3&3&\ell_1\ell_2\ell_3\\
c:c,c&3&3&\ell_1^2\ell_2\\
c:c:c&3&3&\ell_1^3\\
\mu_2&2&1&0,
\endmatrix$$for some non-associate linear forms $\ell_1,\ell_2,\ell_3$ in $\pmb k[\pmb u]$. 

\smallskip\flushpar{\bf (3)} If $\varphi_{\natural}$ is in $M_{\natural}$ and $\Cal C_{\natural}$ is the curve parameterized by   the signed maximal ordered minors of $\varphi_{\natural}$, then the multiplicity $c$ singularities on or infinitely near $\Cal C_{\natural}$ are  given in the following chart. 
$$\alignat 3
&\natural&\qquad&\text{\rm the multiplicity $c$}&\ \ &\text{\rm the number of multiplicity $c$}\\\vspace{-5pt}
&&\ \ &\text{\rm singularities $p$ on $\Cal C_i$}&\ \ &\text{\rm singularities infinitely near to $p$}\\
\noalign{\hrule}
&(\emptyset,\mu_6)&\ \ &\text{\rm none}&\ \  \\\noalign{\hrule}
&(\emptyset,\mu_5)&\ \ &\text{\rm none}&\ \  \\\allowdisplaybreak\noalign{\hrule}
&(c,\mu_5)&\ \ &[0:0:1]&\ \ &0\\\allowdisplaybreak\noalign{\hrule}
&(\emptyset,\mu_4)&\ \ &\text{\rm none}&&\ \ \\\allowdisplaybreak\noalign{\hrule}
&(c,\mu_4)&\ \ &[0:0:1]&\ \ &0\\\allowdisplaybreak\noalign{\hrule}
&c,c&\ \ & &\ \ &\\ 
&&\ \ &[0:0:1]&\ \ &0\\ 
&&\ \ &[0:1:0]&\ \ &0\\\allowdisplaybreak\noalign{\hrule}
&c:c&\ \ &[0:0:1]&\ \ &1\\\allowdisplaybreak\noalign{\hrule}
&c:c:c&\ \ &[0:0:1]&\ \ &2\\\allowdisplaybreak\noalign{\hrule}
&c:c,c&\ \ & &\ \ &\\ 
&&\ \ &[0:0:1]&\ \ &1\\ 
&&\ \ &[1:0:0]&\ \ &0\\\allowdisplaybreak\noalign{\hrule}
&c,c,c&\ \ & &\ \ &\\ 
&&\ \ &[0:0:1]&\ \ &0\\
&&\ \ &[0:1:0]&\ \ &0\\ 
&&\ \ &[1:0:0]&\ \ &0\endalignat$$
\endproclaim
\remark{\bf Remarks \tnum{408}} {\bf(1)}  In (1) and (2), there exists $g\in G$ with $g\varphi_{\natural}\in M_{\natural}^{\text{\rm Bal}}$. Two additional hypotheses have been imposed on the matrix $\varphi_{\natural}$ in part (3). First of all, $\varphi_{\natural}$ is already in $M_{\natural}^{\text{\rm Bal}}$; that is, $g$ may be taken to be $1$. Secondly, the parameterization determined by $\varphi_{\natural}$ is birational; that is, $\varphi_{\natural}\in\Bbb B\Bbb H_d$; thus, 
$\varphi_{\natural}\in M_{\natural}^{\text{\rm Bal}}\cap \Bbb B\Bbb H_d=M_{\natural}$.

\flushpar{\bf(2)} To read the chart of (3), notice that the chart says  for example, that there are $2$ multiplicity $c$ singularities on the  curve $\Cal C_{c:c,c}$; one of these singularities ${([0:0:1])}$ has an infinitely near singularity of multiplicity $c$ and the other singularity ${([1:0:0])}$ does not have any infinitely near singularities of multiplicity $c$.

\smallskip\flushpar{\bf(3)} We did not include $\mu_2$ in the chart of (3) because
the intersection $\operatorname{DO}^{\text{\rm Bal}}_{\mu_2}{}\cap \Bbb B\Bbb H_d$, which is also called $\operatorname{DO}_{\mu_2}$, is empty.

\smallskip\flushpar{\bf(4)} Fix   $\natural$   in $\operatorname{ECP}$ with $\natural\neq \mu_2$. 
We note that $\operatorname{DO}_{\natural}$ is a non-empty  open subset of $\operatorname{DO}^{\text{\rm Bal}}_{\natural}$; see Observation \tref{o5} and Proposition \tref{no-e'}.
  We also note that, according to Observation \tref{2.1}, if $c$ is a prime integer, then $\operatorname{DO}_{\natural}=\operatorname{DO}^{\text{\rm Bal}}_{\natural}$.
\endremark 
 \demo{Proof} By the definition of $\operatorname{DO}^{\text{\rm Bal}}_d$, there is an element $g\in G$ with $g\varphi_{\natural}\in M^{\text{\rm Bal}}_{\natural}$. Remark \tref{305} shows how the matrices $C_{g\natural}$ and $A_{g\natural}$ are obtained from $C_{\natural}$ and $A_{\natural}$. It suffices to prove the result when $\varphi_{\natural}\in M^{\text{\rm Bal}}_{\natural}$. 

We have recorded the matrices $(C_{\natural}',A_{\natural}')$, whose entries are linear forms from $\pmb k[\pmb T]$ and $\pmb k[\pmb u]$, and which satisfy $\pmb T \varphi_{\natural}=\bmatrix Q_1&\cdots&Q_{\mu}\endbmatrix C_{\natural}'$ and $C_{\natural}'\pmb u^{\text{\rm T}}=A_{\natural}'\pmb T^{\text{\rm T}}$, for $\mu=\mu(I_1(\varphi_{\natural}))$. The set of linearly independent  forms $Q_1,\dots,Q_\mu$ from $B_c$ may be extended to a basis $Q_1,\dots,Q_{c+1}$ for $B_c$. The entries of $\rho^{(c)}$ also form a basis for $B_c$; so there is an invertible matrix $\upsilon $, with entries in $\pmb k$, so that $\rho^{(c)}=[Q_1,\dots,Q_{c+1}]\upsilon $.  The matrix $C_{\natural}$ has been defined to satisfy $\pmb T \varphi_{\natural}=\rho^{(c)} C_{\natural}$.
  Thus,
$$\split \bmatrix Q_1&\cdots&Q_{c+1}\endbmatrix\upsilon C_{\natural}&{}=\rho^{(c)}C_{\natural}=\pmb T \varphi_{\natural}=\bmatrix Q_1&\cdots&Q_{\mu}\endbmatrix C_{\natural}'\\
&{}=\bmatrix Q_1&\cdots&Q_{c+1}\endbmatrix \bmatrix I_\mu\\0\endbmatrix C_{\natural}';\endsplit $$
hence,  the matrices $\upsilon C_{\natural}$ and $\left [\smallmatrix I_{\mu} \\0\endsmallmatrix\right ] C_{\natural}'$, of linear forms from $\pmb k[\pmb T]$, are equal,   and $C_{\natural}'$ is obtained from $C_{\natural}$ by applying invertible row operations and suppressing zero rows. It quickly follows that 
$\upsilon A_{\natural}=\left [\smallmatrix I_{\mu} \\0\endsmallmatrix\right ] A_{\natural}'$, and $A_{\natural}'$ is obtained from $A_{\natural}$ by applying invertible row operations and suppressing zero rows.
This completes the proof of (1). To prove (2), one  quickly calculates 
$\mu (I_2(C_{\natural}'))$ for each matrix $C_{\natural}'$ of assertion (1). It is clear that $I_2(C_{\natural}')=I_2(C_{\natural})$. One also calculates $\gcd I_3(A'_{\natural})=\gcd I_3(A_{\natural})$:
$$\matrix \format \c&\ \c&\qquad\qquad \c&\ \c\\
\natural&\gcd I_3(A'_{\natural})&\natural&\gcd I_3(A'_{\natural})\\
(\emptyset,\mu_6)&1&c:c&u_2^2\\
(\emptyset,\mu_5)&1&c:c:c&u_2^3\\
(c,\mu_5)&u_2&c:c,c&u_1u_2^2\\
(\emptyset,\mu_4)&1&c,c,c&(u_1+u_2)u_1u_2\\
(c,\mu_4)&u_2&\mu_2&0\\
c,c&u_2(u_1+u_2).\endmatrix$$Remark \tref{305} shows that the transformation $A_{\natural}$ to $A_{g^{-1}\natural}$, for $g\in G$, replaces $u_1$ and $u_2$ with linearly independent  linear forms $\ell_1$ and $\ell_2$ from $\pmb k[\pmb u]$. 
For (3), 
Theorem \tref{TBA}  guarantees that there are exactly $6-\mu  (I_2(C_{\natural}))$ distinct singularities of multiplicity $c$ on or infinitely near $\Cal C_{\natural}$. We use Corollary \tref{CTLL} to identify these singularities. 
 \qed
\enddemo

The following small calculation provides
 a sufficient condition for the existence of a generalized zero in a matrix $\varphi$. We use this calculation three times as we complete our classification of the matrices $\varphi$ of Lemma \tref{d=2c'}. The most important application of this calculation occurs when the parameters ``$b$'' and ``$N$'' are both taken to be zero. In this case, the matrix ``$\varphi$'' looks like  $\bmatrix A_1&A_2\endbmatrix$. 
\proclaim{Observation \tnum{Ogz}} Let $R$ be an algebra over the field $\pmb k$ and 
$$\varphi=\bmatrix A_1&A_2&A_3\\A_4&A_5&A_6\endbmatrix$$ be a  matrix with entries in $R$, where $A_1$ and $A_2$ are $a\times 1$ matrices, $A_3$ is an $a\times N$ matrix, $A_4$ and $A_5$ are $b\times 1$ matrices, and $A_6$ is an $b\times N$ matrix for some non-negative integers $a$, $b$, and $N$.  Suppose that every polynomial of degree   $a$ in $\pmb k[x]$ has a root in $\pmb k$. Suppose further that each entry of $A_2$ is in the vector space spanned by the entries of $A_1$.   Then   there exist invertible matrices $\chi'$, $\chi$, and $\xi$, with entries in $\pmb k$, such that 
$$\chi \varphi\xi=\bmatrix A_1'&A_2''&A_3'\\A_4&A_5''&A_6\endbmatrix,$$
where at least one entry of $A_2''$ is zero, $A_i'=\chi'A_i$ for $i$ equal to $1$ or $3$, $A_2''= \chi'(A_2-\lambda A_1)$ and 
 ${A_5''=A_5-\lambda A_4}$ for some $\lambda$ in $\pmb k$,  and  the submatrices   $A_4$ and $A_6$ remain unchanged in the transformation from $\varphi$ to $\chi \varphi\xi$.
\endproclaim
\demo{Proof} The hypothesis about the entries of $A_1$ and $A_2$  guarantees the existence of a matrix $M$, with entries in $\pmb k$, such that $A_2=M A_1$. The hypothesis about roots of polynomials ensures that $M$ has an eigenvalue $\lambda$ in $\pmb k$. Let $\xi$ be the elementary matrix which subtracts $\lambda$ times column one from column two.   The second column of $\varphi\xi$ is 
$$\bmatrix (M-\lambda I)A_1\\ A_5-\lambda A_4\endbmatrix; $$the other columns of $\varphi$ are unperturbed under the transformation $\varphi \mapsto \varphi\xi$. 
 The matrix $M-\lambda I$ is singular, so there exists a non-zero row vector $v$, with entries in $\pmb k$, so that $v(M-\lambda I)=0$. Insert $v$ as a row in an invertible $a\times a$ matrix $\chi'$ and let
$$\chi=\bmatrix \chi'&0_{a\times b}\\0_{b\times a}&I_{b\times b}\endbmatrix,$$where $0$ is a zero matrix and $I$ is an identity matrix. The triple $(\chi',\chi,\xi)$ satisfies the required properties.  \qed \enddemo

\proclaim{Lemma \tnum{rest'}}  Let $\pmb k$ be a field. Assume that 
 every quadratic polynomial   in $\pmb k[x]$ has a root in $\pmb k$. 
Let  $R$ be a $\pmb k$-algebra, and $\varphi$ be a $3\times 2$ matrix with entries from $R$. Assume that
the entries of $\varphi$ span a vector space of dimension $4$  and $\operatorname{ht} I_2(\varphi)=2$. Then there exist invertible matrices $\chi$ and $\xi$ over $\pmb k$ so that $\chi\varphi \xi$ has one of the following forms{\rm:}
$$\varphi_{(\emptyset,\mu_4)}=\left [\smallmatrix Q_1&Q_2\\Q_2&Q_3\\Q_3&Q_4\endsmallmatrix\right ],\ \  
 \varphi_{(c,\mu_4)}=\left [\smallmatrix Q_1&Q_2\\Q_3&Q_1\\0&Q_4\endsmallmatrix\right ],\ \  \varphi_{c,c}=\left [\smallmatrix Q_1&Q_2\\Q_3&Q_3\\0&Q_4\endsmallmatrix\right ],\ \ \text{or}\ \ 
\varphi_{c:c}=\left [\smallmatrix Q_1&Q_3\\Q_2&Q_4\\0&Q_2\endsmallmatrix\right ],$$ with $Q_1,Q_2,Q_3,Q_4$ linearly independent.\endproclaim

\demo{Proof} There are two possibilities for the original matrix $\varphi$. In Case 1, the entries in each column of $\varphi$ span a vector space of dimension $2$. In  Case 2, the entries of at least one of the columns of $\varphi$ span a vector space of dimension $3$.  In Case 1, $\varphi$ can be put in the form
$$\bmatrix Q_1&*_1\\Q_2&*_2\\0&*_3\endbmatrix,$$ where $Q_1,Q_2$ are linearly independent  and $*_1$, $*_2$, and $*_3$ span a two dimensional subspace of $R$ which meets the vector space ${<}Q_1,Q_2{>}$ only at $0$.
The grade of $I_2(\varphi)$ is two; so no row of $\varphi$ can be zero. In particular, $*_3$ is not zero. We call $*_3$ by the name $Q_3$ and we have $Q_1,Q_2,Q_3$ linearly independent. At least one of the entries $*_1$ or $*_2$ is not in ${<}Q_1,Q_2,Q_3{>}$. Apply a row exchange and rename $Q_1$ and $Q_2$, if necessary. We have transformed $\varphi$ into the form
$$\bmatrix Q_1&\alpha_3Q_3+\alpha_4Q_4\\Q_2&Q_4\\0&Q_3\endbmatrix,$$with $Q_1,Q_2,Q_3,Q_4$ linearly independent and $\alpha_3,\alpha_4$  in $\pmb k$. Let $\operatorname{Ro}_i$ and $\operatorname{Co}_i$ represent row $i$ and column $i$, respectively. 
Subtract  $\alpha_4\operatorname{Ro}_2+\alpha_3\operatorname{Ro}_3$ from
$\operatorname{Ro}_1$ and rename $Q_1-\alpha_4Q_2$ as $Q_1$. We have transformed $\varphi$ into 
$$\varphi'=\bmatrix Q_1&0\\Q_2&Q_4\\0&Q_3\endbmatrix,$$
which can be transformed into the form of $\varphi_{c,c}$.

In Case 2, the matrix $\varphi$ may be put in the form 
$$\bmatrix Q_1&\alpha_1Q_1+\alpha_2Q_2+\alpha_3Q_3+\alpha_4Q_4\\
Q_2&\beta_1Q_1+\beta_2Q_2+\beta_3Q_3+\beta_4Q_4\\
Q_3&Q_4\endbmatrix,$$ for some constants $\alpha_i,\beta_i$ with $Q_1,Q_2,Q_3,Q_4$ linearly independent. Subtract $\alpha_4\operatorname{Ro}_4$ from  $\operatorname{Ro}_1$   and $\beta_4\operatorname{Ro}_3$ from $\operatorname{Ro}_2$. Rename $Q_1$ and $Q_2$: the old $Q_1-\alpha_4Q_3$ becomes the new $Q_1$ and the old $Q_2-\beta_4Q_3$ becomes the new $Q_2$. Rename the constants $\alpha_i$ and $\beta_i$. We have transformed $\varphi$ into $$\bmatrix Q_1&\alpha_1Q_1+\alpha_2Q_2+\alpha_3Q_3\\Q_2&\beta_1Q_1+\beta_2Q_2+\beta_3Q_3\\Q_3&Q_4\endbmatrix.$$ Subtract $\alpha_1\operatorname{Co}_1$ from $\operatorname{Co}_2$ and rename $Q_4$ and $\beta_2$. The matrix $\varphi$ has become
$$\varphi'=\bmatrix Q_1&\alpha_2Q_2+\alpha_3Q_3\\Q_2&\beta_1Q_1+\beta_2Q_2+\beta_3Q_3\\Q_3&Q_4\endbmatrix.$$ At this point there are three cases. Either $\alpha_3=\beta_3=0$ (Case 2A),  or   $\alpha_3\neq 0$ (Case 2B), or   $\alpha_3= 0$ and $\beta_3\neq 0$ (case 2C).

In case 2A, apply Observation \tref{Ogz} with 
$$ \matrix\format \l&\qquad\l\\ A_1=\bmatrix Q_1\\Q_2\endbmatrix& A_2= \bmatrix  \alpha_2Q_2\\ \beta_1Q_1+\beta_2Q_2\endbmatrix\\\vspace{5pt}
A_4=\bmatrix Q_3\endbmatrix & A_5=\bmatrix Q_4\endbmatrix\endmatrix $$
 and transform $\varphi'$ into
$$\bmatrix Q_1'&0\\Q_2'&*\\Q_3&Q_4'\endbmatrix,$$ where $Q_1',Q_2',Q_3,Q_4'$ are linearly independent and $*$ is a non-zero element of the vector space ${<}Q_1',Q_2'{>}$. If $*\in {<}Q_1'{>}$, then $\varphi$ may be transformed into $\varphi_{c:c}$; otherwise, $\varphi$ may be transformed into $\varphi_{c,c}$. 

In Case 2B, one may quickly transform $\alpha_3$ into $1$. (Indeed, one may multiply $\operatorname{Co}_2$ by $\alpha_3^{-1}$ and rename $Q_4$ and the constants $\alpha_i$ and $\beta_i$.) At this point, $\varphi'$ is
$$\bmatrix Q_1&\alpha_2Q_2+Q_3\\Q_2&\beta_1Q_1+\beta_2Q_2+\beta_3Q_3\\Q_3&Q_4\endbmatrix.$$
Add $\alpha_2\operatorname{Ro}_2$ to $\operatorname{Ro}_3$ and rename $Q_3$, $Q_4$, and the $\beta$'s to obtain
$$\bmatrix Q_1&Q_3\\Q_2&\beta_1Q_1+\beta_2Q_2+\beta_3Q_3\\Q_3&Q_4\endbmatrix.$$Subtract $\beta_3\operatorname{Ro}_1$ from $\operatorname{Ro}_2$ and rename $Q_2$  to obtain
$$\bmatrix Q_1&Q_3\\Q_2&\beta_1Q_1+\beta_2Q_2\\Q_3&Q_4\endbmatrix.$$
Subtract $\beta_2\operatorname{Co}_1$ from $\operatorname{Co}_2$,  $\beta_2\operatorname{Ro}_1$ from $\operatorname{Ro}_3$, and rename $Q_3$ and $Q_4$ to obtain 
$$\bmatrix Q_1&Q_3\\Q_2&\beta_1Q_1\\Q_3&Q_4\endbmatrix.$$If $\beta_1$ is zero, then $\varphi'$ may be transformed into $\varphi_{(c,\mu_4)}$.   If $\beta_1$ is not zero, then $\beta_1$ may be transformed into $1$. (One multiplies  $\operatorname{Ro}_2$ by $\beta_1^{-1}$ and renames $Q_2$.) At this point one uses row and column exchanges to transform $\varphi'$ into the form of $\varphi_{(\emptyset,\mu_4)}$. 

In Case 2C, one starts with 
$$\varphi'=\bmatrix Q_1&\alpha_2Q_2\\Q_2&\beta_1Q_1+\beta_2Q_2+\beta_3Q_3\\Q_3&Q_4\endbmatrix,$$with $\beta_3\neq 0$. Transform $\beta_3$ into $1$ by multiplying $\operatorname{Co}_2$ by $\beta_3^{-1}$ and renaming $Q_4$ and the constants. Add $\beta_1\operatorname{Ro}_1+\beta_2\operatorname{Ro}_2$ to $\operatorname{Ro}_3$ and rename $Q_3$ and $Q_4$. The matrix $\varphi'$ has become
$$\bmatrix Q_1&\alpha_2Q_2\\Q_2&Q_3\\Q_3&Q_4\endbmatrix.$$ If $\alpha_2=0$, then $\varphi'$ may be transformed into the matrix $\varphi_{(c,\mu_4)}$. If $\alpha_2\neq 0$, then $\alpha_2$ may be transformed into $1$ (by multiplying $\operatorname{Ro}_1$ by $\alpha_2^{-1}$ and renaming $Q_1$) and $\varphi'$ may be transformed into the form of $\varphi_{(\emptyset,\mu_4)}$. \qed
\enddemo

\proclaim{Lemma \tnum{L4q}} Let $\pmb k$ be a field. Assume that 
 every polynomial in $\pmb k[x]$ of degree $2$ or $3$   has a root in $\pmb k$. 
Let  $R$ be a $\pmb k$-algebra, and $\varphi$ be a $3\times 2$ matrix with entries from $R$. Assume that
the entries of $\varphi$ span a vector space of dimension $3$  and $\operatorname{ht} I_2(\varphi)=2$. Then there exist invertible matrices $\chi$ and $\xi$ over $\pmb k$ so that $\chi\varphi \xi$ has one of the following forms{\rm:}
$$\varphi_{c:c:c}=\bmatrix Q_1&Q_2\\Q_3&Q_1\\ 0&Q_3\endbmatrix,\quad   \varphi_{c:c,c}=\bmatrix Q_1&0\\Q_2&Q_3\\0&Q_2\endbmatrix,\quad\text{or}\quad \varphi_{{c,c,c}}=\bmatrix Q_1&Q_1\\Q_2&0\\0&Q_3\endbmatrix, $$ with $Q_1,Q_2,Q_3$ linearly independent.\endproclaim

\demo{Proof} First we show that there exist invertible matrices $\chi$ and $\xi$ so that some entry of $\chi\varphi \xi$ is zero. There is nothing to show unless the entries in the first column of $\varphi$ are linearly independent; so, we make this assumption. The hypothesis that
the entries of $\varphi$ span a vector space of dimension $3$  tells us that
every entry in column two of $\varphi$ is contained in the vector space spanned by the entries of column one of $\varphi$. Apply Observation \tref{Ogz} to transform $\varphi$ into a matrix which contains a zero. Further  elementary row and column operations  put $\varphi$ into one of the  forms
$$\bmatrix Q_1&*\\Q_2&*\\0&Q_3\endbmatrix \quad\text{or}\quad \bmatrix Q_1&*\\Q_2&*\\0&Q_1\endbmatrix,\tag\tnum{.q}$$with $Q_1,Q_2,Q_3$ linearly independent. 

We first work on the left hand case of (\tref{.q}). Apply elementary row operations in order to make the entries labeled $*$ be in the vector space ${<}Q_1,Q_2{>}$.  Now we apply Observation \tref{Ogz} 
with 
$$ \matrix\format \l&\qquad\l\\ A_1=\bmatrix Q_1\\Q_2\endbmatrix& A_2= \bmatrix  *\\ *\endbmatrix\\\vspace{5pt}
A_4=\bmatrix 0\endbmatrix & A_5=\bmatrix Q_3\endbmatrix\endmatrix $$
 and transform
$\varphi$ into the form
$$\bmatrix Q_1& aQ_1+bQ_2\\Q_2&0\\0&Q_3\endbmatrix.\tag\tnum{..q}$$(It might be necessary to re-name $Q_1$ and $Q_2$.)
There are two cases.

We first consider the case $a\neq 0$. Add $b/a$ times row 2 to row 1, multiply column 1 by $a$, and rename $Q_1$ to obtain a matrix of the form of $\varphi_{{c,c,c}}$.

Now we consider the case $a=0$ in (\tref{..q}). In this case the hypothesis that $\operatorname{ht} I_2(\varphi)=2$ ensures that $b\neq 0$. It is easy to transform $\varphi$ into a matrix with the form $\varphi_{c:c,c}$. 

Finally, we consider the matrix on the right side of (\tref{.q}). If the entry in position $(2,1)$ is in the vector space ${<}Q_1,Q_2{>}$, we may use elementary row operations to put $\varphi$ into the form
$$\bmatrix Q_1&aQ_2\\Q_2&Q_3\\0&Q_1\endbmatrix.$$ The constant $a$ can not be zero because $\operatorname{ht}(I_2(\varphi))=2$ and $\varphi$ may be transformed into the form of the matrix $\varphi_{c:c:c}$. On the other hand, if the $(2,1)$ entry of $\varphi$ is not in ${<}Q_1,Q_2{>}$, then $\varphi$ may be transformed into a matrix of the form
$$\bmatrix Q_1&Q_3\\Q_2&*\\0&Q_1\endbmatrix.$$ Use a column operation and rename $Q_3$ to put $*$  into ${<}Q_1,Q_3{>}$. Use row operations, and rename $Q_2$, to put $\varphi$ into the form $$\bmatrix Q_1&Q_3\\Q_2&0\\0&Q_1\endbmatrix.$$This matrix may be easily transformed into a matrix of the form $\varphi_{c:c,c}$. \qed\enddemo

When we decompose $\Bbb T_d$ into strata in Section 6, most of our calculations (in particular the verification of irreducibility as well as the calculation of dimension) are made using the closure of a given stratum rather than the stratum itself. To facilitate those calculations, we gather the disjoint orbits $\operatorname{DO}_{\natural}$ for $\natural \in \operatorname{ECP}$, with $\natural <\#$ together to form the combined orbit $\operatorname{CO}_{\#}$  for $\#\in \operatorname{CP}$. Our proofs in Section 6 require that we identify a well understood irreducible variety $N_{\#}$ with $\operatorname{CO}_{\#}=G\cdot N_{\#}$. We lay out our candidates for $N_{\#}$ in Definition \tref{CO} and show that our candidates have the relevant properties in Theorem \tref{orbitz2}. 
\definition{Definition \tnum{CO}} {\bf (1)}  For each $\#\in \operatorname{CP}$, define $\operatorname{CO}^{\text{\rm Bal}}_{\#}$ to be the subset 
$$\operatorname{CO}^{\text{\rm Bal}}_{\#}=\bigcup\limits\limits_{\{\natural\in \operatorname{ECP}\mid \natural\le \# \text{ in } \operatorname{TCP}\}}\operatorname{DO}^{\text{\rm Bal}}_{\natural}$$
of $\operatorname{BalH}_d$.
\medskip\flushpar{\bf (2)} For each $\#\in \operatorname{CP}$,  define $N^{\text{\rm Bal}}_{\#}$ to be the following subset of $\operatorname{BalH}_d$:
$$\allowdisplaybreaks\alignat 1  N^{\text{\rm Bal}}_{c:c:c}&{}=\left\{ \bmatrix Q_1&Q_2\\Q_3&Q_1\\0&Q_3\endbmatrix\in \operatorname{BalH}_d  \right\},\\N^{\text{\rm Bal}}_{c:c,c}&{}=\left\{\left.\bmatrix Q_1&Q_2\\Q_3&Q_4\\0&Q_3\endbmatrix\in \operatorname{BalH}_d\right\vert \dim_{\pmb k}{<}Q_1,Q_2,Q_3,Q_4{>}\le 3\right\},\\
N^{\text{\rm Bal}}_{c,c,c}&{}=\left\{\left.\bmatrix Q_1&Q_2\\Q_3&Q_4\\0&Q_5\endbmatrix\in \operatorname{BalH}_d\right\vert \dim_{\pmb k}{<}Q_1,Q_2,Q_3,Q_4,Q_5{>}\le 3\right\},\\
N^{\text{\rm Bal}}_{c:c}&{}=\left\{\bmatrix Q_1&Q_2\\Q_3&Q_4\\0&Q_3\endbmatrix\in \operatorname{BalH}_d\right\},\\
N^{\text{\rm Bal}}_{c,c}&{}=\left\{\left.\bmatrix Q_1&Q_2\\Q_3&Q_4\\0&Q_5\endbmatrix\in \operatorname{BalH}_d\right\vert \dim_{\pmb k}{<}Q_3,Q_4,Q_5{>}\le2\right\},\\
N^{\text{\rm Bal}}_{c}&{}=\left\{\bmatrix Q_1&Q_2\\Q_3&Q_4\\0&Q_5\endbmatrix\in \operatorname{BalH}_d \right\},\ \text{and}\\
N^{\text{\rm Bal}}_{\emptyset}&{}=  \operatorname{BalH}_d .\endalignat $$
\medskip\flushpar{\bf (3)} For each $\#\in \operatorname{CP}$, let $\operatorname{CO}_{\#}=\operatorname{CO}^{\text{\rm Bal}}_{\#}{}\cap \Bbb B\Bbb H_d$ and $N_{\#}=N^{\text{\rm Bal}}_{\#}\cap \Bbb B\Bbb H_d$.  \enddefinition 
 
\proclaim{Theorem \tnum{orbitz2}}Let $\pmb k$ be a field which satisfies the hypothesis of Theorem \tref{orbitz}. 
\roster
 
\item If $\natural\in \operatorname{ECP}$ and $\#\in \operatorname{CP}$, with $\natural\le \#$ in $\operatorname{TCP}$, then  $M^{\text{\rm Bal}}_{\natural}\subseteq N^{\text{\rm Bal}}_{\#}$ and $M_{\natural}\subseteq N_{\#}$. 

\item If $\#\in \operatorname{CP}$, then $\operatorname{CO}^{\text{\rm Bal}}_{\#}=G\cdot N^{\text{\rm Bal}}_{\#}$ and
$\operatorname{CO}_{\#}=G\cdot N_{\#}$

\item The varieties  $N^{\text{\rm Bal}}_{\#}$ and   $N_{\#}$ are irreducible for all $\#\in \operatorname{CP}$.
\endroster
\endproclaim
\remark{Remark} The definition of $\operatorname{CO}^{\text{\rm Bal}}_{\#}$ ensures that
$$\text{$\#'<\#$ in $\operatorname{CP}$}\implies \operatorname{CO}^{\text{\rm Bal}}_{\#'}\subseteq \operatorname{CO}^{\text{\rm Bal}}_{\#}.$$
However, $$\text{$\#'<\#$ in $\operatorname{CP}$  does not imply }N^{\text{\rm Bal}}_{\#'}\subseteq N^{\text{\rm Bal}}_{\#}.$$ Indeed, 
$\{c,c,c\}<\{c,c\}$ and if 
$$\varphi=\bmatrix Q_1&Q_2\\Q_2&Q_3\\0&Q_1\endbmatrix,$$
with
$Q_1,Q_2,Q_3$  linearly independent elements of $B_c$ and
$Q_1,Q_2$ relatively prime, then $\varphi\in N^{\text{\rm Bal}}_{c,c,c}$ and $\varphi\notin N^{\text{\rm Bal}}_{c,c}$.
\endremark

\demo{Proof} It suffices to establish the assertions in $\operatorname{BalH}_d$. One may then intersect with $\Bbb B\Bbb H_d$ to obtain the comparable result in $\Bbb B\Bbb H_d$. To establish (1), one must verify many inclusions; but each inclusion is completely straightforward. 
We prove (2). The inclusion $\operatorname{CO}^{\text{\rm Bal}}_{\#}\subseteq G\cdot N^{\text{\rm Bal}}_{\#}$ follows from (1). We now prove $G\cdot N^{\text{\rm Bal}}_{\#}\subseteq \operatorname{CO}^{\text{\rm Bal}}_{\#}$.  Fix $\varphi\in N^{\text{\rm Bal}}_{\#}$ for some $\#\in \operatorname{CP}$. We prove that there exists $\natural\in \operatorname{ECP}$ with $\natural\le \#$ and $\varphi\in \operatorname{DO}^{\text{\rm Bal}}_{\natural}$. 

If $\#$ is $c:c:c$, then either $\mu(I_1(\varphi))=3$ and $\varphi$ is already in $M^{\text{\rm Bal}}_{c:c:c}$ or  $\mu(I_1(\varphi))=2$  and $\varphi\in \operatorname{DO}^{\text{\rm Bal}}_{\mu_2}$ by the proof of Theorem \tref{orbitz}. 

If $\#$ is $c:c,c$, then 
$$\varphi=\bmatrix Q_1&*\\Q_3&*\\0&Q_3\endbmatrix$$ with $Q_1,Q_3$ linearly independent and $\mu(I_1(\varphi))\le 3$. If $\mu(I_1(\varphi))=2$, then $\varphi$ is in $\operatorname{DO}^{\text{\rm Bal}}_{\mu_2}$, as above. If $\mu(I_1(\varphi))= 3$, then $\varphi$ appears on the right side of (\tref{.q}) and the proof of Lemma \tref{L4q} shows that $\varphi$ is in $\operatorname{DO}_{c:c:c}$ or $\operatorname{DO}_{c:c,c}$.

If $\#$ is $c,c,c$, then $\mu(I_1(\varphi))\le 3$ and Theorem \tref{orbitz}, together with the chart of (2) in Lemma \tref{d=2c'} shows that $\varphi\in \operatorname{DO}^{\text{\rm Bal}}_{\natural}$ for $\natural$ equal to $\mu_2$, $\{c:c:c\}$, $\{c:c,c\}$, or $\{c,c,c\}$. In any event, $\natural\le \{c,c,c\}$.  

If $\#$ is $c:c$, then either $\mu(I_1(\varphi))= 4$ and $\varphi$ is in $M^{\text{\rm Bal}}_{c:c}$; or else, $\mu(I_1(\varphi))\le 3$ and $\varphi\in N^{\text{\rm Bal}}_{c:c,c}$. We have already shown that $N^{\text{\rm Bal}}_{c:c,c}\subseteq \operatorname{CO}^{\text{\rm Bal}}_{c:c,c}$. The definition of $\operatorname{CO}^{\text{\rm Bal}}_{\#}$ shows $\operatorname{CO}^{\text{\rm Bal}}_{c:c,c}\subseteq \operatorname{CO}^{\text{\rm Bal}}_{c:c}$, since $\{c:c,c\}<\{c:c\}$. 

If $\#$ is $c,c$, then $$\varphi=\bmatrix Q_1&Q_2\\Q_3&Q_4\\0&Q_5\endbmatrix,$$with $\dim_{\pmb k}{<}Q_3,Q_4,Q_5{>}\le 2$. If $\mu(I_1(\varphi))\le 3$, then $\varphi\in N^{\text{\rm Bal}}_{c,c,c}\subseteq \operatorname{CO}^{\text{\rm Bal}}_{c,c,c}\subseteq \operatorname{CO}^{\text{\rm Bal}}_{c,c}$. Henceforth, we assume $\mu(I_1(\varphi))=4$. We look at the proof of Lemma \tref{rest'}. If Case 1 is in effect, then $\varphi\in \operatorname{DO}^{\text{\rm Bal}}_{c,c}$. We assume that Case 2 is in effect. This forces $Q_2,Q_4,Q_5$ to be linearly independent. The hypothesis $\dim_{\pmb k}{<}Q_3,Q_4,Q_5{>}\le 2$ may be re-written as $Q_3\in{<}Q_4,Q_5{>}$. The hypothesis $\mu(I_1(\varphi))=4$ now forces 
$Q_1,Q_2,Q_4,Q_5$ to be linearly independent. We permute the rows and columns of $\varphi$ to obtain 
$$ g\varphi=\bmatrix Q_5&0\\Q_4&Q_3\\Q_2&Q_1\endbmatrix,$$ with $Q_3=\beta_1Q_5+\beta_2Q_4$. This is Case 2A from the proof of Lemma \tref{rest'}. We conclude $\varphi\in \operatorname{DO}^{\text{\rm Bal}}_{c,c}\cup \operatorname{DO}^{\text{\rm Bal}}_{c:c}$.

Take $\#$ to be $c$. If $\mu(I_1(\varphi))=5$, then $\varphi\in M^{\text{\rm Bal}}_{(c,\mu_5)}$ and $(c,\mu_5)<c$. If $\mu(I_1(\varphi))\le 4$, then the chart in (2) of Lemma \tref{d=2c'} shows that either $\varphi\in M^{\text{\rm Bal}}_{(\emptyset,\mu_4)}$; or else, $\varphi\in M^{\text{\rm Bal}}_{\natural}$ for some $\natural$ in $\operatorname{ECP}$ with $\natural <c$. On the other hand,
$$M^{\text{\rm Bal}}_{(\emptyset,\mu_4)}=\{\theta\in \operatorname{BalH}_d\mid \mu(I_1(\theta))=4\ \text{and $\theta$ does not have a generalized zero}\}.$$ One entry of $\varphi$ is zero; so, $\varphi\notin M^{\text{\rm Bal}}_{(\emptyset,\mu_4)}$ and $\varphi$ is in $\operatorname{CO}^{\text{\rm Bal}}_c$. (One could also argue that $\varphi\notin M^{\text{\rm Bal}}_{(\emptyset,\mu_4)}$ because if $C$ is the companion to $\varphi$ in the sense of (\tref{no1}), then $\mu(I_2(C))\le 5$ and this is too small for $\varphi$ to be in $ M^{\text{\rm Bal}}_{(\emptyset,\mu_4)}$.) 

Finally, if $\#$ is $\emptyset$, then Theorem \tref{orbitz} shows that $\varphi\in \operatorname{DO}^{\text{\rm Bal}}_{\natural}$ for some $\natural \in \operatorname{ECP}$. This $\natural$ automatically satisfies $\natural<\emptyset$.

The proof of (2) is complete. For (3), notice that each set $N^{\text{\rm Bal}}_{\#}$ and $N_{\#}$ has the form $V{}\cap {}\operatorname{BalH}_d$ or $V{}\cap{} \Bbb B\Bbb H_d$, where $V$ is a closed irreducible subset of the affine space $\Bbb H_d$. Indeed, for example, the condition $\dim_{\pmb k}{<}Q_3,Q_4,Q_5{>}\le 2$ is defined by the prime ideal generated by the maximal order minors of a generic $3\times (c+1)$ matrix. Furthermore, $\operatorname{BalH}_d$ and $\Bbb B\Bbb H_d$ are open subsets of $\Bbb H_d$, see Observation \tref{o5}.
\qed
\enddemo

Assume that the field $\pmb k$ is infinite  throughout the rest of the section. Proposition \tref{no-e'} shows that  $\operatorname{DO}_{\natural}$ is non-empty whenever it has a chance of being non-empty. In other words, $\operatorname{DO}_{\mu_2}$ is always empty; and, if an element of $\operatorname{DO}_{\natural}$ requires more linearly independent entries than are available in $B_c$, then $\operatorname{DO}_{\natural}$ is empty. Otherwise, $\operatorname{DO}_{\natural}$ is non-empty. To prove Proposition \tref{no-e'} we must show that the signed maximal order minors of a particular  $3\times 2$ matrix determine a birational parameterization of a curve. Lemma \tref{DC} provides an explicit  sufficient condition for establishing this birationality; recall that the spaces $\Bbb H_d$ and $\Bbb B\Bbb H_d$ are defined in Definition \tref{H}. We use the Avoidance Lemma \tref{avoid} repeatedly in the proof of Proposition \tref{no-e'}.

\proclaim{Lemma \tnum{DC}} Let $\pmb k$ be an infinite field and let $\varphi$ be an element of $\Bbb H_d$ with $\operatorname{ht} I_2(\varphi)=2$ and  $\mu(I_1(\varphi))\ge 3$. If two of the entries of $\varphi$ are $x^c$ and $y^{c-1}(x+y)$,  then $\varphi\in \Bbb B\Bbb H_d$.  \endproclaim 

\demo{Proof} We must show that the morphism, $\Psi\:\Bbb P^1\to \Bbb P^2$, given by the signed maximal order minors of $\varphi$, is   birational.  Let $\bar{\pmb k}$ be the algebraic closure of $\pmb k$. The field $\pmb k$ is infinite; so $\Psi$ is birational if and only if the induced morphism $\Psi_{\bar{\pmb k}}\:\Bbb P_{\bar{\pmb k}}^1\to \Bbb P_{\bar{\pmb k}}^2$  is birational. Throughout the rest of this proof, we assume that $\pmb k$ is algebraically closed. 

In the language of (\tref{PHI}), write $\Phi(\varphi)$ as $(g_1,g_2,g_3)$. Adopt the notation of Theorem \tref{2}. In particular, let $r$ be the degree of the field extension 
$$r=[\operatorname{Quot} \pmb k[B_d]:\operatorname{Quot} \pmb k[g_1,g_2,g_3]],$$ and $e$ be the multiplicity of the standard graded $\pmb k$-algebra $\pmb k[g_1,g_2,g_3]$. Theorem \tref{2} guarantees that $re=d$,  and $\Psi$ is birational if and only if $r=1$. Furthermore, according to Theorem \tref{2} and Observation \tref{2.1}, there exist homogeneous forms $f_1$ and $f_2$ in $B$ of degree $r$ so that every entry of $\varphi$ is in $\pmb k[f_1,f_2]$. Let $$\theta\:\pmb k[z_1,z_2]\to \pmb k[f_1,f_2]$$ be the $\pmb k$-algebra homomorphism with $\theta(z_i)=f_i$, as described in the proof of Observation \tref{2.1}. Select homogeneous forms $Q_1'$ and $Q_2'$ in $\pmb k[z_1,z_2]$  with  $\theta(Q_1')=x^c$ and $\theta(Q_2')=y^{c-1}(x+y)$.
Let $\varepsilon$ be the degree of the $Q_i'$. It follows that $c=\varepsilon r$. The equations
$$er=d=2c=2
\varepsilon r$$ yield $e=2\varepsilon$.
 Write $Q_i'=\prod_{j=1}^{\varepsilon}\ell_{i,j}$, for linear forms $\ell_{i,j}$ in $\pmb k[z_1,z_2]$. We notice that $\varepsilon\ge 2$. Indeed, if $\varepsilon$ were equal to $1$, then $r$ would equal $c$ and the entries of $\varphi$ would all live in the two-dimensional vector space ${<}f_1,f_2{>}$ and this would violate the hypothesis that $\mu(I_1(\varphi))\ge 3$. The  equations
$$x^c=\theta(Q_1')=\prod \ell_{1,j}(f_1,f_2)\quad\text{and}\quad y^{c-1}(x+y)=\theta(Q_2')=\prod \ell_{2,j}(f_1,f_2),$$which take place in the Unique Factorization Domain $B$, tell us that, after re-numbering the linear factors and adjusting their constants, we have
$$\matrix \format \r&\c&\l&\l\\x^r&{}={}&\ell_{1,j}(f_1,f_2),&\quad \text{for all $j$,}\\y^{r-1}(x+y)&{}={}&\ell_{2,1}(f_1,f_2),&\quad\text{and}\\y^r&{}={}&\ell_{2,j}(f_1,f_2),&\quad \text{for all $j\ge 2$}.\endmatrix $$The linear forms $\ell_{1,1}$, $\ell_{2,1}$, and $\ell_{2,2}$ in $\pmb k[z_1,z_2]$ look like 
$$\ell_{1,1}=\alpha z_1+\beta z_2,\quad \ell_{2,1}=\alpha' z_1+\beta' z_2,\quad\text{and}\quad \ell_{2,2}=\alpha'' z_1+\beta'' z_2$$ for constants $\alpha$, $\beta$, $\alpha'$, $\beta'$, $\alpha''$, and $\beta''$ in $\pmb k$. It follows that 
$$\bmatrix x^r\\ y^r\endbmatrix =\bmatrix \alpha &\beta \\\alpha'' &\beta''\endbmatrix \bmatrix f_1\\ f_2\endbmatrix.$$ The $2\times 2$ matrix in the above equation    is necessarily invertible. We also have that $y^{r-1}(x+y)$ is equal to 
$$\ell_{2,1}(f_1,f_2)= \bmatrix \alpha'&\beta'\endbmatrix  \bmatrix f_1\\f_2\endbmatrix =\bmatrix \alpha'&\beta'\endbmatrix \bmatrix \alpha &\beta \\\alpha'' &\beta''\endbmatrix ^{-1} \bmatrix x^r\\ y^r\endbmatrix =ax^r+by^r,$$ 
where $a$ and $b$ are the elements of $\pmb k$ which are defined by
$$\bmatrix a&b\endbmatrix= \bmatrix \alpha'&\beta'\endbmatrix \bmatrix \alpha &\beta \\\alpha'' &\beta''\endbmatrix ^{-1}.$$
The equation $y^{r-1}(x+y)=ax^r+by^r$  is impossible in $\pmb k[x,y]$, unless $r=1$. Thus, $r=1$, $\Psi$ is birational, and the proof is complete.  \qed \enddemo

\proclaim{Lemma \tnum{avoid}} Let $B$ be the polynomial ring $B=\pmb k[x,y]$, where $\pmb k$ is an infinite field, and $c$ be a positive integer. If $V$ is  proper subspace of $B_c$ and 
 $f$ is a non-zero homogeneous polynomial in $B$, then there exists a polynomial $Q$ in $B_c \setminus V$ with $Q$ and $f$ relatively prime.
\endproclaim

\demo{Proof} Let $\prod f_i$ be the factorization of $f$ into homogeneous irreducible factors in $B$.
  Each irreducible factor $f_i$  gives rise to the
proper subspace $V_i = B_{c-\deg f_i}f_i$ of $B_c$. (If $\deg f_i>\deg Q$, then the vector space $V_i$ is automatically equal to zero.) The field $\pmb k$ is infinite; so,
$V \cup \bigcup_{i} V_i$ is a proper subset of $B_c$. 
 Any $Q$ in the
complement of $V \cup \bigcup_{i} V_i$ has the desired
property. \qed\enddemo

\proclaim{Proposition \tnum{no-e'}} Let $\pmb k$ be an infinite field and $d=2c$ be an even integer. Fix $\natural\in \operatorname{ECP}\setminus\{\mu_2\}$.
\flushpar{\bf (1)} If $d=2$, then $\operatorname{DO}_{\natural}$ is empty.
\flushpar{\bf (2)} If $d=4$, then $\operatorname{DO}_{\natural}$ is non-empty if and only if $\natural\le  \{c,c,c\}$.
\flushpar{\bf (3)} If $d=6$, then $\operatorname{DO}_{\natural}$ is non-empty if and only if $ \natural\le (c,\mu_4)$ or $\natural=(\emptyset,\mu_4)$.
\flushpar{\bf (4)} If $d=8$, then $\operatorname{DO}_{\natural}$ is non-empty if and only if $\natural\le (\emptyset,\mu_5)$ or $\natural=(\emptyset,\mu_4)$
\flushpar{\bf (5)} If $10\le d$, then $\operatorname{DO}_{\natural}$ is non-empty.
\endproclaim
\demo{Proof} We first consider $\natural\in \operatorname{ECP}$ with $\natural<\emptyset$.  We show that there exist linearly independent $Q_1,Q_2,\dots $ in $B_c$   so that $\varphi_{\natural}$ is in $\operatorname{DO}_{\natural}$, where
$$\alignat 3\varphi_{c:c:c}&{}=\bmatrix Q_1&Q_2\\Q_3&Q_1\\0&Q_3\endbmatrix,\ 
\varphi_{c:c,c}&{}={}&\bmatrix Q_1&0\\Q_2&Q_3\\0&Q_2\endbmatrix,\ 
\varphi_{c,c,c}&{}={}&\bmatrix Q_1&Q_1\\Q_2&0\\0&Q_3\endbmatrix,\\ 
\varphi_{c:c}&{}=\bmatrix Q_1&Q_2\\Q_3&Q_4\\0&Q_3\endbmatrix,\ 
\varphi_{c,c}&{}={}&\bmatrix Q_1&Q_2\\Q_3&Q_3\\0&Q_4\endbmatrix,\ 
\varphi_{(c,\mu_4)}&{}={}&\bmatrix Q_1&Q_2\\Q_3&Q_1\\0&Q_4\endbmatrix,\\ 
\varphi_{(c,\mu_5)}&{}=\bmatrix Q_1&Q_2\\Q_3&Q_4\\0&Q_5\endbmatrix.
\endalignat$$
 For each $\natural$, take $Q_1=x^c$ and $Q_2=y^{c-1}(x+y)$. We show how to pick the rest of the $Q_i$. According to Lemma \tref{DC} we need only verify that $\operatorname{ht} (I_2(\varphi_{\natural}))=2$. We apply the Avoidance Lemma \tref{avoid} repeatedly. 

Take $\natural$ to be $\{c:c:c\}$, $\{c:c,c\}$, or $\{c,c,c\}$ . Pick $Q_3\in B_c$ so that $Q_3$ is not in  ${<}Q_1,Q_2{>}$  and $Q_3$ and $Q_1Q_2$ are relatively prime. Observe $\operatorname{ht} (I_2(\varphi_{\natural}))=2$.

Take $\natural$ to be $\{c:c\}$. Pick $Q_3\in B_c$ so that $Q_3\notin{<}Q_1,Q_2{>}$ and $Q_3$ and $Q_1$ are relatively prime. Pick $Q_4\in B_c$ so that $Q_4\notin{<}Q_1,Q_2,Q_3{>}$ and $Q_4$ and $Q_3$ are relatively prime.
Observe $\operatorname{ht} (I_2(\varphi_{\natural}))=2$.

Take $\natural$ to be $\{c,c\}$. Pick $Q_3\in B_c$ so $Q_3\notin{<}Q_1,Q_2{>}$ and $Q_3$ is relatively prime to $Q_1$. Pick $Q_4\in B_c$ so that $Q_4\notin{<}Q_1,Q_2,Q_3{>}$ and $Q_4$ is relatively prime to $Q_3(Q_1-Q_2)$.
Observe $\operatorname{ht} (I_2(\varphi_{\natural}))=2$.

Take $\natural$ to be $(c,\mu_4)$. Pick $Q_3\in B_c$ so $Q_3\notin{<}Q_1,Q_2{>}$ and $Q_3$ is relatively prime to $Q_1$. Pick $Q_4\in B_c$ so that $Q_4\notin{<}Q_1,Q_2,Q_3{>}$ and $Q_4$ is relatively prime to $Q_1^2-Q_2Q_3$.
Observe $\operatorname{ht} (I_2(\varphi_{\natural}))=2$.

Take $\natural$ to be $(c,\mu_5)$.  Pick $Q_3\in B_c$ so $Q_3\notin{<}Q_1,Q_2{>}$ and $Q_3$ is relatively prime to $Q_1$. Pick $Q_4\in B_c$ so that $Q_4\notin{<}Q_1,Q_2,Q_3{>}$ and $Q_4$ is relatively prime to $Q_1Q_3$. Pick
$Q_5\in B_c$ so that $Q_5\notin{<}Q_1,Q_2,Q_3,Q_4{>}$ and $Q_5$ is relatively prime to $Q_1(Q_1Q_4-Q_2Q_3)$. 
Observe $\operatorname{ht} (I_2(\varphi_{\natural}))=2$.

For the final three elements of $\operatorname{ECP}\setminus\{\operatorname{DO}_{\mu_2}\}$ we verify directly that the matrix $\varphi_{\natural}$, given below,  is in $\operatorname{DO}_{\natural}$:

$$\eightpoint \matrix\format \r&\c&\l\ &\r&\c&\l\\
\varphi_{(\emptyset,\mu_4)}&=&\bmatrix yx^{c-1}&x^c\\x^c&y^c\\y^c&y^{c-1}(x+y)\endbmatrix,&
\varphi_{(\emptyset,\mu_5)}&=&\bmatrix x^c&y^{c-2}(x^2+y^2)\\y^c&yx^{c-1}\\y^{c-1}(x+y)&x^c\endbmatrix,\ \text{and}\\\vspace{5pt} 
\varphi_{(\emptyset,\mu_6)}&=&\bmatrix x^{c-2}(x^2+y^2)&y^{c-2}(x^2+y^2)\\y^c&yx^{c-1}\\y^{c-1}(x+y)&x^c\endbmatrix.\endmatrix$$
The polynomials  $x^c$ and $y^{c-1}(x+y)$ each appear as entries  in each matrix. 
 Each matrix $\varphi_{\natural}$ has the correct form to be in $\operatorname{DO}_{\natural}$. It is not difficult to see that $\operatorname{ht} I_2(\varphi_{\natural})=2$ for each $\natural$. Notice when calculating $I_2(\varphi_{(\emptyset,\mu_i)})$ that one of the $2\times 2$ minors is equal to $\pm x^{c-1}y^{c+1}$. Indeed, when $i$ is $5$ or $6$, then the minor obtained by deleting row $1$ of $\varphi_{(\emptyset,\mu_i)}$ is
$$\left\vert\matrix y^c&yx^{c-1}\\y^{c-1}(x+y)&x^c\endmatrix\right\vert=x^{c-1}y^{c-1}\left\vert\matrix y&y\\x+y&x\endmatrix\right\vert=-x^{c-1}y^{c+1}.$$The comparable minor for $\natural=(\emptyset,\mu_4)$ involves rows 1 and 3.  If  $\natural=(\emptyset,\mu_4)$ then the appropriate entries of $\varphi_{\natural}$ are linearly independent provided $3\le c$. If $\natural= (\emptyset,\mu_5)$, then the appropriate entries of $\varphi_{\natural}$ are linearly independent provided $4\le c$. If
$5\le c$, then all entries of $\varphi_{\emptyset,\mu_6}$ are linearly independent. 
Apply Lemma \tref{DC}. \qed \enddemo

\bigpagebreak
\SectionNumber=5\tNumber=1
\heading Section \number\SectionNumber.  \quad The space of true triples of forms of degree $d$: the base point free locus, the birational locus, and the generic   Hilbert-Burch matrix.
\endheading 

In the previous sections we considered one rational curve at a time. At this point our attention turns to the  family of all rational plane curves of degree $d$. This family is parameterized by the space $\Bbb T_d$ of true triples of forms of degree $d$. The space $\Bbb T_d$ sits naturally in the affine space $\Bbb A_d$.

\definition{Definition \tnum{.'Org}}Let $\pmb k$ be a  field, $B$ be the polynomial ring  $B=\pmb k[x,y]$, and $d$ be a positive integer. Define 
$$\Bbb A_d=B_d\times B_d\times B_d.$$\enddefinition
\remark{\bf Remark \tnum{R5}} Each element of $\Bbb A_d$ is an ordered triple $\pmb g=(g_1,g_2,g_3)$ of homogeneous forms of degree $d$ from  the polynomial ring $B$. The  space $\Bbb A_d$ is isomorphic to affine space $\Bbb A^{3d+3}$. If $\pmb g\in \Bbb A_d$, then the corresponding element of $\Bbb A^{3d+3}$  is denoted $\lambda_{\pmb g}$; and if $\lambda\in \Bbb A^{3d+3}$, then the corresponding element of $\Bbb A_d$ is denoted $\pmb g_{\lambda}$. The correspondence is given as follows. If $\pmb g=(g_1,g_2,g_3)\in \Bbb A_d$, with $g_j=\sum\limits_{i=0}^{d}\lambda_{i,j}x^iy^{d-i}$, then
$$\lambda_{\pmb g}=(\lambda_{0,1},\dots,\lambda_{d,1},\lambda_{0,2},\dots,\lambda_{d,2},\lambda_{0,3},\dots,\lambda_{d,3}).$$
\endremark 

\definition{Definition \tnum{.'D27.13}}Fix an element $\pmb g=(g_1,g_2,g_3)\in \Bbb A_d$. Define      $I_{\pmb g}$ to be the ideal $(g_1,g_2,g_3)B$ of $B$ and
 $\Psi_{\pmb g}$ to be the morphism $$\Psi_{\pmb g}\:\Bbb P^1\setminus V(I_{\pmb g})\to \Bbb P^2,$$ which is given by
$$\Psi_{\pmb g}(q)=[g_1(q):g_2(q):g_3(q)],$$ for each point $q$ in $\Bbb P^1$,  where $V(I_{\pmb g})$ is the zero locus in $\Bbb P^1$ of the ideal $I_{\pmb g}$. Define the curve $\Cal C_{\pmb g}$ to be the closure of the image of $\Psi_{\pmb g}$, and define  $\pmb d_1(\pmb g)$ to be  the row vector $\bmatrix g_1&g_2&g_3\endbmatrix$. \enddefinition

\definition{Definition \tnum{.'D27.14}}Define  subsets $\operatorname{BPF}_d$, $\operatorname{Bir}_d$, and $\Bbb T_d$   of $\Bbb A_d$ as follows:
$$\matrix \format  \r&\c&\l\\ 
\operatorname{BPF}_d&{}={}&\{\pmb g\in \Bbb A_d\mid \text{the rational map $\Psi_{\pmb g}$ is base point free}\}\\
\operatorname{Bir}_d&{}={}&\{\pmb g\in \Bbb A_d\mid \text{the rational map $\Psi_{\pmb g}$ is birational}\}\\
\Bbb T_d&{}={}&\operatorname{BPF}_d\cap \operatorname{Bir}_d
=\left \{\pmb g\in \Bbb A_d\left\vert \matrix\format\l\\\text{$\Psi_{\pmb g}$ is a birational morphism with}\\\text{no base points}\endmatrix\right.\right\}.
\endmatrix$$
If $d$ is even, then define subsets $\operatorname{Bal}_d$, $\Bbb B_d$, and $\Bbb U\Bbb B_d$ of $\Bbb A_d$ as follows:
$$\matrix \format  \r&\c&\l\\
\operatorname{Bal}_d&{}={}&\left\{\pmb g\in \Bbb A_d\left\vert \matrix\format\l\\\text{every entry in a homogeneous  Hilbert-Burch}\\\text{matrix for $\pmb d_1(\pmb g)$ has degree $d/2$}\endmatrix\right.\right\}\\
\Bbb B_d&{}={}&\operatorname{Bal}_d{}\cap {}\Bbb T_d\\
\Bbb U\Bbb B_d&{}={}&\Bbb T_d\setminus \Bbb B_d.\endmatrix$$
\enddefinition

\remark{\bf Remark \tnum{UB}}We call $\Bbb T_d$ the space of {\bf true} triples of forms of degree $d$, $\Bbb B_d$ the space of  {\bf balanced} true triples of  forms of degree $d$, and 
$\Bbb U\Bbb B_d$   the space of {\bf unbalanced} true triples of forms of degree $d$.  The true triples of $\Bbb B_d$ are called ``balanced'' because the corresponding  Hilbert-Burch matrices are balanced in the sense that the column degrees $d_1$ and $d_2$ (in the language of Data \tref{34.1}) are equal; an unbalanced Hilbert-Burch matrix has $d_1<d_2$. Keep in mind that if $\pmb g$ is in $\Bbb A_d$, then
$$\pmb g\in \Bbb B_d \iff \matrix\format\l\\\text{the morphsim $\Psi_{\pmb g}$ is {\bf B}irational, {\bf B}ase point free and the}\\\text{Hilbert-Burch matrix for $\pmb d_1(\pmb g)$ is {\bf B}alanced}.\endmatrix$$
\endremark

In Theorems \tref{XXXX} and \tref{YYY}  we prove that  $\Bbb T_d$ and $\Bbb B_d$ are open subsets of $\Bbb A_d$.

In practice we are only interested in the open subset $\Bbb T_d$ of $\Bbb A_d$. Every element $\pmb g$ of $\Bbb A_d$ which is not in $\Bbb T_d$ corresponds to an unsuitable parameterization of the curve $\Cal C_{\pmb g}$.  Some of these unsuitable parameterizations have base points; others  are not birational. One can remove  base points by factoring out and removing the greatest common factor of the parameterizing forms. Also, if the parameterization is not birational, then one can reparameterize to find a birational parameterization; see, for example, Theorem \tref{2} or \cite{\rref{SWP}, Section 6.1}. Notice that if $\pmb g$ is in $\Bbb A_d$, then the data $\Psi_{\pmb g}$, $I_{\pmb g}$, and $\Cal C_{\pmb g}$ are uniquely determined by $\pmb g$; furthermore, 
$$\matrix\format\l\\ \text{the data $(\Psi_{\pmb g}$, $I_{\pmb g}$, $\Cal C_{\pmb g})$ satisfy the}\\\text{conditions and hypotheses of Data \tref{34.1}}\endmatrix \iff \pmb g\in \Bbb T_{d}.\tag\tnum{dta}$$
Let $\Bbb P\Bbb T_d$ represent the space of true parameterizations of plane curves of degree $d$; that is, $$\Bbb P\Bbb T_d=\left\{\text{morphisms  $\Psi\:\Bbb P^1\to \Bbb P^2$}\left\vert \matrix\format\l\\ \text{$\Psi$ is birational, base point free, and has}\\ \text{degree $d$}\endmatrix\right.\right\}.$$ We have established the following statement, where $\pmb k^*$ means $\pmb k\setminus\{0\}$ and $u(g_1,g_2,g_3)$ means $(ug_1,ug_2,ug_3)$.
\proclaim{Observation \tnum{O40}} The function $\Bbb T_d\to \Bbb P\Bbb T_d$, which is given by $\pmb g\mapsto \Psi_{\pmb g}$, is surjective and the fiber over $\Psi_{\pmb g}$ is $\{u\pmb g\mid u\in \pmb k^*\}$. \endproclaim

If $d=2c$ is an even integer, then the open subset $\Bbb B_d$ of $\Bbb T_d$ is of particular interest to us because 
one of the main topics of study in this paper is singularities of multiplicity $c$ on,  or infinitely near, curves of degree $d$. 
We proved in Corollary \tref{CTLL} that if such a singularity exists for the curve $\Cal C_{\pmb g}$, with $\pmb g\in \Bbb T_d$, then $\pmb g$ must be in $\Bbb B_d$. 

In the second half of this section we return  to the idea of (\tref{dta}). In order to have all of the data of \tref{34.1} one also must identify a Hilbert-Burch matrix for $\pmb d_1(\pmb g)$. It is not possible to define a unique Hilbert-Burch matrix for $\pmb d_1(\pmb g)$ as a function of $\pmb g$ over all of $\Bbb T_d$; however, if we restrict our attention to $\Bbb B_d$, then this is almost possible. In Corollary \tref{C5.30} we identify a complex $\Bbb F$ of free $(\pmb k[\pmb z])[x,y]$-modules, where $\pmb z$ is the $1\times (3d+3)$ matrix 
$$\pmb z=[z_{0,1},\dots,z_{d,1},z_{0,2},\dots,z_{d,2},z_{0,3},\dots,z_{d,3}]\tag\tnum{714}$$ 
of indeterminates and we prove that for $\lambda\in \Bbb A^{3d+3}$, 
$$\tsize \Bbb F\otimes_{\pmb k[\pmb z]}\frac {\pmb k[\pmb z]}{(\pmb z -\lambda)}\text{ is a resolution of $I_{\pmb g_{\lambda}}$} \iff \pmb g_{\lambda}\in \operatorname{Bal}_d.\tag\tnum{698}$$ The complex $\Bbb F$ would furnish a generic Hilbert-Burch matrix for $\pmb d_1(\pmb g)$ for ${\pmb g\in \operatorname{Bal}_d}$, except, unfortunately, the rank of $\Bbb F_2$ is three instead of two. There are three ways to interpret what we do get. First of all,  in Theorem \tref{T1}, we identify an open cover $\cup_{i=1}^3 \operatorname{Bal}_{d}^{(i)}$ of $\operatorname{Bal}_d$ and matrices $\pmb d_2^{(i)}$, for $1\le i\le 3$, in $(\pmb k[\pmb z])[x,y]$ such that if ${\pmb g\in \operatorname{Bal}_{d}^{(i)}}$, then $\pmb d_2^{(i)}\otimes_{\pmb k[\pmb z]} \frac {\pmb k[\pmb z]}{(\pmb z -\lambda_{\pmb g})}$ is a Hilbert-Burch matrix for $\pmb d_1(\pmb g)$. Secondly, in Corollary \tref{UPR} we identify a universal projective resolution $\Bbb U\Bbb P\Bbb R_{\Bbb Z}$ for all resolutions with graded Betti numbers
$$0\to \pmb k[x,y](-3c)^2\to \pmb k[x,y](-2c)^3\to \pmb k[x,y].$$
Do notice that the module in position $2$ of $\Bbb U\Bbb P\Bbb R_{\Bbb Z}$ is a projective module which is not necessarily free.
 Finally, in Section 4 we studied the morphism $\Phi\:\Bbb H_d \to \Bbb A_d$ which sends a $3\times 2$ matrix with entries from $\pmb k[x,y]_c$ to a triple $\pmb g$ of $\Bbb A_d$. One may restrict $\Phi$ to become
$$\Phi|\: \Phi^{-1}(\operatorname{Bal}_d)\to \operatorname{Bal}_d.\tag\tnum{Phi|}$$Corollary \tref{sectn} gives a local section of the morphism (\tref{Phi|}). 

Before we get to work we make one small observation about how the subsets 
of Definition \tref{.'D27.14}  fit together.

\proclaim{Observation \tnum{O5.22}} If $d=2c$ is an even integer, then $\operatorname{Bal}_d\subseteq \operatorname{BPF}_d$; however, in general, $\operatorname{Bal}_d\not\subseteq \operatorname{Bir}_d$.\endproclaim

\demo{Proof}The triple $\pmb g=(x^d,x^cy^c,y^d)$ of $\Bbb A_d$ is in $\operatorname{Bal}_d$, but does not correspond to a birational parameterization of the curve $T_2^2=T_1T_3$ unless $d=2$. On the other hand, if $\pmb g=(g_1,g_2,g_3)$ is any non-zero element of $\Bbb A_d$, then  the Hilbert-Burch Theorem shows that  the graded Betti numbers of $B/I$ are 
$$0\to B(-d-d_2)\oplus B(-d-d_1) \to B(-d)^3\to B,$$ for some non-negative integers $d_1$ and $d_2$ with $d_1+d_2+\deg\gcd(g_1,g_2,g_3)=d$. If  $\pmb g\in \operatorname{Bal}_d$, then $d_1+d_2=c+c=d$, the $\gcd$ of $(g_1,g_2,g_3)$ is a unit,  the height of $I_{\pmb g}$ is $2$, and $\Psi_{\pmb g}$ has no base points. \qed
\enddemo

Theorem \tref{L37.7} is our main tool for proving that subsets of $\Bbb A_d$ are closed. 
 We apply Theorem \tref{L37.7} in Theorem \tref{XXXX} when we prove that $\Bbb T_d$ is an open subset of $\Bbb A_d$; we also apply Theorem \tref{L37.7} in throughout  Section 6 when we decompose $\Bbb B_d$ into locally closed strata. 
The first assertion in the statement below  is well-known. The proof we offer for Theorem \tref{L37.7} is based on the Generic Freeness Lemma and is  inspired by the proof of assertion (1) which is given by Eisenbud \cite{\rref{E95}, Thm. 14.8}. We prove both assertions simultaneously.  If $p\in \operatorname{Spec} R$, then 
$$\text{$\pmb k(p)$ denotes the residue class field $R_p/pR_p$ of the local ring $R_p$.}\tag\tnum{rcf}$$

\proclaim{Theorem \tnum{L37.7}} Let $S$ be a standard graded Noetherian algebra over $S_0=R$. For non-negative integers $d$ and $e$,
define
$$\alignat 2 X(\ge d)&{}=&{}X(S;\ge d){}&=\{p\in \operatorname{Spec} R\mid \dim S\otimes_R\pmb k(p)\ge d\}\\
 X(= d)&{}=&{}X(S;= d){}&=\{p\in \operatorname{Spec} R\mid \dim S\otimes_R\pmb k(p)= d\}\\
X(=d,\ge e)&{}={}&X(S;=d,\ge e){}&=\left \{p\in \operatorname{Spec} R\left\vert \matrix \format\l\\\dim S\otimes_R\pmb k(p)= d \quad\text{and}\\ e(S\otimes_R\pmb k(p))\ge e\endmatrix\right.\right\}.\endalignat$$
Then 
\roster 
\item $X(S;\ge d)$ is a closed subset of $\operatorname{Spec} R$, and
\item $X(S;=d,\ge e)$ is a closed subset of $X(S;=d)$.\endroster
\endproclaim
\remark{Note} Once one has identified the ring $S$, then the degree zero component $S_0=R$ is automatically determined. For this reason we have denoted with $X(S;\ge d)$ the subset $X(\ge d)$ of $\operatorname{Spec} R$.\endremark 
\demo{Proof}  The proof proceeds by induction on the dimension of $R$. If $\dim R=0$, then $\operatorname{Spec} R$ is finite and every subset of $\operatorname{Spec} R$ is closed. Henceforth, we assume that $\dim R$ is positive. 

We next reduce to the case where $R$ is a domain. Let $\{p_1,\dots,p_s\}$ be the set of minimal prime ideals of $R$. Observe that
$$\alignat 2  X(\ge d)&{}={}&\bigcup\limits_{i=1}^s(X(\ge d)\cap V(p_i))&{}=\bigcup\limits_{i=1}^s(X(S\otimes_RR/p_i;\ge d), \ \text{and}\\
 X(=d,\ge e)&{}={}&\bigcup\limits_{i=1}^s(X(=d,\ge e)\cap V(p_i)){}&=\bigcup\limits_{i=1}^s(X(S\otimes_RR/p_i;=d,\ge e)).\endalignat$$
Once we prove the result for each $R/p_i\subseteq S\otimes_RR/p_i$, then we have also established the result for $R\subseteq S$.
Henceforth, we also assume that $R$ is a domain. 

Now we apply the Generic Freeness Lemma. The ring $R$ is a domain and $S$ is a finitely generated $R$-algebra; so there exists a non-zero element $a$ of $R$ such that $S_a=S\otimes_RR_a$ is a free $R_a$-module. We notice that $S_a$ is still a standard graded $R_a$-algebra with degree zero component equal to $R_a$; furthermore, each component $[S\otimes_RR_a]_i$ is finitely generated free $R_a$-module.

Let $K=\operatorname{Quot} (R)$. Take $p$ in the open subset $D(a)=\operatorname{Spec} R\setminus V(a)$, and take $q$ in $\operatorname{Spec} R$. (In this context, $V(a)$ is the closed subset ${\{p\in \operatorname{Spec} R\mid a\in p\}}$ of $\operatorname{Spec} R$.) We make the following observations about the Hilbert functions $\operatorname{H}_{S\otimes_R \pmb k(p)}(i)$ and $\operatorname{H}_{S\otimes_R \pmb k(q)}(i)$:
$$\split \operatorname{H}_{S\otimes_R \pmb k(p)}(i)&{}=\lambda_{R_p}[S\otimes_R \pmb k(p)]_i=\mu_{R_p}[S\otimes_RR_p]_i=\dim_K[S\otimes_RK]_i\\&{}=\mu_{K}[S\otimes_RK]_i\le \mu_{R_q}[S\otimes_RR_q]_i=\lambda_{R_q}[S\otimes \pmb k(q)]_i=\operatorname{H}_{S\otimes_R \pmb k(q)}(i).\endsplit$$The first and last equalities are the definition of Hilbert function. The second equality and second from last equalities follow from Nakayama's Lemma. The third equality is due to the fact that $[S\otimes_RR_p]_i$ is a free  $R_p$-module. The fourth equality and the inequality are obvious. We have shown that  

$$\matrix\format\l\\  \operatorname{H}_{S\otimes_R \pmb k(p)}(i)=\operatorname{H}_{S\otimes_R K}(i),\text{ for all }p\in D(a),\text{ and} \\\vspace{5pt}
\operatorname{H}_{S\otimes_R \pmb k(p)}(i)\le\operatorname{H}_{S\otimes_R \pmb k(q)}(i)\text{ for all }p\in D(a)\text{ and }q\in \operatorname{Spec} R. \endmatrix\tag{\tnum{up}} $$
Let $d_0=\dim S\otimes_RK$, and $e_0=e(S\otimes_RK)$. The Krull dimension of a standard graded algebra over a field may be read from its Hilbert function; hence, (\tref{up}) shows that
$$d_0=\dim S\otimes_R\pmb k(p)\le \dim S\otimes_R\pmb k(q)\text{ for all }p\in D(a)\text{ and }q\in \operatorname{Spec} R.$$
Furthermore, once the dimension of a standard graded algebra is determined, then its multiplicity may also be read from its Hilbert function. It follows that 
$$e_0=e(S\otimes_R\pmb k(p))\le e(S\otimes_R\pmb k(q))$$   for all  $p\in D(a)$ and $q\in \operatorname{Spec} R$  with $\dim S\otimes_R\pmb k(q)=d_0$.

The Krull dimension of  $R/(a)$ is less than the  Krull dimension of $R$ and so 
$$\split &V(a)\cap X(S;\ge d)=X(S\otimes_R R/(a);\ge d) \quad\text{and}\\ &V(a)\cap X(S;=d,\ge e)=X(S\otimes_R R/(a);=d,\ge e)\endsplit$$ are closed subsets of $\operatorname{Spec} R/(a)$ and $X(S\otimes_R R/(a);= d)$, respectively,  by the induction hypothesis. 

We first consider the case $d>d_0$. In this case, $X(S;\ge d)=X(S;\ge d)\cap V(a)$  is a closed subset of $\operatorname{Spec} R/(a)$; hence, a closed subset of $\operatorname{Spec} R$. Also in this case,  $X(S;=d,\ge e)=X(S;=d, \ge e)\cap V(a)$
is a closed subset of $$X(S\otimes_R R/(a);= d)=X(S;= d).$$

Finally, we consider the case $d=d_0$. In this case, $X(S;\ge d_0)=\operatorname{Spec} R$. If $e>e_0$, then $X(S;=d_0,\ge e)=X(S;=d_0,\ge e)\cap V(a)$ is a closed subset of
$$ X(S\otimes_R R/(a);= d_0)=
V(a)\cap (X(S;= d_0); $$and therefore $X(S;=d_0,\ge e)$ is a closed subset of
$X(S;= d_0)$. If $e=e_0$, then $X(S;=d_0,\ge e_0)=X(S;= d_0)$. \qed
 \enddemo

\definition{Conventions \tnum{CP}} Let $\pmb k$ be a field.\vphantom{\tnum{Gj}}   

\smallskip\flushpar{\bf(1)} We denote the coordinate ring of $\Bbb A^{3d+3}$ by $\pmb R=\pmb k[\pmb z]$, where $\pmb z$ is a $1\times (3d+3)$ matrix of indeterminates given in (\tref{714}).
 We continue to write $B$ for the standard graded polynomial ring $\pmb k[x,y]$. Let $\pmb S$ be the bi-graded polynomial ring  $\pmb S=\pmb k[x,y,\pmb z]$, where $x$ and $y$ each have bi-degree $(1,0)$, and $z_{i,j}$ has bi-degree $(0,1)$.

\smallskip\flushpar{\bf(2)} For each $\lambda\in \Bbb A^{3d+3}$, let $\pmb R_{\lambda}$ denote the ring $\pmb R/\frak m_\lambda=\pmb k(\frak m_\lambda)=\pmb k$ where $\frak m_\lambda$ is the maximal ideal $(\{z_{i,j}-\lambda_{i,j}\})$ of $\pmb R$. (Recall the notation $\pmb k(p)$ for residue class field from (\tref{rcf}).) If $S$ is an $\pmb R$-algebra, then  let $S_{\lambda}$ denote $S\otimes_{\pmb R}\pmb R_{\lambda}$ and if $G\in S$, then let $G|_{\lambda}$ (read as ``$G$ evaluated at $\lambda$'')  be the image of $G$ in $S_{\lambda}$.

\smallskip\flushpar{\bf(3)}  For each index $j$, with $1\le j\le 3$, let $G_j$ be the polynomial 
$$G_j=\sum_{i=0}^d z_{i,j}x^iy^{d-i} \tag\tref{Gj}$$ in $\pmb S$. If 
$\pmb g=(g_1,g_2,g_3)\in \Bbb A_d$, then, according to Convention (2) and Remark \tref{R5}, $G_j|_{\lambda_{\pmb g}}$ is equal to $g_j$ in $\pmb S_{\lambda_{\pmb g}}=B$.

\smallskip\flushpar{\bf(4)} We often apply Theorem  \tref{L37.7} to the affine space $\Bbb A_d$ by way of the identification of $\Bbb A_d$ and $\Bbb A^{3d+3}$ which is given in Remark \tref{R5}. Of course, the Zariski topology on $\Bbb A^{3d+3}$ is the same as the subspace topology that $\Bbb A^{3d+3}$ inherits as subset of  $\operatorname{MaxSpec} \pmb R\subseteq \operatorname{Spec} \pmb R$.  In this language,   Theorem  \tref{L37.7} shows that if $S$ is a standard graded Noetherian algebra over $S_0=\pmb R$ and $a$ and $b$ are non-negative integers then  
$$X(S;\ge a)\cap \Bbb A_{d}=\{\pmb g\in \Bbb A_d\mid \dim S_{\lambda_{\pmb g}}\ge a\}$$ is a closed subset of $\Bbb A_d$ and 
$$X(S; =a,\ge b)\cap \Bbb A_d=\{\pmb g\in \Bbb A_d\mid \dim S_{\lambda_{\pmb g}}= a\quad\text{and}\quad e(S_{\lambda_{\pmb g}})\ge b\}$$ is a closed subset of 
$X(S;=a)\cap \Bbb A_d=\{\pmb g\in \Bbb A_d\mid \dim S_{\lambda_{\pmb g}}= a\}$.
\enddefinition

\proclaim{Proposition \tnum{prop1}} The subset $\operatorname{BPF}_d$  of 
  $\Bbb A_d$ is open.\endproclaim
\demo{Proof} Let $\pmb R=\pmb k[\pmb z]$ be the coordinate ring of $\Bbb A^{3d+3}$, as described in item (1) of  Conventions  \tref{CP}, and $S$ be the $\pmb R$-algebra $\pmb R[x,y]/(G_1,G_2,G_3)$, where the polynomials $G_j$ are described in (\tref{Gj}). Fix an element $\pmb g\in \Bbb A_d$. Observe that $S_{\lambda_{\pmb g}}$ is equal to $\pmb k[x,y]/I_{\pmb g}$; furthermore,
$$\split \pmb g\in \Bbb A_d\setminus \operatorname{BPF}_d &{}\iff \text{$\Psi_{\pmb g}$ has base points}\iff \operatorname{ht} I_{\pmb g}\le 1\\&{}\iff \dim \pmb k[x,y]/I_{\pmb g}\ge 1\iff \dim S_{\lambda_{\pmb g}}\ge 1.\endsplit $$ Thus,
$$\Bbb A_d\setminus \operatorname{BPF}_d=
\{\pmb g\in \Bbb A_d\mid \dim S_{\lambda_{\pmb g}}\ge 1\}=X(S;\ge 1)\cap \Bbb A_d,
$$which is a closed subset of  $\Bbb A_d$ by Theorem \tref{L37.7} by way of item (4) of Conventions \tref{CP}.
 \qed
\enddemo

\proclaim{Proposition \tnum{prop2}} The subset $\operatorname{Bir}_d$  of 
  $\Bbb A_d$ is open.\endproclaim

\demo{Proof} Let $\pmb R=\pmb k[\pmb z]$ be the coordinate ring of $\Bbb A^{3d+3}$ as described in Conventions \tref{CP}. Form polynomials $G_1,G_2,G_3$ in $\pmb R[x,y]$ also as described in Conventions \tref{CP}. Consider the set   $$\Cal G=\{G_1^{i}G_2^jG_3^k\mid i+j+k=d-1\}$$ 
 of $\binom{d+1}2$ polynomials in $\pmb R[x,y]$.
Each element of $\Cal G$ is a polynomial of degree ${d(d-1)}$ in the variables $x,y$ with coefficients coming from the ring $\pmb R$. Let $Z$ be the ${(d(d-1)+1) \times  \binom{d+1}2}$ matrix which expresses the elements of $\Cal G$ in terms of the usual monomial basis for $\pmb k[x,y]_{d(d-1)}$, and let 
$\frak a$ be the ideal in $\pmb R$ generated by the maximal minors of $Z$; that is $\frak a=I_{\binom{d+1}2}(Z)$ . 

Take $\pmb g=(g_1,g_2,g_3)\in \Bbb A_d$. Item (3) of Conventions \tref{CP} shows that
 $G_1^iG_2^jG_3^k|_{\lambda_{\pmb g}}$ is equal to $g_1^ig_2^jg_3^k$. Observe that the polynomials
$$g_1^ig^j_2g_3^k,\quad \text{with $i+j+k=d-1$},$$ are linearly independent in $\pmb k[g_1,g_2,g_3]\subseteq \pmb k[x,y]$ if and only if $\lambda_{\pmb g}$ is not in  $V(\frak a)$.

  The ring $\pmb k[g_1,g_2,g_3]$ is the coordinate ring of the curve $\Cal C_{\pmb g}$. Furthermore, the ring homomorphism $\pmb k[T_1,T_2,T_3]\to \pmb k[g_1,g_2,g_3]$, which sends $T_i$ to $g_i$, induces an isomorphism
$$\frac {\pmb k[T_1,T_2,T_3]}{(f_{\pmb g})}\cong \pmb k[g_1,g_2,g_3],\tag\tnum{506}$$ where $f_{\pmb g}$ is the defining equation of the curve $\Cal C_{\pmb g}$. The degree of $f_{\pmb g}$ is equal to the multiplicity $e$ of $\pmb k[g_1,g_2,g_3]$ and $\Psi_{\pmb g}$ is a birational morphism if and only if $e=d$; see, for example, Theorem \tref{2}.  
Thus, 
$$\alignat 1 \text{$\Psi_{\pmb g}$ is birational}&{}\iff \deg f_{\pmb g}=d\\ 
&{}\iff \dim_{\pmb k} \left (\pmb k[T_1,T_2,T_3]/(f_{\pmb g})\right )_{d-1}\ge \dim_{\pmb k} \pmb k[T_1,T_2,T_3]_{d-1}\\
&{}\iff \dim_{\pmb k} \pmb k[g_1,g_2,g_3]_{d-1}\ge \tsize\binom{d+1}2\\
&{}\iff \text{the elements of $\Cal G|_{\lambda_{\pmb g}}$ are linearly independent}\\&\phantom{{}\iff{}}\text{in $\pmb k[g_1,g_2,g_3]\subseteq \pmb k[x,y]$}\\ 
&{}\iff \lambda_{\pmb g}\notin V(\frak a). \qed\endalignat$$ \enddemo

\proclaim{Theorem \tnum{XXXX}} The subset $\Bbb T_d$ of $\Bbb A_d$ is open.\endproclaim
\demo{Proof} The set $\Bbb T_d$ is equal to $\operatorname{BPF}_d\cap \operatorname{Bir}_d$.  
Apply Propositions \tref{prop1} and \tref{prop2}.  \qed\enddemo

Most of our work takes place over a field $\pmb k$. However, Theorem \tref{ZZ}, Corollary \tref{UPR}, and Corollary \tref{pop1} are   about  $3$ generic forms of the same even degree   in the polynomial ring $\Bbb Z[x,y]$. In order to facilitate the transition from working over a field to working over the integers, we set up   our data over $\Bbb Z$ right from the beginning. Every polynomial with coefficients in $\Bbb Z$ automatically represents a unique polynomial with coefficients in $\pmb k$.

\definition{Definition \tnum{data5}} Let $d=2c$ be a positive even integer and $\pmb S_{\Bbb Z}$
be the bi-graded polynomial ring  $\pmb S_{\Bbb Z}=\Bbb Z[x,y,\pmb z]$, where $\pmb z$   is the $1\times (3d+3)$ matrix of indeterminates given in (\tref{714}),  $\deg x=\deg y=(1,0)$, and $\deg z_{i,j}=(0,1)$.  For each index $j$, with $1\le j\le 3$, let $G_j$ be the bi-homogeneous polynomial of (\tref{Gj}) 
  in $\pmb S_{\Bbb Z}$ of degree ${(d,1)}$. Let $\Delta_{\Bbb Z}$ be the ring $\Delta_{\Bbb Z}= \pmb S_{\Bbb Z}/(G_1,G_2,G_3)$. 
Let $\pmb R_{\Bbb Z}$ be the standard graded polynomial ring   $\pmb R_{\Bbb Z}=\Bbb Z[\pmb z]$.
\vphantom{\tnum{Ai}}
\smallskip\flushpar{\bf (1)}   For each positive integer $i$, let $A_{\Bbb Z}^{(i)}$ be the $(d+i+1)\times 3(i+1)$ matrix 
$$A_{\Bbb Z}^{(i)}= \left [\smallmatrix 
z_{0,1}&      0&\cdots&      0&z_{0,2}&      0&\cdots&0      &z_{0,3}&0      &\cdots&0\\
z_{1,1}&z_{0,1}&\cdots&0      &z_{1,2}&z_{0,2}&\cdots&0      &z_{1,3}&z_{0,3}&\cdots&0\\
z_{2,1}&z_{1,1}&\cdots&0&z_{2,2}&z_{1,2}&\cdots&0&z_{2,3}&z_{1,3}&\cdots&0\\
\vdots&\vdots&\cdots&\vdots&\vdots&\vdots&\cdots&\vdots&\vdots&\vdots&\cdots&\vdots\\
z_{d-1,1}&z_{d-2,1}&\cdots&\vdots&z_{d-1,2}&z_{d-2,2}&\cdots&\vdots&z_{d-1,3}&z_{d-2,3}&\cdots&\vdots\\
z_{d,1}&z_{d-1,1}&\cdots&\vdots&z_{d,2}&z_{d-1,2}&\cdots&\vdots&z_{d,3}&z_{d-1,3}&\cdots&\vdots\\
0      &z_{d,1}&\cdots&\vdots&0      &z_{d,2}&\cdots&\vdots&0      &z_{d,3}&\cdots&\vdots\\
\vdots&\vdots&\cdots&\vdots&\vdots&\vdots&\cdots&\vdots&\vdots&\vdots&\cdots&\vdots\\
0      &0      &\cdots&z_{d,1}&0        &0      &\cdots&z_{d,2}&0      &0      &\cdots&z_{d,3} \endsmallmatrix\right ].\tag\tref{Ai}$$ Each polynomial $G_j$ contributes exactly $i+1$ columns to  $A_{\Bbb Z}^{(i)}$.

\smallskip\flushpar{\bf (2)} Let  $w_{\Bbb Z}$ be the determinant of   the $3c\times 3c$ matrix $A_{\Bbb Z}^{(c-1)}$; so, $w_{\Bbb Z}$ is a bi-homogeneous element of $\pmb S_{\Bbb Z}$ of degree $(0,3c)$. 

\smallskip\flushpar{\bf (3)} Let
$A_{\Bbb Z}$  be the $(3c+1)\times (3c+3)$ matrix $A_{\Bbb Z}^{(c)}$ of (\tref{Ai}). One may obtain $3c+3$ ``Eagon-Northcott'' relations on $A_{\Bbb Z}$, by crossing one column of $A_{\Bbb Z}$ at a time and computing the signed maximal minors of the resulting $(3c+1)\times (3c+2)$ matrix. In particular, when one crosses out columns 
$1$, $c+2$, or $2c+3$ of $A_{\Bbb Z}$, one obtains the relations 
$\pmb b^{(1)}_{\Bbb Z}$, $\pmb b^{(c+2)}_{\Bbb Z}$, and $\pmb b^{(2c+3)}_{\Bbb Z}$, on  $A_{\Bbb Z}$, which are given in Table 1,  where $A_{\Bbb Z}(i,j)$ is the determinant of the submatrix of $A_{\Bbb Z}$ which is obtained by deleting columns $i$ and $j$. Each $A_{\Bbb Z}(i,j)$ is a bi-homogeneous element of $\pmb S_{\Bbb Z}$  of degree $(0,3c+1)$.
 \topinsert  
$$ \bmatrix 0\\A_{\Bbb Z}(1,2)\\-A_{\Bbb Z}(1,3)\\\vdots\\(-1)^{c+1}A_{\Bbb Z}(1,c+1)\\\noalign{\hrule}\\(-1)^cA_{\Bbb Z}(1,c+2)\\(-1)^{c+1}A_{\Bbb Z}(1,c+3)\\\vdots\\ A_{\Bbb Z}(1,2c+2)\\\noalign{\hrule}\\-A_{\Bbb Z}(1,2c+3)\\A_{\Bbb Z}(1,2c+4)\\\vdots \\(-1)^{c+1}A_{\Bbb Z}(1,3c+3)\endbmatrix\hskip-2.51pt
 \bmatrix A_{\Bbb Z}(1,c+2)\\-A_{\Bbb Z}(2,c+2)\\A_{\Bbb Z}(3,c+2)\\\vdots\\ (-1)^cA_{\Bbb Z}(c+1,c+2)\\\noalign{\hrule}\\0\\(-1)^{c+1}A_{\Bbb Z}(c+2,c+3)\\
\vdots\\A_{\Bbb Z}(c+2,2c+2)\\\noalign{\hrule}\\-A_{\Bbb Z}(c+2,2c+3)\\A_{\Bbb Z}(c+2,2c+4)\\\vdots\\ (-1)^{c+1}A_{\Bbb Z}(c+2,3c+3)
\endbmatrix\hskip-2.51pt
  \bmatrix A_{\Bbb Z}(1,2c+3)\\-A_{\Bbb Z}(2,2c+3)\\A_{\Bbb Z}(3,2c+3)\\\vdots\\(-1)^cA_{\Bbb Z}(c+1,2c+3)\\\noalign{\hrule}\\(-1)^{c+1}A_{\Bbb Z}(c+2,2c+3)\\
(-1)^{c}A_{\Bbb Z}(c+3,2c+3)\\\vdots\\
 -A_{\Bbb Z}(2c+2,2c+3)\\\noalign{\hrule}\\0\\A_{\Bbb Z}(2c+3,2c+4)\\\vdots\\(-1)^{c+1}A_{\Bbb Z}(2c+3,3c+3)\endbmatrix $$
\centerline{{\bf Table 1:} The relations 
$\pmb b_{\Bbb Z}^{(1)}$, $\pmb b_{\Bbb Z}^{(c+2)}$, and $\pmb b_{\Bbb Z}^{(2c+3)}$, on  $A_{\Bbb Z}^{(c)}$ from (3) of Definition \tref{data5}.}\endinsert   

\smallskip\flushpar{\bf (4)} For each positive integer $i$, recall the $1\times (i+1)$ matrix  $\rho_{\Bbb Z}^{(i)}=[y^i,xy^{i-1},\dots,x^i]$ from (\tref{rho}) and let
$N^{(i)}_{\Bbb Z}$ be the $3\times 3(i+1)$ matrix 
$$N^{(i)}_{\Bbb Z}=\left [\matrix \rho_{\Bbb Z}^{(i)}&0&0\\0&\rho_{\Bbb Z}^{(i)}&0\\0&0&\rho_{\Bbb Z}^{(i)}\endbmatrix.$$Define $\pmb q_{1,\Bbb Z}=N_{\Bbb Z}^{(c)}\pmb b^{(1)}_{\Bbb Z}$, $\pmb q_{2,\Bbb Z}=N^{(c)}_{\Bbb Z}\pmb b^{(c+2)}_{\Bbb Z}$, and 
$\pmb q_{3,\Bbb Z}=N^{(c)}_{\Bbb Z}\pmb b^{(2c+3)}_{\Bbb Z}$. 
Each $\pmb q_{j,\Bbb Z}$ is a column vector with three entries and each entry of each $\pmb q_{j,\Bbb Z}$ is a bi-homogeneous element of $\pmb S_{\Bbb Z}$ of degree $(c,3c+1)$. 

\smallskip\flushpar{\bf (5)} Define the matrix $\pmb d_{2,\Bbb Z}$ to be  $\pmb d_{2,\Bbb Z}=\bmatrix \pmb q_{1,\Bbb Z}&\pmb q_{2,\Bbb Z}&\pmb q_{3,\Bbb Z}\endbmatrix$ and for each $j$, with $1\le j\le 3$, define $\pmb d^{(j)}_{2,\Bbb Z}$ to be $\pmb d_{2,\Bbb Z}$ with column $j$ removed. 

\smallskip\flushpar{\bf (6)} Define $\Bbb F_{\Bbb Z}$ to be the maps and modules 
   $$\Bbb F_{\Bbb Z}:\quad 0\to \pmb S_{\Bbb Z}(-3c,-3c-3)@>\pmb d_{3,\Bbb Z}>> \pmb S_{\Bbb Z}(-3c,-3c-2)^3 @> \pmb d_{2,\Bbb Z}>> \pmb S_{\Bbb Z}(-2c,-1)^3 @> \pmb d_{1,\Bbb Z}>>\pmb S_{\Bbb Z},$$
where $\pmb d_{1,\Bbb Z}=\bmatrix G_1&G_2&G_3\endbmatrix$, $\pmb d_{2,\Bbb Z}$ is given in (5), and $$\pmb d_{3,\Bbb Z}=\bmatrix z_{0,1}\\(-1)^{c+1}z_{0,2}\\z_{0,3}\endbmatrix.$$  
\enddefinition

\remark{Remark} If $\pmb k$ is a field and $\pmb S_{\Bbb Z}$ is the polynomial ring of Definition \tref{data5}, then $\pmb k\otimes_{\Bbb Z}\pmb S_{\Bbb Z}$ is equal to the polynomial ring $\pmb S=\pmb k[x,y,\pmb z]$ of Conventions \tref{CP}.\endremark 

\definition{Definition \tnum{seg}}Fix  a field $\pmb k$. Let $\Delta$ be the ring $\pmb k\otimes_{\Bbb Z}\Delta_{\Bbb Z}=\pmb S/(G_1,G_2,G_3)$.
If $s_{\Bbb Z}$ is an element of $\pmb S_{\Bbb Z}$, then write $s$ for the image   of $s_{\Bbb Z}$ under the natural homomorphism $\pmb S_{\Bbb Z}\to \pmb S$; and if $M_{\Bbb Z}$ is a matrix with entries in $\pmb S_{\Bbb Z}$, then write $M$ for the corresponding matrix with entries in $\pmb S$. In particular, the elements $w$ and $A(i,j)$ of $\pmb S$ are obtained from the elements $w_{\Bbb Z}$ and $A_{\Bbb Z}(i,j)$ of $\pmb S_{\Bbb Z}$ in this manner,  the matrices $A^{(i)}$, $\pmb d_i$, and $\pmb d_2^{(j)}$, with entries in $\pmb S$, are obtained from the matrices $A_{\Bbb Z}^{(i)}$, $\pmb d_{i,\Bbb Z}$, and $\pmb d^{(j)}_{2,\Bbb Z}$, with entries in $\pmb S_{\Bbb Z}$, in this manner, and the maps and modules $\pmb k\otimes_{\Bbb Z}\Bbb F_{\Bbb Z}$ are written as 
$$  \Bbb F: \quad 0\to   \pmb S(-3c,-3c-3)@>\pmb d_3>> \pmb S(-3c,-3c-2)^3 @> \pmb d_2>> \pmb S(-2c,-1)^3 @> \pmb d_1>>\pmb S.  $$ 
\enddefinition

\proclaim{Theorem \tnum{YYY}} Retain the notation and Conventions of {\rm\tref{CP}}. If $d=2c$ is an even integer and $w$ in $\pmb R$ is given in {\rm(3)} of Definition {\rm\tref{data5}} by way of Definition {\rm\tref{seg}}, then  $\Bbb A_d\setminus \operatorname{Bal} _d$ is the closed subset 
$\{\pmb g\in \Bbb A_d\mid \lambda_{\pmb g}\in V(w)\}$ of $\Bbb A_d$.
\endproclaim

\remark{Remark} It is an immediate consequence of Theorems \tref{YYY} and \tref{XXXX} that $\operatorname{Bal}_d$ and $\Bbb B_d$ are open subsets of $\Bbb A_d$ and  that $\Bbb U\Bbb B_d$ is a hypersurface section of $\Bbb T_d$.\endremark

\demo{Proof}   
Recall $\rho_{\Bbb Z}^{(i)}$, $A_{\Bbb Z}^{(i)}$, $\pmb d_{1,\Bbb Z}$, and  $N_{\Bbb Z}^{(i)}$ from Definition \tref{data5}. 
Matrix multiplication yields that 
$$\rho_{\Bbb Z}^{(d+i)}A_{\Bbb Z}^{(i)}=\pmb d_{1,\Bbb Z} N_{\Bbb Z}^{(i)},\tag\tnum{key}$$
over $\pmb S_{\Bbb Z}$, for each $i$.
Let $\pmb q$ be a $3\times 1$ matrix  of forms from $B=\pmb k[x,y]$ of degree $i$. Then $\pmb q=N^{(i)} \pmb b$ for some  $(3i+3)\times 1$ matrix of scalars $\pmb b$.   Take the image of (\tref{key}) under the homomorphism $\pmb S_{\Bbb Z}\to \pmb S$ to
 see that  $$\pmb d_1 \pmb q=\pmb d_1 N^{(i)} \pmb b=\rho^{(d+i)}A^{(i)}\pmb b\tag\tnum{oneup}$$  in $\pmb S$. Fix $\pmb g=(g_1,g_2,g_3)\in \Bbb A_d$. Apply the homomorphism $\pmb S\to \pmb S_{\lambda_{\pmb g}}$ to  (\tref{oneup}) to see that  $$\bmatrix g_1&g_2&g_3\endbmatrix\pmb q=0\text{ in }B\iff A^{(i)}|_{\lambda_{\pmb g}}\pmb b=0\text{ in } \pmb k.\tag\tnum{ky}$$ Recall that 
 $w$ is the determinant of   the $3c\times 3c$ matrix $A^{(c-1)}$.   Apply (\tref{ky}), with $i=c-1$, to see that
$$  \split \pmb g\in \Bbb A_d\setminus \operatorname{Bal}_d\iff{}&\text{there exists a non-zero $3\times 1$ matrix $\pmb q$ of   forms of degree $c-1$}\\ &\text{from $B$ with $\bmatrix g_1&g_2&g_3\endbmatrix\pmb q=0$}\hfill\\ {}\iff{}&\text{there exists a non-zero $3c\times 1$ matrix $\pmb b$ of constants with}\\ &A^{(c-1)}|_{\lambda_{\pmb g}} \pmb b=0\hfill\\
{}\iff{}&w|_{\lambda_{\pmb g}}=0.  \qed \hfill \endsplit$$
\enddemo

Now that we have shown $\operatorname{Bal}_d$ and $\Bbb B_d$ to be open subsets of $\Bbb A_d$ it is time to show that the subsets $\operatorname{BalH}_d$ and $\Bbb B\Bbb H_d$ of $\Bbb H_d$ are also open. These subsets were  introduced in  Definition \tref{H}.

\definition{Definition \tnum{Phi!}} Let $\pmb k$ be a field and $d=2c$ be an even integer. Define $$\Phi\:\Bbb H_d\to \Bbb A_d$$ to be the morphism of affine varieties which sends a matrix $\varphi$ in $\Bbb H_d$ to the ordered triple $\Phi(\varphi)$ of  signed maximal order minors of $\varphi$ as given  in (\tref{PHI}).\enddefinition
\proclaim{Observation \tnum{o5}} The subsets  $\operatorname{BalH}_d$ and $\Bbb B\Bbb H_d$ of $\Bbb H_d$ are   open. \endproclaim
\demo{Proof} Observe that $\Phi^{-1}(\operatorname{Bal}_d)=\operatorname{BalH}_d$ and $\Phi^{-1}(\Bbb B_d)=\Bbb B\Bbb H_d$. The subsets $\operatorname{Bal}_d$ and $\Bbb B_d$ of $\Bbb A_d$ are open by Theorem \tref{YYY}. The proof is complete because the morphism $\Phi$ is continuous. \qed\enddemo 

\proclaim{Theorem \tnum{T1}} Retain the notation of Conventions {\rm\tref{CP}} and let  $\pmb d^{(1)}_2,\pmb d^{(2)}_2,\pmb d^{(3)}_2$ be the matrices in $\pmb k[\pmb z,x,y]$ from {\rm(5)} of Definition {\rm\tref{data5}}, by way of Definition {\rm\tref{seg}}. Then there exists an open cover $\operatorname{Bal}_{d}^{(1)}\cup \operatorname{Bal}_{d}^{(2)}\cup\operatorname{Bal}_{d}^{(3)}$ of $\operatorname{Bal}_d$   such that if $\pmb g\in \operatorname{Bal}_{d}^{(j)}$, then $\pmb d_2^{(j)}|_{\lambda_{\pmb g}}$ is a Hilbert-Burch matrix for $\pmb d_1(\pmb g)$.\endproclaim

\remark{Remark} Define $\Bbb B_{d}^{(j)}= \operatorname{Bal}_d^{(j)}{}\cap{} \Bbb B_d$. Theorem \tref{T1} also asserts that 
there exists an open cover $\Bbb B_{d}^{(1)}\cup \Bbb B_{d}^{(2)}\cup\Bbb B_{d}^{(3)}$ of $\Bbb B_d$ such that if $\pmb g\in \Bbb B_{d}^{(j)}$, then $\pmb d_2^{(j)}|_{\lambda_{\pmb g}}$ is a Hilbert-Burch matrix for $\pmb d_1(\pmb g)$.\endremark

\demo{Proof} The set $U_j=\{\lambda\in \Bbb A^{3d+3}\mid \lambda_{0,j}\neq 0\}$ is open in $\Bbb A^{3d+3}$ for $1\le j\le 3$; and therefore,
$$\operatorname{Bal}_{d}^{(j)}=\operatorname{Bal}_d{}\cap {}U_j=\{\pmb g\in \operatorname{Bal}_d\mid \lambda_{\pmb g}\in U_j\}$$ is open in $\operatorname{Bal}_d$. Furthermore, if $\pmb g\in \operatorname{Bal}_d$, then the ideal $I_{\pmb g}$ has height $2$ by Observation \tref{O5.22}; hence $I_{\pmb g}\not\subseteq (x)$ and $\pmb g$ is in $\operatorname{Bal}_{d}^{(j)}$ for some $j$. We have established that $\operatorname{Bal}_{d}^{(1)}\cup \operatorname{Bal}_{d}^{(2)}\cup\operatorname{Bal}_{d}^{(3)}$ is an open cover of $\operatorname{Bal}_d$.

The construction of $\pmb d_{2,\Bbb Z}$ in Definition \tref{data5} guarantees that $$\pmb d_{1,\Bbb Z}\pmb d_{2,\Bbb Z}=0.\tag\tnum{cx}$$ Indeed, the  definition of $\pmb q_{j,\Bbb Z}$ from (4) of Definition \tref{data5}, together with (\tref{key}) and the definition of  $\pmb b_{\Bbb Z}^{(j+(j-1)c)}$   from (3) of Definition \tref{data5}, shows that
$$ \pmb d_{1,\Bbb Z} \pmb q_{j,\Bbb Z}= \pmb d_{1,\Bbb Z} N_{\Bbb Z}^{(c)} \pmb b_{\Bbb Z}^{(j+(j-1)c)}=\rho_{\Bbb Z}^{(d+c)}A_{\Bbb Z}^{(c)}\pmb b_{\Bbb Z}^{(j+(j-1)c)}=0.$$The equation (\tref{cx}) is preserved under the homomorphism $\pmb S_{\Bbb Z}\to\pmb S$; hence, the composition $\pmb d_{1}\pmb d_{2}$ is zero. 
Let $\pmb g\in \Bbb A_d$. Evaluate $\pmb d_1\pmb d_2=0$  at $\lambda_{\pmb g}$ to see that $\pmb d_1(\pmb g)\pmb d_2|_{\lambda_{\pmb g}}=0$    and therefore,
 $\pmb d_1(\pmb g)\pmb d_2^{(j)}|_{\lambda_{\pmb g}}=0$   for all  $\pmb g\in \Bbb A_d$ and all $j$.

Fix $\pmb g\in \operatorname{Bal}_{d}^{(j)}$. We show that $\pmb d_2^{(j)}|_{\lambda_{\pmb g}}$ is a Hilbert-Burch matrix for $\pmb d_1(\pmb g)$. 
The fact that $\pmb g\in \operatorname{Bal}_{d}^{(j)}\subseteq \operatorname{Bal}_d$ ensures that a minimal resolution of $\pmb R/I_{\pmb g}$ looks like
$$0\to \pmb R(-3c)^2 @>\varphi>> \pmb R(-2c)^3 @> \pmb d_1(\pmb g) >> \pmb R.$$To identify $\varphi$, we need only find two linearly independent degree $c$ relations on $\pmb d_1(\pmb g)$.
We already have seen that each of the  two columns of $\pmb d_2^{(j)}|_{\lambda_{\pmb g}}$ is a  degree $c$ relation on $\pmb d_1(\pmb g)$. To show that  $\pmb d_2^{(j)}|_{\lambda_{\pmb g}}$ is a Hilbert-Burch matrix for $\pmb d_1(\pmb g)$, we need only see that the two columns of  $\pmb d_2^{(j)}|_{\lambda_{\pmb g}}$ are linearly independent. Look at $\pmb d_{2,\Bbb Z}$ after $x$ has been set equal to zero:
$$\pmb d_{2,\Bbb Z}\equiv \bmatrix 0&A_{\Bbb Z}(1,c+2)y^c&A_{\Bbb Z}(1,2c+3)y^c\\
(-1)^c A_{\Bbb Z}(1,c+2)y^c&0&(-1)^{c+1} A_{\Bbb Z}(c+2,2c+3)y^c\\
- A_{\Bbb Z}(1,2c+3)y^c&- A_{\Bbb Z}(c+2,2c+3)y^c&0\endbmatrix$$ mod $(x)$.
Expand the minors $A_{\Bbb Z}(c+2,2c+3)$, $A_{\Bbb Z}(1,2c+3)$, and $A_{\Bbb Z}(1,c+2)$ of $A_{\Bbb Z}$ across the first row to see that
$$A_{\Bbb Z}(c+2,2c+3)=z_{0,1}w_{\Bbb Z},\ A_{\Bbb Z}(1,2c+3)=(-1)^cz_{0,2}w_{\Bbb Z},\ \text{and}\ A_{\Bbb Z}(1,c+2)=z_{0,3}w_{\Bbb Z}.$$ (Recall that $A_{\Bbb Z}$ is given in (\tref{Ai}) as $A_{\Bbb Z}^{(c)}$ and $w_{\Bbb Z}=\det A_{\Bbb Z}^{(c-1)}$ is defined in (2) of Definition \tref{data5}.) Thus, 
$$\pmb d_{2,\Bbb Z}\equiv y^c\bmatrix 0&w_{\Bbb Z}z_{0,3}&(-1)^cw_{\Bbb Z}z_{0,2}\\
(-1)^c w_{\Bbb Z}z_{0,3}&0&(-1)^{c+1} w_{\Bbb Z}z_{0,1}\\
(-1)^{c+1} w_{\Bbb Z}z_{0,2}&- w_{\Bbb Z}z_{0,1}&0\endbmatrix \mod (x).\tag\tnum{511}$$
We know from Theorem \tref{YYY} and the definition of $\operatorname{Bal}_{d}^{(j)}$ that $(wz_{0,j})|_{\lambda_{\pmb g}}\neq 0$. 
The matrix $\pmb d_2^{(j)}$ is defined to be $\pmb d_2$ with column $j$ removed. It is now clear that the columns of 
$\pmb d_2^{(j)}|_{\lambda_{\pmb g}}\mod(x)$ are linearly independent; hence, the columns of $\pmb d_2^{(j)}|_{\lambda_{\pmb g}}$ are linearly independent and $\pmb d_2^{(j)}|_{\lambda_{\pmb g}}$ is a Hilbert-Burch matrix for $\pmb d_1(\pmb g)$. \qed \enddemo

\remark{\bf Remark \tnum{R5.44}}Define $C^{(j)}$ and $A^{(j)}$ to be  the bi-homogeneous matrices which satisfy (\tref{no1}) and  (\tref{no2})  for $\pmb d_2^{(j)}$. In other words, $C^{(j)}$ is a $(c+1)\times 2$ matrix and each entry is a bi-homogeneous element of $\pmb k[\pmb z,\pmb T]$ of degree $3c+1$ in $\pmb z$ and degree $1$ in $\pmb T$. Furthermore, $A^{(j)}$ is a $(c+1)\times 3$ matrix and each entry is bi-homogeneous element in $\pmb k[\pmb z,\pmb u]$  of degree $3c+1$ in $\pmb z$ and degree $1$ in $\pmb u$. These matrices satisfy 
$$\pmb T\pmb d_2^{(j)}= \rho^{(c)}C^{(j)}\quad\text{and}\quad C^{(j)}\pmb u^{\text{\rm T}}=A^{(j)}\pmb T^{\text{\rm T}}.$$ Theorem \tref{T1} shows that if $\pmb g\in \Bbb B_d^{(j)}$, then $\pmb d_2^{(j)}|_{\lambda_{\pmb g}}$ is a Hilbert-Burch matrix for $\pmb d_1(\pmb g)$; hence, the pair of matrices $(C^{(j)}|_{\lambda_{\pmb g}},A^{(j)}|_{\lambda_{\pmb g}})$ complete the requirements of Data \tref{biData} for $\pmb d_2^{(j)}|_{\lambda_{\pmb g}}$ and Corollary \tref{CorND}
may be applied to the ring $\pmb k[\pmb u]/I_3(A^{(j)}|_{\lambda_{\pmb g}})$ in order to determine 
the configuration of multiplicity $c$ singularities on, or infinitely near, the  curve $\Cal C_{\pmb g}$. This analysis is carried out in Theorem \tref{strata}.\endremark 

\medskip
The rest of this section is devoted to giving alternate interpretations of Theorem~\tref{T1}.   In Corollary~\tref{C5.30} we prove that the maps and modules of $\Bbb F$, from Definition \tref{seg}, behave as promised in (\tref{698}). The ideas of flat families of algebras then allow us to prove, in Corollary \tref{C5.31}, that the localization $\Bbb F_w$ is a resolution, at least when the ambient field is algebraically closed. In Theorem \tref{ZZ}, we use Corollary \tref{C5.31}, which holds  over an algebraically closed field, together with an ``evaluate the constant'' argument, to identify the {\bf exact} relationship, over $\Bbb Z$, between the maximal order minors of $\pmb d_{2,\Bbb Z}^{(j)}$ and the ``generic $d$-forms'' $G_1$, $G_2$, $G_3$ in $\pmb S_{\Bbb Z}=\Bbb Z[x,y,\pmb z]$. One could    obtain this result using identities involving determinants. In Theorem \tref{ZZ}, we also prove that $(\Bbb F_{\Bbb Z})_{w_{\Bbb Z}}$ is a resolution and that the ideal $(G_1,G_2,G_3)(\pmb S_{\Bbb Z})_{w_{\Bbb Z}}$, which is defined over $\Bbb Z$, is perfect of grade $2$. At this point, one could recover Corollary \tref{C5.31} from Theorem \tref{ZZ} if one gives a direct argument for assertion (1) in Theorem \tref{ZZ} in place of the ``evaluate the constant'' argument that we give. 

\proclaim{Observation \tnum{O5.99}} The maps and modules of $\Bbb F_{\Bbb Z}$ from {\rm(6)} of Definition {\rm\tref{data5}} are  a bi-homogeneous complex of $\pmb S_{\Bbb Z}$-modules.\endproclaim
\remark{Remark} It follows that $\Bbb F$ from Definition \tref{seg} is a bi-homogeneous complex of free $\pmb S$-modules.  \endremark
\demo{Proof}We saw in (\tref{cx}) that the composition $\pmb d_{1,\Bbb Z}\pmb d_{2,\Bbb Z}$ is zero. The Eagon-Northcott complex which is associated to the map $A_{\Bbb Z}=A_{\Bbb Z}^{(c)}$ of (\tref{Ai}) looks like :
$$0\to \pmb R_{\Bbb Z}^{3c+1}\otimes \bigwedge^{3c+3} \pmb R_{\Bbb Z}^{3c+3} \to \bigwedge^{3c+2}\pmb R_{\Bbb Z}^{3c+3} \to \pmb R_{\Bbb Z}^{3c+3} @> A>> \pmb R_{\Bbb Z}^{3c+1}.$$ Each row of $A_{\Bbb Z}$ gives rise to a relation on the $3c+3$ Eagon-Northcott relations on $A_{\Bbb Z}$. In particular, the top row of $A_{\Bbb Z}$ gives rise to the relation
$$z_{0,1} \pmb b_{\Bbb Z}^{(1)}+(-1)^{c+1}z_{0,2}\pmb b_{\Bbb Z}^{(c+2)}+z_{0,3} \pmb b_{\Bbb Z}^{(2c+3)}=0,\tag\tnum{zb}$$ for $\pmb b_{\Bbb Z}^{(j)}$ as defined in (3) of Definition \tref{data5}.  Multiply (\tref{zb}) by $N_{\Bbb Z}^{(c)}$ and use (4) of Definition \tref{data5} to see that    $\pmb d_{3,\Bbb Z}$ is a non-trivial relation on $\pmb d_{2,\Bbb Z}$. \qed \enddemo

\proclaim{Corollary \tnum{C5.30}} Retain the notation and hypotheses of Conventions {\rm \tref{CP}} and let $\Bbb F$ be the complex of free $S$-modules which is defined in Definition {\rm \tref{seg}}. Then  $\Bbb F\otimes_{\pmb R}\pmb R_{\lambda_{\pmb g}}$ is a resolution of $B/I_{\pmb g}$ for all
$\pmb g\in \operatorname{Bal}_d$.\endproclaim

\remark{Remark} If $\pmb g\in \Bbb A_d\setminus \operatorname{Bal}_d$, then the complex $\Bbb F\otimes_{\pmb R}\pmb R_{\lambda_g}$ is not a resolution of $B/I_{\pmb g}$ because it  has the wrong graded Betti numbers. Thus, this Corollary  establishes the claim made in (\tref{698}).\endremark

\demo{Proof}   We saw in Theorem \tref{T1} that $\{\operatorname{Bal}_{d}^{(j)}\}$ is an open cover of $\operatorname{Bal}_d$. Fix an element $\pmb g$ of $\operatorname{Bal}_{d}^{(j)}$. Theorem \tref{T1} also shows that 
$$0\to B(-3c)^2@> \pmb d_2^{(j)}|_{\lambda_{\pmb g}} >> B(-2c)^3 @> \pmb d_1(\pmb g)>>B$$ is a resolution of $B/I_{\pmb g}$. Observe that $\pmb d_1\otimes_{\pmb R} \pmb R_{\lambda_{\pmb g}}=\pmb d_1(\pmb g)$, two of the columns of $\pmb d_2\otimes_{\pmb R} \pmb R_{\lambda_{\pmb g}}$ form $\pmb d_2^{(j)}|_{\lambda_{\pmb g}}$, and $\pmb d_3\otimes_{\pmb R} \pmb R_{\lambda_{\pmb g}}$ is the relation 
on the columns of $\pmb d_2\otimes_{\pmb R} \pmb R_{\lambda_{\pmb g}}$ which expresses  
  the redundant column   in terms of the columns of $\pmb d_2^{(j)}|_{\lambda_{\pmb g}}$. In other words, $\Bbb F\otimes_{\pmb R}\pmb R_{\lambda_{\pmb g}}$ is a (non-minimal) resolution of $B/I_{\pmb g}$. \qed \enddemo

\proclaim{Corollary \tnum{C5.31}} Retain the notation and hypotheses of Conventions {\rm \tref{CP}} with $\pmb k$ an algebraically closed field. Let $w$, $\Bbb F$, and $\Delta$ be the element of $\pmb S$, the complex of free $\pmb S$-modules, and the quotient ring of $\pmb S$ which are defined in Definition {\rm\tref{seg}}. 
Then the  following statements hold.
\roster
\item"{(1)}"  The ring homomorphism $\pmb R_w\to \Delta_w$,  which is the composition $$\pmb R_w\to\pmb S_w\to  \Delta_w$$ of inclusion followed by the natural quotient map, is flat.
\item"{(2)}"   The ring $\Delta_w$ is  Cohen-Macaulay. 
\item"{(3)}" The complex $\Bbb F_w$ is a resolution of $\Delta_w$ by free $\pmb S_w$-modules.
\endroster
\endproclaim 

\demo{Proof}
Recall the maximal ideal $\frak m_{\lambda}=(\{z_{i,j}-\lambda_{i,j}\})$ of $\pmb R$ which is defined in (2) of Conventions \tref{CP}. The field $\pmb k$ is algebraically closed; so $\operatorname{MaxSpec} \pmb R_w$ is identified with $\Bbb A^{3d+3}\setminus V(w)$. To show assertion (1) it suffices to show that the ring homomorphism $\pmb R_{\frak m_{\lambda}}\to \Delta_{\frak m_{\lambda}}$ is flat for all 
$\frak m_{\lambda}$ in $\operatorname{MaxSpec} \pmb R_w$. Fix one such $\frak m_{\lambda}$. The ring $\Delta_{\frak m_{\lambda}}$ is equal to 
$\pmb R_{\frak m_{\lambda}}[x,y]/(G_1,G_2,G_3)$; and therefore, $\frak m_{\lambda} \Delta_{\frak m_{\lambda}}$ is contained in the Jacobson radical of $\Delta_{\frak m_{\lambda}}$.  We apply the local criterion for flatness, (see, for example, \cite{\rref{M80}, Thm.~49}).  To show that $\pmb R_{\frak m_{\lambda}}\to \Delta_{\frak m_{\lambda}}$ is flat it suffices to prove that
$\operatorname{Tor}_1^{\pmb R_{\frak m_{\lambda}}}(\pmb R_{\lambda},\Delta_{\frak m_{\lambda}})=0$ (since $\pmb R_{\lambda}$ means $\pmb R_{\frak m_{\lambda}}/\frak m_{\lambda}\pmb R_{\frak m_{\lambda}}$). Consider a resolution $\Bbb G$ of $\Delta_{\frak m_{\lambda}}$ by free $\pmb R_{\frak m_{\lambda}}[x,y]$-modules. One can create $\Bbb G$ by starting with $\Bbb F\otimes_{\pmb R}\pmb R_{\frak m_{\lambda}}$ and adjoining extra summands in positions $2$ and higher. We saw in Corollary \tref{C5.30} that $\operatorname{H}_1(\Bbb F\otimes_{\pmb R} \pmb R_{\frak m_{\lambda}}\otimes_{\pmb R_{\frak m_{\lambda}}}\pmb R_{\lambda})=0$. It follows that
$$0=\operatorname{H}_1(\Bbb G\otimes_{\pmb R_{\frak m_{\lambda}}}\pmb R_{\lambda})=\operatorname{Tor}_1^{\pmb R_{\frak m_{\lambda}}}(\pmb R_{\lambda},\Delta_{\frak m_{\lambda}}),$$
and (1) is established. 

The field $\pmb k$ remains algebraically closed; so, the maximal ideals of $\Delta$ continue to correspond to points in $\Bbb A^{3d+3}\times \Bbb P^1$. If $\frak M$ is a maximal ideal of $\Delta$, then $\frak M\cap \pmb R_w$ is a maximal ideal  $\frak m_{\lambda}$ of $\pmb R$ and $\pmb R_{\frak m_{\lambda}}\to \Delta_{\frak M}$ is a flat local homomorphism. The base ring $\pmb R_{\frak m_{\lambda}}$ and the local ring of the fiber $\Delta_{\frak M}\otimes_{\pmb R_{\frak m_{\lambda}}}\pmb R_{\lambda}=B/I_{\pmb g_{\lambda}}$ are both Cohen-Macaulay. It follows (see, for example, \cite{\rref{M80}, Thm. 50}) that $\Delta_{\frak M}$ is a Cohen-Macaulay ring of dimension 
$$\dim \pmb R_{\frak m_{\lambda}}+\dim B/I_{\pmb g_{\lambda}}=3d+3,$$ and (2) is established.

Consider the augmented complex $\Bbb F_{w}^{\text{aug}}\: 0\to \Bbb F_w\to \Delta\to 0$. Fix a maximal ideal $\frak M$ of  $\pmb S_w$. Let
$\Bbb F_{\frak M}^{\text{aug}}$ denote the localization of $\Bbb F_{w}^{\text{aug}}$ at $\frak M$ and let $\frak m_{\lambda}$ denote $\frak M\cap \pmb R_w$. The generators $\{z_{i,j}-\lambda_{i,j}\}$ of $\frak m_{\lambda}$ are a regular sequence on $\pmb S_w$ and they form a system of parameters on the Cohen-Macaulay ring $\Delta$; hence, $\{z_{i,j}-\lambda_{i,j}\}$ is a regular sequence on each non-zero module in $\Bbb F_{\frak M}^{\text{aug}}$. Apply Lemma \tref{L5.40} to see that
$$\text{$\Bbb F_{\frak M}^{\text{aug}}$ is exact} \iff \text{$\Bbb F_{\frak M}^{\text{aug}}\otimes_{\pmb S_{\frak M}}\pmb S_{\frak M}/\frak m_{\lambda}$ is exact}.\tag\tnum{*}$$
The augmented complex on the right hand side of (\tref{*}) is the localization $$\left (\Bbb F\otimes_{\pmb R}\pmb R_{\lambda_{\pmb g}}\to B/I_{\pmb g_{\lambda}}\right )_{\frak M}$$
of $\Bbb F\otimes_{\pmb R}\pmb R_{\lambda_{\pmb g}}\to B/I_{\pmb g_{\lambda}}$, which is exact by (2), and the proof of (3) is complete. 
\qed \enddemo

We used one direction of the following result in the course of proving Corollary \tref{C5.31}. 
\proclaim{Lemma \tnum{L5.40}} Let $R$ be a local ring and $$\Bbb F:\quad  \cdots @>d_2>>F_1@>d_1 >> F_0\to 0$$ be a complex of $R$-modules. Suppose that $z_1,\dots,z_r$ are elements of $R$ which form a regular sequence on each non-zero module in $\Bbb F$. Then $\Bbb F$ is exact if and only if $\Bbb F\otimes_RR/(z_1,\dots,z_r)$ is exact\endproclaim
\demo{Proof} By induction it suffices to prove the result when $r=1$. Write $z$ for $z_1$. 
\flushpar$(\Rightarrow)$ Assume that $\Bbb F$ is exact. Decompose $\Bbb F$ into a collect of of short exact sequences:
$$ 
0\to B_1\to F_1\to F_0\to 0, \qquad
0\to B_2\to F_2\to B_1 \to 0,\qquad
\dots,$$
where $B_i=\operatorname{im} d_{i+1}$. 
Observe that $z$ is regular on each $B_i$. Apply ${\underline{\phantom{X}}\otimes_RR/(z)}$ and glue the resulting short exact sequences back together to form the exact sequence $\Bbb F\otimes_RR/(z)$.
\flushpar$(\Leftarrow)$ Let $x$ be a cycle in $F_i$. The hypothesis that $\Bbb F\otimes_RR/(z)$ is exact ensures that $x=b+zx_2$ for some $b\in B_{i}$ and $x_2\in F_i$. The hypothesis that $z$ is regular on $F_i$ yields that $x_2$ is also a cycle. Repeat this process to see that $z$ is in $\cap_n (B_i+ z^nF_i)$, and this module, by the Krull Intersection Theorem, is $B_i$. \qed\enddemo

On at least three occasions, we apply the Buchsbaum-Eisenbud acyclicity criterion; furthermore, we apply the criterion in both directions. Roughly speaking, this criterion states that if 
$$F_{\bullet}:\quad 0\to F_s@> \varphi_s>> \cdots @> \varphi_1>> F_0$$ is a complex of finitely generated free modules over the Noetherian ring $R$, then $F_{\bullet}$ is acyclic if and only if $\operatorname{grade} I_{r_i}(\varphi_i)\ge i$ (or $I_{r_i}(\varphi_i)=R$) for $1\le i\le s$, where $r_i$ is the expected rank of $\varphi_i$. More details may be found in 
\cite{\rref{BE}}, \cite{\rref{E95}, Thm.~20.9}, or \cite{\rref{BH}, Thm.~1.4.12}, and many other places.

Recall the function $\Phi:$ 
$$\bmatrix Q_{1,1}&Q_{1,2}\\Q_{2,1}&Q_{2,2}\\Q_{3,1}&Q_{3,2}\endbmatrix@>\Phi>>\left [\,\left\vert \matrix  Q_{2,1}&Q_{2,2}\\Q_{3,1}&Q_{3,2}\endmatrix\right\vert, - \left\vert \matrix Q_{1,1}&Q_{1,2}\\ Q_{3,1}&Q_{3,2}\endmatrix\right\vert, \left\vert \matrix Q_{1,1}&Q_{1,2}\\Q_{2,1}&Q_{2,2}\endmatrix\right\vert\,\right ],\tag\tnum{Phi}$$which
sends a $3\times 2$ matrix  to a $1\times 3$ matrix   of signed  $2\times 2$ minors. 

\proclaim{Theorem \tnum{ZZ}} Adopt the notation of Definition {\rm\tref{data5}}. Let $\Phi$ be the function defined in {\rm(\tref{Phi})}. Then the following statements hold.
\roster\item The row vector   $\Phi (\pmb d_{2,\Bbb Z}^{(j)})+(-1)^{jc}z_{0,j}w_{\Bbb Z}^2 \pmb d_{1,\Bbb Z}$, with entries in $\pmb S_{\Bbb Z}$, is identically zero, for $1\le j\le 3$.
\item The complex $(\Bbb F_{\Bbb Z})_{w_{\Bbb Z}}$ is a resolution of $(\Delta_{\Bbb Z})_{w_{\Bbb Z}}$ by free $(\pmb S_{\Bbb Z})_{w_{\Bbb Z}}$-modules.
\item The ideal $(G_1,G_2,G_3)$ of $(\pmb S_{\Bbb Z})_{w_{\Bbb Z}}$ is perfect of grade two.
\item The complex $$\Bbb F_{\Bbb Z}^{(j)}:\quad 0\to (\pmb S_{\Bbb Z})_{w_{\Bbb Z}z_{0,j}}^2 @> \pmb d_{2,\Bbb Z}^{(j)}>> (\pmb S_{\Bbb Z})_{w_{\Bbb Z}z_{0,j}}^3@>\pmb d_{1,\Bbb Z}>> (\pmb S_{\Bbb Z})_{w_{\Bbb Z}z_{0,j}}$$is a resolution of $\Delta_{wz_{0,j}}$ by free $(\pmb S_{\Bbb Z})_{w_{\Bbb Z}z_{0,j}}$-modules for   $1\le j\le 3$.
\endroster
\endproclaim
\demo{Proof} Let $\pmb k$ be an algebraically closed field of characteristic zero. The natural ring homomorphism $\pmb S_{\Bbb Z}\to \pmb k\otimes_{\Bbb Z}\pmb S_{\Bbb Z}=\pmb S$ is an injection. Each element in each matrix in (1) is an element of $\pmb S_{\Bbb Z}$. If the left side of (1) is sent to zero in $\pmb S$, then the left side of (1) is already zero in $\pmb S_{\Bbb Z}$. We make our calculation in $\pmb S$. Recall the convention of Definition \tref{seg}: if $s_{\Bbb Z}$ is an element of $\pmb S_{\Bbb Z}$, then the image of $s_{\Bbb Z}$ in $\pmb S$ is denoted by $s$. We saw in Corollary \tref{C5.31} that $\Bbb F_{w}$ is a resolution of $\Delta_w$ by free $\pmb S_w$-modules. Fix $j$. We localize further to see that $\Bbb F_{wz_{0,j}}$ is a resolution of $\Delta_{wz_{0,j}}$ by free $\pmb S_{wz_{0,j}}$-modules. One entry of $\pmb d_3$ is a unit of $S_{wz_{0,j}}$. We split off a rank one summand from positions $2$ and $3$ of $\Bbb F_{wz_{0,j}}$ to obtain the resolution 
$$0\to \pmb S_{wz_{0,j}}^2@> \pmb d_2^{(j)} >> \pmb S_{wz_{0,j}}^3@> \pmb d_1 >> \pmb S_{wz_{0,j}}@>>> \Delta_{wz_{0,j}}\to 0.$$The Hilbert-Burch Theorem (or the Buchsbaum-Eisenbud acyclicity criterion) ensures that the maximal order minors of $\pmb d_2^{(j)}$ generate a grade two ideal of $S_{wz_{0,j}}$; hence, 
$$0\to \pmb S_{wz_{0,j}}@> \Phi(\pmb d_2^{(j)})^{\text{\rm T}} >> \pmb S_{wz_{0,j}}^3 @>{\pmb d_2^{(j)}}^{\text{\rm T}}>> 
\pmb S_{wz_{0,j}}^2$$ is acyclic. Of course, $\pmb d_1^{\text{\rm T}}$ is in the kernel of ${\pmb d_2^{(j)}}^{\text{\rm T}}$.   It follows that 
$\pmb d_1=\theta \Phi(\pmb d_2^{(j)})$, for some $\theta$ in $\pmb S_{wz_{0,j}}$. View $\pmb S_{wz_{0,j}}$ as the polynomial ring $\pmb R_{wz_{0,j}}[x,y]$ where each element of  $\pmb R_{wz_{0,j}}$ has degree zero and $x$ and $y$ have degree $1$. The entries of $\pmb d_1$ and the entries of  $\Phi(\pmb d_2^{(j)})$ are homogeneous forms of degree $d$; hence, $\theta$ is in $\pmb R_{wz_{0,j}}$. We identify $\theta$ by looking at the coefficient of $y^d$ in position $j$. In $\pmb d_1$, this coefficient is $z_{0,j}$. In $\Phi(\pmb d_2^{(j)})$, this coefficient may easily be read from (\tref{511}); it is 
$$ \cases \phantom{-}\left\vert\smallmatrix  
0&(-1)^{c+1} wz_{0,1}\\
- wz_{0,1}&0\endsmallmatrix\right\vert&\text{if $j=1$}\\
-\left\vert\smallmatrix 0&(-1)^cwz_{0,2}\\
(-1)^{c+1} wz_{0,2}&0\endsmallmatrix\right\vert&\text{if $j=2$}\\
\phantom{-}\left\vert\smallmatrix 0&wz_{0,3}\\
(-1)^c wz_{0,3}&0\\
\endsmallmatrix\right\vert&\text{if $j=3$}.\endcases$$ In other words, in $\Phi(\pmb d_2^{(j)})$, this coefficient is $(-1)^{jc+1}w^2z_{0,j}^2$.  It follows that $\theta$ is equal to $\frac{(-1)^{jc+1}}{w^2z_{0,j}}$ and   the proof of (1) is complete.

We   make two observations before proving (2). First of all,
the polynomials $G_1,G_2$ in $\pmb S_{\Bbb Z}$ are generic linear combinations of the generators of the grade $2$ ideal $(x,y)^d$ of $\pmb k[x,y]$; and therefore, $$\text{$G_1,G_2$ form a regular sequence on $\pmb S_{\Bbb Z}$};\tag\tnum{Ho}$$ see, for example, \cite{\rref{Ho}}. Secondly, 
 $w_{\Bbb Z}$ is the determinant of the matrix $A_{\Bbb Z}^{(c-1)}$ of (\tref{Ai}). One can expand $\det A_{\Bbb Z}^{(c-1)}$ across the first row and write   $$w_{\Bbb Z}=z_{0,1}Z_{1,\Bbb Z}+z_{0,2}Z_{2,\Bbb Z}+z_{0,3}Z_{3,\Bbb Z},\tag\tnum{zZ}$$ where $Z_{1,\Bbb Z}$, $Z_{2,\Bbb Z}$, and $Z_{3,\Bbb Z}$ are the appropriate signed minors of $A_{\Bbb Z}^{(c-1)}$.

To prove (2) it suffices to show that the localization of $\Bbb F_{\Bbb Z}$ at $P$ is exact for all prime ideals $P$ of $(\pmb S_{\Bbb Z})_{w_{\Bbb Z}}$.  If $P$ is such a prime ideal, then   $w_{\Bbb Z}\notin P$ and (\tref{zZ}) shows that   some $z_{0,j}\notin P$. Consequently, to prove (2), it suffices to show that $(\Bbb F_{\Bbb Z})_{w_{\Bbb Z}z_{0,j}}$ is exact for all $j$. We employ the Buchsbaum-Eisenbud acyclicity criterion once again. We see that $I_1((\pmb d_{3,\Bbb Z})_{w_{\Bbb Z}z_{0,j}})$ is the entire ring $(S_{\Bbb Z})_{w_{\Bbb Z}z_{0,j}}$. Assertion (1) yields that $G_1$ and $G_2$ are in $I_2((d_{2,\Bbb Z})_{w_{\Bbb Z}z_{0,j}})$; and therefore   (\tref{Ho}) shows that the grade of $I_2((d_{2,\Bbb Z})_{w_{\Bbb Z}z_{0,j}})$ is at least $2$. 
The proof of (2) is complete. 

The inequality 
$$\operatorname{grade}(G_1,G_2,G_3)(\pmb S_{\Bbb Z})_{w_{\Bbb Z}}\le \operatorname{pd} _{(\pmb S_{\Bbb Z})_{w_{\Bbb Z}}} (\Delta_{\Bbb Z})_{w_{\Bbb Z}}$$ holds automatically. In (3) we assert that both numbers are equal to $2$. See, for example, \cite{\rref{BH}, pg. 25} for more information about perfect modules. The observation (\tref{Ho}) gives
  $2\le \operatorname{grade} (G_1,G_2,G_3)(\pmb S_{\Bbb Z})_{w_{\Bbb Z}}$. On the other hand, we saw in (2) that 
$(\Bbb F_{\Bbb Z})_{w_{\Bbb Z}}$ is a resolution of $(\Delta_{\Bbb Z})_{w_{\Bbb Z}}$; 
 hence, $\operatorname{pd} _{(\pmb S_{\Bbb Z})_{w_{\Bbb Z}}} (\Delta_{\Bbb Z})_{w_{\Bbb Z}}\le 3$.
 Furthermore, we know an explicit splitting map for $$(\pmb S_{\Bbb Z})_{w_{\Bbb Z}}@> \pmb d_{3,\Bbb Z}>> (\pmb S_{\Bbb Z})_{w_{\Bbb Z}}^3.$$ Define $$\sigma_{\Bbb Z}\:(\pmb S_{\Bbb Z})_{w_{\Bbb Z}}^3\to (\pmb S_{\Bbb Z})_{w_{\Bbb Z}}\text{ to be }\tsize\frac 1{w_{\Bbb Z}}\bmatrix Z_{1,\Bbb Z}&Z_{2,\Bbb Z}&Z_{3,\Bbb Z}\endbmatrix.$$ We see that $\sigma_{\Bbb Z}\circ \pmb d_{3,\Bbb Z}$ is the identity on $(\pmb S_{\Bbb Z})_{w_{\Bbb Z}}$. It follows that $(\pmb S_{\Bbb Z})_{w_{\Bbb Z}}^3$ is equal to $\operatorname{im} \pmb d_{3,\Bbb Z}\oplus\ker \sigma_{\Bbb Z}$ and that 
$$\Bbb U\Bbb P\Bbb R_{\Bbb Z}:\quad 0\to \ker \sigma_{\Bbb Z} @> \pmb d_{2,\Bbb Z}|\ker \sigma_{\Bbb Z} >> (\pmb S_{\Bbb Z})_{w_{\Bbb Z}}^3 @> \pmb d_{1,\Bbb Z} >> (\pmb S_{\Bbb Z})_{w_{\Bbb Z}}@>>> (\Delta_{\Bbb Z})_{w_{\Bbb Z}}\to 0\tag\tnum{UR}$$
 is a length $2$ projective resolution of $(\Delta_{\Bbb Z})_{w_{\Bbb Z}}$ and the proof of (3) is complete. 

Assertion (1) shows that $(G_1,G_2)\subseteq I_2(\pmb d_{2,\Bbb Z}^{(j)})(\pmb S_{\Bbb Z})_{w_{\Bbb Z}z_{0,j}}$; so, (\tref{Ho}) gives $$2\le \operatorname{grade} I_2(\pmb d_{2,\Bbb Z}^{(j)})(\pmb S_{\Bbb Z})_{w_{\Bbb Z}z_{0,j}}$$ and (4) follows from the Buchsbaum-Eisenbud acyclicity criterion (or the Hilbert-Burch Theorem). 
\qed\enddemo

\proclaim{Corollary \tnum{sectn}} 
Let $\pmb k$ be a field, $d=2c$ be an even integer,  $\Phi\:\Bbb H_d \to \Bbb A_d$ be the morphism which sends a $3\times 2$ matrix with entries from $\pmb k[x,y]_c$ to a triple $\pmb g$ of $\Bbb A_d$, and let $\Phi|$ be the restriction
$$\Phi|\: \Phi^{-1}(\operatorname{Bal}_d)\to \operatorname{Bal}_d\tag\tref{Phi|}$$of $\Phi$ to $\Phi^{-1}(\operatorname{Bal}_d)$. Then there exists a local section of the morphism $\Phi|$ of {\rm(\tref{Phi|})}.\endproclaim

\demo{Proof} We saw in Theorem \tref{T1} that $\cup_{j=1}^3\operatorname{Bal}_d^{(j)}$ is an open cover of $\operatorname{Bal}_d$. Consider $w$ and $\pmb d_{2}^{(j)}$ as found in Definition \tref{seg}. Let $D^{(j)}$ be the matrix whose first column is $\frac{(-1)^{jc+1}}{w^2z_{0,j}}$ times the first column of $\pmb d_{2}^{(j)}$ and whose second column is the second column of $\pmb d_{2}^{(j)}$.
Define $\sigma_j\:\operatorname{Bal}_d^{(j)}\to \Bbb H_d$ by $$\pmb g\mapsto \left [\tsize D^{(j)}\right ]|_{\lambda_{\pmb g}}$$
for $\pmb g\in \operatorname{Bal}_d^{(j)}$. Assertion (1) of Theorem \tref{ZZ} shows that $\Phi|\circ \sigma_j$ is the identity morphism on $\operatorname{Bal}_d^{(j)}$. \qed\enddemo

\proclaim{Corollary \tnum{UPR}} Let $\pmb k$ be a  field, $B$ be the standard graded polynomial ring $\pmb k[x,y]$ and $c$ be a positive integer. Then the projective resolution $\Bbb U\Bbb P\Bbb R_{\Bbb Z}$ of {\rm(\tref{UR})} is a universal projective resolution for the graded Betti numbers 
$$0\to B(-3c)^2@>>> B(-2c)^3@>> >B\to B/I\to 0.\tag\tnum{gbn}$$
\endproclaim

\demo{Proof} Let $I$ be a homogeneous ideal of $B$ with the property that the graded Betti numbers of $B/I$ are given by {\rm(\tref{gbn})}. In particular, $I$ is generated by three forms $g_1,g_2,g_3$ of degree $2c$. Let $d$ be the integer $2c$ and $\pmb g$ be the triple $(g_1,g_2,g_3)$ of $\Bbb A_d$. In the language of Definition \tref{.'D27.13}, we have  $I_{\pmb g}=I$. The hypothesis (\tref{gbn}) also shows that $\pmb g\in \operatorname{Bal}_d$.  It follows from Theorem \tref{YYY} that $\lambda_{\pmb g}\notin
V(w)$; and therefore, $\pmb R_{\lambda_{\pmb g}}$ is an $\pmb R_w$-algebra.
Consider the  homomorphism $\alpha\:(\pmb R_{\Bbb Z})_{w_{\Bbb Z}}\to \pmb R_{\lambda_{\pmb g}}$ which is the composition of the natural maps
$$(\pmb R_{\Bbb Z})_{w_{\Bbb Z}}\to (\pmb R_{\Bbb Z})_{w_{\Bbb Z}}\otimes_{\Bbb Z}\pmb k=\pmb R_w\quad\text{and}\quad \pmb R_w\to \tsize \frac {\pmb R_w}{(\frak m_{\lambda_{\pmb g}})\pmb R_w}=\pmb R_{\lambda_{\pmb g}}.$$
We establish the result by showing that 
$\Bbb U\Bbb P\Bbb R_{\Bbb Z}\otimes_{(\pmb R_{\Bbb Z})_{w_{\Bbb Z}}}\pmb R_{\lambda_{\pmb g}}$ is a resolution of $B/I$ by free $B$-modules with graded Betti numbers given by (\tref{gbn}). Corollary \tref{C5.30} shows that $\Bbb F\otimes_{\pmb R} \pmb R_{\lambda_{\pmb g}}$ is a resolution of $B/I_{\pmb g}=B/I$. On the other hand, $\pmb R_{\lambda_{\pmb g}}$ is an $\pmb R_w$-algebra; so, 
$$\Bbb F\otimes_{\pmb R} \pmb R_{\lambda_{\pmb g}}=\Bbb F\otimes_{\pmb R} \left (\pmb R_w\otimes_{\pmb R_w}\pmb R_{\lambda_{\pmb g}}\right )
=\Bbb F_w\otimes_{\pmb R_w}\pmb R_{\lambda_{\pmb g}}=\left ((\Bbb F_{\Bbb Z})_{w_{\Bbb Z}}\otimes_{\Bbb Z}\pmb k\right )\otimes_{\pmb R_w}\pmb R_{\lambda_{\pmb g}}.$$The definition of $\alpha\:(\pmb R_{\Bbb Z})_{w_{\Bbb Z}}\to \pmb R_{\lambda_{\pmb g}}$ shows that
$$\left ((\Bbb F_{\Bbb Z})_{w_{\Bbb Z}}\otimes_{\Bbb Z}\pmb k\right )\otimes_{\pmb R_w}\pmb R_{\lambda_{\pmb g}}=(\Bbb F_{\Bbb Z})_{w_{\Bbb Z}}\otimes_{(\pmb R_{\Bbb Z})_{w_{\Bbb Z}}} \pmb R_{\lambda_{\pmb g}}.$$Thus, $(\Bbb F_{\Bbb Z})_{w_{\Bbb Z}}\otimes_{(\pmb R_{\Bbb Z})_{w_{\Bbb Z}}} \pmb R_{\lambda_{\pmb g}}$ is a resolution of $B/I$ by free $B$-modules.
The construction of $\Bbb U\Bbb P\Bbb R_{\Bbb Z}$ in the proof of  Theorem \tref{ZZ} shows that $(\Bbb F_{\Bbb Z})_{w_{\Bbb Z}}\otimes_{(\pmb R_{\Bbb Z})_{w_{\Bbb Z}}} \pmb R_{\lambda_{\pmb g}}$ is the direct sum of the subcomplexes $\Bbb U\Bbb P\Bbb R_{\Bbb Z}\otimes_{(\pmb R_{\Bbb Z})_{w_{\Bbb Z}}}\pmb R_{\lambda_{\pmb g}}$  and $\Bbb S$, where $\Bbb S$ is the  exact complex 
$$0\to [(\pmb S_{\Bbb Z})_{w_{\Bbb Z}}\otimes_{(\pmb R_{\Bbb Z})_{w_{\Bbb Z}}}\pmb R_{\lambda_{\pmb g}}](-3c)\to [(\pmb S_{\Bbb Z})_{w_{\Bbb Z}}\otimes_{(\pmb R_{\Bbb Z})_{w_{\Bbb Z}}}\pmb R_{\lambda_{\pmb g}}](-3c)\to 0,$$ 
concentrated in positions two and three. Keep in mind that $$(\pmb S_{\Bbb Z})_{w_{\Bbb Z}}\otimes_{(\pmb R_{\Bbb Z})_{w_{\Bbb Z}}}\pmb R_{\lambda_{\pmb g}}=\pmb R_{\lambda_{\pmb g}}[x,y]\cong \pmb k[x,y]=B.$$
The Quillen-Suslin theorem guarantees that every finitely generated projective $B$-module is free. We conclude that
 $\Bbb U\Bbb P\Bbb R_{\Bbb Z}\otimes_{(\pmb R_{\Bbb Z})_{w_{\Bbb Z}}}\pmb R_{\lambda_{\pmb g}}$ is a resolution of $B/I$ by free $B$-modules and   the graded Betti numbers of $\Bbb U\Bbb P\Bbb R_{\Bbb Z}\otimes_{(\pmb R_{\Bbb Z})_{w_{\Bbb Z}}}\pmb R_{\lambda_{\pmb g}}$ are given in (\tref{gbn}).   \qed \enddemo

The conclusion in Theorem \tref{ZZ} that $(G_1,G_2,G_3)(\pmb S_{\Bbb Z})_{w_{\Bbb Z}}$ is a perfect ideal is particularly useful because of the   ``Persistence of Perfection Principle'', which is also known as  the ``transfer of perfection'' (see \cite{\rref{Ho-T}, Prop. 6.14} or \cite{\rref{BV}, Thm. 3.5}). \proclaim{The Persistence of Perfection Principle}Let $R\to S$ be a homomorphism of Noetherian rings, $I$ be a perfect ideal of  $R$ of grade $g$, and   $\Bbb P$ be a resolution of $R/I$ by projective $R$-modules. If  $IS$ is a proper ideal of $S$ with grade at least $g$, then $IS$ is a perfect ideal of $S$ of grade $g$ and $\Bbb P\otimes_RS$ is a  resolution of $S/IS$ by projective $S$-modules. \endproclaim
In the context of Definition \tref{data5}, we have identified two projective resolutions for $(\Delta_{\Bbb Z})_{w_{\Bbb Z}}$: the length three free  resolution $(\Bbb F_{\Bbb Z})_{w_{\Bbb Z}}$ of part (2) of Theorem \tref{ZZ}   and the length two projective resolution $\Bbb U\Bbb P\Bbb R_{\Bbb Z}$ of (\tref{UR}).  In part (4) of Theorem \tref{ZZ} we have identified a free resolution of length two for $(\Delta_{\Bbb Z})_{w_{\Bbb Z}z_{0,j}}$.  The  Persistence of Perfection Principle yields the following result.
\proclaim{Corollary \tnum{pop1}}Adopt the notation of Definition {\rm\tref{data5}} and let $S$ be an arbitrary Noetherian ring. \roster\item If $(\pmb S_{\Bbb Z})_{w_{\Bbb Z}}\to S$ is a ring homomorphism and $(G_1,G_2,G_3)S$ is a proper ideal of grade at least $2$, then $(G_1,G_2,G_3)S$ is a perfect ideal of grade equal to $2$ and
$$(\Bbb F_{\Bbb Z})_{w_{\Bbb Z}}\otimes_{(\pmb S_{\Bbb Z})_{w_{\Bbb Z}}}S\quad\text{and}\quad \Bbb U\Bbb P\Bbb R_{\Bbb Z}\otimes_{(\pmb S_{\Bbb Z})_{w_{\Bbb Z}}}S$$ both are projective resolutions of $S/(G_1,G_2,G_3)S$. 
\item If $(\pmb S_{\Bbb Z})_{w_{\Bbb Z}z_{0,j}}\to S$ is a ring homomorphism and $(G_1,G_2,G_3)S$ is a proper ideal of grade at least $2$, then $(G_1,G_2,G_3)S$ is a perfect ideal of grade equal to $2$ and
$\Bbb F_{\Bbb Z}^{(j)}\otimes_{(\pmb S_{\Bbb Z})_{w_{\Bbb Z}z_{0,j}}}S$ is free resolution of $S/(G_1,G_2,G_3)S$ of length $2$. \endroster
 \endproclaim

\bigpagebreak
\SectionNumber=6\tNumber=1
\heading Section \number\SectionNumber. \quad Decomposition of the space of true triples.
\endheading 

Let $\pmb k$ be an algebraically closed field and $d=2c$ be an even integer. 
Recall the open subsets  $\Bbb B_d\subseteq\Bbb T_d\subseteq \Bbb A_d$ of  triples of  forms of degree $d$ from Definition \tref{.'D27.14}.
In Observation \tref{O40} we saw that $\pmb g\mapsto \Psi_{\pmb g}$ gives a surjection from $\Bbb T_d$ to 
the space of true (i.e., birational and  base-point free) parameterizations of plane curves of degree $d$.
In this section we decompose $\Bbb T_d$ into locally closed subsets which depend on the configuration of multiplicity $c$ singularities on, or infinitely near, the corresponding curve. Recall the Configuration Poset $(\operatorname{CP},\le)$ from Definition \tref{po}. 
\definition{Definition \tnum{D39.1}} For each $\#$ in $\operatorname{CP}$ define
$$S_{\#}=\left\{\pmb g\in \Bbb T_d\left\vert \matrix\format\l\\\text{the configuration of multiplicity $c$ singularities on or infinitely}\\\text{near $\Cal C_{\pmb g}$ are described by \#}\endmatrix\right.\right\}$$ and define $T_{\#}=\bigcup\limits_{\#'\le \#} S_{\#'}$.
\enddefinition

\bigskip It is clear that $\Bbb T_d$ is the disjoint union $$\Bbb T_d=\bigcup\limits_{\#\in \operatorname{CP}}S_\#=T_{\emptyset}.$$   
It is shown in Corollary \tref{CTLL} that $S_{\#}\subseteq \Bbb B_d$ for all $\#$ in $\operatorname{CP}$, except $\#=\emptyset$;  thus
$$S_{\emptyset}=(S_{\emptyset} \cap \Bbb B_d)\cup \Bbb U\Bbb B_d,\tag\tnum{notn}$$  where $\Bbb U\Bbb B_d=\Bbb T_d\setminus \Bbb B_d$, as was defined in Remark \tref{UB}. It follows that the space 
$\Bbb T_d$ is equal to the disjoint union of the eight subsets 
$$\tsize \bigcup\limits_{\#\in \operatorname{CP}\setminus\{\emptyset\}}S_{\#}\cup (S_{\emptyset}\cap \Bbb B_d)\cup \Bbb U\Bbb B_d.$$
 The set $\Bbb U\Bbb B_d$ is a closed subset of $\Bbb T_d$; indeed it is shown in Theorem \tref{YYY}  that $\Bbb U\Bbb B_d$ is a hypersurface section of $\Bbb T_d$. Once $\Bbb U\Bbb B_d$ is removed from $\Bbb T_d$, then one is left with
the space $\Bbb B_d$ which is the disjoint union of the seven subsets $\bigcup\limits_{\#\in \operatorname{CP}\setminus\{\emptyset\}}S_{\#}\cup (S_{\emptyset}\cap \Bbb B_d)$. We call each  of these seven subsets of $\Bbb B_d$  a stratum of $\Bbb B_d$. We consider the picture:
$$\xymatrix{&&S_{c:c}\ar[dr]\\
S_{c:c:c} \ar[r]  &S_{c:c,c} \ar[ru] \ar[rd]& &S_{c,c} \ar[r] &S_{c}\ar[r] &S_{\emptyset}\cap \Bbb B_d.\\&&S_{c,c,c}\ar[ur]}\tag\tnum{pictt}$$
 In this section, we  prove that each stratum $S_\#\cap \Bbb B_d$ in the above picture  is a locally closed subset of $\Bbb B_d$ and wherever we have  drawn $S_{\#'}\to S_{\#}\cap \Bbb B_d$, then $S_{\#'}$ is contained in the closure $\overline{S_{\#}\cap \Bbb B_d}$ of $S_{\#}\cap \Bbb B_d$ in $\Bbb B_d$; indeed, we prove that $\overline{S_{\#}\cap \Bbb B_d}=T_{\#}\cap \Bbb B_d$.  Furthermore, we prove that  each $T_{\#}\cap \Bbb B_d$ is a closed irreducible subset of $\Bbb B_d$ and we calculate its dimension.
We emphasize that if $\#\in \operatorname{CP}$, then
$$\left\{\matrix\format\r&\c&\c&\c&\l&\c&\l&\c&\l\\ 
\#<\emptyset&{}\implies{}& S_\# &{}={}&S_\#\cap \Bbb B_d \ &\text{ and }& T_{\#}&{}={}&T_{\#} \cap \Bbb B_d, \text{ whereas}\\
\#=\emptyset&{}\implies{}&
 S_{\emptyset} &{}={}& (S_{\emptyset}\cap \Bbb B_d) \cup \Bbb U\Bbb B_d\ &\text{ and }&
T_{\emptyset}&{}={}&(T_{\emptyset}\cap \Bbb B_d)\cup \Bbb U\Bbb B_d=\Bbb T_d,\\
\endmatrix\right.\tag\tnum{emp}$$
with $T_{\emptyset}\cap \Bbb B_d=\Bbb B_d$.

The set $\Bbb B_d$ is an open subset of the affine space $\Bbb A_d$ of Data \tref{.'Org} and the topology on $\Bbb B_d$ is the subspace topology: a subset of $Y$ of $\Bbb B_d$ is closed in $\Bbb B_d$ if and only $Y=\Bbb B_d\cap X$ for some closed subset $X$ of $\Bbb A_d$.

We now connect the subsets $S_{\#}$ and $T_{\#}$ of $\Bbb T_d\subseteq \Bbb A_d$ with the subsets $M_{\natural}$ and $N_{\#}$ of $\Bbb B\Bbb H_d\subseteq \Bbb H_d$. Recall the morphisms 
$$G\times \Bbb H_d@> \Upsilon>> \Bbb H_d \quad\text{and}\quad \Bbb H_d @> \Phi >> \Bbb A_d$$
from Remark (6) of \tref{401} and Definition \tref{Phi!}. The definitions of the sets $\Bbb T_d$, $\Bbb B_d$, and $\Bbb B\Bbb H_d$ ensures 
$$\split \Phi^{-1}(\Bbb A_d\setminus \Bbb T_d)&{}=\Bbb H_d\setminus \Bbb B\Bbb H_d,\\ \Phi^{-1}(\Bbb T_d\setminus \Bbb B_d)&{}\text{ is empty, and}\\ \Phi^{-1}(\Bbb B_d)&{}=\Bbb B\Bbb H_d.\endsplit$$Furthermore, the restriction of $\Phi$ to $\Bbb B\Bbb H_d$ is a surjection onto $\Bbb B_d$:
$$\Phi\:\Bbb B\Bbb H_d\twoheadrightarrow \Bbb B_d.$$
In Definition \tref{D39.1}, we decomposed $\Bbb B_d$ as $\cup_{\#\in \operatorname{CP}} S_{\#}\cap \Bbb B_d$ and in Theorem \tref{orbitz} we decomposed $\Bbb B\Bbb H_d$ as $\cup_{\natural \in \operatorname{ECP}}\operatorname{DO}_{\natural}$. Furthermore, 
$$\Phi^{-1}(S_{\#}\cap \Bbb B_d)=\cases \operatorname{DO}_{\#}&\text{if $\#<c$}\\ \operatorname{DO}_{c,\mu_4}\cup \operatorname{DO}_{c,\mu_5}&\text{if $\#=c$}\\
\operatorname{DO}_{\emptyset,\mu_4}\cup \operatorname{DO}_{\emptyset,\mu_5}\cup \operatorname{DO}_{\emptyset,\mu_6}&\text{if $\#=\emptyset$}.\endcases$$
We have accounted for all eleven elements of $\operatorname{ECP}$ because, as was noted in Remark (3) of \tref{408}, $\operatorname{DO}_{\mu_2}$ is empty. 

Each orbit $\operatorname{DO}_{\natural}$ is defined (in Definition \tref{Mnat}) to be $\operatorname{DO}_{\natural}=G\cdot M_{\natural}$. Recall from (4) of Definition \tref{Mnat} that $$M_c= M_{c,\mu_4}\cup M_{c,\mu_5}.$$
Using  this definition, we have
$$\Phi^{-1}(S_{\#}\cap \Bbb B_d)=G\cdot M_{\#}\quad\text{for all $\#\in \operatorname{CP}\setminus \{\emptyset\}$}.\tag\tnum{abr}$$
Fix $\#\in \operatorname{CP}$, with $\#<\emptyset$. Recall from (\tref{emp}) that $S_{\#}\cap\Bbb B_d=S_{\#}$. We have
$$\xymatrix{M_{\#}\subseteq G\cdot M_{\#}=\Phi^{-1}(S_{\#}\cap \Bbb B_d)= \Phi^{-1}(S_{\#}) \ar@{>>}[r]^-{\Phi}&S_{\#}}\tag\tnum{M29.1}$$

Keep $\#\in \operatorname{CP}$, with $\#<\emptyset$. Apply $\Phi^{-1}$ to 
$$T_{\#}=\bigcup\limits_{\#'\in \operatorname{CP}\atop{ \#'\le\#}}S_{\#}$$ to see that
$$\Phi^{-1}(T_{\#})=\bigcup\limits_{\#'\in \operatorname{CP}\atop{ \#'\le\#}}\Phi^{-1}(S_{\#})=\bigcup\limits_{\#'\in \operatorname{CP}\atop{\#'\le\#}}G\cdot M_{\#'}= \bigcup\limits_{\natural\in \operatorname{ECP}\atop{ \natural\le\#}}\operatorname{DO}_{\natural}.$$Definition \tref{CO} gives
$$\bigcup\limits_{\natural\in \operatorname{ECP}\atop{ \natural\le\#}}\operatorname{DO}_{\natural}=\operatorname{CO}_{\#}$$and assertion (2) of Theorem \tref{orbitz2} gives $\operatorname{CO}_{\#}=G\cdot N_{\#}$. Thus 
$$\xymatrix{N_{\#}\subseteq G\cdot N_{\#}= \Phi^{-1}(T_{\#}) \ar@{>>}[r]^-{\Phi}&T_{\#}}\tag\tnum{M29.2}$$
Combine (\tref{M29.1}) and (\tref{M29.2}) to obtain the following picture for all $\#\in \operatorname{CP}$ with $\#<\emptyset$.
$$ \xymatrix{ 
&&\Bbb H_d\ar[r]^{\Phi}&\Bbb A_d\\
&\Bbb B\Bbb H_d\ar@{=}[r]&\Phi^{-1}(\Bbb B_d)\ar[r]^{\Phi}\ar@{^{(}->}[u]&\Bbb B_d\vphantom{E^{E^{E}}_{E_{E}}}\ar@{^{(}->}[u]\\
N_{\#}\ar@{^{(}->}[r]&GN_{\#}\ar@{=}[r]\ar@{^{(}->}[u]&\Phi^{-1}(T_{\#})\ar[r]^{\Phi}\ar@{^{(}->}[u]&T_{\#}\ar@{^{(}->}[u]\\
M_{\#}\ar@{^{(}->}[r]\ar@{^{(}->}[u]&GM_{\#}\ar@{=}[r]\ar@{^{(}->}[u]&\Phi^{-1}(S_{\#})\ar[r]^{\Phi}\ar@{^{(}->}[u]&S_{\#}.\ar@{^{(}->}[u]\\
}\tag\tnum{P3}
$$
\remark{\bf Remark \tnum{fws}} It follows from Proposition \tref{no-e'} that each stratum $S_{\#}$ is non-empty for $c\ge 3$; and if $c=2$, then $S_{\#}$ is non-empty if and only if $\#\le \{c,c,c\}$.\endremark

\proclaim{Theorem \tnum{strata}} For every $\#\in \operatorname{CP}$,  the set
$T_{\#}\cap \Bbb B_d$   is a closed and irreducible subset of $\Bbb B_d$.
\endproclaim

\demo{Proof} The statement is true, but not particularly interesting, for $\#=\emptyset$ because $T_{\emptyset}=\Bbb T_d$; and therefore $T_\emptyset\cap \Bbb B_d=\Bbb B_d$ which is a closed subset of $\Bbb B_d$. Furthermore, $\Bbb B_d$ is an open subset of the affine space $\Bbb A_d$; so $\Bbb B_d$ is irreducible. 

Henceforth, we take $\#\neq \emptyset$. In this case, $T_{\#}\subseteq \Bbb B_d$; hence,   $T_{\#}\cap \Bbb B_d=T_{\#}$. 
We first prove that $T_{\#}$ is a closed subset of $\Bbb B_d$. Recall the open cover $\cup_{i=1}^3 \Bbb B_{d}^{(i)}$ of $\Bbb B_d$ which is given in the Remark after Theorem \tref{T1}. It suffices to prove that $T_{\#}\cap \Bbb B_{d}^{(i)}$ is a closed subset of $\Bbb B_{d}^{(i)}$
for each $i$. Adopt the notation of Conventions \tref{CP}; in particular,  $\pmb R=\pmb k[\pmb z]$ is the coordinate ring of $\Bbb A^{3d+3}$. Let $A^{(i)}$ be the matrix of Remark \tref{R5.44} with entries from $\pmb R[\pmb u]$, and let  $S_i$ and $D_i$ be the $\pmb R$-algebras $S_i=\pmb R[\pmb u]/I_3(A^{(i)})$ and $D_i=S_i/\operatorname{Jac}(S_i/\pmb R)$. Corollary \tref{CorND} shows that if $\pmb g\in \Bbb B_{d}^{(i)}$, then the configuration of multiplicity $c$ singularities on or infinitely near $\Cal C_{\pmb g}$ may be read from the rings $(S_i)_{\lambda_{\pmb g}}$ and $(D_i)_{\lambda_{\pmb g}}$.  
We proved in Corollary \tref{CorND}  that
$$\alignat 1 T_c\cap \Bbb B_{d}^{(i)}&{}= X(S_i;\ge 1) \cap \Bbb B_{d}^{(i)}\\\allowdisplaybreak
T_{c,c}\cap \Bbb B_{d}^{(i)}&{}=  X(S_i;=1,\ge 2) \cap \Bbb B_{d}^{(i)}\\\allowdisplaybreak
T_{c:c}\cap \Bbb B_{d}^{(i)} &{}= X(S_i; \ge 1)\cap X(D_i,\ge 1)\cap \Bbb B_{d}^{(i)}= X(D_i, \ge 1) \cap \Bbb B_{d}^{(i)}\\\allowdisplaybreak
T_{c,c,c}\cap \Bbb B_{d}^{(i)}&{}=X(S_i;=1,\ge 3)\cap \Bbb B_{d}^{(i)}\\\allowdisplaybreak
T_{c:c,c}\cap \Bbb B_{d}^{(i)}&{}=X(S_i;=1,\ge 3)\cap X(D_i;\ge 1)\cap \Bbb B_{d}^{(i)}\\\allowdisplaybreak
T_{c:c:c}\cap \Bbb B_{d}^{(i)}&{}= X(S_i; =1, \ge 3)\cap X(D_i;=1,\ge 2)\cap \Bbb B_{d}^{(i)}=X(D_i;=1,\ge 2)\cap \Bbb B_{d}^{(i)}\endalignat $$ 
Apply Theorem \tref{L37.7}, by way of Convention (4) of \tref{CP}. The set   $X(S_i;\ge 1)\cap \Bbb A_d$ is closed in $\Bbb A_d$;
hence, $T_c\cap \Bbb B_{d}^{(i)}$ is a closed subset of $\Bbb B_{d}^{(i)}$. The set ${X(S_i;=1,\ge 2)\cap \Bbb B_{d}^{(i)}}$ is closed in 
 $$X(S_i,=1)\cap \Bbb B_{d}^{(i)}=X(S_i;\ge 1) \cap \Bbb B_{d}^{(i)};$$ but $X(S_i;\ge 1) \cap \Bbb B_{d}^{(i)}$ is already closed in $\Bbb B_{d}^{(i)}$, 
so, $T_{c,c}\cap \Bbb B_{d}^{(i)}$ is closed in $\Bbb B_{d}^{(i)}$. The set ${X(D_i;\ge 1) \cap \Bbb B_{d}^{(i)}}$ is closed in $\Bbb B_{d}^{(i)}$; so,  $T_{c:c}\cap \Bbb B_{d}^{(i)}$ is closed in 
 $\Bbb B_{d}^{(i)}$. The set $X(S_i;=1,\ge 3)\cap \Bbb B_{d}^{(i)}$ is closed in ${X(S_i;=1)\cap \Bbb B_{d}^{(i)}}$. The argument from $T_{c,c}$ yields that  $T_{c,c,c}\cap \Bbb B_{d}^{(i)}$ is closed in $\Bbb B_{d}^{(i)}$. The set $T_{c:c,c}\cap \Bbb B_{d}^{(i)}$ is the intersection of the closed subsets $T_{c,c,c}\cap \Bbb B_{d}^{(i)}$ and $T_{c:c}\cap \Bbb B_{d}^{(i)}$; so,  $T_{c:c,c}\cap \Bbb B_{d}^{(i)}$ is also a closed subset of $\Bbb B_{d}^{(i)}$. The set ${X(D_i;=1,\ge 2)\cap \Bbb B_{d}^{(i)}}$ is closed in 
the closed set $X(D_i;\ge 1)\cap \Bbb B_{d}^{(i)}$; so, $T_{c:c:c}\cap \Bbb B_{d}^{(i)}$ is a closed set in $\Bbb B_{d}^{(i)}$. 

Fix $\#\in \operatorname{CP}$ with $\#\neq \emptyset$. We have established that each set $T_{\#}$    is a closed subset of $\Bbb B_d$. Now we show that $T_{\#}$ is irreducible. 
Consider the surjective morphisms 
$$\xymatrix{G\times N_{\#}\ar@{>>}^{\Upsilon}[r]&\Phi^{-1}(T_{\#})\ar@{>>}^{\Phi}[r]&T_{\#}}.\tag\tnum{surj}$$
We have shown in  Theorem \tref{orbitz2} that $N_{\#}$ is an irreducible variety. It is clear that $G$ is an irreducible variety; hence,   $G\times N_\#$  is irreducible. It follows that $T_{\#}$ and $\Phi^{-1}(T_{\#})$ both are also irreducible. \qed
\enddemo

\proclaim{Corollary \tnum{May2}}If $\#\in \operatorname{CP}$ and, either $c\ge 3${\rm;} or else, $c=2$ and $\#\le \{c,c,c\}$, then \roster
\item $S_{\#}=T_{\#}\setminus (\bigcup\limits_{\#'<\#}T_{\#'})$,
\item $S_{\emptyset}\cap \Bbb B_d$ is open in $\Bbb T_d$ and $S_{\emptyset}=(S_{\emptyset}\cap \Bbb B_d)\cup\Bbb U\Bbb B_d$ is the union of an open subset of $\Bbb T_d$ and a closed subset of $\Bbb T_d$,
\item $S_{\#}\cap \Bbb B_d$ is open in $T_{\#}\cap \Bbb B_d$ and is locally closed in $\Bbb T_d$,
\item $T_{\#}\cap \Bbb B_d$ is the closure of $S_{\#}\cap \Bbb B_d$ in $\Bbb B_d$
\item $\Bbb T_d$ is the closure of $S_{\emptyset}$ in $\Bbb T_d$,
\item $S_{c:c:c}=T_{c:c:c}$ is closed in $\Bbb B_d$, and
\item $S_{\#}\cap \Bbb B_d$ is irreducible.
\endroster
\endproclaim

\demo{Proof} Assertion (1) follows from the definition of $T_{\#}$. We next prove (2).
The subset $T_c$ of $\Bbb B_d$ is closed by Theorem \tref{strata}; so, $S_{\emptyset}\cap \Bbb B_d$, which is equal to $\Bbb B_d\setminus T_c$, is open in $\Bbb B_d$. But, $\Bbb B_d$ is open in $\Bbb T_d$; so $S_{\emptyset}\cap \Bbb B_d$ is also open in $\Bbb T_d$. The second part of (2) is established in (\tref{notn}) and  
Theorem \tref{YYY}. 
 The first assertion  of (3) follows from (1) and the fact that $\Bbb T_{\#}\cap \Bbb B_d$ is closed in $\Bbb B_d$. The second assertion of (3) also uses the fact that $\Bbb B_d$ is open in $\Bbb T_d$. 
We prove (4). The subset $S_{\#}\cap \Bbb B_d$ of $\Bbb T_{\#}\cap \Bbb B_d$ is open by (3) and non-empty by Remark \tref{fws}. Theorem \tref{strata} shows that   $\Bbb T_{\#}\cap \Bbb B_d$ is irreducible. 
For (5), use the fact that $\Bbb T_d$ is irreducible and that $S_{\emptyset}$ contains a non-empty open subset $S_{\emptyset}\cap \Bbb B_d$ of $\Bbb T_d$; see (2).
Assertion (6) is clear and (7) 
follows from (4) and the fact that $T_{\#}\cap \Bbb B_d$ is irreducible.  \qed
\enddemo

\proclaim{Theorem \tnum{stratif}} Assume $c\ge 3$. The picture {\rm(\tref{pictt})} gives a stratification of $\Bbb B_d$. In other words, $\Bbb B_d$ is the disjoint union of $\{S_{\#}\cap \Bbb B_d\mid \#\in \operatorname{CP}\}$ and if $\#'\le \#$ in $\operatorname{CP}$, then $S_{\#'}\cap \Bbb B_d$ is contained in   the closure $\overline{S_{\#}\cap \Bbb B_d}$ of $S_{\#}\cap \Bbb B_d$ in $\Bbb B_d$. \endproclaim

\demo{Proof} The fact that $\Bbb B_d$ is the disjoint union of $\{S_{\#}\cap \Bbb B_d\mid \#\in \operatorname{CP}\}$ follows from the geometric description of $S_{\#}$: any given curve has exactly one configuration of multiplicity $c$ singularities. The poset $\operatorname{CP}$ contains all possible configurations as was discussed at the beginning of Section 4. The second assertion follows from Corollary \tref{May2}~(4) and the definition of $T_{\#}$. \qed \enddemo

\proclaim{Theorem \tnum{dim}} Assume $c\ge 3$. For each fixed $\#\in \operatorname{CP}$, the sets $T_{\#}$ and $S_{\#}$ have the same dimension and this dimension is 
 given in the following chart{\rm:}
$$\matrix 
\#&c:c:c&c:c,c&c,c,c&c:c&c,c&c&\emptyset\\
\dim S_{\#}&3c+7&3c+8&3c+9&4c+6&4c+7&5c+5&6c+3.
\endmatrix$$
\endproclaim
\demo{Proof} The closure of $S_{\emptyset}$ in $\Bbb T_d$ is $T_{\emptyset}=\Bbb T_d$. It follows that
$$\dim S_{\emptyset}=\dim T_{\emptyset}=\Bbb T_d.$$ On the other hand, $\Bbb T_d$ is a non-empty open subset of $\Bbb A_d$; so, its dimension is $6c+3$. 

Henceforth, we consider $\#\in \operatorname{CP}$ with $\#\neq \emptyset$. We have seen that $ T_{\#}$ is irreducible and that $S_{\#}$ is a non-empty open subset of $T_{\#}$. Hence, $\dim  T_{\#}=\dim S_{\#}$. We now compute this dimension. 
Recall from (\tref{surj}) 
  that $G\times N_{\#}$, $T_{\#}$, and $\Phi^{-1}(T_{\#})$ all are irreducible.

Notice that for every $\pmb g\in T_{\#}$ with $\pmb g= \Phi(\varphi)$, we have $$\Phi^{-1}(\pmb g)=(\{I_{3\times 3}\}\times \operatorname{SL}_2(\pmb k))\varphi\cong \operatorname{SL}_2(\pmb k),$$ which has dimension 3. Thus, since $T_{\#}$ is irreducible, the surjective morphism on the right hand side of (\tref{surj}) yields  that $$\dim T_{\#}=\dim\Phi^{-1}(T_{\#})-3.\tag\tnum{form1}$$ 
Indeed, once we pass to the closure of $T_{\#}$ and $\Phi^{-1}(T_{\#})$, in their respective affine spaces, then the map $\Phi$ corresponds to an embedding of affine $k$-domains and the given formula follows from the additivity of transcendence degrees; see, for example, \cite{\rref{E95}, Cor~14.6}. We apply the same technique to the surjective morphism on the left hand side of (\tref{surj}) to see that 
$$\dim \Phi^{-1}(T_\#)=\dim (G\times N_{\#})-\dim \Upsilon^{-1}(\varphi)\tag\tnum{form2}$$ for general $\varphi$. Combine (\tref{form1}) and (\tref{form2}) to see $$\dim T_{\#}=\dim \Phi^{-1}(T_{\#})-3=\dim(G\times N_{\#})-\dim \Upsilon^{-1}(\varphi)-3;$$ and therefore,
$$\dim T_{\#}=\dim N_{\#}-\dim \Upsilon^{-1}(\varphi)+10.\tag\tnum{***}$$ 

We compute the dimensions of the fibers $\Upsilon^{-1}(\varphi)$, for $\varphi$ general in $\Phi^{-1}(T_{\#})$. 

Observe that $\Phi^{-1}(S_{\#})$ is a dense open subset of $\Phi^{-1}(T_{\#})$. Indeed, 
the set $\Phi^{-1}(T_{\#})$ is irreducible, the subset $S_{\#}$ of $T_{\#}$ is open, and the morphism $\Phi$ is continuous. Thus, we may assume that the general element $\varphi$ of $\Phi^{-1}(T_{\#})$ is actually in  $\varphi\in \Phi^{-1}(S_{\#})$ and we can write $\varphi=h\varphi_M$ for some $h\in G$ and $\varphi_M\in M_{\#}$, see (\tref{P3}).
Now $$\alignat 1\Upsilon^{-1}(\varphi)&{}=\{(g,\psi)\in G\times N_{\#}\mid g\psi=\varphi\}\\\allowdisplaybreak
&{}=\{(g,\psi)\in G\times N_{\#}\mid g\psi=h\varphi_M\}\\\allowdisplaybreak
&{}=\{(g,\psi)\in G\times N_{\#}\mid \psi=g^{-1}h\varphi_M\}\\\allowdisplaybreak
&{}=\{(g,g^{-1}h\varphi_M)\in G\times N_{\#}\}\\\allowdisplaybreak
&{}\cong \{g\in G\mid g^{-1}h\varphi_M\in N_{\#}\}\\\allowdisplaybreak
&{}= \{hk\in G\mid k^{-1}\varphi_M\in N_{\#}\}\\\allowdisplaybreak
&{}\cong \{k\in G\mid k^{-1}\varphi_M\in N_{\#}\}.\endalignat$$
Define $$F_{\#}(\varphi_M)=\{k\in G\mid k^{-1}\varphi_M\in N_{\#}\}.$$
We compute the dimension of $F_{\#}(\varphi_M)$ for various choices of $\#$.

First let  $\#=c$.   Recall that
$$N_c=\left\{\left (\matrix Q_1&Q_2\\Q_3&Q_4\\0&Q_5\endmatrix\right )\right\}\cap \Phi^{-1}(\Bbb B_d).$$
We show that$$F_c(\varphi_M)=\left\{\left (\bmatrix x_1&x_2&x_3\\x_4&x_5&x_6\\0&0&x_7\endbmatrix ,\left (\matrix y_1&y_2\\0&y_3\endmatrix\right )\right )\in G\right\},\tag\tnum{case1}$$which has dimension $10$.  Once (\tref{case1}) is established, then 
we use $\dim N_c=5c+5$ and (\tref{***}) to see that $\dim T_{c}=5c+5-10+10=5c+5$.

It is clear that the inclusion $(\supseteq)$ holds in (\tref{case1}). To show the other inclusion, let $k\in F_c(\varphi_M)$. Now $\varphi_M$, as well as $k^{-1}\varphi_M$ are in $N_c$. As $\varphi$ is general in $\Phi^{-1}(T_c)$ and this set contains matrices with five linearly independent entries, it follows that $\mu(I_1(\varphi))\ge 5$; therefore, $\mu(I_1(\varphi_M))\ge 5$ and $\mu(I_1(k^{-1}\varphi_M))\ge 5$. Thus,
$$\varphi_M=\left (\matrix Q_1&Q_2\\Q_3&Q_4\\0&Q_5\endmatrix\right )\quad\text{and}\quad k^{-1}\varphi_M=\left (\matrix Q_1'&Q_2'\\Q_3'&Q_4'\\0&Q_5'\endmatrix\right ),
$$ where $\{Q_1,\dots,Q_5\}$ and $\{Q_1',\dots,Q_5'\}$ are linearly independent sets. 

Let $\Cal C_1$ and $\Cal C_2$ be the curves parameterized by $\Phi(\varphi_M)$ and $\Phi(k^{-1}\varphi_M)$, respectively. Then $\Cal C_1$ and $\Cal C_2$ have exactly one singularity of multiplicity $c$, and this singularity is the point $[0:0:1]$. (See Theorem \tref{d=2c}). Write $k=(\chi^{-1},\xi)\in \operatorname{GL}_3(\pmb k)\times \operatorname{GL}_2(\pmb k)$. Notice that $\chi$ induces a linear automorphism of $\Bbb P^2$ that maps $\Cal C_1$ to $\Cal C_2$ (Remark \tref{sigh}) and hence leaves the point $[0:0:1]$ fixed. It follows that the third row of $\chi$ is $[0,0,a]$ for some $a\in \pmb k^*$. Since $$(\chi\varphi_M\xi)_{(3,1)}=(k^{-1}\varphi_M)_{(3,1)}=0,$$we obtain 
$[0,0,a]\varphi_M\xi_1=0$ when $\xi_1$ is the first column of $\xi$. Clearly, $$[0,0,a]\varphi_M\bmatrix 1\\0\endbmatrix=0.$$ Hence Remark \tref{CCC} gives $\xi_1=\bmatrix b\\0\endbmatrix$, for some $b\in \pmb k^*$. It follows that $k=(\chi^{-1},\xi)$
has the form of the elements in the right hand side of (\tref{case1}).

 For $\#=c,c$, recall that
$$N_{c,c}=\left\{\left. \bmatrix Q_1&Q_2\\Q_3&Q_4\\0&Q_5\endbmatrix \right \vert \dim(Q_3,Q_4,Q_5)\le 2\right\}\cap \Phi^{-1}(\Bbb B_d),$$ which has dimension $4c+4$. In this case, we have 
$$F_{c,c}(\varphi_M)=\{k\in G\mid k^{-1}\varphi_M\in N_{c,c}\}.$$ 
We now show that $$F_{c,c}(\varphi_M)=\{k\in G\mid k^{-1}\varphi_M\in M_{c,c}\}.$$ 
Indeed, $\varphi_M\in M_{c,c}\subseteq \Phi^{-1}(S_{c,c})$. Thus, $k^{-1}\varphi_M\in \Phi^{-1}(S_{c,c})$. Hence, if $k^{-1}\varphi_M\in N_{c,c}$, then $$k^{-1}\varphi_M\in N_{c,c}\cap \Phi^{-1}(S_{c,c})=M_{c,c}.$$ Recall that 
$$M_{c,c}=\left\{\left.\bmatrix Q_1&Q_2\\Q_3&Q_3\\0&Q_4\endbmatrix\right\vert \dim _{\pmb k}{<}Q_1,Q_2,Q_3,Q_4{>}=4\right\}\cap \Phi^{-1}(\Bbb B_d).$$ We now claim that
$$F_{c,c}(\varphi_M)=\left\{\left (\left [\smallmatrix x_1&x_2&x_3\\0&x_4&0\\0&0&x_5\endsmallmatrix\right ],\left [\smallmatrix y_1&y_2\\0&y_1-y_2\endsmallmatrix\right ]\right ) \right\} \cup \left\{\left (\left [\smallmatrix x_1&x_2&x_3\\0&0&x_4\\0&x_5&0\endsmallmatrix\right ],\left [\smallmatrix y_1&y_2\\-y_1&-y_1\endsmallmatrix\right ]\right ) \right\},\tag\tnum{2box}$$ which has dimension $7$. Thus, this claim, together with (\tref{***}) gives $\dim T_{c,c}=4c+7$. Clearly, the left hand side of (\tref{2box}) contains the right hand side. To show the other containment, let $k=(\chi^{-1},\xi)$ be an element of $F_{c,c}(\varphi_M)$. As before, let $\Cal C_1$ and $\Cal C_2$ be the curves parameterized by $\Phi(\varphi_M)$ and $\Phi(k^{-1}\varphi_M)$. Then $\Cal C_1$ and $\Cal C_2$ have exactly two singularities of multiplicity $c$: one at the point $[0:1:0]$ and the other at $[0:0:1]$. 
Since $\chi$ induces a linear automorphism of $\Bbb P^2$ mapping $\Cal C_1$ to $\Cal C_2$ and hence a permutation of the points  $[0:1:0]$ and   $[0:0:1]$, we have that $\chi$ is of the form 
$$\left [\smallmatrix x_1&x_2&x_3\\0&x_4&0\\0&0&x_5\endsmallmatrix\right ]\quad  \text{or}\quad \left [\smallmatrix x_1&x_2&x_3\\0&0&x_4\\0&x_5&0\endsmallmatrix\right ].$$ As noted before, the BiProj Lemma (Remark \tref{CCC}) now determines $\xi$ in both cases. This establishes the claim.

For $\#=c:c$, recall that $$N_{c:c}=\left\{ \bmatrix Q_1&Q_2\\Q_3&Q_4\\0&Q_3\endbmatrix\right\}\cap \Phi^{-1}(\Bbb B_d),$$ which has dimension $4c+4$. As before, we have
$$F_{c:c}(\varphi_M)=\{k\in G\mid k^{-1}\varphi_M\in N_{c:c}\}=\{k\in G\mid k^{-1}\varphi_M\in M_{c:c}\}.$$ Now recall that
$$M_{c:c}=\left\{\left. \bmatrix Q_1&Q_2\\Q_3&Q_4\\0&Q_3\endbmatrix\right\vert  \dim_{\pmb k}{<}Q_1,Q_2,Q_3,Q_4{>}=4 \right\}\cap \Phi^{-1}(\Bbb B_d).$$  We now claim that
$$F_{c:c}(\varphi_M)=\left\{\left.\left (\left [\smallmatrix x_1&x_2&x_3\\0&x_4&x_5\\0&0&x_6\endsmallmatrix\right ]^{-1},\left [\smallmatrix y_1&y_2\\0&y_3\endsmallmatrix\right ]\right )\right \vert x_4y_1=x_6y_3 \right\},\tag\tnum{bs}$$ which has dimension $8$. Thus, this claim, together with (\tref{***}), gives $\dim T_{c:c}=4c+6$. 

One easily sees that the left hand side of (\tref{bs}) contains the right hand side. To show the other containment, let $k=(\chi^{-1},\xi)$ be in $F_{c:c}(\varphi_M)$. As before, one proves that the last row of $\chi$ is of the form $[0,0,a]$ and the first column of $\xi$ is of the form $\xi_1=\left [\smallmatrix b\\0\endsmallmatrix\right ]$. Therefore, we can write
$$\chi=\bmatrix a_1&a_2&a_3\\a_4&a_5&a_6\\0&0&a_7\endbmatrix \quad\text{and}\quad \xi=\bmatrix b_1&b_2\\0&b_3\endbmatrix.$$ 
Writing $\varphi_M=\left [\smallmatrix Q_1&Q_2\\Q_3&Q_4\\0&Q_3\endsmallmatrix\right ] $ one obtains
$$[\chi\varphi_M\xi]_{(2,1)}=a_4b_1Q_1+a_5b_1Q_3\quad\text{and}\quad [\chi\varphi_M\xi]_{(3,2)}=a_7b_3Q_3.$$ Since $\chi \varphi_M\xi=k^{-1}\varphi_M\in M_{c:c}$, we conclude that $a_4b_1=0$ and $a_5b_1=a_7b_3$. As $\xi\in \operatorname{GL}_2(\pmb k)$, one has $b_1\neq 0$; therefore, $a_4=0$ and the claim now follows. 

For $\#=c,c,c$, recall that
$$N_{c,c,c}=\left\{\left. \left [ \smallmatrix Q_1&Q_2\\Q_3&Q_4\\0&Q_5\endsmallmatrix\right ]\right\vert\dim_{\pmb k}{<}Q_1,Q_2,Q_3,Q_4,Q_5{>}\le 3\right\}\cap \Phi^{-1}(\Bbb B_d),$$ which has dimension $3c+3$. As before,
$$F_{c,c,c}(\varphi_M)=\{k\in G\mid k^{-1}\varphi_M\in N_{c,c,c}\}=\{k\in G\mid k^{-1}\varphi_M\in M_{c,c,c}\}.$$ Now,
$$M_{c,c,c}=\left\{\left. \left [ \smallmatrix Q_1&Q_1\\Q_2&0\\0&Q_3\endsmallmatrix\right ]\right\vert\dim_{\pmb k}{<}Q_1,Q_2,Q_3{>}= 3\right\}\cap \Phi^{-1}(\Bbb B_d);$$ thus, we claim that $F_{c,c,c}(\varphi_M)$ is equal to  \phantom{\tnum{3dot}}
$$\alignat 2 &{}\phantom{{}\cup{}}\left\{\left (\left [\smallmatrix x_1&0&0\\0&x_2&0\\0&0&x_3\endsmallmatrix\right ],\left [\smallmatrix y&0\\0&y\endsmallmatrix\right ]\right )\in G\right\}&{}\cup{}& 
\left\{\left (\left [\smallmatrix 0&0&x_1\\0&x_2&0\\x_3&0&0\endsmallmatrix\right ],\left [\smallmatrix -y&0\\y&y\endsmallmatrix\right ]\right )\in G\right\}\\
&{}\cup\left\{\left (\left [\smallmatrix 0&x_1&0\\x_2&0&0\\0&0&x_3\endsmallmatrix\right ],\left [\smallmatrix y&y\\0&-y\endsmallmatrix\right ]\right )\in G\right\}&{}\cup{}& \left\{\left (\left [\smallmatrix 0&0&x_1\\x_2&0&0\\0&x_3&0\endsmallmatrix\right ]^{-1},\left [\smallmatrix 0&-y\\y&y\endsmallmatrix\right ]\right )\in G\right\}\tag{\tref{3dot}}\\
&{}\cup\left\{\left (\left [\smallmatrix 0&x_1&0\\0&0&x_2\\x_3&0&0\endsmallmatrix\right ]^{-1},\left [\smallmatrix y&y\\-y&0\endsmallmatrix\right ]\right )\in G\right\}&{}\cup{}&\left\{\left (\left [\smallmatrix x_1&0&0\\0&0&x_2\\0&x_3&0\endsmallmatrix\right ],\left [\smallmatrix 0&y\\y&0\endsmallmatrix\right ]\right )\in G\right\} ,\endalignat$$ 
which has dimension $4$. Thus, $\dim T_{c,c,c}=3c+9$. It is easy to see that subset of $G$ listed in (\tref{3dot}) is contained in $F_{c,c,c}(\varphi_M)$. To see the reverse inclusion, observe that the curves $\Cal C_1$ and $\Cal C_2$, parameterized by $\Phi(\varphi_M)$ and $\Phi(k^{-1}\varphi_M)$ have exactly three singularities of multiplicity $c$ and they are at points $[1:0:0]$, $[0:1:0]$, and $[0:0:1]$. Thus, if $k=(\chi^{-1},\xi)$ is in $F_{c,c,c}(\varphi_M)$, then the linear automorphism of $\Bbb P^2$ induced by $\chi$ gives a permutation of these $3$ points. Hence, $\chi$ is of the form 
$\left [\smallmatrix x_1&0&0\\0&x_2&0\\0&0&x_3\endsmallmatrix\right ]\sigma$, where $\sigma$ is one of the six permutation matrices. Once again the BiProj Lemma determines the shape of $\xi$, and the claim follows. 

If $\#=c:c:c$, recall that $$N_{c:c:c}= \left\{\left. \left [ \smallmatrix Q_1&Q_2\\Q_3&Q_1\\0&Q_3\endsmallmatrix\right ]\right\vert\dim_{\pmb k}{<}Q_1,Q_2,Q_3{>}= 3\right\}\cap \Phi^{-1}(\Bbb B_d)=M_{c:c:c},$$ which has dimension $3c+3$.
We claim that
$$F_{c:c:c}(\varphi_M)=\left\{\left.\left (\left [\smallmatrix x_1&x_2&x_3\\0&x_4&x_5\\0&0&x_6\endsmallmatrix\right ]^{-1},\left [\smallmatrix y_1&y_2\\0&y_3\endsmallmatrix\right ]\right )\in G\right \vert \smallmatrix x_6y_3=x_4y_1\\x_5y_3=x_2y_1-x_4y_2\\x_4y_3=x_1y_1\endsmallmatrix \right\}.\tag\tnum{dcirc}$$ Since $y_3\neq 0$, the dimension of the fiber is $6$; thus, $\dim T_{c:c:c}=3c+7$. A matrix calculation gives immediately the containment of the right hand side in the left hand side of (\tref{dcirc}). To prove the other containment, one sees as before, that for $k=(\chi^{-1},\xi)$ in $F_{c:c:c}(\varphi_M)$, one has $\chi=\bmatrix a_1&a_2&a_3\\a_4&a_5&a_6\\0&0&a_7\endbmatrix$ and $\xi=\bmatrix b_1&b_2\\0&b_3\endbmatrix$. Writing $\varphi_M= \bmatrix Q_1&Q_2\\Q_3&Q_1\\0&Q_3\endbmatrix$, we obtain
$$\chi \varphi_M\xi= \bmatrix a_1b_1Q_1+a_2b_1Q_3&*\\a_4b_1Q_1+a_5b_1Q_3&(a_4b_2+a_5b_3)Q_1+a_4b_3Q_2+(a_5b_2+a_6b_3)Q_3\\0&a_7b_3Q_3\endbmatrix.$$Since this matrix is in $M_{c:c:c}$, we conclude that $a_4b_1$ and $a_4b_3$ are both zero, which gives $a_4=0$ since neither $b_1$ nor $b_3$ can be zero. Using this, we also obtain 
$$a_2b_1=a_5b_2+a_6b_3,\quad a_1b_1=a_5b_3\quad\text{and}\quad a_5b_1=a_7b_3,$$ which are the required conditions.

Finally, to compute the dimension of $T_{c:c,c}$, notice that one has proper containments of closed irreducible subsets of $\Bbb B_d$: 
$$T_{c:c:c}\subsetneq T_{c:c,c}\subsetneq T_{c,c,c},$$ where $\dim T_{c:c:c}=3c+7$ while $\dim T_{c,c,c}=3c+9$. It follows that $\dim T_{c:c,c}=3c+8$. \qed
\enddemo
\bigpagebreak
\SectionNumber=7\tNumber=1
\heading Section \number\SectionNumber. \quad The Jacobian matrix and the ramification locus.
\endheading

Remark \tref{brn} provides a method of parameterizing the   branches of a parameterized curve. Theorem \tref{,T4.1} shows that the Jacobian matrix associated to the parameterization identifies the 
non-smooth branches of the curve as well as the multiplicity of each branch.  The starting point for this line of reasoning  is the result  that if $D$ is an algebra which is essentially of finite type over the ring $C$, then the ramification locus of $D$ over $C$ is equal to the support of the module of K\"ahler differentials $\Omega_{D/C}$. (See, for example,   \cite{\rref{Kz86}, Cor.~6.10}.)  In our ultimate application of this result, we consider the module of differentials 
$\Omega$ for the ring extension 
 $$\tsize \frac{\hat R}{J_i}\to \hat S_{\frak M_i}\tag\tnum{re}$$from the proof of Lemma \tref{.gl}. We have two presentations of the $\hat S_{\frak M_i}$-module $\Omega$. One presentation comes from the Jacobian matrix of the parameterization of the curve $\Cal C$. The other presentation comes from the geometry which gives rise to the ring extension (\tref{re}). The Fitting ideal $\operatorname{Fitt}_0(\Omega)$ may be computed using either presentation. 

In addition to \cite{\rref{Kz86}} one may consult \cite{\rref{E95}, Chapt.~16} or \cite{\rref{M80}, Sect.~26} for elementary facts and notation pertaining to K\"ahler differentials.

\proclaim{Theorem \tnum{,T4.1}} Adopt the Data of \tref{34.1} with $\pmb k$ an algebraically closed field of characteristic zero. Consider the inclusion map $\pmb k[I_d]\subseteq B$ of homogeneous coordinate rings which is induced by the morphism $\Psi\:\Bbb P^1\to \Cal C$. The $\gcd$ of the zeroth Fitting ideal of $\Omega_{B/\pmb k[I_d]}$ is a polynomial in $B$. Let 
$$\gcd \operatorname{Fitt}_0(\Omega_{B/\pmb k[I_d]})=\prod_{i=1}^s\ell_i^{f_i},$$ where $(\ell_1),\dots,(\ell_s)$ are distinct height one linear ideals of $B$. If $(\ell)$ is an arbitrary height one linear ideal of $B$ and $\Cal C(\ell)$ is the branch of $\Cal C$ which corresponds to $\ell$, in the sense of Remark \tref{brn}, then  the multiplicity of  the branch $\Cal C(\ell)$ is $$\cases f_i+1&\text{if $(\ell)=(\ell_i)$ for some $i$}\\1&\text{otherwise}.\endcases$$
Furthermore, the Fitting ideal $\operatorname{Fitt}_0(\Omega_{B/\pmb k[I_d]})$ is equal to the ideal $I_2(N)$ of $B$, where
$N$ is  the $2 \times n$ transposed Jacobian matrix $$N=\bmatrix \frac{\partial g_1}{\partial x}&\dots & \frac{\partial g_n}{\partial x}\\
 \frac{\partial g_1}{\partial y}&\dots & \frac{\partial g_n}{\partial y}\endbmatrix.$$
\endproclaim
Before proving  Theorem \tref{,T4.1}, we describe various special cases.  Corollary \tref{cl} follows immediately from Theorem \tref{,T4.1} without any further proof. Also, Theorem \tref{Cor2'} 
requires only a small amount of additional proof. 

\proclaim{Corollary \tnum{cl}} Retain the notation and hypotheses of Theorem \tref{,T4.1}.
\roster
\item All of the branches through all of the points of $\Cal C$ are smooth if and only if $\operatorname{ht} I_2(N)=2$.
\item The multiplicity of each branch of $\Cal C$ is at most two if and only if the $\gcd$ of  $I_2(N)$ decomposes into a product of pairwise non-associated   linear factors.\endroster \endproclaim 

\proclaim{Theorem \tnum{Cor2'}} 
Adopt Data \tref{34.1} with $\pmb k$ an algebraically closed field of characteristic zero.
Let  $p_1,\dots,p_z$ be the   singularities of $\Cal C$.
For each singular point $p_j$, let  $m_j$ be the multiplicity of $\Cal C$ at $p_j$  and  
$s_j$ be the number of branches of $\Cal C$ at $p_j$.
 For each index $j$, with ${1\le j\le z}$,  let $\dsize \gcd I_1( p_j\varphi)=\prod_{u=1}^{s_j}\ell_{u,j}^{e_{u,j}}$, where the $\ell_{u,j}$ are pairwise non-associated linear factors and the exponents $e_{u,j}$ are positive.
Let $$N=\bmatrix \frac{\partial g_1}{\partial x}&\dots & \frac{\partial g_n}{\partial x}\\
 \frac{\partial g_1}{\partial y}&\dots & \frac{\partial g_n}{\partial y}\endbmatrix$$ be the $2 \times n$ Jacobian matrix of the parametrization.
Then  
\roster
\item    
$\gcd I_2(N)=\dsize \prod_{j=1}^{z}\prod_{u=1}^{s_j}\ell_{u,j}^{e_{u,j}-1}$, and  
\item the degree of $\gcd I_2(N)$ is equal to $\dsize \sum_{j=1}^z(m_j-s_j)$.
\endroster\endproclaim 

\demo{Proof}We are given the singular points $p_1,\dots,p_z$ on $\Cal C$ and the factorizations $\gcd I_1(p_j \varphi)=\prod \ell_{u,j}^{e_{u,j}}$. Lemma \tref{.gl} tells us that the multiplicity of the branch $\Cal C(\ell_{u,j})$ of $\Cal C$ is $e_{u,j}$. Assertion (1) is now an immediate consequence Theorem \tref{,T4.1}. We prove (2). Theorem \tref{Cor1} shows that for each $j$, with $1\le j\le z$, the polynomial 
$\gcd I_1( p_j\varphi)$ has degree $m_j$ and $s_j$ pairwise non-associated linear factors; so, one has 
$$\dsize \deg \prod_{u=1}^{s_j}\ell_{u,j}^{e_{u,j}-1}=\deg \prod_{u=1}^{s_j}\ell_{u,j}^{e_{u,j}}-s_j=m_j-s_j;$$
hence, (1) gives 
$$\deg \gcd I_2(N)=\deg \prod_{j=1}^z\prod_{u=1}^{s_j}\ell_{u,j}^{e_{u,j}-1}=\sum_{j=1}^z(m_j-s_j). \qed$$
\enddemo

\demo{\bf Proof of Theorem \tref{,T4.1}}The relative cotangent complex
$$ \Omega_{\pmb k[I_d]/\pmb k}\otimes_{\pmb k}B\to \Omega_{B/\pmb k}\to \Omega_{B/\pmb k[I_d]}\to 0$$
gives rise to the presentation
$$B^n@> N>> B^2 @>>>\Omega_{B/\pmb k[I_d]}\to 0\tag\tnum{,.prez}$$
of  $\Omega_{B/\pmb k[I_d]}$ as a $B$-module. It follows immediately that the Fitting ideal of the $B$-module $\Omega_{B/\pmb k[I_d]}$ is 
$$\operatorname{Fitt}_0(\Omega_{B/\pmb k[I_d]})=I_2(N)B.$$
Fix a point $q\in \Bbb P^1$ and  a non-zero linear form $\ell\in B$ with $\ell(q)=0$. Let $p$ be the point $(g_1(q),\dots,g_n(q))$  on $\Cal C$ in $\Bbb P^{n-1}$ and $e$ be the multiplicity of the branch $\Cal C(\ell)$. Define   $f$ to be the exponent with   $$\gcd I_2(N)=\ell^f\cdot \theta,\tag\tnum{,exp}$$ where $\theta$ is a  polynomial in $B$ which is   relatively prime to $\ell$. 
We prove $e=f+1$.
Let $\frak p$ be the ideal $$I_2\pmatrix g_1(q)&\dots&g_n(q)\\ g_1&\dots&g_n\endpmatrix$$of $\pmb k[I_d]$.

The formation of $\Omega$ is preserved under base change. That is, if $C'$ and $D$ are $C$-algebras, then
$$C'\otimes_C\Omega_{D/C}=\Omega_{(C'\otimes_C D)/C'}.$$ In particular, any presentation of $\Omega_{D/C}$ as an $D$-module:
$$D^a@> \sigma >> D^b @> \tau >> \Omega_{D/C}\to 0$$
gives rise to a presentation of $\Omega_{(C'\otimes_C D)/C'}$ as an $C'\otimes_C D$-module:
 $$(C'\otimes_C D)^a@> \sigma >> (C'\otimes_C D)^b @> \tau >> C'\otimes_C \Omega_{D/C}=\Omega_{(C'\otimes_C D)/C'} \to 0.$$
For example, if $C'=U^{-1}C$ for some multiplicatively closed subset $U$ of $C$, then we  write $U^{-1}D$ in place of $U^{-1}C\otimes_CD$; so we have the presentation
$$(U^{-1}D)^a@> \sigma >> (U^{-1}D)^b @> \tau >> \Omega_{U^{-1}D/U^{-1}C} \to 0.$$
In our situation, we  localize at $U=\pmb k[I_d]\setminus \frak p$. We write $\pmb k[I_d]_{\frak p}$ in place of $U^{-1}\pmb k[I_d]$ and $B_{\frak p}$ in place of $U^{-1}(B)$. Apply the base change 
$\pmb k[I_d]_{\frak p}\otimes_{\pmb k[I_d]}\underline{\phantom{X}}$ to (\tref{,.prez}) to obtain the following presentation by free $B_{\frak p}$-modules
$$B_{\frak p}^n@> N>> B_{\frak p}^2 @>>>\Omega_{B_{\frak p}/\pmb k[I_d]_{\frak p}}\to 0.\tag\tnum{,.nonum}$$
In the statement of Lemma \tref{.gl} we have called $\pmb k[I_d]_{\frak p}=R\subseteq S=B_{\frak p}$. 
The ring $R$ is local with maximal ideal $\frak m_R=\frak pR_{\frak p}$. In this new language, (\tref{,.nonum}) becomes the exact sequence of $S$-modules:
$$S^n@> N>> S^2 @>>>\Omega_{S/R}\to 0.\tag\tnum{,.non1}$$
Complete both $R\subseteq S$ in the $\frak m_R$-adic topology to obtain the rings $\hat R\subseteq \hat S$. One of the maximal ideals of $\hat S$ is $(\ell)\hat S$, which we denote by $\frak M$. Let $J$ be the kernel of $\hat R\ \to \hat S_{\frak M}$. It is shown at (\tref{shw}) that the multiplicity $e(\hat R/J)$ is equal to $e$, which is the multiplicity of  the branch $\Cal C(\ell)$.

Apply the base change $\hat R\otimes_R\underline{\phantom{X}}$ to (\tref{,.non1}) to obtain the exact sequence of $\hat R\otimes_R S=\hat S$ modules:
$$\hat S^n@> N>> \hat S^2 @>>>\Omega_{\hat S/\hat R}\to 0.$$
Localize at the multiplicatively closed set $\hat S\setminus \frak M$ of $\hat S$ to obtain an exact sequence of $\hat S_{\frak M}$-modules:
$$\hat S_{\frak M}^n@> N>> \hat S_{\frak M}^2 @>>>\Omega_{\hat S_{\frak M}/\hat R}\to 0.$$ Apply the base change $\hat R/J\otimes_{\hat R}\underline{\phantom{X}}$. Keep in mind that $\hat R/J\otimes_{\hat R}\hat S_{\frak M}=\hat S_{\frak M}$. Obtain an exact sequence of $\hat S_{\frak M}$-modules 
$$\hat S_{\frak M}^n@> N>> \hat S_{\frak M}^2 @>>>\Omega_{\hat S_{\frak M}/\frac{\hat R}J}\to 0.$$ The zeroth Fitting ideal of $\Omega_{\hat S_{\frak M}/\frac{\hat R}J}$ is $$ \operatorname{Fitt}_0(\Omega_{\hat S_{\frak M}/\frac{\hat R}J})=I_2(N)\hat S_{\frak M}=(\ell^f)\hat S_{\frak M}=\frak M^f\hat S_{\frak M}.\tag\tnum{,Fitt1}$$ 

Now we calculate $\Omega_{\hat S_{\frak M}/\frac{\hat R}J}$ in a completely different manner. Recall the Veronese ring $\pmb k[B_d]$ and the ring  $T=\pmb k[B_d]_\frak p$ from Lemma \tref{.gl}. The completion of $T$ in the $\frak m_R$-adic topology is denoted $\hat T$. We have $$\hat R/J\subseteq \hat T_{\frak M \cap \hat T}\subseteq \hat S_{\frak M},$$ with $\hat T_{\frak M \cap \hat T}$ equal to the integral closure of $\hat R/J$. The rings $\hat T_{\frak M \cap \hat T}$ and $\hat R/J$ share the same residue class field, which, in the language   of the proof of Lemma \tref{.gl},  was called $\pmb k(g_n')$. Furthermore, $\pmb k(g_n')\subseteq \hat R/J$. The rings $\hat T_{\frak M \cap \hat T}$ and $\hat S_{\frak M}$ are complete DVRs with the same uniformizing parameter $t=\frac \ell  m$, where $m$ is any linear form in $B$ for which $\ell,m$ is a basis for the vector space $B_1$. The ring $\hat T_{\frak M \cap \hat T}$ is equal to $\pmb k(m^d)[[t]]$ and the ring $\hat S_{\frak M}$ is equal to $\pmb k(m)[[t]]$. It was observed above that  $e(\hat R/J)=e$. Proposition \tref{.P33.16} shows that $\Omega_{\hat S_{\frak M}/\frac{\hat R}J}$ is isomorphic to $\hat S_{\frak M}/{(t^{e-1})}$. We conclude that the Fitting ideal of $\Omega_{\hat S_{\frak M}/\frac{\hat R}J}$ is also equal to $$\operatorname{Fitt}_0(\Omega_{\hat S_{\frak M}/\frac{\hat R}J})=(t^{e-1})\hat S_{\frak M}=\frak M^{e-1}\hat S_{\frak M}.\tag\tnum{,Fitt2}$$ Compare (\tref{,Fitt1}) and (\tref{,Fitt2}) to see that  $f=e-1$.
\qed\enddemo

\proclaim{Proposition \tnum{.P33.16}} Let $K\subseteq A\subseteq B\subseteq C$ be local rings. Assume that 
\roster \item $B=L[[t]]$ and $C=M[[t]]$ are formal power series rings in one variable where $L\subseteq M$ are fields of characteristic zero with $\dim_LM$ finite,
\item $K$ is a field and   the natural maps 
$$K\to A\to A/\frak m_A\quad\text{and}\quad K\to B\to B/(t)=L$$ are isomorphisms, 
\item  $B$   finitely generated as an $A$-module, and 
\item $B\subseteq \operatorname{Quot}(A)$. 
\endroster
Then the $C$-modules $\Omega_{C/A}$ and $\frac{C}{(t^{e-1})C}$ are isomorphic, where $e=e(A)$ is the multiplicity of the local ring $A$.\endproclaim 

\demo{Proof}The relative cotangent sequence gives  an   exact sequence of $C$-modules:
   $$\Omega_{A/K}\otimes_AC@> \alpha>> \Omega_{C/K}@> \beta>>\Omega_{C/A}\to 0,$$ where
$\alpha(da\otimes c)=(da)c$ and $\beta(dc)=dc$. The ring $C$ is  generated as an $A$-algebra by  $t$ together with a finite generating set for $M$ over $K$; so, $\Omega_{C/A}$ is finitely generated as a $C$-module. The ring $C$ is local; so the Krull Intersection Theorem guarantees that $\beta$ sends $\cap(t^n)\Omega_{C/K}$ to zero. Thus, the above exact sequence induces the  exact sequence of $C$-modules
$$\Omega_{A/K}\otimes_AC@> \bar\alpha>> \overline{\Omega}_{C/K}@> \bar\beta>>\Omega_{C/A}\to 0,$$ where
$$\overline{\Omega}_{C/K}=\frac{\Omega_{C/K}}{\cap(t^n)\Omega_{C/K}}$$ and $\bar\alpha$ and $\bar \beta$ are induced by $\alpha$ and $\beta$. 
If $\omega\in \Omega_{C/K}$, then we write $\bar\omega$ for the image of $\omega$ in $\overline{\Omega}_{C/K}$. 
The field extension $K\subseteq M$ is separable and algebraic so the universal derivation $d=d_{C/K}\:C\to \Omega_{C/K}$ sends $M$ to zero. Therefore, if $f\in C$, then the elements $df$ and $f'dt$ of ${\Omega}_{C/K}$ represent the same class in $\overline{\Omega}_{C/K}$. 
It follows that $\overline{\Omega}_{C/K}$ is generated as a $C$-module by $\overline{dt}$. 
To complete the proof we show that \roster \item $\overline{\Omega}_{C/K}$ is a free $C$-module and \item the image of $\bar \alpha$ is  equal to $(t^{e-1})C\overline{dt}$.\endroster
We prove (2) first. The one-dimensional domains  $A\subseteq B$ are local and $B$ is the integral closure of $A$; so, Lemma \tref{.O33.10'} guarantees the existence of an element $z$ in $\frak m_A$ such that 
$z=t^{e}+$ higher order terms and  $\frak m_AB=zB$. It follows that $\frak m_AC=zC$. The image of $\bar \alpha$ is the $C$-submodule of $\overline{\Omega}_{C/K}=C\overline{dt}$ which is generated by $\bar\alpha(dz)$ and this is the equal to $(t^{e-1})C\overline{dt}$ since the field $K$ has characteristic zero.

 Now we prove (1). Suppose that $\theta\in C$ and that the element $\theta dt$ of $\Omega_{C/K}$ is in $\bigcap\limits_{i=0}^{\infty}(t^i)\Omega_{C/K}$. We prove that $\theta$ is zero in $C$. Let $n$ be a positive integer. Consider   the conormal exact sequence 
$$(t^n)/(t^{2n})\to \Omega_{C/K}\otimes_CC/(t^n)\to \Omega_{\frac{C}{(t^n)}/K}\to 0\tag\tnum{.con}$$
of $C/(t^n)$-modules associated to the $K$-algebra homomorphism $C\to C/(t^n)$. The sequence (\tref{.con}) induces an isomorphism
$$\frac{\Omega_{C/K}}{(t^{n-1})Cdt+(t^n)C\Omega_{C/K}}@>\cong>> \Omega_{\frac{C}{(t^n)}/K},\tag\tnum{.osi}$$ which is given by
$\text{class of $df$}\mapsto d(\text{class of $f$})$, for all $f\in C$. The element $\theta dt$ of $\Omega_{C/K}$ represents the class of  zero in the module on the left side of (\tref{.osi}); so $\bar \theta dt$ is zero in $\Omega_{\frac{C}{(t^n)}/K}$, where $\bar\theta$ is the image of $\theta$ in $C/(t^n)$. On the other hand, $C/(t^n)=M[t]/(t^n)$ and it is well-known that $1\mapsto dt$ gives an isomorphism $$\tsize\frac{M[t]}{(t^{n-1})}@>\cong>> \Omega_{\frac{M[t]}{(t^n)}/K}.$$ Thus, the image $\bar\theta$ of $\theta$ in  $C/(t^{n-1})$ is zero. This process may be repeated for all $n$. We conclude that $\theta$ is zero in $C$.  \qed \enddemo

\bigpagebreak
\SectionNumber=8\tNumber=1
\heading Section \number\SectionNumber. \quad  The conductor and the branches of a rational plane curve.
\endheading 

Let $\pmb g$ be an element of $\Bbb T_d$ and $\Cal C_{\pmb g}$ be the corresponding parameterized plane curve. We produce a homogeneous polynomial $c_{\pmb g}$ in $B=\pmb k[x,y]$. The linear factors of $c_{\pmb g}$ correspond to the branches of $\Cal C_{\pmb g}$ at its singularities. (Recall from Remark \tref{brn} that the given parameterization of $\Cal C_{\pmb g}$ induces a one-to-one correspondence between the set of height one linear ideals of $B$ and the branches of $\Cal C_{\pmb g}$.) Moreover, the exponents that appear on the linear factors of $c_{\pmb g}$ give information about the singularity degree $\delta_p$ at the corresponding singularity. See Theorem \tref{52.6} for the general statement and Theorem \tref{52.7} for the form of $c_{\pmb g}$ for each of the $13$ possible configurations of singularities on a quartic. We are able to capture some information about $\delta_p$ by way of  $c_{\pmb g}$ in a polynomial manner from the coefficients of the entries of a Hilbert-Burch matrix for $\pmb g$. If $d$ is even and $\pmb g\in \Bbb B_d$, then a Hilbert-Burch matrix for $\pmb g$ may be built in a polynomial manner from the coefficients of $\pmb g$, see Theorem \tref{T1}. Therefore,   we can use features of the factorization of $c_{\pmb g}$ to separate, by way of open and closed subsets of $\Bbb B_d$, various  configurations of singularities; see Theorem \tref{52.9} and Proposition \tref{52.10}.  

There are $13$ possible configurations of singularities on a rational plane quartic, see (\tref{intro}). The techniques of Sections 1, 2, and 3 (which involve the multiplicity  $m_p$ at each singular point $p$) and Section 7 (which involve the  number of branches $s_p$ at each singular point $p$) are able to separate, using open and closed sets in $\Bbb T_4$, the 13 possible configurations of singularities into 12 subsets. The techniques of the previous  sections are not able to separate 
those curves whose singularities are described by $(2:2:1),(2:1,1)$ from those curves whose singularities are described by $(2:2:1,1),(2:1)$. (This notation is explained at (\tref{osc}).)
We were motivated to prove Theorem \tref{52.9} and Proposition \tref{52.10} in order to separate these two configurations of singularities.  

The polynomial $c_{\pmb g}$ is defined in Theorem \tref{52.5}. Let $\frak c_{\pmb g}$ be the conductor of the $d^{\text{th}}$ Veronese subring $\pmb V$ of $B=\pmb k[x,y]$ into the coordinate ring $A_{\pmb g}=\pmb k[g_1,g_2,g_3]$ of $\Cal C_{\pmb g}$. 
In Corollary \tref{52.2} we show that $\pmb V/\frak c_{\pmb g}$ is Cohen-Macaulay and we calculate its multiplicity. 
The polynomial $c_{\pmb g}$ is the $\gcd$ of the extension $\frak c_{\pmb g}\cdot B$ of the ideal $\frak c_{\pmb g}$ in $\pmb V$ to the larger ring $B$. The extension of $\frak c_{\pmb g}$ to $B$ does not define a Cohen-Macaulay quotient; however, the ease of computation in $B$ (as compared to $\pmb V$) compensates for the inconvenience of having to saturate $\frak c_{\pmb g}$ in order to produce $c_{\pmb g}$.
Corollary \tref{52.4} expresses $\frak c_{\pmb g}$ in terms of the Hilbert-Burch matrix for $\pmb g$. Theorem \tref{52.6} explains the geometric significance of $c_{\pmb g}$. 

Laurent Bus\'e has shown us that polynomials similar to $c_{\pmb g}$ have been studied elsewhere in the literature. Abhyankar \cite{\rref{Ab}, Pg.~153} called his version of this polynomial a Taylor resultant. Bus\'e further developed the concept of Taylor resultants in paragraph 4.4 of \cite{\rref{B09}}. The Taylor resultant is essentially the determinant of a square matrix over a Principal Ideal Domain. Recently Bus\'e and D'Andrea \cite{\rref{BD'A}} examined all of the invariant factors of the square matrix, not only the determinant, in order to gain much more detailed information about the blow-up history of the corresponding singularity. Even more recently, Bus\'e and Luu Ba \cite{\rref{BB}, Sect.~5} have applied these ideas to the study of singularities on rational space curves. In this context the matrix in question is no longer square; so, they use Fitting ideals in place of invariant factors.

\definition{Data \tnum{data52}}In this section $\pmb k$ continues to be a field, $B$ continues to be the polynomial ring $B=\pmb k[x,y]$, and $d\ge 3$ is a fixed integer. Consider an ordered triple $\pmb g=(g_1,g_2,g_3)$ in the set $\Bbb T_d$ of Definition \tref{.'D27.14}. Let $A_{\pmb g}$ be the subring $\pmb k[g_1,g_2,g_3]$ of $B$ and $T$ be the polynomial ring $\pmb k[T_1,T_2,T_3]$. The ring $A_{\pmb g}$ is a $T$-module by way of the $\pmb k$-algebra homomorphism $T\twoheadrightarrow A_{\pmb g}$ which sends $T_i$ to $g_i$. Indeed, $A_{\pmb g}\cong T/(f_{\pmb g})$ as described in (\tref{506}). Let $\pmb V$ be the $d^{\text{th}}$ Veronese subring of $B$; so $\pmb V$ is the standard graded $\pmb k$-algebra $\pmb V=\bigoplus_{i} \pmb V_i$, with $\pmb V_i=B_{id}$. 
We have
$A_{\pmb g}\subseteq \pmb V\subseteq B$. 
The assumption that $\pmb g$ is in $\Bbb T_d$ ensures that    $(g_1,g_2,g_3)B$ is an ideal in  $B$   of height $2$ and that $A_{\pmb g}\subseteq \pmb V$ is a birational extension.   
Define the conductor $\frak c_{\pmb g}$ by $\frak c_{\pmb g}=A_{\pmb g}\:\!\!_Q\pmb V$, where $Q=\operatorname{Quot} (A_{\pmb g})$. Notice that $\frak c_{\pmb g}$ is an ideal of $\pmb V$ and also is an ideal of $A_{\pmb g}$. Let $L_{\pmb g}$ be the preimage of $\frak c_{\pmb g}$ in $T$. It follows that $T/L_{\pmb g}\cong A_{\pmb g}/\frak c_{\pmb g}$. \enddefinition

\proclaim{Lemma \tnum{52.1}}Adopt Data {\rm\tref{data52}}. The homogeneous minimal free resolution  of $T/L_{\pmb g}\cong A_{\pmb g}/\frak c_{\pmb g}$ by free $T$-modules has the form
$$0\to T^{d-2}(-d+1)\to T^{d-1}(-d+2)\to T\to T/L_{\pmb g}\to 0.$$Furthermore, $T/L_{\pmb g}$ is Cohen-Macaulay and $\omega_{T/L_{\pmb g}}\cong \omega_{A_{\pmb g}/\frak c_{\pmb g}}\cong \pmb V/A_{\pmb g}$.\endproclaim

\demo{Proof} The $a$-invariant of $B$ is $-2$; hence, $a(\pmb V)=\lfloor\frac {-2}{d}\rfloor=-1$. Therefore, $$\operatorname{reg}(\pmb V)=-1+\dim \pmb V=1$$ because $\pmb V$ is Cohen-Macaulay. 

Since $\pmb V$ is a finitely generated $T$-module, we have 
$$\operatorname{reg}_T(\pmb V)=\operatorname{reg}(\pmb V)=1.$$ Hence,  $\pmb V$ is generated in degree at most $1$ as a $T$-module. Therefore, $\pmb V/A_{\pmb g}$ is generated in degree at most $1$ as a $T$-module. Since $[\pmb V/A_{\pmb g}]_0=0$, it follows that $\pmb V/A_{\pmb g}$ is generated as a $T$-module in degree $1$. But $[\pmb V/A_{\pmb g}]_1=\pmb V_1/(A_{\pmb g})_1=\pmb k^{d+1-3}=\pmb k^{d-2}$. Therefore, $\pmb V/A_{\pmb g}$ is minimally generated by $d-2$ elements of degree $1$ as a graded $T$-module. 

The rings $A_{\pmb g}$ and $\pmb V$ are Cohen-Macaulay of dimension $2$; hence, $A_{\pmb g}$ and $\pmb V$ are Cohen-Macaulay $T$-modules of dimension $2$. It follows that each of these $T$-modules is perfect of projective dimension equal to $1$. Consider the $T$-free resolutions of these modules:
$$\CD 0@>>> T(-d)@> f_{\pmb g} >> T@>>>A_{\pmb g}\\@. @V VV @V VV \\ 0@>>> G_1 @>>> G_0@> >> \pmb V,\endCD \tag\tnum{52.*}$$for some graded free $T$-modules $G_0$ and $G_1$. The mapping cone of (\tref{52.*}) gives a $T$-free resolution of $\pmb V/A_{\pmb g}$ of length $2$. On the other hand, $\operatorname{grade} \operatorname{ann}_T(\pmb V/A_{\pmb g})$ is at least two. Indeed,  $f_{\pmb g}$ is a prime element of $T$ which is in $  \operatorname{ann}_T(\pmb V/A_{\pmb g})$. Furthermore,
the extension $A_{\pmb g}\subseteq \pmb V$ is birational and $\pmb V$ is finitely generated as an $A_{\pmb g}$-module; so there is a non-zero homogeneous element of $A_{\pmb g}$  which is in $\operatorname{ann}_T(\pmb V/A_{\pmb g})$. Lift this element back to $T$ to complete  a regular sequence of length $2$ on $T$   in $\operatorname{ann}_T (\pmb V/A_{\pmb g})$.  It follows that $\pmb V/A_{\pmb g}$ is a perfect $T$-module with $\operatorname{pd}_T (\pmb V/A_{\pmb g})=\operatorname{grade}\operatorname{ann}_T(\pmb V/A_{\pmb g})=2$. In particular, $\pmb V/A_{\pmb g}$ is a Cohen-Macualay $T$-module of dimension $1$ and the $T$-free resolution of $\pmb V/A_{\pmb g}$ has the form
$$0\to T(-d)\to F_1\to T^{d-2}(-1)\to \pmb V/A_{\pmb g},$$ for some graded free $T$-module $F_1$. Dualize the above exact sequence and shift by $d$ in order to obtain the exact sequence:
$$0\to T^{d-2}(-d+1)\to F_1^*(-d)\to T\to T/H\to 0,$$ where $H$ is a perfect ideal of height $2$. Since $\operatorname{rank} F_1^*=d-1$, we have $\mu(H)=d-1$; and hence, the initial degree of $H$ is at least $d-2$ by the Hilbert-Burch Theorem and therefore, 
$$F_1^*(-d)=T^{d-1}(-d+2).$$ Finally, $\pmb V/A_{\pmb g}\cong \omega_{T/H}$; and therefore, $\pmb V/A_{\pmb g}$ is a faithful $T/H$-module. Thus, $H=\operatorname{ann}_T(\pmb V/A_{\pmb g})$. But the later is $L_{\pmb g}$ by the definition of $L_{\pmb g}$; therefore, $H=L_{\pmb g}$. \qed \enddemo

\proclaim{Corollary \tnum{52.2}} Adopt Data {\rm\tref{data52}}. The following statements hold.
\item{\rm (a)} The ring $A_{\pmb g}/\frak c_{\pmb g}$ is $1$-dimensional and Cohen-Macaulay of multiplicity $\binom {d-1}2$.
\item{\rm(b)} The ring $\pmb V/\frak c_{\pmb g}$ is $1$-dimensional and Cohen-Macaulay of multiplicity $(d-1)(d-2)$.
\endproclaim
\demo{Proof} Assertion (a) follows from the resolution in Lemma \tref{52.1}, using for instance \cite{\rref{HM}}. We prove (b). By Lemma \tref{52.1} we have the short exact sequence
$$0\to A_{\pmb g}/\frak c_{\pmb g}\to \pmb V/\frak c_{\pmb g}\to \pmb V/A_{\pmb g}\cong \omega_{A_{\pmb g}/\frak c_{\pmb g}}\to 0.$$Since $\omega_{A_{\pmb g}/\frak c_{\pmb g}}$ is a Cohen-Macaulay module of dimension one and multiplicity $\binom{d-1}{2}$ by part (a), it follows that $\pmb V/\frak c_{\pmb g}$ is a Cohen-Macaulay $1$-dimensional ring of multiplicity $2\binom{d-1}2=(d-1)(d-2)$. 
\qed \enddemo

We define the matrix $M'_{\pmb g}$ which appears in the next result.  Let $\varphi_{\pmb g}$ be a homogeneous Hilbert-Burch matrix for the row vector $[g_1,g_2,g_3]$ of degree $(d_1,d_2)$ and let $$N=\bmatrix \rho^{(d_2)}&0\\0&\rho^{(d_1)}\endbmatrix,$$ for $\rho^{(i)}$ as given in (\tref{rho}). Each entry in the row vector $[T_1,T_2,T_3]\varphi_{\pmb g} N$ is a homogeneous form of degree $1$ in the $T$'s and degree $d$ in $x,y$.  Define $M'_{\pmb g}$ to be the matrix with homogeneous linear entries from $T$ which satisfies the equation
$$[T_1,T_2,T_3]\varphi_{\pmb g} N=\rho^{(d)}M'_{\pmb g}.\tag\tnum{52.3-}$$

\proclaim{Proposition \tnum{52.3}} Adopt Data {\rm\tref{data52}}. A homogeneous $T$-presentation of $\pmb V_+$ is given by 
$$T^{d+2}(-2)@>M'_{\pmb g}>> T^{d+1}(-1)@>>>\pmb V_+\to 0$$\endproclaim

\demo{Proof}The sequence $0\to \pmb V_+\to \pmb V\to \pmb k\to 0$ shows that $\operatorname{reg}_T\pmb V_+=1$; so, $\pmb V_+$ as a $T$-module is minimally generated by the monomials $y^d,\dots, x^d$ and the relations between those generators have degree $1$. Let
$$\pmb \ell=\bmatrix \ell_0\\\vdots\\\ell_d\endbmatrix$$ be a matrix of linear forms from $T$, and let $\Lambda$ be a matrix of scalars with $$\pmb \ell=\Lambda \bmatrix T_1\\T_2\\T_3\endbmatrix.$$
The polynomial $T_i$ of $T$ acts like $g_i$ on $A_{\pmb g}$; so
$$\rho^{(d)}\pmb \ell=0 \text{ in }\pmb V_+\iff \rho^{(d)} \Lambda \bmatrix g_1\\g_2\\g_3\endbmatrix =0\text{ in } B.$$The matrix $\varphi_{\pmb g}$ is a Hilbert-Burch matrix for the row vector $[g_1,g_2,g_3]$; so 
$\rho^{(d)}\pmb \ell=0$   in $\pmb V_+$ if and only if there exist homogeneous forms $F_1$ and $F_2$ in $B$ with $\deg F_1=d_2$, $\deg F_2=d_1$, and $$\Lambda^{\text{\rm T}} (\rho^{(d)})^{\text{\rm T}}=\varphi_{\pmb g}\bmatrix F_1\\F_2\endbmatrix.$$ The vector $\left [\smallmatrix F_1\\F_2\endsmallmatrix\right ]$ is equal to $N\pmb F$ for some columns vector of scalars $\pmb F$. We have  
$$\split \rho^{(d)}\pmb \ell=0 \text{ in }\pmb V_+&{}\iff \exists \pmb F \text{ with } \Lambda^{\text{\rm T}} (\rho^{(d)})^{\text{\rm T}}=\varphi_{\pmb g} N\pmb F\\
&{}\iff \exists \pmb F \text{ with }[T_1,T_2,T_3] \Lambda^{\text{\rm T}} (\rho^{(d)})^{\text{\rm T}}=[T_1,T_2,T_3] \varphi_{\pmb g} N\pmb F\\
&{}\iff \exists \pmb F \text{ with } \rho^{(d)}\Lambda\bmatrix T_1\\T_2\\T_3\endbmatrix = \rho^{(d)}M'_{\pmb g} \pmb F\\
&{}\iff \exists \pmb F \text{ with } \pmb \ell= M'_{\pmb g} \pmb F. \qed\endsplit$$
\enddemo

\definition{Notation set up} Complete $g_1,g_2,g_3$ to a $\pmb k$-basis $g_1,g_2,g_3,\dots,g_{d+1}$ of $B_d=\pmb V_1$ and let $E'_{\pmb g}$ be an element of $\operatorname{GL}_{d+1}(\pmb k)$ with
$$[g_1,g_2,g_3,\dots,g_{d+1}]=\rho^{(d)}E'_{\pmb g}.\tag\tnum{in?}$$Let $E_{\pmb g}=\operatorname{Adj}(E'_{\pmb g})$ and $M_{\pmb g}$ be the matrix obtained from $E_{\pmb g}M'_{\pmb g}$ by deleting the first three rows. \enddefinition

\proclaim{Corollary \tnum{52.4}}Adopt Data {\rm\tref{data52}}. The following statements hold.
\item{\rm (a)}The matrix $M_{\pmb g}$ is a presentation matrix for the $T$-module $\pmb V/A_{\pmb g}$. 
\item{\rm (b)} The ideals $L_{\pmb g}$ and $\frak c_{\pmb g}$ of $T$ and $A_{\pmb g}$ are  equal to $I_{d-2}(M_{\pmb g})$ and $I_{d-2}(M_{\pmb g}) A_{\pmb g}$, respectively.  
\endproclaim
\remark{Remark} To obtain $\frak c_{\pmb g}$ from $L_{\pmb g}$, we substitute $g_i$ in place of $T_i$.\endremark
\demo{Proof}(a) Notice that $E_{\pmb g}M'_{\pmb g}$ is a matrix of relations on the row vector $[g_1,\dots,g_{d+1}]$; hence, $M_{\pmb g}$ is a presentation matrix for $$\frac{T g_1+\dots +T g_{d+1}}{Tg_1+Tg_2+Tg_3}=\pmb V/A_{\pmb g}.$$ (b) The proof of Lemma \tref{52.1} shows that $L_{\pmb g}=I_{d-2}(M_{\pmb g})$ since $M_{\pmb g}$ is a presentation matrix for $\pmb V/A_{\pmb g}$. \qed\enddemo

\proclaim{Theorem \tnum{52.5}} Adopt Data {\rm\tref{data52}}. The ideal $\frak c_{\pmb g} B$ is equal to $c_{\pmb g}(x,y)^{d-2}$, for some form $c_{\pmb g}$ in $B$ of degree $(d-1)(d-2)$.\endproclaim

\demo{Proof}Since the  $\pmb V$-module  $B$ contains $\pmb V$ as a   direct summand, we have $\frak c_{\pmb g} B\cap \pmb V=\frak c_{\pmb g}$; and hence, $\pmb V/\frak c_{\pmb g}\hookrightarrow B/\frak c_{\pmb g} B$. Therefore, $$\pmb V/\frak c_{\pmb g}=(B/\frak c_{\pmb g} B)^{(d)}; \tag\tnum{vero}$$ that is, $\pmb V/\frak c_{\pmb g}$ is the $d^{\text{th}}$ Veronese subring of $B/\frak c_{\pmb g} B$. It follows that $B/\frak c_{\pmb g} B$ is also a one-dimensional ring; furthermore, the Hilbert functions of $B/\frak c_{\pmb g} B$ and $\pmb V/\frak c_{\pmb g}$ eventually take the same value. It follows that $B/\frak c_{\pmb g} B$ and $\pmb V/\frak c_{\pmb g}$ have the same multiplicity and, according to Corollary \tref{52.2}~(b), this multiplicity is $(d-1)(d-2)$. The ideal $\frak c_{\pmb g} B$ has height one; so we can write $\frak c_{\pmb g} B=c_{\pmb g}J$ for some form $c_{\pmb g}$ of $B$ and some ideal $J$ of $B$ of height at least $2$. We have
$$(d-1)(d-2)=e(B/\frak c_{\pmb g} B)=e(B/(c_{\pmb g}))=\deg c_{\pmb g}.$$The ideal $\frak c_{\pmb g} B$ is generated by $I_{d-2}(M_{\pmb g}) B$, see Corollary \tref{52.4}; and therefore, the initial degree of $\frak c_{\pmb g} B$ is $d(d-2)$. Compare initial degrees to see that $J$ is generated by forms of degree $(d-2)$. 

We have $\frak c_{\pmb g}=L_{\pmb g}/(f_{\pmb g})$. 
Recall  that  $L_{\pmb g}$ is minimally generated as a $T$-module by $d-1$ forms of degree $d-2$ (see   Lemma \tref{52.1}) and $f_{\pmb g}$ has degree $d$. Hence, as an $A_{\pmb g}$-module, $\frak c_{\pmb g}$ is minimally generated by $d-1$ forms of the same degree; and therefore, as a $B$-module $\frak c_{\pmb g}$ is minimally generated by $d-1$ forms because $B_0=(A_{\pmb g})_0$. Thus, $J$ is minimally generated by $d-1$ forms of degree $d-2$; and therefore, $J=(x,y)^{d-2}$. \qed
\enddemo

\proclaim{Theorem \tnum{52.6}}Adopt Data {\rm\tref{data52}} with $\pmb k$ an algebraically closed field. Factor the polynomial  $c_{\pmb g}$ of Theorem \tref{52.5} as $\prod_{j=1}^{n}\ell_j^{t_j}$, where the  $\ell_j$ are pairwise linearly independent linear forms in $B$ and the $t_j$ are positive integers.
\item{\rm(a)}The sets   $V(\ell_1)\cup\dots\cup V(\ell_n)$ and  $\Psi_{\pmb g}^{-1}(\operatorname{Sing} (\Cal C_{\pmb g}))$ are equal, where   the curve $\Cal C_{\pmb g}$ and the parameterization $\Psi_{\pmb g}\:\Bbb P^1\to \Bbb P^2$ are described in Definition {\rm\tref{.'D27.13}}, and $\operatorname{Sing} (\Cal C_{\pmb g})$ is the singular locus of $\Cal C_{\pmb g}$.
\item{\rm(b)}The sets $V(\ell_1B\cap A_{\pmb g})\cup\dots\cup V(\ell_nB \cap A_{\pmb g})$ and $\operatorname{Sing} (\Cal C_{\pmb g})$ are equal.
\item{\rm(c)}If $P\in \operatorname{Sing} (\Cal C_{\pmb g})$, then $$\delta_P={\tsize\frac 12} \sum\limits_{\{j\mid P=V(\ell_j B\cap A_{\pmb g})\}}t_j.$$
\item{\rm(d)} The equation $t_j=2\delta_{\ell_j}+i_{\ell_j}$ holds, where $\delta_{\ell_j}$ is the $\delta$-invariant of the branch $\ell_j$ and $i_{\ell_j}$ is the intersection multiplicity of $\ell_j$ with the other branches. {\rm(}In the language of Lemma \tref{.gl}, $\delta_{\ell_j}=\delta(\hat{R}/J_j)=\lambda(\overline{(\hat{R}/J_j)}/(\hat{R}/J_j))$ and $i_{\ell_j}=\sum\limits_{i \not=j} \lambda(\hat{R}/(J_i,J_j))$, where $\lambda=\lambda_{\hat{R}}$ is length.{\rm)}
\endproclaim

\demo{Proof} For (b), notice that $\operatorname{Sing} (\Cal C_{\pmb g}) =V(\frak c_{\pmb g})\subseteq \operatorname{Proj}(A_{\pmb g})=\Cal C_{\pmb g}$; otherwise, (a) and (b) are clear. We prove (c). Let $p$ be the prime ideal of $A_{\pmb g}$ which corresponds to  the point $P$ on $\Cal C_{\pmb g}$. We have have
$$\delta_P=\lambda (\overline{\Cal O_{\Cal C_{\pmb g},p}}/\Cal O_{\Cal C_{\pmb g},p})=\lambda (\overline{(A_{\pmb g})_p}/(A_{\pmb g})_p)={\tsize\frac 12} \lambda(\pmb V_p/(\frak c_{\pmb g})_p).$$ The right most equality holds because the ring $(A_{\pmb g})_p$ is Gorenstein and $V_p$ is the integral closure of $(A_{\pmb g})_p$. Observe that 
$\lambda(\pmb V_p/(\frak c_{\pmb g})_p)= \lambda(B_p/\frak c_{\pmb g} B_p)$ since, according to (\tref{vero}), the rings $\pmb V/\frak c_{\pmb g}\subseteq B/\frak c_{\pmb g} B$  are equal on the punctured spectrum.  We have $$\split \delta(\Cal O_{\Cal C_{\pmb g},p})&{}={\tsize \frac 12}\lambda(B_p/\frak c_{\pmb g} B_p)={\tsize \frac 12}\lambda(B_p/  c_{\pmb g} B_p)\\&{}={\tsize\frac 12} \sum\limits_{\{j\mid p=\ell_j B\cap A_{\pmb g}\}}t_j={\tsize\frac 12} \sum\limits_{\{j\mid P=V(\ell_j B\cap A_{\pmb g})\}}t_j.\endsplit $$The proof of (d) may be found in \cite{\rref{Kunz}, Thm.~17.12}. \qed
\enddemo

\proclaim{Theorem \tnum{52.7}} Adopt Data {\rm\tref{data52}} with $d=4$ and $\pmb k$ an algebraically closed field. Let   $\%$ be one of the $13$ possible configurations of singularities on a rational plane quartic as listed in {\rm(\tref{intro})}.  If $\pmb g\in \Bbb T_d$ and the singularities on $\Cal C_{\pmb g}$ are described by $\%$,  then $c_{\pmb g}$ has the form described in the following chart. In each case, the linear forms $\{\ell_j\}$ in $\pmb k[x,y]$ are pairwise linearly independent. 
$$\allowdisplaybreaks \alignat 2
\%\phantom{(2:1,1)^3}&\qquad&c_{\pmb g}\phantom{\ell_4\ell_5\ell_6}&\\ 
(2:1,1)^3&\qquad&\ell_1\ell_2\ell_3\ell_4\ell_5\ell_6&\\
(2:2:1,1),(2:1,1)&\qquad&\ell_1^2\ell_2^2\ell_3\ell_4&\\
(2:1,1)^2,(2:1)&\qquad&\ell_1^2\ell_2\ell_3\ell_4\ell_5&\\
(2:2:2:1,1)&\qquad&\ell_1^3\ell_2^3&\\
(2:2:1),(2:1,1)&\qquad&\ell_1^4\ell_2\ell_3&\\
(2:2:1,1),(2:1)&\qquad&\ell_1^2\ell_2^2\ell_3^2&\\
(2:1,1),(2:1)^2&\qquad&\ell_1^2\ell_2^2\ell_3\ell_4&\\
(2:2:2:1)&\qquad&\ell_1^6&\\
(2:2:1),(2:1)&\qquad&\ell_1^4\ell_2^2&\\
(2:1)^3&\qquad&\ell_1^2\ell_2^2\ell_3^2&\\
(3:1,1,1)&\qquad&\ell_1^2\ell_2^2\ell_3^2&\\
(3:1,1)&\qquad&\ell_1^4\ell_2^2&\\
(3:1)&\qquad&\ell_1^6&\\
\endalignat$$ 
\endproclaim

\demo{Proof}Throughout this discussion,   $\Cal C=\Cal C_{\pmb g}$, $c=c_{\pmb g}$, and the singularities on $\Cal C_{\pmb g}$ are described by $\%$. The polynomial  $c$ is a homogeneous form of degree $(d-2)(d-1)=6$.

If $\%=(2:2:2:1)$, then there is only one singular point on $\Cal C$ and there is only one branch at this singular point; thus, $c=\ell_1^6$. 

If $\%=(2:2:2:1,1)$, then there only one singular point $P$ on $\Cal C$. This singularity has two branches. Each branch is smooth (i.e., $\delta_{\ell_j}=0$) because the multiplicity of the singularity is $2$. Thus, $t_j=i_{\ell_j}$ and $i_{\ell_1}=i_{\ell_2}$. Notice that $$i_{\ell_1}=\lambda (\widehat{\Cal O}_{\Cal C,P}/(q_1+q_2))=i_{\ell_2},$$ where $q_j$ is the minimal prime ideal of $\widehat{\Cal O}_{\Cal C,P}$ which corresponds to $\ell_j$. Thus, $c=\ell_1^3\ell_2^3$. 

If $\%=(2:2:1),(2:1)$, then there are two singular points, $P_1,P_2$, on $\Cal C$, each with one branch: $q_1 \leftrightarrow \ell_1$ and $q_2\leftrightarrow \ell_2$. Theorem \tref{52.6}~(c) shows that $t_1=4$ and $t_2=2$; and therefore, $c=\ell_1^4\ell_2^2$. 
 
If $\%=(2:2:1),(2:1,1)$, then there are two singular points, $P_1,P_2$, on $\Cal C$. The singular point $P_1$ has one branch $q_1\leftrightarrow \ell_1$ and $\delta_{P_1}=2$. The corresponding factor of $c$ is $\ell_1^4$ by part (c) of Theorem \tref{52.6}. The singular point $P_2$ has two branches: $q_2\leftrightarrow \ell_2$ and $q_3\leftrightarrow \ell_3$. The factors $\ell_2$ and $\ell_3$ of $c$ appear with exponent $1$ because $c$ is a form of degree $6$ and there is no room for any higher exponent. Thus, $c=\ell_1^4\ell_2\ell_3$.  

If $\%=(2:2:1,1),(2:1)$, then there are two singular points, $P_1,P_2$, on $\Cal C$. The singular point $P_1$ has two branches:  $q_1\leftrightarrow \ell_1$ and $q_2\leftrightarrow \ell_2$, with $\delta_{P_1}=2$. The exponents $t_1$ and $t_2$ of $\ell_1$ and $\ell_2$ in $c$ satisfy $t_1+t_2=4$   by part (c) of Theorem \tref{52.6}. The two branches are smooth; so $t_1=i_1=i_2=t_2$; hence, $t_1=t_2=2$. The singular point $P_2$ has one branch: $q_3\leftrightarrow \ell_3$. The fact that $\deg c=6$ forces $t_3$ to be $2$. In this case, we have $c=\ell_1^2\ell_2^2\ell_3^2$.

One may proceed through the other eight cases using the same techniques. \qed
\enddemo

\definition{Definition \tnum{jr}}Let $d$ be an even integer. Recall the open subset $\Bbb B_d$ of $\Bbb A_d$ from Definition \tref{.'D27.14}. For each pair of indices $\pmb r=(r_1,r_2)$ with $1\le r_1\le r_2\le d$, define 
$\Bbb B_d^{(j,\pmb r)}$ to be the set of all $\pmb g\in \Bbb B_d$ with
$$\lambda_{0,j}\det\bmatrix \lambda _{r_1,c_1} &\lambda _{r_1,c_2}\\\lambda _{r_2,c_1} &\lambda _{r_2,c_2}\endbmatrix \neq 0,$$where $\{j,c_1,c_2\}=\{1,2,3\}$ and $\lambda_{\pmb g}$ is given in Remark \tref{R5}. \enddefinition 
\remark{Remark \tnum{jroc}} The sets 
$\cup_{j,\pmb r} \Bbb B_d^{(j,\pmb r)}$ form an open cover of $\Bbb B_d$. 
\endremark

\definition{Notation \tnum{52.7.5}}
Recall $\pmb z$, $\pmb R$, and $\pmb S$ as given in Convention (1) of \tref{CP} and view $\pmb V[\pmb z]$ as a subring of $\pmb S$. Recall also  Convention (2) of \tref{CP}. 
Let $\lambda$ be in $\Bbb A^{3d+3}$.
The ring    $\pmb R_{\lambda}$ is defined to be $\pmb R/(\{z_{i,j}-\lambda_{i,j}\})$.   If $S$ is any $\pmb R$-algebra, then $S_{\lambda}$ is defined to be $S\otimes_{\pmb R}\pmb R_{\lambda}$. Furthermore, if $G$ is an element of $S$, then $G|_{\lambda}$ is the image of $G$ in $S_{\lambda}$. We use the same ``evaluation'' notation for ideals of $S$ or matrices with entries in $S$. In other words, if $I$ is an ideal of $S$, then $I|_\lambda$ is the image of $I$ in $S_{\lambda}$ and if $E$ is a matrix with entries in $S$, then $E|_{\lambda}$ is the image of $E$ with entries in $S_{\lambda}$. 
\enddefinition

\proclaim{Lemma \tnum{52.8}} Adopt Data \tref{data52} with $d$ even. Fix one of the  pairs $(j,\pmb r)$ as described in Remark {\rm\tref{jroc}}. Then  there exists  an ideal $C\subseteq \pmb V[\pmb z]$ so that if  $\pmb g$ is in $\Bbb B_{d}^{(j,\pmb r)}$, then  
$C|_{\lambda_{\pmb g}}$   is equal to the conductor $\frak c_{\pmb g}$.  \endproclaim
\remark{Remark} We recall that $C|_{\lambda_{\pmb g}}$ is the image
 of $C$ in $\pmb V[\pmb z]_{\lambda_{\pmb g}}=\pmb V$. \endremark
\demo{Proof} Let $G_k$, for $1\le k\le 3$, be the polynomials in $\pmb V[\pmb z]$ which are defined in (\tref{Gj}), and let $G_4,\dots,G_{d+1}$ be the polynomials 
$$\{y^d,xy^{d-1},\dots,x^d\}\setminus\{y^d,x^{r_1}y^{d-r_1},x^{r_2}y^{d-r_2}\}$$in $\pmb V$. If $\pmb g\in \Bbb B_d^{(j,\pmb r)}$, then the definition of $\Bbb B_d^{(j,\pmb r)}$ shows that 
$$G_1|_{\lambda_{\pmb g}},G_2|_{\lambda_{\pmb g}},G_3|_{\lambda_{\pmb g}},G_4,\dots,G_{d+1}$$ is a basis for $B_d=[(\pmb V[\pmb z])_{\lambda_{\pmb g}}]_d$. Let $E'$ be the following $(d+1)\times (d+1)$ matrix over $\pmb V[\pmb z]$

\def\fakerow{&&&\vrule height 2pt\\ }
\def\hline{\fakerow\noalign{\hrule}\fakerow}
\def\nohline{\fakerow\fakerow}
$$
\normalbaselineskip=0pt\normallineskip=0pt
E'=\bmatrix z_{1,1}&z_{1,2}&z_{1,3}&\vrule\\
\nohline
  z_{2,1}&z_{2,2}&z_{3,3}&\vrule\\
\vdots&\vdots&\vdots&\vrule&\quad*\quad\\
  \vdots&\vdots&\vdots&\vrule&\quad \\
\nohline
 z_{d,1}&z_{d,2}&z_{d,3}&\vrule\\
\endbmatrix$$with $$[G_1,\dots,G_{d+1}]=\rho^{(d)}E'\tag\tnum{in!}$$ and each entry of $*$ is either $0$ or $1$. (The matrix $\rho^{(d)}$ is defined in (\tref{rho}).) If $\pmb g\in \Bbb B_d^{(j,\pmb r)}$, then the image, $E'|_{\lambda_{\pmb g}}$, of $E'$ in $\pmb V[\pmb z]_{\lambda_{\pmb g}}$ is in $\operatorname{GL}_{d+1}(\pmb k)$ and the map $\pmb V[\pmb z]\to \pmb V[\pmb z]_{\lambda_{\pmb g}}=\pmb V$ sends (\tref{in!}) to    (\tref{in?}). 

Recall the matrix $\pmb d_2^{(j)}$ of Theorem \tref{T1}. If $\pmb g\in \Bbb B_d^{(j,\pmb r)}$, then $\pmb d_2^{(j)}|_{\lambda_{\pmb g}}$ is a Hilbert-Burch matrix for the row vector $[g_1,g_2,g_3]$. Form the matrix $M'$ which satisfies 
$$\bmatrix T_1&T_2&T_3\endbmatrix \pmb d_2^{(j)}\bmatrix \rho^{(d/2)}&0\\0&\rho^{(d/2)}\endbmatrix =\rho^{(d)}M',$$exactly as was done in (\tref{52.3-}). Each entry of $M'$ is a homogeneous form in $\pmb k[\{T_k\},\pmb z]$ of degree $1$ in the $T$'s and degree $3c+1$ in the $z$'s. 
The matrix $M$ is obtained from $\operatorname{Adj}(E') M'$ by deleting the first three rows and substituting $G_k$ for $T_k$, for $1\le k\le 3$. Let $C$ be the ideal $I_{d-2}(M)$ in $\pmb V[\pmb z]$.  If $\pmb g\in \Bbb B_d^{(j,\pmb r)}$, then according to Corollary \tref{52.4}, we have
$$C|_{\lambda_{\pmb g}}=I_{d-2}(M_{\pmb g})\cdot A_{\pmb g}=\frak c_{\pmb g}A_{\pmb g}=\frak c_{\pmb g}\pmb V. \qed$$
\enddemo

\proclaim{Theorem \tnum{52.9}} Adopt Data {\rm\tref{data52}} with  $\pmb k$   an algebraically closed field of characteristic zero and $d$ even. 
Define 
$$X=\left\{\pmb g\in \Bbb B_d\left \vert \matrix \format\l\\ c_{\pmb g}=\prod \ell_j^{t_j} \text{with the $\ell_j$ pairwise linearly independent}\\\text{and $t_j\ge 2$ for all $j$}\endmatrix \right.\right\}.$$
Then the subset $X$ of $\Bbb B_d$ is closed.
\endproclaim

\demo{Proof} According to Remark \tref{jroc}, the sets $\cup_{j,\pmb r} \Bbb B_d^{(j,\pmb r)}$ form an open cover of $\Bbb B_d$. Therefore, it suffices to prove that $X\cap \Bbb B_d^{(j,\pmb r)}$ is a closed subset of $\Bbb B_d^{(j,\pmb r)}$ for each pair $(j,\pmb r)$ as described in Definition \tref{jr}. Fix such a pair $(j,\pmb r)$. 
By Lemma \tref{52.8} there exists an ideal $C\subseteq \pmb V[\pmb z]$ so that $C|_{\lambda_{\pmb g}}$ is equal to $\frak c_{\pmb g}$ for each $\pmb g\in \Bbb B_d^{(j,\pmb r)}$. Let $D=\pmb V[\pmb z]/C$. Recall from Notation \tref{52.7.5}  that $D_{\lambda_{\pmb g}}$  is defined    to be $D\otimes_{\pmb R}\pmb R_{\lambda_g}$. It follows that  $D_{\lambda_{\pmb g}}=\pmb V/(C|_{\lambda_{\pmb g}})=\pmb V/\frak c_{\pmb g}$.

Theorem \tref{52.5} ensures that $\frak c_{\pmb g}B=\ell_1^{t_1}\cdots \ell_n^{t_n}(x,y)^{d-2}$
for some $n$ and some pairwise linearly independent forms $\{\ell_j\}$. We see that
$$\text{$t_j\ge 2$ for all $j$}\iff   \operatorname{Reg} (B/\frak c_{\pmb g}B)=\emptyset,$$  where $\operatorname{Reg}(R)$ is the regular locus of the ring $R$. Recall from (\tref{vero}) that $\pmb V/\frak c_{\pmb g}$ is equal to $(B/\frak c_{\pmb g} B)^{(d)}$; hence, 
$\operatorname{Proj} (\pmb V/\frak c_{\pmb g})\cong \operatorname{Proj} (B/\frak c_{\pmb g}B)$, and 
$$\operatorname{Reg} (B/\frak c_{\pmb g}B)=\emptyset\iff \operatorname{Reg} (\pmb V/\frak c_{\pmb g})=\emptyset.$$ Apply the Jacobian criterion to see  that $\operatorname{Reg} (\pmb V/\frak c_{\pmb g})=\emptyset$ if and only if Jacobian ideal $\operatorname{Jac}(\pmb V/\frak c_{\pmb g})$ is nilpotent. 

Let $J=\operatorname{Fitt}_1(\Omega_{D/\pmb k[z]})\subseteq D$. Notice that $\operatorname{Jac}(\pmb V/\frak c_{\pmb g})$ is the image 
$J|_{\lambda_{\pmb g}}$ of $J$ in  $D_{\lambda_{\pmb g}}=\pmb V/\frak c_{\pmb g}$.
Recall from Corollary \tref{52.2}~(b) that $D_{\lambda_{\pmb g}}$ is Cohen-Macaulay and has multiplicity $(d-1)(d-2)$ for all $\pmb g\in \Bbb B_d^{(j,\pmb r)}$. 

We claim that $J|_{\lambda_{\pmb g}}$ is nilpotent if and only if $J|_{\lambda_{\pmb g}}^{(d-1)(d-2)}=0$. Indeed, if $J|_{\lambda_{\pmb g}}$ is nilpotent, then $J|_{\lambda_{\pmb g}}\subseteq q$ for all $q\in \min(D_{\lambda_{\pmb g}})$. However, the length of $(D_{\lambda_{\pmb g}})_q$ is at most 
$  e(D_{\lambda_{\pmb g}})=(d-1)(d-2)$; hence, $q_q^{(d-1)(d-2)}=0$ which forces $(J|_{\lambda_{\pmb g}})_q^{(d-1)(d-2)}$ to be zero. Since this holds for every minimal prime and $D_{\lambda_{\pmb g}}$ is Cohen-Macaulay, we have $J|_{\lambda_{\pmb g}}^{(d-1)(d-2)}=0$. Write $\frak a=J^{(d-1)(d-2)}$. It remains to show that the set 
$$\{\pmb g\in \Bbb B_d^{(j,\pmb r)}\mid \frak a|_{\lambda_{\pmb g}}=0\}\tag\tnum{52.set}$$ is closed in $\Bbb B_d^{(j,\pmb r)}$.

Since the $D_{\lambda_{\pmb g}}$ are Cohen-Macaulay rings, we have that $\frak a|_{\lambda_{\pmb g}}=0$ if and only if $\dim( (D/\frak a)_{\lambda_{\pmb g}})=\dim D_{\lambda_{\pmb g}}$ and $e((D/\frak a)_{\lambda_{\pmb g}})=e(D_{\lambda_{\pmb g}})$. Since $\dim D_{\lambda_{\pmb g}}=1$  and $e(D_{\lambda_{\pmb g}})=(d-1)(d-2)$ for all $\pmb g\in \Bbb B_d^{(j,\pmb r)}$, we conclude that set described in (\tref{52.set}) is equal to 
$$\{\pmb g\in \Bbb B_d^{(j,\pmb r)}\mid \dim ( (D/\frak a)_{\lambda_{\pmb g}})\ge 1\}\cap \{\pmb g\in \Bbb B_d^{(j,\pmb r)}\mid e( (D/\frak a)_{\lambda_{\pmb g}})\ge (d-1)(d-2)\}.$$
But the latter set is closed in $\Bbb B_d^{(j,\pmb r)}$ by the upper semi-continuity of dimension and multiplicity; in other words, we  apply Theorem \tref{L37.7} to $\pmb R\to D/\frak a$.    
\qed 
\enddemo

\proclaim{Proposition \tnum{52.10}} Adopt Data {\rm\tref{data52}} with  $\pmb k$   an algebraically closed field of characteristic zero and $d$ even. For each positive integer $i$, define 
$$X_i=\left\{\pmb g\in \Bbb B_d\left \vert \matrix \format\l\\ c_{\pmb g}=\prod \ell_h^{t_h} \text{with the $\ell_h$ pairwise linearly independent}\\\text{and $t_k\ge i$ for some $k$}\endmatrix \right.\right\}.$$
Then each subset $X_i$ of $\Bbb B_d$ is closed.
\endproclaim

\demo{Proof} Fix $(j,\pmb r)$ as in the proof of Theorem \tref{52.9}. We prove that $X_i\cap \Bbb B_d^{(j,\pmb r)}$ is closed in   $\Bbb B_d^{(j,\pmb r)}$.  Keep the ideal $C$ of $\pmb V[\pmb z]$ from the proof of  Theorem \tref{52.9}. We continue to have $C|_{\lambda_{\pmb g}}$ equal to the ideal $\frak c_g$ of $\pmb V[\pmb z]$ for all $\pmb g$ in $\Bbb B_d^{(j,\pmb r)}$. Let $C_1,\dots,C_{s}$, for some $s$,  be the generators for $C$ in $\pmb V[\pmb z]$ as described in Lemma \tref{52.8}. Recall that $\pmb V\subseteq B$. Let $\frak b_i$ be the ideal of $B[\pmb z]$ which is generated by $C_1,\dots,C_{s}$ together with all partial derivatives
$\frac{\partial ^k C_w}{\partial^{k_1}x\partial^{k_2}y}$ with $1\le w\le s$, $1\le k\le i-1$, and $k_1+k_2=k$. Define the ring $\frak B_i=B[\pmb z]/\frak b_i$. 

Let $\pmb g\in \Bbb B_d ^{(j,\pmb r)}$. We have $(\frak B_1)_{\lambda_{\pmb g}}= B/\frak c_{\pmb g}B$. We saw in Theorem \tref{52.5} that   
$\frak c_{\pmb g}B=\prod \ell_h^{t_h}(x,y)^{d-2}B$ where $\sum t_h=(d-1)(d-2)$ and the $\ell_h$ are pairwise linearly independent linear forms in $B$.   It is clear that $\dim(\frak B_1)_{\lambda_{\pmb g}}=1$. The definition of the ring  $\frak B_i$ shows that  
$(\frak B_i)_{\lambda_{\pmb g}}= B/(\frak c_{\pmb g}B+\frak dB)$, where $\frak d$ is the ideal of $B$ generated by all  partial derivatives
$\frac{\partial ^k \Delta}{\partial^{k_1}x\partial^{k_2}y}$ as $\Delta$ varies over all listed generators of $\frak c_{\pmb g}B$ and the parameters satisfy   $1\le k_1+k_2=k\le i-1$. We see that  if there exists $h$ with $t_h\ge i-1$, then 
$0\le \dim(\frak B_i)_{\lambda_{\pmb g}}\le 1$  and 
$$\text{there exists $k$ with $t_k\ge i$} \iff \dim(\frak B_i)_{\lambda_{\pmb g}}=1.$$
 Theorem \tref{L37.7} shows that for each index $i$, $\{\lambda\in \Bbb A^{3d+3}\mid \dim (\frak B_i)_{\lambda}\ge 1\}$ is a closed subset of $\Bbb A^{3d+3}$. It follows that 
$X_i\cap \Bbb B_d ^{(j,\pmb r)}$, which is equal to $$\{\pmb g\in \Bbb B_d ^{(j,\pmb r)}\mid 
\dim (\frak B_i)_{\lambda_{\pmb g}}\ge 1\},$$is a closed subset of $\Bbb B_d ^{(j,\pmb r)}$. \qed
\enddemo

\bigpagebreak
\SectionNumber=9\tNumber=1
\heading Section \number\SectionNumber. \quad Rational plane quartics: a stratification   and the  correspondence between the Hilbert-Burch matrices and the configuration of singularities.
\endheading

 The configuration of all singularities that can appear on, or infinitely near, a  rational plane quartic are completely determined by two classical formulas:
$$g=\binom {d-1}2-\sum_{q} \binom{m_q}2\qquad\quad\text{and}\quad\qquad
\sum_{q'} m_{q'} \le m_p,$$
The formula on the left was known by Max Noether; see (\tref{MNF}). It gives the genus $g$ of the irreducible plane curve  $\Cal C$ of degree $d$, where $q$ varies over all singularities  on, and infinitely near,    $\Cal C$, and  $m_q$ is the multiplicity at $q$. The formula on the right holds for any point $p$ on any curve $\Cal C$; the points $q'$ vary over all of the points in the first neighborhood of $p$. 
The above
formulas   permit   $9$ possible  singularities on a rational plane curve of degree 4:
$$\matrix
 \text{Classical}&\text{modern}&\text{multiplicity}&(m,\delta,s)\\
\text{name}&\text{name}&\text{sequence}&\\\vspace{5pt}
\text{Node}&A_1&(2:1,1)&(2,1,2)\\
\text{Cusp}&A_2&(2:1)&(2,1,1)\\
\text{Tacnode}&A_3&(2:2:1,1)&(2,2,2)\\
\text{Ramphoid Cusp}&A_4&(2:2:1)&(2,2,1)\\
\text{Oscnode}&A_5&(2:2:2:1,1)&(2,3,2)\\
\text{$A_6$-Cusp}&A_6&(2:2:2:1)&(2,3,1)\\
\text{Ordinary Triple Point}&D_4&(3:1,1,1)&(3,3,3)\\
\text{Tacnode Cusp}&D_5&(3:1,1)&(3,3,2)\\
\text{Multiplicity 3 Cusp}&E_6&(3:1)&(3,3,1)\\
\endmatrix$$   The singularity ``$A_6$'' does not have a consistent classical name and for that reason we have introduced the name ``$A_6$-cusp''. Beware that Namba's terminology is not completely consistent with the terminology used above: he uses ``double cusp'' (respectively, ``ramphoid cusp'') for what we call a ramphoid cusp (respectively, ``$A_6$-cusp'').  The  thirteen possible ways to configure the above singularities on a rational plane quartic are given in (\tref{intro}).  For each  configuration $\%$ from (\tref{intro}) define
$$S_{\%} =\{\pmb g\in \Bbb T_4|\text{the configuration of singularities on $\Cal C_{\pmb g}$ are described by \%}\}.$$
Thus, $\Bbb T_4$ is the disjoint union of $S_{\%}$ as $\%$ roams over the $13$ configurations listed in (\tref{intro}). We use open and closed subsets in $\Bbb T_4$ and $\Bbb B_4$ to separate  the subsets $S_{\%}$. 
We describe the open and closed sets in $\Bbb T_4$ by making use of partial orders on the set of $13$ possible configurations of singularities, exactly as was done in Section 6. In Definition \tref{D8.1} we introduce two partial orders: $\operatorname{QCP}$ and $\operatorname{BQP}$.
The closed sets that we identify in the proof of Proposition \tref{T8.1} come from the techniques of Section 6. The closed sets of Lemma \tref{next} arise from the techniques of Section 7. The techniques of Sections 6 and 7 are not able to separate  
$$S_{(2:2:1),(2:1,1)} \quad \text{from} \quad S_{(2:2:1,1),(2:1)};\tag\tnum{diff}$$ and for that reason we set
$$T_{(2:2:1),(2:1,1)}^{\operatorname{QCP}}=T_{(2:2:1,1),(2:1)}^{\operatorname{QCP}}.$$
The sets of (\tref{diff}) are separated in Proposition \tref{sep}, which uses  the techniques of Section 8. This separation is reflected in the poset $\operatorname{BQP}$. 

The main result of the present section is 
 Corollary \tref{MC9} where  we exhibit a stratification of $\Bbb B_4$ in which every curve associated to a given stratum has the same configuration of singularities. The stratification of Corollary \tref{MC9} is based on the poset $\operatorname{BQP}$. A second major result in this section is Theorem \tref{cid} where we identify a large collection of closed irreducible subsets of $\Bbb T_4$; these subsets are parameterized by the elements of of the poset $\operatorname{QCP}$.  In both Theorems the dimension of each closed irreducible set is given.

 The third major result in this section is Theorem \tref{Clda} where we extend Theorem \tref{d=2c} in order to describe a pretty parameterization of each rational plane quartic $\Cal C$ as a function of the singularity configuration of $\Cal C$. In addition to being a complete classification of parameterizations, Theorem \tref{Clda} is used extensively in our proof of irreducibility in Proposition \tref{P8.3}.

\definition{Definition \tnum{D8.1}} We define two posets: the ``Quartic Configuration Poset'' $\operatorname{QCP}$ and 
the  ``Balanced Quartic Poset'' $\operatorname{BQP}$. 
The elements of  $\operatorname{QCP}$  are the $13$ configurations of (\tref{intro})   once $(2:2:1),(2:1,1)$ has been set equal to $(2:2:1,1),(2:1)$. The order in $\operatorname{QCP}$ is   given in Table 2, where we have drawn $\%'\to \%$ to mean $\%'\le \%$.  For each $\%\in \operatorname{QCP}$, define 
$$T^{\operatorname{QCP}}_{\%} =\bigcup_{\{\%'\in\operatorname{QCP}\mid \%'\le \%\}}S_{\%'}.$$
\midinsert$$\eightpoint \xymatrix{
&&&\boxed{(2:1,1)^3}\\
&&\boxed{\tsize(2:2:1,1),\atop\tsize(2:1,1)}\ar[ru]&\boxed{\tsize(2:1,1)^2,\atop\tsize(2:1)}\ar[u]\\
&\boxed{(2:2:2:1,1)}\ar[ru]&\boxed{{(2:2:1),\atop(2:1,1)}}=\boxed{{(2:2:1,1),\atop(2:1)}}\ar[u]\ar[ru]&\boxed{\tsize(2:1,1),\atop\tsize(2:1)^2}\ar[u]\\
\boxed{(3:1,1,1)}\ar[ru]&\boxed{(2:2:2:1)}\ar[u]\ar[ru]&\boxed{\tsize(2:2:1),\atop\tsize(2:1)}\ar[u]\ar[ru]&\boxed{(2:1)^3}\ar[u]\\
\boxed{(3:1,1)}\ar[u]\ar[ru]\\
\boxed{(3:1)}\ar[u]\ar[rruu]
}$$
\centerline{{\bf Table 2:} The poset $\operatorname{QCP}$.} \endinsert
The elements of $\operatorname{BQP}$ are the $10$ balanced singularity configurations for a rational plane quartic. (These are the configurations which involve singularities of multiplicity $2$.) The order in $\operatorname{BQP}$ is given in Table 3 where we have drawn $\%'\to \%$ in $\operatorname{BQP}$ to mean $\%'\le \%$ in $\operatorname{BQP}$.
For each $\%$ in $\operatorname{BQP}$, define $$T^{\operatorname{BQP}}_{\%}=\bigcup_{\{\%'\in \operatorname{BQP}\mid \%'\le \%\}}S_{\%'}.$$

\midinsert$$\eightpoint \xymatrix{&&\boxed{\tsize(2:1,1)^3}\\
&\boxed{\tsize(2:2:1,1)\atop\tsize(2:1,1)}\ar[ur]&&\boxed{\tsize(2:1,1)^2\atop\tsize(2:1)}\ar[ul]\\
\boxed{\tsize(2:2:2:1,1)}\ar[ur]&\boxed{\tsize(2:2:1)\atop\tsize(2:1,1)}\ar[u]\ar[urr]&&\boxed{\tsize(2:2:1,1)\atop\tsize(2:1)}\ar[ull]\ar[u]&\boxed{\tsize(2:1,1)\atop\tsize(2:1)^2}\ar[ul]\\&\boxed{\tsize (2:2:2:1)}\ar[ul]\ar[u]\ar[urr]&\boxed{\tsize(2:2:1)\atop\tsize(2:1)}\ar[ul]\ar[ur]\ar[urr]&&
\boxed{\tsize(2:1)^3}\ar[u]
}$$\centerline{{\bf Table 3:} The poset $\operatorname{BQP}$.} \endinsert 

\enddefinition \remark{Remark} Let $\%$ and $\%'$ be singularity configurations   $(2:2:1),(2:1,1)$ and $(2:2:1,1),(2:1)$, respectively.  Please notice that $S_{\%}$ and $S_{\%'}$ are disjoint subsets of $\Bbb T_4$; but the subsets $T^{\operatorname{QCP}}_{\%}$ and $T^{\operatorname{QCP}}_{\%'}$ of $\Bbb T_4$ are equal.\endremark

\proclaim{Proposition \tnum{T8.1}} For each $\%\in \operatorname{QCP}$, $T^{\operatorname{QCP}}_{\%}$ is a closed subset of $\Bbb T_4$.
\endproclaim

\demo{Proof} The subset $T^{\operatorname{QCP}}_{(3:1,1,1)}$ of $\Bbb T_4$ is the same as the subset $\Bbb U\Bbb B_4$ and this is closed in $\Bbb T_4$ by Theorem \tref{YYY}. The subset $S_{(2:2:2:1,1)}\cup S_{(2:2:2:1)}$ in the present notation is the same as the subset $T_{2:2:2}$ of $\Bbb B_4$ in the notation of Section 6. Theorem  \tref{strata} shows that $T_{2:2:2}$ is closed in $\Bbb B_4$; hence, $$T^{\operatorname{QCP}}_{(2:2:2:1,1)}=\Bbb U\Bbb B_4\cup T_{2:2:2}$$ is closed in $\Bbb T_4$. Repeat this argument for the closed subsets $T_{2:2,2}$ and $T_{2,2,2}$ of $\Bbb B_4$ to conclude that $T^{\operatorname{QCP}}_{\%}$ is closed in $\Bbb T_4$ for $\%=(2:2:1,1),(2:1,1)$ and for $\%=(2:1,1)^3$. We have shown that  $T^{\operatorname{QCP}}_{\%}$ is closed in $\Bbb T_4$ for each $\%$ across the top row of Table 2.

Lemma \tref{next} shows that $T^{\operatorname{QCP}}_{\%}$ is a closed subset of $\Bbb T_4$ for $\%$
equal to 
$$ (2:1,1)^2,(2:1),\qquad (2:1,1),(2:1)^2,\qquad\text{or}\qquad(2:1)^3.$$That is;   $T^{\operatorname{QCP}}_{\%}$ is closed in $\Bbb T_4$ for each $\%$ down the right hand column of Table 2. 

If $\%$ is an element of $\operatorname{QCP}$ and $\%$ does not live in the top row or the right hand column of Table 2, then    $T^{\operatorname{QCP}}_{\%}$ is the intersection of two closed subsets of $\Bbb T_4$.  \qed 
\enddemo

\proclaim{Lemma \tnum{next}}Take $\pmb R$, $\pmb S$, and $G_j$ from Conventions \tref{CP}. Let $N$ be the matrix 
$$N=\bmatrix \frac{\partial G_1}{\partial x}&\frac{\partial G_2}{\partial x}&\frac{\partial G_3}{\partial x}\\
\frac{\partial G_1}{\partial y}&\frac{\partial G_2}{\partial y}&\frac{\partial G_3}{\partial y}\endbmatrix$$
and $S$ be the ring $\pmb S/I_2(N)$. Let $\pmb R\to S$ be the natural map induced by the inclusion $\pmb R\subseteq \pmb S$. Recall the language of Convention {\rm(4)} from  {\rm\tref{CP}}.
The following statements hold.

\smallskip \flushpar{\bf (1)} The set $X(S;=1,\ge 1)\cap \Bbb T_4=T^{\operatorname{QCP}}_{(2:1,1)^2,(2:1)}$ is a closed subset of $\Bbb T_4$.

\smallskip \flushpar{\bf (2)} The set $X(S;=1,\ge 2)\cap \Bbb T_4=T^{\operatorname{QCP}}_{(2:1,1),(2:1)^2}$ is a closed subset of $\Bbb T_4$.

\smallskip \flushpar{\bf (3)} The set $X(S;=1,\ge 3)\cap \Bbb T_4=T^{\operatorname{QCP}}_{(2:1)^3}$ is a closed subset of $\Bbb T_4$.
\endproclaim

\demo{Proof} Observe first that if $\pmb g\in \Bbb A_d$, then $S_{\lambda_{\pmb g}}$ is equal to the ring 
$B/(\operatorname{Jac} (B/\pmb k[I_{\pmb g}]))$. If $\pmb g\in \Bbb T_4$, then $\operatorname{ht} I_{\pmb g}=2$ and this guarantees that 
$\operatorname{Jac} (B/\pmb k[I_{\pmb g}])$ is not the zero ideal. It follows that the dimension of $S_{\lambda_{\pmb g}}$ is either $0$ or $1$. Thus, $X(S;\ge 1)\cap\Bbb T_4$ is equal to $X(S;= 1)\cap\Bbb T_4$. According to  Theorem \tref{L37.7}, $X(S;\ge 1)\cap\Bbb T_4$ is closed in $\Bbb T_4$. It follows that $X(S;= 1)\cap\Bbb T_4$ is closed in $\Bbb T_4$.
Theorem \tref{L37.7} also shows that 
 $X(S;= 1,\ge i)\cap\Bbb T_4$ is closed in $X(S;=1)\cap\Bbb T_4$; hence,
$X(S;= 1,\ge i)\cap\Bbb T_4$ is closed in $\Bbb T_4$.
Theorem \tref{Cor2'} shows $$X(S;\ge 1,\ge i)\cap \Bbb T_d=\{\pmb g\in \Bbb T_d\mid \sum\limits_{p\in \operatorname{Sing} \Cal C_{\pmb g}} m_p-s_p\ge i\}.$$ It is now easy easy to compute $\sum  m_p-s_p$ for each $\Cal C$ in   $S_{\%}$ for each $\%$ from (\tref{intro}). The order in $\operatorname{QCP}$ was chosen in order to make   $$\split
X(S;=1,\ge 1)\cap \Bbb T_4&{}=T^{\operatorname{QCP}}_{(2:1,1)^2,(2:1)},\\
X(S;=1,\ge 2)\cap \Bbb T_4&{}=T^{\operatorname{QCP}}_{(2:1,1),(2:1)^2}, \text{ and}\\ 
X(S;=1,\ge 3)\cap \Bbb T_4&{}=T^{\operatorname{QCP}}_{(2:1)^3},\endsplit$$and this completes the proof. \qed \enddemo

\proclaim{Proposition \tnum{sep}} 
 For each $\%\in \operatorname{BQP}$, $T^{\operatorname{BQP}}_{\%}$ is a closed subset of $\Bbb B_4$.
 \endproclaim
 \midinsert$$\eightpoint \xymatrix{
&&\boxed{(2:1,1)^3\atop c_{\pmb g}=\ell_1\ell_2\ell_3\ell_4\ell_5\ell_6}\\
&\boxed{\tsize(2:2:1,1),\atop{\tsize(2:1,1) \atop{\tsize c_{\pmb g}=\ell_1^2\ell_2^2\ell_3\ell_4}}}\ar[ru]&\boxed{\tsize(2:1,1)^2,\atop{\tsize(2:1)\atop{\tsize c_{\pmb g}=\ell_1^2\ell_2\ell_3\ell_4\ell_5}}}\ar[u]\\
\boxed{(2:2:2:1,1)\atop{\tsize c_{\pmb g}=\ell_1^3\ell_2^3}}\ar[ru]&\boxed{{(2:2:1),\atop{\tsize(2:1,1)\atop{\tsize c_{\pmb g}=\ell_1^4\ell_2\ell_3}}}}=\boxed{{(2:2:1,1),\atop{\tsize(2:1)\atop{\tsize c_{\pmb g}=\ell_1^2\ell_2^2\ell_3^3}}}}\ar[u]\ar[ru]&\boxed{\tsize(2:1,1),\atop{\tsize(2:1)^2\atop{\tsize c_{\pmb g}=\ell_1^2\ell_2^2\ell_3\ell_4}}}\ar[u]\\
\boxed{(2:2:2:1)\atop{\tsize c_{\pmb g}=\ell_1^6}}\ar[u]\ar[ru]&\boxed{\tsize(2:2:1),\atop{\tsize(2:1)\atop{\tsize c_{\pmb g}=\ell_1^4\ell_2^2}}}\ar[u]\ar[ru]&\boxed{(2:1)^3\atop{\tsize c_{\pmb g}=\ell_1^2\ell_2^2\ell_3^2}}\ar[u]\\
}$$
\centerline{{\bf Table 4:} The poset  ${\operatorname{QCP}}\cap \Bbb B_4$ with $c_{\pmb g}$.} \endinsert
\demo{Proof} Let $\%$ be an element of $\operatorname{BQP}$. Assume first that $\%$ is not equal to  
$$(2:2:1),(2:1,1)\text{ or }(2:2:1,1),(2:1).$$ 
In this case $T^{\operatorname{QCP}}_{\%}\cap \Bbb B_4=T^{\operatorname{BQP}}_{\%}$. We proved in Proposition \tref{T8.1} that $T^{\operatorname{QCP}}_{\%}$
is closed in $\Bbb T_4$. It follows that $T^{\operatorname{QCP}}_{\%}\cap \Bbb B_4$ is closed in the subspace $\Bbb B_4$ of $\Bbb T_4$.    

We now prove that 
the sets $T^{\operatorname{BQP}}_{(2:2:1),(2:1,1)}$ and $T^{\operatorname{BQP}}_{(2:2:1,1),(2:1)}$ are closed in $\Bbb B_4$.
In Table 4, we recorded the poset ${\operatorname{QCP}}\cap {\Bbb B_4}$ together with the form of the factorization of the polynomial $c_{\pmb g}$  for each $\pmb g\in S_{\%}$ for each $\%$ in ${\operatorname{QCP}}\cap \Bbb B_4$. The factorizations of $c_{\pmb g}$ were calculated in Theorem \tref{52.7}. 
Recall the  sets $$X_{\tref{52.9}}=\{\pmb g\in \Bbb B_d\mid c_{\pmb g}=\prod \ell_j^{t_j}\text{ with all $t_j\ge 2$}\}$$ and
$$X_{4,\tref{52.10}}=\{\pmb g\in \Bbb B_d\mid c_{\pmb g}=\prod \ell_j^{t_j}\text{ with at least one $t_j\ge 4$}\}$$  from Theorems \tref{52.9} and \tref{52.10}. A quick look at Table 4 shows 
 that $$X_{\tref{52.9}}\cap (T^{\operatorname{QCP}}_{(2:2:1,1),(2:1)}\cap \Bbb B_4)=
S_{(2:2:1,1),(2:1)}\cup S_{(2:2:1),(2:1)}\cup S_{(2:2:2:1)}=T^{\operatorname{BQP}}_{(2:2:1,1),(2:1)}$$and
$$X_{4,\tref{52.10}}\cap (T^{\operatorname{QCP}}_{(2:2:1),(2:1,1)}\cap \Bbb B_4)=
S_{(2:2:1),(2:1,1)}\cup S_{(2:2:1),(2:1)}\cup S_{(2:2:2:1)}=T^{\operatorname{BQP}}_{(2:2:1),(2:1,1)}.$$

Theorems \tref{52.9} and \tref{52.10} show that $X_{\tref{52.9}}$ and $X_{4,\tref{52.10}}$ are closed subsets of $\Bbb B_4$ and Proposition \tref{T8.1} shows that $$ T^{\operatorname{QCP}}_{(2:2:1),(2:1,1)}\cap \Bbb B_4=T^{\operatorname{QCP}}_{(2:2:1,1),(2:1)}\cap \Bbb B_4$$
is a closed subset of $\Bbb B_4$. We conclude that $T^{\operatorname{BQP}}_{(2:2:1,1),(2:1)}$ and $T^{\operatorname{BQP}}_{(2:2:1),(2:1,1)}$ are closed subsets of $\Bbb B_4$. \qed
\enddemo

Theorem \tref{Clda} gives the canonical form for a parameterization of some representative of the right equivalence class of each  rational plane quartic 
 $\Cal C$. If  $\Cal C$ has any singularities of multiplicity $2$, then all of its singularities have multiplicity $2$ and the techniques of Section 4 apply to all of the singularities. 
For this reason, the  description of the Hilbert-Burch matrix 
in the present result is much more detailed 
than the description of       
Theorem \tref{d=2c}.
In addition to being a very pretty classification, Theorem \tref{Clda} is the starting point in the proof of Proposition \tref{P8.3} about the irreducibility of many subsets of $\Bbb T_4$. 
(We first used the phrase ``signed maximal order minors''   in Remark \tref{401}~(2).)

\proclaim{Theorem \tnum{Clda}} Let $\Cal C$ be a rational plane quartic over an algebraically closed field. Then there exists a linear automorphism $\Lambda$ of $\Bbb P^2$ so that $\Lambda \Cal C$ is parameterized by the signed maximal order minors of $\varphi_{\%}$, where $\varphi_{\%}$ is given below. In each case, the entries of $\varphi_{\%}$ must be chosen so that $\operatorname{ht} I_2(\varphi_{\%})=2$ and $\mu(I_1(\varphi_{\%}))\ge 3$. We write $\ell_i$ for linear form, $Q_i$ for quadratic form, and $C_i$ for cubic form. 

\smallskip\flushpar{\bf(1)} If the singularity configuration of $\Cal C$ is given by $\%$ is equal to one of
$$(2:1,1)^3,\qquad  (2:1,1)^2,(2:1),\qquad  {(2:1,1),(2:1)^2},\qquad \text{or}\qquad   (2:1)^3,$$ then
$$\varphi_{\%}= \bmatrix Q_1&Q_1\\Q_2&0\\0&Q_3\endbmatrix.$$
\itemitem {\rm(a)}If $\%=(2:1,1)^3$, then each $Q_i$ is the product of $2$ non-associate linear forms in $B$. The singularities of $\Lambda \Cal C$ are $[0:0:1]$, $[0:1:0]$, and $[1:0:0]$. Each singularity is an ordinary node.

\itemitem {\rm(b)}If $\%=(2:1,1)^2,(2:1)$, then each of the polynomials  $Q_1$ and $Q_2$ is the product of $2$ non-associate linear forms in $B$ and $Q_3$ is a perfect square. The singularities of $\Lambda \Cal C$ are $[0:0:1]$, $[0:1:0]$, and $[1:0:0]$. The    singularities $[1:0:0]$, $[0:1:0]$ are ordinary nodes; and $[0:0:1]$ is a cusp.
\itemitem {\rm(c)}If $\%=(2:1,1),(2:1)^2$, then    $Q_1$ is the product of $2$ non-associate linear forms in $B$ and each of the polynomials $Q_2,Q_3$ is a perfect square. The singularities of $\Lambda \Cal C$ are $[0:0:1]$, $[0:1:0]$, and $[1:0:0]$. The    singularity $[1:0:0]$ is an ordinary node; and   $[0:1:0]$ and  $[0:0:1]$ are cusps.
 \itemitem {\rm(d)}If $\%= (2:1)^3$, then    all three polynomials $Q_1$, $Q_2$, and $Q_3$ are perfect squares.    The singularities of $\Lambda \Cal C$ are $[0:0:1]$, $[0:1:0]$, and $[1:0:0]$. Each     singularity   is a  cusp.
 \smallskip\flushpar{\bf(2)} If the singularity configuration of $\Cal C$ is given by $\%$ is equal to
$$\matrix\format \l\quad&\l\quad\\(2:2:1,1),(2:1,1),&(2:2:1,1),(2:1),\\ (2:2:1),(2:1,1),\ \text{or}&(2:2:1),(2:1),\endmatrix $$ then 
$$\varphi_{\%}=\bmatrix Q_1&0\\Q_2&Q_3\\0&Q_2\endbmatrix.$$
\itemitem {\rm(a)}If $\%=(2:2:1,1),(2:1,1)$, then $Q_1$ and $Q_2$ are   the product of $2$ non-associate linear forms in $B$. The singularities of $\Lambda \Cal C$ are $[0:0:1]$,   and $[1:0:0]$. The singularity $[1:0:0]$ is an  ordinary node and $[0:0:1]$ is a tacnode. 
\itemitem {\rm(b)}If $\%=(2:2:1,1),(2:1)$, then $Q_1$ is a perfect square and   $Q_2$   is  the product of $2$ non-associate linear forms in $B$. The singularities of $\Lambda \Cal C$ are $[0:0:1]$,   and $[1:0:0]$. The singularity $[1:0:0]$ is an  cusp and $[0:0:1]$ is a tacnode. 
\itemitem {\rm(c)}If $\%=(2:2:1),(2:1,1)$, then $Q_2$ is a perfect square and   $Q_1$   is  the product of $2$ non-associate linear forms in $B$. The singularities of $\Lambda \Cal C$ are $[0:0:1]$,   and $[1:0:0]$. The singularity $[1:0:0]$ is an  ordinary node and $[0:0:1]$ is a ramphoid cusp. 
\itemitem {\rm(d)}If $\%=(2:2:1),(2:1)$, then $Q_1$ and $Q_2$ are both  perfect squares.  The singularities of $\Lambda \Cal C$ are $[0:0:1]$,   and $[1:0:0]$. The singularity $[1:0:0]$ is a cusp and $[0:0:1]$ is a ramphoid cusp. 
\smallskip\flushpar{\bf(3)} If the singularity configuration of $\Cal C$ is given by $\%$ is equal to
$$(2:2:2:1,1)\qquad  \text{or}\qquad (2:2:2:1) $$ then 
$$\varphi_{\%}=\bmatrix Q_1&Q_2\\Q_3&Q_1\\0&Q_3\endbmatrix.$$
\itemitem {\rm(a)}If $\%=(2:2:2:1,1)$, then  $Q_3$ is   the product of $2$ non-associate linear forms in $B$. The only singularity on $\Lambda \Cal C$ is $[0:0:1]$ and this  singularity  is an  oscnode.
\itemitem {\rm(b)}If $\%=(2:2:2:1)$, then  $Q_3$ is   a perfect square. The only singularity on $\Lambda \Cal C$ is $[0:0:1]$ and this  singularity  is an  $A_6$-cusp.
\smallskip\flushpar{\bf(4)} If the singularity configuration of $\Cal C$ is given by $\%$ is equal to
$$(3:1,1,1),\qquad    (3:1,1),\qquad\text{or}\qquad (3:1), $$ then 
$$\varphi_{\%}=\bmatrix \ell_1&C_1\\\ell_2&C_2\\0&C_3\endbmatrix.$$
\itemitem {\rm(a)}If $\%=(3:1,1,1)$, then  $C_3$ is   the product of $3$ non-associate linear forms in $B$. The only singularity on $\Lambda \Cal C$ is $[0:0:1]$ and this  singularity  is an  ordinary triple point.
\itemitem {\rm(b)}If $\%=(3:1,1)$, then  $C_3$ has the form $\ell_3^2\ell_4$, where $\ell_3$ and $\ell_4$    non-associate linear forms from  $B$. The only singularity on $\Lambda \Cal C$ is $[0:0:1]$ and this  singularity  is a tacnode cusp.
\itemitem {\rm(c)}If $\%=(3:1)$, then  $C_3$ is a perfect cube.   The only singularity on $\Lambda \Cal C$ is $[0:0:1]$ and this  singularity  is a multiplicity three cusp.
\endproclaim

\demo{Proof} The basic form of $\varphi_{\%}$ is described by Theorem \tref{d=2c} when $\Cal C$ has a multiplicity $2$ singularity and by Corollary \tref{ctgl} when $\Cal C$ has a multiplicity $3$ singularity. The information about the number of non-associate linear factors is an immediate consequence of Theorem \tref{Cor1}. \qed \enddemo

\proclaim{Proposition \tnum{P8.3}} 

\flushpar{\bf (a)} The subsets $T^{\operatorname{QCP}}_{(3:1)}$,  $T^{\operatorname{QCP}}_{(3:1,1)}$,  $T^{\operatorname{QCP}}_{(3:1,1,1)}$, $T^{\operatorname{QCP}}_{(2:1,1)^3}$, and $T^{\operatorname{QCP}}_{(2:1,1)^2,(2:1)}$   of $\Bbb T_4$ are irreducible.

\smallskip\flushpar{\bf (b)} The subset $T^{\operatorname{BQP}}_{\%}$   of  $\Bbb B_4$ is irreducible for all ten $\%$ in the poset $\operatorname{BQP}$.
\endproclaim
\demo{Proof} The subset $T^{\operatorname{QCP}}_{(2:1,1)^3}$ of $\Bbb T_4$ is by definition equal to all of $\Bbb T_4$. Theorem \tref{XXXX} shows that $\Bbb T_4$ is an open subset of the irreducible space $\Bbb A_4$. It follows that $\Bbb T_4$ is also irreducible.
 The subsets $T^{\operatorname{BQP}}_{(2:2:2:1,1)}$ and   $T^{\operatorname{BQP}}_{(2:2:1,1),(2:1,1)}$ of $\Bbb B_4$ are called $S_{2:2:2}$ and $S_{2:2,2}$ in Section 6. Theorem \tref{strata} shows that these sets are irreducible. We use the following trick repeatedly throughout the rest of the argument. 

\medskip\flushpar (\tnum{trk}) If $Y$ is an irreducible variety and $\theta\:Y\to \Bbb A_d$ is a morphism, then $\operatorname{im} \theta\cap \Bbb T_d$ is an irreducible set.

\demo{Proof of {\rm(\tref{trk})}} The set $\Bbb T_d$ is open in $\Bbb A_d$; so, $\theta^{-1}(\Bbb T_d)$ is an open subset of  the irreducible set $Y$. It follows that $Y\cap \theta^{-1}(\Bbb T_d)$ is an irreducible set. The image of the irreducible set $Y\cap \theta^{-1}(\Bbb T_d)$ under the morphism is necessarily irreducible; and this image is equal to $\theta(Y)\cap \Bbb T_d$. The proof of  (\tref{trk}) is complete.
\enddemo

We consider a situation which is very similar to Definition \tref{H} and Remarks \tref{401}. Let $\Cal G$ be the multiplicative group of matrices 
$$\Cal G=\left\{\left.\bmatrix u_1&Q\\0&u_2\endbmatrix \right\vert u_1,u_2\in \pmb k^*, Q\in B_2\right\}$$ and $G'$ be the group $\operatorname{GL}_3(\pmb k)\times \Cal G$. The group $G'$ acts on the space of matrices
$$\Bbb H'=\left\{\left.\bmatrix \ell_1&C_1\\\ell_2&C_2\\\ell_3&C_3\endbmatrix\right\vert \ell_i\in B_1\ \text{and}\ C_i\in B_3\right\}$$ by $(\chi,\xi)\cdot \varphi=\chi \varphi\xi^{-1}$ for $\chi\in \operatorname{GL}_3(\pmb k)$, $\xi\in \Cal G$, and $\varphi\in \Bbb H'$. Let $\Phi$ be the morphism which is defined in (\tref{PHI}).

Consider the morphism $\theta\: G'\times B_1^3\times B_3^2\to \Bbb A_4$, which is given by
$$\theta(g,\ell_1,\ell_2,\ell_3,C_1,C_2)=\Phi\left (g\cdot \bmatrix \ell_1&C_1\\\ell_2&C_2\\0&\ell_3^3\endbmatrix\right ).$$
Theorem \tref{Clda} shows that $\operatorname{im} \theta\cap \Bbb T_4=T^{\operatorname{QCP}}_{(3:1)}$. The domain of $\theta$ is an irreducible variety; so (\tref{trk}) shows that $T^{\operatorname{QCP}}_{(3:1)}$ is an irreducible subset of $\Bbb T_4$. Let 
$$\Delta=\{C\in B_3\mid \text{the discriminant of $C$ is zero}\}.$$  Define $\theta\: G'\times B_1^2\times B_3^2\times \Delta\to \Bbb A_4$ by 
$$\theta(g,\ell_1,\ell_2,C_1,C_2,C)=\Phi\left (g\cdot \bmatrix \ell_1&C_1\\\ell_2&C_2\\0&C\endbmatrix\right ).$$The domain of $\theta$ is an irreducible variety and Theorem \tref{Clda} shows that $\operatorname{im} \theta\cap \Bbb T_4$ is $T^{\operatorname{QCP}}_{(3:1,1)}$. It follows from (\tref{trk})   that $T^{\operatorname{QCP}}_{(3:1,1)}$ is an irreducible subset of $\Bbb T_4$. We abbreviate the argument. Consider the morphisms 
$$\matrix\format\l&\quad\l&\quad\l\\ 
\theta_1\:G'\times B_1^2\times B_3^3\to \Bbb A_4, 
&\theta_2\:G\times B_1\times B_2^3\to \Bbb A_4,
&
\theta_3\:G\times B_1\times B_2^2\to\Bbb A_4, \\\vspace{4pt}
\theta_4\:G\times B_1^2\times B_2\to\Bbb A_4,
&\theta_5\:G\times B_1^3\to\Bbb A_4,
&\theta_6\:G\times B_1^2\times B_2^2\to\Bbb A_4,\\\vspace{4pt}
\theta_7\:G\times B_1\times B_2^4\to\Bbb A_4, &\theta_8\:G\times \Bbb A^3\times B_1\times B_2^2\to \Bbb A_4
\endmatrix$$ which are given by
$$ \theta_1(g,\ell_1,\ell_2,C_1,C_2,C_3)=\Phi\left (g\cdot \bmatrix \ell_1&C_1\\\ell_2&C_2\\0&C_3\endbmatrix\right ),$$
$$\theta_2(g,\ell,Q_1,Q_2,Q_3)=\Phi\left (g\cdot \bmatrix Q_1&Q_2\\\ell^2&Q_3\\0&\ell^2\endbmatrix\right ),  $$

$$\theta_3(g,\ell,Q_1,Q_2)=\Phi\left (g\cdot \bmatrix Q_1&Q_2\\\ell^2&Q_1\\0&\ell^2\endbmatrix\right ),$$
$$\theta_4(g,\ell_1,\ell_2,Q)=\Phi\left (g\cdot \bmatrix \ell_1^2&0\\\ell_2^2&Q\\0&\ell_2^2\endbmatrix\right ),$$
$$\theta_5(g,\ell_1,\ell_2,\ell_3)=\Phi\left (g\cdot \bmatrix \ell_1^2&\ell_1^2\\\ell_2^2&0\\0&\ell_3^2\endbmatrix\right ),$$
$$\phantom{\text{ and}}\theta_6(g,\ell_1,\ell_2,Q_3,Q_4)=\Phi\left (g\cdot \bmatrix Q_3&Q_4\\\ell_2^2&0\\0&\ell_1^2\endbmatrix\right ),$$
$$\theta_7(g,\ell_1,Q_1,Q_2,Q_3,Q_4)=\Phi\left (g\cdot \bmatrix Q_1&Q_2\\Q_3&Q_4\\0&\ell_1^2\endbmatrix\right )\text{ and}$$
$$\theta_8( g,a,b,c,\ell,Q_1,Q_2)=\Phi\left (g \bmatrix Q_1&Q_2\\a\ell^2&Q_1\\b\ell^2&c\ell^2\endbmatrix\right ).$$
The domain of each of these morphisms is irreducible. Theorem \tref{Clda} shows that 
$$\operatorname{im} \theta_i\cap \Bbb T_4=\cases 
T^{\operatorname{QCP}}_{(3:1,1,1)}&\text{if $i=1$}\\
T_{(2:2:1),(2:1,1)}^{\operatorname{BQP}}&\text{if $i=2$}\\
T^{\operatorname{BQP}}_{(2:2:2:1)}&\text{if $i=3$}\\
T^{\operatorname{BQP}}_{(2:2:1),(2:1)}&\text{if $i=4$}\\
T^{\operatorname{BQP}}_{(2:1)^3}&\text{if $i=5$}\\
T^{\operatorname{BQP}}_{(2:1,1),(2:1)^2}&\text{if $i=6$}\\
T^{\operatorname{BQP}}_{(2:1,1)^2,(2:1)}&\text{if $i=7$}\\
T^{\operatorname{BQP}}_{(2:2:1,1),(2:1)}&\text{if $i=8$}.\endcases \tag\tnum{tnim}$$
We offer some details which are involved in the proof of (\tref{tnim}). First we show that $$\operatorname{im} \theta_6\cap \Bbb T_4\subseteq T^{\operatorname{QCP}}_{(2:1,1),(2:1)^2}\cap \Bbb B_4.$$ Let $$\varphi=\bmatrix Q_3&Q_4\\\ell_2^2&0\\0&\ell_1^2\endbmatrix,$$ with $\Phi(\varphi)\in \Bbb T_4$. There are two cases. Either $Q_3\in {<}\ell_1^2,\ell_2^2{>}$ or $Q_3\not\in {<}\ell_1^2,\ell_2^2{>}$. In the first case, there is an element $g$ of $G$ with 
$$g\varphi=\bmatrix \ell_2^2&0\\\ell_1^2&Q_4\\0&\ell_1^2\endbmatrix;$$ hence, $\Phi(\varphi)\in S_{(2:2:1),(2:1)}$. In the second case, $\ell_1^2,\ell_2^2,Q_3$ is a basis for $B_2$ and $Q_4=\alpha_1\ell_1^2+\alpha_2\ell_2^2+\alpha_3Q_3$, for some constants $\alpha_1,\alpha_2,\alpha_3$. Once again there are two cases: either $\alpha_3=0$ or $\alpha_3\neq 0$. If $\alpha_3=0$, then there exists $g\in G$ with $$g\varphi=\bmatrix \ell_1^2&0\\\ell_2^2&Q_3\\0&\ell_2^2\endbmatrix;$$hence, $\Phi(\varphi)\in S_{(2:2:1),(2:1)}$. If $\alpha_3\neq 0$, then there exists $g\in G$ and $Q_3'\in B_2$ such that
$$g\varphi=\bmatrix Q_3'&Q_3'\\\ell_2^2&0\\0&\ell_1^2\endbmatrix;$$hence, $\Phi(\varphi)\in S_{(2:1,1),(2:1)^2}\cup S_{(2:1)^3}$. 

It is obvious that $T^{\operatorname{QCP}}_{(2,1,1)^2,(2:1)}\cap \Bbb B_4 \subseteq \operatorname{im} \theta_7\cap \Bbb T_4$. We show the converse. Suppose that $\pmb g\in \operatorname{im} \theta_7\cap \Bbb T_4$. It is clear that $\pmb g\in \Bbb B_4$. Also, the curve $\Cal C_{\pmb g}$ which is parameterized by $\pmb g$ has at least one singularity with exactly one branch. A quick look at Table 2 shows that $\pmb g\in T^{\operatorname{QCP}}_{(2,1,1)^2,(2:1)}$.

We show $\operatorname{im} \theta_8\cap \Bbb T_4= T^{\operatorname{BQP}}_{(2:2:1,1),(2:1)}$. We first prove the inclusion  ``$  \supseteq$''. Set $a=c=0$, $b=1$, and $g=1$ to get 
$$\Phi\left (\bmatrix Q_1&Q_2\\ 0&Q_1\\ \ell^2&0\endbmatrix\right ).$$Thus, $S_{(2:2:1,1),(2:1)}\cup S_{(2:2:1),(2:1)}$ is contained in $\operatorname{im} \theta_8\cap \Bbb T_4$.  Set $a=c=1$, $b=0$, and $g=1$ to get 
$$\Phi\left ( \bmatrix Q_1&Q_2\\\ell^2&Q_1\\0&\ell^2\endbmatrix\right ).$$Thus, $S_{(2:2:2:1)}$ is contained in $\operatorname{im} \theta_8\cap \Bbb T_4$. The set $T^{\operatorname{BQP}}_{(2:2:1,1),(2:1)}$ is equal  
$$S_{(2:2:1,1),(2:1)}\cup S_{(2:2:1),(2:1)}\cup S_{(2:2:2:1)}.$$ We have shown that 
$\operatorname{im} \theta_8\cap \Bbb T_4\supseteq T^{\operatorname{BQP}}_{(2:2:1,1),(2:1)}$. Now we prove the direction `` $ \subseteq$''. Let
$$\varphi=\bmatrix Q_1&Q_2\\a\ell^2&Q_1\\b\ell^2&c\ell^2\endbmatrix.$$ Assume that $\Phi(\varphi)\in \Bbb T_4$. 
There are two cases: either $b\neq 0$ or $b=0$. If $b\neq 0$, then $\varphi$  may be quickly transformed  into
$$\bmatrix Q_1&Q_2-\frac cb Q_1\\a\ell^2&Q_1-\frac {ac}b\ell^2 \\b\ell^2&0\endbmatrix;$$thus, $\varphi$ may be transformed into
$$\bmatrix Q_1-\frac {ac}b\ell^2 &Q_2-\frac cb Q_1\\0&Q_1-\frac {ac}b\ell^2 \\\ell^2&0\endbmatrix.$$We see that  $\Phi(\varphi)$  is in $S_{(2:2:1,1),(2:1)}\cup S_{(2:2:1),(2:1)}$. If $b=0$, then the hypothesis $\operatorname{ht} \Phi(\varphi)=2$ forces $ac$ to be non-zero. One may transform $\varphi$ into $$\bmatrix Q_1&aQ_2\\\ell^2&Q_1\\0&\ell^2\endbmatrix.$$
We see that  $\Phi(\varphi)$  is in $S_{(2:2:2:1)}$.

Now that (\tref{tnim}) is established, we apply   (\tref{trk}) to see that 
most of the subsets $T^{\operatorname{QCP}}_{\%}$ from (a) and all of the subsets $T^{\operatorname{BQP}}_{\%}$ from (b) are irreducible. It remains to show that $T^{\operatorname{QCP}}_{(2:1,1)^2,(2:1)}$ is irreducible. Consider the morphism 
$$\theta\: B_1\times B^2_2\times B_4\times \operatorname{GL}_3(\pmb k)\to \Bbb A_4,$$which is given by
$$\theta(\ell,Q_1,Q_2,F,\chi)=(\ell^2Q_1,\ell^2Q_2,F)\chi.$$In light of (\tref{trk}) it suffices to show that $\operatorname{im} \theta \cap \Bbb T_4= T^{\operatorname{QCP}}_{(2:1,1)^2,(2:1)}$. The inclusion ``$\supseteq$'' may be read from Theorem \tref{Clda}. To prove ``$\subseteq$'', we recall from Lemma \tref{next} that
$$T^{\operatorname{QCP}}_{(2:1,1)^2,(2:1)}=\{g\in \Bbb T_4\mid \deg \gcd I_2(N|_{\lambda_{\pmb g}})\ge 1\},$$where
$$N|_{\lambda_{\pmb g}}= \bmatrix \frac{\partial g_1}{\partial x}&\frac{\partial g_2}{\partial x}&\frac{\partial g_3}{\partial x}\\
\frac{\partial g_1}{\partial y}&\frac{\partial g_2}{\partial y}&\frac{\partial g_3}{\partial y}\endbmatrix.$$One can quickly check that if $\pmb g=(x^2Q_1,x^2Q_2,F)$, then $x$ divides $\gcd I_2(N|_{\lambda_{\pmb g}})$. \qed
\enddemo

In Corollary \tref{MC9} we exhibit a stratification of $\Bbb B_4$ where every curve associated to a given stratum has the same configuration of singularities; this is the main result in the present section. The hard work in the proof of this result is carried out in the proof of Theorem \tref{T53.3}.

\proclaim{Theorem \tnum{T53.3}}If  $\%$ is in $\operatorname{BQP}$, then $T^{\operatorname{BQP}}_{\%}$ is a closed irreducible subset of $\Bbb B_4$ and the dimension of $T^{\operatorname{BQP}}_{\%}$ is given in the following table{\rm:}
$$\eightpoint \matrix
\dim T^{\operatorname{BQP}}_{\%}& \%\\
15&(2:1,1)^3\\
14&(2:2:1,1)(2:1,1); & (2:1,1)^2,(2:1)\\
13&(2:2:2:1,1);&(2:2:1),(2:1,1);&(2:2:1,1),(2:1);& (2:1,1),(2:1)^2  \\
12&(2:2:2:1);& (2:2:1),(2:1);& (2:1)^3.\\\endmatrix$$
\endproclaim

\demo{Proof} Proposition \tref{sep} shows that $T_{\%}^{\operatorname{BQP}}$ is closed in $\Bbb B_4$; Proposition \tref{P8.3} shows that $T_{\%}^{\operatorname{BQP}}$ is  irreducible; and Proposition   \tref{P8.9} calculates its dimension. \qed\enddemo
 
\proclaim{Corollary \tnum{MC9}} The poset $\operatorname{BQP}$ gives a stratification of $\Bbb B_4$. In other words,
\roster
\item $\Bbb B_4$ is the disjoint union of $\{S_{\%}\mid \%\in \operatorname{BQP}\}$, 
\item if $\%'\le \%$ in $\operatorname{BQP}$, then $S_{\%'}$ is contained in the closure $\overline{S_{\%}}$ of $S_{\%}$ in $\Bbb B_4$, and 
\item $T_{\%}^{\operatorname{BQP}}$ is the closure of $S_{\%}$ in $\Bbb B_4$ for all $\%$ in $\operatorname{BQP}$.
\endroster\endproclaim 
\demo{Proof} This result follows immediately from Theorem \tref{T53.3} using the technique of the proof of Corollary \tref{May2}~(d) and Theorem \tref{stratif}. \qed\enddemo 

Theorem \tref{cid} is our best result about closed irreducible subsets of $\Bbb T_4$. 
\proclaim{Theorem \tnum{cid}} 
\flushpar{\bf (1)} The subsets $S_{(3:1)}$, $S_{(2:2:2:1)}$, $S_{(2:2:1),(2:1)}$, and $S_{(2:1)^3}$ of $\Bbb T_4$ are irreducible, closed in $\Bbb T_4$, and have   dimension $12$.
\smallskip\flushpar{\bf (2)} The subset $T^{\operatorname{QCP}}_{(3:1,1)}$ of $\Bbb T_4$ is irreducible, closed in $\Bbb T_4$, and has  dimension $13$.

 \smallskip\flushpar{\bf (3)} The subsets $T^{\operatorname{QCP}}_{(3:1,1,1)}$ and $T^{\operatorname{QCP}}_{(2:1,1)^2,(2:1)}$ of $\Bbb T_4$ are irreducible, closed in $\Bbb T_4$, and have   dimension $14$.

\smallskip\flushpar{\bf (4)} The set $T^{\operatorname{QCP}}_{(2:1,1)^3}$ is equal to $\Bbb T_4$. This set is irreducible, closed in $\Bbb T_4$, and has dimension $15$.
\endproclaim

\demo{Proof} The sets $S_{(2:2:2:1)}$ and $S_{(2:2:1),(2:1)}$ are shown to be closed in $\Bbb T_4$ in Proposition \tref{Prom}. The other listed subsets  are shown to be closed in   $\Bbb T_4$ in Proposition \tref{T8.1}. (Recall that  $S_{(2:1)^3}=T^{QCP}_{(2:1)^3}$ and $S_{(3:1)}=T^{QCP}_{(3:1)}$.) All of the sets are shown to be irreducible in Proposition  \tref{P8.3}.  The dimensions are calculated in Proposition \tref{P8.9}. \qed\enddemo

\proclaim{Observation \tnum{O9.10}} Each subset $S_{\%}$, with $\%$ from {\rm(\tref{intro})}, has dimension at least $12$.
\endproclaim

\demo{Proof}Recall  the group $\Sigma=\operatorname{GL}_1(\pmb k)\times \operatorname{SL}_2(\pmb k)\times \operatorname{SL}_3(\pmb k)$ 
from (\tref{sm}). We define an action of $\Sigma$ on $\Bbb T_d$. 
If $\left [\smallmatrix a&c\\b&d\endsmallmatrix\right ]$ is an element of $\operatorname{SL}_2 (\pmb k)$, then $ax+by,cx+dy$ is a basis for vector space $B_1$ of linear forms in the polynomial ring $B=\pmb k[x,y]$. If $\sigma=(u,(\ell_1,\ell_2),\chi)$ is an element of $\Sigma$, where the pair of linear forms $(\ell_1,\ell_2)$ represents an element of $\operatorname{SL}_2(\pmb k)$,
and
  $\pmb g=(g_1,g_2,g_3)$ is an element of $\Bbb T_d$, then define $$\sigma\cdot \pmb g=u[\pmb g(\ell_1,\ell_2)]\chi,$$ where  $\pmb g(\ell_1,\ell_2)$ is the ordered triple 
$$(g_1(\ell_1,\ell_2),g_2(\ell_1,\ell_2),g_3(\ell_1,\ell_2)),\tag\tnum{ot}$$  $[\pmb g(\ell_1,\ell_2)]\chi$ is the product of the ordered triple (\tref{ot}) with the matrix $\chi$, and $u\pmb h$ is $(uh_1,uh_2,uh_3)$ for any scalar $u$ and any ordered triple $\pmb h=(h_1,h_2,h_3)$ in $\Bbb A_d$.
It is clear that the curves $\Cal C_{\pmb g}$ and $\Cal C_{\sigma\cdot \pmb g}$ have the same configuration of singularities; and therefore, $\Sigma \cdot S_{\%}=S_{\%}$ for all $\%$ from (\tref{intro}). 

Fix $\%$. We have seen that the set $S_{\%}$ is non-empty. Fix $\pmb g\in S_{\%}$ and define
 $$\theta_{\pmb g}: \Sigma\to S_{\%}$$ by $\theta_{\pmb g}(\sigma)=\sigma\cdot \pmb g$.
Let $\overline{\operatorname{im} \theta_{\pmb g}}$  and $\overline{S_{\%}}$ be the closure of $\operatorname{im} \theta_{\pmb g}$ and $S_{\%}$, respectively,  in $\Bbb T_d$.
The group $\Sigma$ is an irreducible closed algebraic variety of dimension $12$. It follows that $\operatorname{im} \theta_{\pmb g}$; and therefore, $\overline{\operatorname{im} \theta_{\pmb g}}$ are irreducible subsets of $\Bbb T_d$.
Thus,  $\theta_{\pmb g}\: \Sigma \to \overline{\operatorname{im} \theta_{\pmb g}}$ is a  dominate morphism of closed irreducible affine algebraic varieties over an algebraically closed field. A quick calculation shows that if $\pmb h$ is in $\operatorname{im} \theta_{\pmb g}$, then
the fiber $\theta_{\pmb g}^{-1}(\pmb h)$ of $\theta_{\pmb g}$ over $\pmb h$ is a finite set of points; so $\dim \theta_{\pmb g}^{-1}(\pmb h)=0$. The Fiber Dimension Theorem applies, exactly as it did in (\tref{form1}) and (\tref{form2}), to yield that
$$\dim (\overline{\operatorname{im} \theta_{\pmb g}})=\dim (\Sigma)- \dim (\text{a generic fiber})=12-0=12.$$We have
$\overline{\operatorname{im} \theta_{\pmb g}}$ is a closed irreducible subset of $\overline{S_{\%}}$ of dimension $12$.
We conclude that $\dim (S_{\%})=\dim (\overline{S_{\%}})\ge   \dim (\overline{\operatorname{im} \theta_{\pmb g}})=12$.\qed\enddemo

\remark{\bf Remark \tnum{R9.11}} The proof of Observation \tref{O9.10} shows that if $S$ is any non-empty subset of $\Bbb T_d$ which is closed under the action of $\Sigma$, and $\overline {S}$ is the closure of $S$ in $\Bbb T_d$, then every irreducible component of $\overline{S}$ has dimension at least $12$.\endremark
\proclaim{Proposition \tnum{P8.9}} The statements about dimension in Theorems {\rm\tref{T53.3}} and {\rm\tref{cid}} are correct. \endproclaim
\demo{Proof} Consider
$$S_{(3:1)}=T^{\operatorname{QCP}}_{(3:1)}\subsetneq T^{\operatorname{QCP}}_{(3:1,1)}\subsetneq T^{\operatorname{QCP}}_{(3:1,1,1)}\subsetneq T^{\operatorname{QCP}}_{(2:1,1)^3}.$$
Each set is irreducible by Proposition \tref{P8.3} and closed in $\Bbb T_4$ by Proposition \tref{T8.1}. One has
$$12 \le \dim T^{\operatorname{QCP}}_{(3:1)}< \dim T^{\operatorname{QCP}}_{(3:1,1)}< \dim T^{\operatorname{QCP}}_{(3:1,1,1)}<\dim T^{\operatorname{QCP}}_{(2:1,1)^3}=15.$$
One concludes that 
$$\dim T^{\operatorname{QCP}}_{(3:1)}=12,\quad \dim T^{\operatorname{QCP}}_{(3:1,1)}=13,\quad  \text{and}\quad \dim T^{\operatorname{QCP}}_{(3:1,1,1)}=14.$$
Consider 
$$S_{(2:2:2:1)}=T^{\operatorname{BQP}}_{(2:2:2:1)} \subsetneq T^{\operatorname{BQP}}_{(2:2:2:1,1)}\subsetneq T^{\operatorname{BQP}}_{(2:2:1,1),(2:1,1)}
\subsetneq T^{\operatorname{BQP}}_{(2:1,1)^3}=\Bbb B_4.$$ Each set is irreducible by Proposition \tref{P8.3} and the fact that $\Bbb B_4$ is an open subset of the irreducible space $\Bbb A_4$; hence, $\Bbb B_4$ is also an irreducible space. Each set is closed in $\Bbb B_4$ by Proposition  \tref{sep}. One has
$$  12\le \dim \left (T^{\operatorname{BQP}}_{(2:2:2:1)}\right ) <\dim  \left (T^{\operatorname{BQP}}_{(2:2:2:1,1)}\right ) <\dim \left (T^{\operatorname{BQP}}_{(2:2:1,1),(2:1,1)}\right ) <\dim \Bbb B_4=15; $$ 
and therefore,
$$   \dim \left (T^{\operatorname{BQP}}_{(2:2:2:1)}\right )=12,\   \dim  \left (T^{\operatorname{BQP}}_{(2:2:2:1,1)}\right )=13,\ \   \dim \left (T^{\operatorname{BQP}}_{(2:2:1,1),(2:1,1)}\right )=14. \tag {\tnum{two}}$$
 Consider 
$$T^{\operatorname{BQP}}_{(2:1)^3}\subsetneq T^{\operatorname{BQP}}_{(2:1,1),(2:1)^2}\subsetneq T^{\operatorname{BQP}}_{(2:1,1)^2,(2:1)}\subseteq T^{\operatorname{BQP}}_{(2:1,1)^3}=\Bbb B_4.$$All four sets are closed in $\Bbb B_4$ and irreducible. It follows that 
$$  12\le \dim T^{\operatorname{BQP}}_{(2:1)^3}< \dim \left (T^{\operatorname{BQP}}_{(2:1,1),(2:1)^2}\right ) <\dim\left (T^{\operatorname{BQP}}_{(2:1,1)^2,(2:1)}\right ) <\dim \Bbb B_4=15, $$
and $$  \dim T^{\operatorname{BQP}}_{(2:1)^3}=12,\quad  \dim \left (T^{\operatorname{BQP}}_{(2:1,1),(2:1)^2}\right )= 13,\quad \dim\left (T^{\operatorname{BQP}}_{(2:1,1)^2,(2:1)}\right )=14.$$ One may also consider the following chain of closed irreducible subsets of $\Bbb T_4$:
$$T^{\operatorname{QCP}}_{(3:1,1)}\subsetneq T^{\operatorname{QCP}}_{(2,1,1)^2,(2:1)}\subsetneq T^{\operatorname{QCP}}_{(2:1,1)^3}$$ to conclude $\dim T^{\operatorname{QCP}}_{(2,1,1)^2,(2:1)}=14$.

Consider $$S_{(2:2:1),(2:1)}=T^{\operatorname{BQP}}_{(2:2:1),(2:1)}\subsetneq T^{\operatorname{BQP}}_{(2:2:1),(2:1,1)} ,T^{\operatorname{BQP}}_{(2:2:1,1),(2:1)} \subsetneq T^{\operatorname{BQP}}_{(2:2:1,1),(2:1,1)}.$$All four   sets   are closed in $\Bbb B_4$ and irreducible by Proposition   \tref{P8.3}.   
The dimension of the set on the right was calculated in (\tref{two}). We have
$$  12\le \dim T^{\operatorname{BQP}}_{(2:2:1),(2:1)}<\dim T^{\operatorname{BQP}}_{(2:2:1),(2:1,1)}, \dim  T^{\operatorname{BQP}}_{(2:2:1,1),(2:1)}  <\dim  T^{\operatorname{BQP}}_{(2:2:1,1),\atop{(2:1,1)}} =14; $$and we conclude that
$$\dim T^{\operatorname{BQP}}_{(2:2:1),(2:1)}=12\quad\text{and}\quad \dim T^{\operatorname{BQP}}_{(2:2:1),(2:1,1)}= T^{\operatorname{BQP}}_{(2:2:1,1),(2:1)}=13. \qed$$  
\enddemo

\proclaim{Proposition \tnum{Prom}}
The sets $S_{(2:2:2:1)}$ and $S_{(2:2:1),(2:1)}$ are  closed in $\Bbb T_4$.\endproclaim
\demo{Proof} Let $\%$ equal $(2:2:2:1)$ or $(2:2:1),(2:1)$. We saw in Proposition \tref{P8.9} that $\dim S_{\%}=12$. Let $\overline{S_{\%}}$ be the closure of $S_{\%}$ in $\Bbb T_4$. We saw in Proposition \tref{T8.1} that $T_{\%}^{QCP}$ is a closed subset of $\Bbb T_4$. It follows that $\overline{S_{\%}}\subseteq T_{\%}^{QCP}$. On the other hand, we know from Table 2 that $T_{\%}^{QCP}\cap \Bbb B_4=S_{\%}$. To prove this result we must show that $\overline{S_{\%}}\cap \Bbb U\Bbb B_4=\emptyset$. (Recall from Definition \tref{.'D27.14} that $\Bbb U\Bbb B_4$ is the closed subset $\Bbb T_4\setminus \Bbb B_4$ of $\Bbb T_4$.) Suppose that  $\overline{S_{\%}}\cap \Bbb U\Bbb B_4$ is non-empty. It is clear that $\overline{S_{\%}}\cap \Bbb U\Bbb B_4$ is closed under the action of the group $\Sigma$. Let $X$ be a closed irreducible component of $\overline{S_{\%}}\cap \Bbb U\Bbb B_4$. It follows from Remark \tref{R9.11} that $\dim X\ge 12$. It is not possible to have two closed irreducible sets $X\subsetneq \overline{S_{\%}}$ with $\dim X\ge 12$ and $\dim \overline{S_{\%}}=12$. \qed
\enddemo

\medskip

\flushpar {\bf Acknowledgment.} Critical conversations which lead to this work occurred when the authors were all in the same city. The authors are especially    appreciative of the travel support    which brought them together at the Pan-American Advanced Study Institute (PASI) conference on
Commutative Algebra and its Connections to Geometry, in honor of 
Wolmer Vasconcelos,   in Olinda, Brazil in August, 2009 and the Midwest Algebra, Geometry and their Interactions Conference (MAGIC'10)    in Notre Dame, Indiana in April, 2010.

\bigpagebreak


\enddocument